The University of Sheffield

D Singh

# The Moduli Space Of Stable $N$ Pointed Curves Of Genus Zero

Submitted for the degree of PhD

in the Department of Pure Mathematics

2004

## Summary

In this thesis I give a new description for the moduli space of stable n pointed curves of genus zero and explicitly specify a natural isomorphism and inverse between them that preserves many important properties. I also give a natural description for the universal curve of this space. These descriptions are explicit and defined in a straight forward way. I also compute the tangent bundle of this space. In the second part of the thesis I compute the ordinary integral cohomology ring from the above description and specify a basis for it.

## Declaration

No portion of this work has been submitted in support of an application for another degree or qualification of this or any other University or other institute of learning.

# Contents

CHAPTER 1

# Introduction

The theory of moduli spaces is a well established tool in modern algebraic geometry. Here we study the moduli space of n-pointed stable curves of genus zero denoted $\overline{\mathcal{X}}_n$ and give a new very concrete description of it. A stable curve of genus zero is a complex algebraic curve with some marked points, satisfying conditions to be described later. $\overline{\mathcal{X}}_n$ is the set of isomorphism classes of stable curves with $n$ marked points. It has a natural structure as a smooth projective variety over $\mathbb{C}$. It was introduced by Grothendieck and has been widely studied. The cohomology of $\overline{\mathcal{X}}_n$ has been determined by Keel [9], but his answer is not very explicit. The point of the thesis is to introduce a new, more explicit model for $\overline{\mathcal{X}}_n$ and to deduce a new, more explicit description of its cohomology. There are many other descriptions of $\overline{\mathcal{X}}_n$, however they are non-explicit and difficult to analyze. For example it has been shown that $\overline{\mathcal{X}}_n$ may be constructed as certain iterated blowups, see Keel [9] or Kapranov [8], for example, who give two different such representations in this form. Keel uses his description to compute its Chow ring, prove we have an homology isomorphism and to give a recursive definition of its Poincaré series. We will show that our approach is closer to Kapranov's construction but be more explicit.

We will denote our space $\overline{\mathcal{M}}_S$ and define it as a subspace of a product of projective spaces $\overline{\mathcal{M}}_S \subseteq \prod PV_T$ where the product is taken over each $T \subseteq S$ with $|T|$ at least 2. $V_T$ is a complex vector space of dimension $|T| - 1$ . Thus our space comes equipped with projections to each $PV_T$. There are natural projection maps $\pi_U^T : V_T \to V_U$ whenever $U \subseteq T \subseteq S$ which have a natural composition rule. We will use these to impose the conditions on our subspace in a straightforward way. From this it will immediately follow that $\overline{\mathcal{M}}_S$ is a projective variety. We then proceed to analyze $\overline{\mathcal{M}}_S$ and prove it is the stated space. We also calculate its cohomology ring $H^*(\overline{\mathcal{M}}_S, \mathbb{Z})$ whose presentation is computationally easy to work with. In particular our work will give an independent proof that it is finitely generated by its elements of degree 2 and free as an abelian group, the rank of $H^2(\overline{\mathcal{M}}_S)$ is $2^n - 1 - n - \frac{n(n-1)}{2}$. That is we have one generator for each subset $T$ of $S$ of size at least 3, for subsets $T$ of size 2 $PV_T$ is just a point. These classes are just the pullback of the standard generators of $H^*(PV_T, \mathbb{Z})$ under the projection maps $\pi_T : \overline{\mathcal{M}}_S \to PV_T$. We also give a basis $\overline{B}[S]$ for this ring which uses the combinatorics of trees.

Let $S$ be a finite set of size $n$ at least 2 and write $S_+ = S \amalg \{0\}$. By a generic S-curve we will mean a pair $(C, x)$ where $C$ is an algebraic curve isomorphic to $\mathbb{C}P^1$ and $x$ is an injective map $x : S_+ \to C$, see example 1 below for a pictorial representation of such a curve.

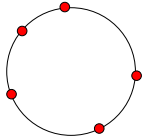

We write $[C, x]$ for the isomorphism class of the S-curve $(C, x)$. In this situation two such curves $(C, x)$ and $(D, y)$ are isomorphic if there is an algebraic isomorphism $\theta : C \to D$ with $\theta \circ x = y$, that is the fixed points are sent to the fixed points and their order is preserved. For $C = D = \mathbb{C}P^1$ we know any such map is a rational homogenous polynomial of degree 1. We define $\mathcal{X}_S$ to be the set of isomorphism classes of all generic S-curves. It is straightforward to identify $\mathcal{X}_S$ with a Zariski open dense subset of $\mathbb{C}P^{|S|-2}$ that consists of the complement of a finite union of hyperplanes, all the diagonals. In this document we will study a certain compactification $\overline{\mathcal{X}}_S$ of $\mathcal{X}_S$ that was first constructed by Deligne and Mumford in the 1960's. These spaces have been studied extensively and may be constructed in several ways, using geometric invariant theory, Chow quotients, or iterated blowups. They form an operad in the category of schemes, which is important in quantum cohomology and conformal field theory. The complex points of $\overline{\mathcal{X}}_S$ parameterize the isomorphism classes of stable S-curves. Knudsen and Mumford proved in the late 1970's that $\overline{\mathcal{X}}_S$ is a projective variety. The proof was developed over 3 papers and relies heavily on the theory of sheaves. From our definition this will be immediate.

An important combinatorial invariant attached to any S-curve is that of its tree type. For $\overline{\mathcal{M}}_S$ we will define the notion of a tree that is right for our construction. More generally we define a forest on $S$ as a collection $\mathcal{F}$ of subsets of $S$ such that for every $U, V \in \mathcal{F}$ either $U \subseteq V$ or $V \subseteq U$ or $U \cap V$ is empty and an $S$-tree $\mathcal{T}$ to be a forest on $S$ with $S \in \mathcal{T}$. The combinatorics of these trees and forests will feature heavily throughout this document in the structure of $\overline{\mathcal{M}}_S$ and later in its cohomology ring when we construct a basis for it, although in a different context.

In chapter 2 we begin by defining the vector spaces $V_T$, one for each $T \subseteq S$ of complex dimension $|T| - 1$ and maps between them $\pi_U^T : V_T \to V_U$ which are required in our definition. We then define our space $\overline{\mathcal{M}}_S$ as a subspace of $\prod PV_T$. It will be immediate from our definition that $\overline{\mathcal{M}}_S$ is a projective variety. We then generalize our construction to a space $\overline{\mathcal{M}}_\mathcal{L}$ defined in an analogous way as a subspace of $\prod_{T \in \mathcal{L}} PV_T$ where $\mathcal{L}$ is a collection of subsets of $S$. It is again immediate that this is a projective variety. It is this approach that is important for our analysis of $\overline{\mathcal{M}}_S$. Of particular importance to us is the study of the space $\overline{\mathcal{M}}_\mathcal{L}$ when $\mathcal{L}$ is an $S$-tree $\mathcal{T}$. In this case we will show that $\overline{\mathcal{M}}_\mathcal{T}$ is the total space of an iterated projective bundle. We will use this case, which is the easiest, to analyze more complicated examples by comparing them with the evident projection map onto these tree cases. In particular this approach will make it easy to prove that everything we consider is smooth and irreducible.

In chapter 3 we define and develop some simple combinatorial objects such as forests and trees and attach to them various numbers that we use throughout this thesis in both the algebraic and topological aspects. A notion of particular use is that of the length of a tree, this is equivalent to its usual combinatorial counterpart. We will use this inductively when studying the space $\overline{\mathcal{M}}_\mathcal{T}$. Many of the results presented here are trivial and included only for completeness. We also gather some one-off results that fit more neatly into this chapter.

In chapter 4 we will develop a catalogue of results that will form the basis for the understanding of our space $\overline{\mathcal{M}}_S$. In particular to any $S$-tree $\mathcal{T}$ we define a space $\overline{\mathcal{M}}_\mathcal{T}$ as a subset of $\prod PV_T$ with the product taken over $\mathcal{T}$, the space $\overline{\mathcal{M}}_\mathcal{T}$ is a natural generalization to that of $\overline{\mathcal{M}}_S$. Our analysis of this space is important for many later results. We also show that to study $\overline{\mathcal{M}}_\mathcal{F}$ when $\mathcal{F}$ is a forest easily reduces to the case of trees. We will show that $\overline{\mathcal{M}}_\mathcal{T}$ can be seen as the total space to a tower of projective bundles, this forms part of a more general result regarding fibre bundles. From this it will immediately follow that $\overline{\mathcal{M}}_\mathcal{T}$ is a smooth irreducible projective variety of complex dimension $|S| - 2$. We will then use the embedding $i : \overline{\mathcal{M}}_\mathcal{T} \to \prod PV_T$ to compute the tangent space of $\overline{\mathcal{M}}_\mathcal{T}$ as a sub-bundle of $T \prod PV_T = \prod \hom(M_T, V_T/M_T)$. This computation will easily generalize to more complicated collections other than trees and in particular we will be able to compute the tangent bundle of $\overline{\mathcal{M}}_S$. This we consider in chapter 6. Lastly in this chapter we will show that $\overline{\mathcal{M}}_\mathcal{T}$ can also be obtained from the iterated blowup of $PV_S$ along certain projective subspaces and their strict transforms. This result will be easily obtained and we will use it later to prove a more general result but in the same spirit. In particular we will be able to establish the existence of an isomorphism from $\overline{\mathcal{X}}_S$ to $\overline{\mathcal{M}}_S$.

In chapter 5 we are mainly concerned with the computation of the cohomology ring of $\overline{\mathcal{M}}_{\mathcal{F}}$ for forests $\mathcal{F}$, we will denote this by $R_{\mathcal{F}}$. This ring will help us to simplify our picture of the cohomology ring of $\overline{\mathcal{M}}_S$ which we denote by $R_S$. We will use the iterated projective bundle description of $\overline{\mathcal{M}}_{\mathcal{T}}$ for trees $\mathcal{T}$ given in chapter 4 to inductively construct the cohomology ring of $\overline{\mathcal{M}}_{\mathcal{T}}$, this calculation will be straightforward. We then deduce the general case for forests $\mathcal{F}$, this will just be an application of the Künneth theorem. In particular we will show that $R_{\mathcal{F}}$ finitely generated by its elements of degree 2 and free as an abelian group, these generators will be the Euler classes over pullbacks of tautological line bundles. The rank of $R_{\mathcal{F}}^2$ will be the number of elements of $\mathcal{F}$ of size at least 3. We will compute two different bases for this ring, the second description will be explicit in its use of forests and will extend naturally to our ring $R_S$. To finish this section we specify necessary and sufficient conditions for a monomial $x$ of $R_{\mathcal{F}}$ to be zero. In particular we will show that any monic monomial of maximum degree that is non-zero is equal to the top class, again we will see how to generalize these results to the ring $R_S$. By the end of this chapter we will have a good algebraic understanding of the ring $R_{\mathcal{F}}$ as well as a topological impression. This will serve as a good model when we come to consider the ring $R_S$ in chapter 9.

In chapter 6 we examine the topological structure of the space $\overline{\mathcal{M}}_{\mathcal{L}}$ for a quite general class of $\mathcal{L}$. These sets we call thickets whose elements consists of subsets of $S$. The advantages of doing this is that we obtain many general results without much extra work and in particular we will be able to use a result by Kapranov [**8**] to prove that $\overline{\mathcal{X}}_S$ and $\overline{\mathcal{M}}_S$ are in fact isomorphic. Although this will not be as satisfactory as we would like, this we explain in a moment. To start with we attach to every element $\underline{M} \in \overline{\mathcal{M}}_{\mathcal{L}}$ a tree $\mathcal{T}$, trees will encode the minimal amount of data that one requires to describe an element of $\overline{\mathcal{M}}_{\mathcal{L}}$ from its image in $\overline{\mathcal{M}}_{\mathcal{T}}$. Then under the projection map $\pi : \overline{\mathcal{M}}_{\mathcal{L}} \to \overline{\mathcal{M}}_{\mathcal{T}}$ we will be able to compare a Zariski open set of $\overline{\mathcal{M}}_{\mathcal{L}}$ with one of $\overline{\mathcal{M}}_{\mathcal{T}}$, a space we are already familiar with, moreover the open set of $\overline{\mathcal{M}}_{\mathcal{L}}$ will be saturated. These open sets will cover $\overline{\mathcal{M}}_{\mathcal{L}}$ and enable us to show that $\overline{\mathcal{M}}_{\mathcal{L}}$ is a smooth irreducible projective variety of dimension $|S|-2$. We will then proceed, again using results developed in chapter 4, to compute the tangent bundle of $\overline{\mathcal{M}}_{\mathcal{L}}$. In particular all of these result apply to $\overline{\mathcal{M}}_S$. We will then prove that the cohomology ring of $\overline{\mathcal{M}}_{\mathcal{L}}$ denoted $R_{\mathcal{L}}$ is well behaved, that is finitely generated by its elements of degree 2 and free as an abelian group. We will also compute the total degree of our ring in a combinatorial way using forests. This we will use in the last chapter to compute a basis for $R_S$. The final section in this chapter is a generalization of the blowup description from chapter 4. More precisely we will show that the space $\overline{\mathcal{M}}_\mathcal{L}$ may be seen as a blowup of $PV_S$ over linear spaces and their strict transforms. This in particular will enable us to prove that $\overline{\mathcal{X}}_S$ and $\overline{\mathcal{M}}_S$ are isomorphic, while this is good to know at this stage it only stipulates the existence of an isomorphism without specifying it. This is of course not a satisfactory answer to our question but is useful in at least that it provides limits to what form the functions intuitively may take. We will explain this at the end of this section. We give an explicit isomorphism in chapter 8 where we apply a more careful analysis to $\overline{\mathcal{X}}_S$ and morphisms from it to the projective space $PV_S$.

In chapter 7 we only consider the space $\overline{\mathcal{M}}_S$. We start by examining the various notions of trees that can arise, we then explains how we consider these different representations equivalent by specifying various bijections between them. These bijections will preserve several important properties. There are well known stratifications of the space $\overline{\mathcal{X}}_S$ by the combinatorics of trees. Several important results are known about these subspaces which we state in the next chapter when we examine in more detail the space $\overline{\mathcal{X}}_S$. We will prove the analogue of these results for our space here, this will mean considering subspaces of $\overline{\mathcal{M}}_S$ that consists of elements of like tree type. We denote these spaces by $\mathcal{M}_S(\mathcal{T})$ and also consider there closures $\overline{\mathcal{M}}_S(\mathcal{T})$ where $\mathcal{T}$ is a tree. These results will be required in chapter 8 to prove neatly that a certain morphism from $\overline{\mathcal{X}}_S$ to $\overline{\mathcal{M}}_S$ is bijective and thus an isomorphism as everything is smooth and complete. However these results are useful in their own right to understand the structure of the space $\overline{\mathcal{M}}_S$. In the last section of this chapter we apply a more local analysis and consider the morphism $\pi : \overline{\mathcal{M}}_{S+} \to \overline{\mathcal{M}}_S$ together with its structure sections. In particular we will show that the restricted map $\pi : \pi^{-1}(\mathcal{M}_S) \to \mathcal{M}_S$ is a $\mathbb{C}P^1$ bundle. We also explain in this section how to construct a map of sets from $\overline{\mathcal{M}}_S$ to $\overline{\mathcal{X}}_S$, although we will not see that it is a morphism of varieties here. This map it will turn out is the inverse to the morphism already mentioned.

In chapter 8 we consider in more detail the space $\overline{\mathcal{X}}_S$. Most of this chapter is already well known. Our contribution will appear towards the end. Parts of chapter 8 are closely related to unpublished notes by Professor N.P. Strickland. We begin by defining the points of $\overline{\mathcal{X}}_S$, that is we define what a stable $S$-curve is. Firstly we define the generic elements of our space, a representative of which is a $S$-marked copy of $\mathbb{C}P^1$. We then proceed to define the general $S$-curve. Each irreducible component of this is a generic $T$-curve for some collection $T$. We continue by defining the graph associated to each element $[C, x] \in \overline{\mathcal{X}}_S$ and prove that it is a tree and this is equivalent to having $H^1(C, \mathcal{O}_C) = 0$ which is part of the definition of an $S$-curve. After defining families of $S$-curves and stating Mumford's crucial theorem, that proves we have a moduli space which can be given the structure of a variety, we give a catalogue of results that are already well established. We will require these later. We will define a regular map $\theta_S : \overline{\mathcal{X}}_S \to PV_S$ this will be analogous to a construction by Kapranov [7]. It is well known that maps into projective spaces are characterized by locally free sheaves of rank one over the domain scheme that are generated by global sections. We use the universal properties of $\overline{\mathcal{X}}_S$ to analyze this map. We will then use the contraction maps defined by Knudsen $\pi_T^S : \overline{\mathcal{X}}_S \to \overline{\mathcal{X}}_T$ to construct regular morphisms $\theta_T^S : \overline{\mathcal{X}}_S \to PV_T$ given by $\theta_T^S = \theta_T \pi_T^S$ which will then give us an induced regular morphism of varieties $\theta_S : \overline{\mathcal{X}}_S \to \prod PV_T$. We will prove that the image actually lies in $\overline{\mathcal{M}}_S$ and that it is an isomorphism. Then we show that the projection map $\pi_S : \overline{\mathcal{M}}_{S+} \to \overline{\mathcal{M}}_S$ is the universal curve whose induced map $p_S : \overline{\mathcal{M}}_S \to \overline{\mathcal{X}}_S$ specified by $p_S(\underline{M}) = \pi_S^{-1}(\underline{M})$ is the inverse for $\theta_S$.

In chapter 9, the final part of this thesis, we are concerned with the computation of the integral cohomology ring of $\overline{\mathcal{M}}_S$. We will define a ring $R_S$ as the quotient of the polynomial ring $\mathbb{Z}_S$ over $\mathbb{Z}$ by a homogeneous ideal $I_S$. Because $\overline{\mathcal{M}}_S \subseteq \prod PV_T$ there is one generator $y_T$ for each $T$ contained in $S$ with size $|T| > 2$. When the size of $T$ is 2 the space $PV_T$ is just a point. We will show that this gives us the rank of $H^2\overline{\mathcal{M}}_S$ as $2^n - 1 - n - \binom{n}{2}$ where $n = |S|$. One of the relations in $I_S$ is then $y_T^{|T|-1} = 0$ inherited from the cohomology of the space $PV_T$. We have shown in chapter 6 that $H^2(\overline{\mathcal{M}}_S)$ is generated by elements of degree two. These are the Euler classes of a natural collection of line bundles, specifically the pullback of the tautological vector bundles over $PV_T$ by the projections $\pi : \overline{\mathcal{M}}_S \to PV_T$. We specify a natural ring map $r_S : R_S \to H^*(\overline{\mathcal{M}}_S)$ and prove this is an isomorphism after analyzing the ring $R_S$ further. We also compute a basis $\overline{B}[S]$ for $R_S$ that uses the combinatorics of forests and is a natural extension to the work of chapter 5. We next specify necessary and sufficient conditions for elements of $R_S$ to be zero, this requires us to further develop some combinatorial ideas that we started in chapter 5. In the next section we define a smaller set of relations that defines an ideal $J_S$ and we will show that $J_S = I_S$ this result is of combinatorial interest and easy to prove. We then explicitly compute the natural map $\text{cl} : A^*(\overline{\mathcal{M}}_S) \to H^*(\overline{\mathcal{M}}_S)$ from the Chow ring to the cohomology ring that we proved in chapter 6 was an isomorphism. This result should be compared with Keel in [**9**] who also proves we have an homology isomorphism. In the final section of this chapter we offer an alternative approach to analyzing the ring $R_S$ that does not require us to compare it with the cohomology ring $H^*(\overline{\mathcal{M}}_S)$ but we do not prove any of the claims. This is where we leave our analysis of the space $\overline{\mathcal{M}}_S$.

We aim in this thesis to convince the reader that our approach to this moduli space is useful for geometric intuition and computation and hope that our description will be a useful alternative to other representations.

# CHAPTER 2

# The definition of $\overline{\mathcal{M}}_S$ and generalizations

In this chapter we define the main objects that are required in the definition of $\overline{\mathcal{M}}_S$ and give some of the notation we will be using throughout this thesis.

DEFINITION 2.0.1. Let $n$ be an integer greater than 1 and $S$ a finite set of size $n$. Then for any subset $T$ of $S$ whose size $|T|$ is at least 1 we define the following spaces.

$$\begin{aligned} F(T,\mathbb{C}) &= \{\, f : T \to \mathbb{C} \mid f \text{ a function}\,\} \\ C_T &= \{\, f \in F(T,\mathbb{C}) \mid f \text{ is constant}\,\} \\ V_T &= F(T,\mathbb{C})/C_T \\ PV_T &= \text{ associated projective space} \end{aligned}$$

REMARK 2.0.2. For each subset $T$ of $S$, $F(T,\mathbb{C})$ is a complex vector space of dimension $|T|$ and we give it the usual topology by identifying it with $\mathbb{C}^{|T|}$ in the standard way. For each $T \subseteq S$ of size at least one $C_T$ is a 1 dimensional subspace of $F(T,\mathbb{C})$ consisting of the constant functions. We write $q_T : F(T,\mathbb{C}) \to V_T$ for the quotient map, $V_T$ is then a complex vector space of dimension $|T| - 1$. Throughout this document we shall use the notation $\leq$ to denote "is a subspace of"

DEFINITION 2.0.3. For any sets $U \subseteq T \subseteq S$ we have the natural projection maps $p_U^T : F(T,\mathbb{C}) \to F(U,\mathbb{C})$ given by $p_U^T(f) = f|_U$. Since $p_U^T(C_T) = C_U$ we then have maps $\pi_U^T : V_T \to V_U$ given by $\pi_U^T(f + C_T) = f|_U + C_U$. These in turn induce partial maps $\rho_U^T : PV_T \dashrightarrow PV_U$ given by $\rho_U^T([f + C_T]) = [\pi_U^T(f + C_T)]$ where $[x] := \mathrm{span}_{\mathbb{C}}\{x\}$. These maps are undefined for elements $[v]$ with $v \in \ker(\pi_U^T)$. We will use the partial arrow $\dashrightarrow$ to denote a map being partial.

LEMMA 2.0.4. *For every $U \subseteq V \subseteq T$ we have $\pi_U^T = \pi_U^V \circ \pi_V^T$ and $\rho_U^T = \rho_U^V \circ \rho_V^T$ whenever defined.* □

Given our previous definitions we are now in a position to define our main object of interest, the space $\overline{\mathcal{M}}_S$.

DEFINITION 2.0.5. $\overline{\mathcal{M}}_S$ is the subspace of $\prod_{\substack{T \subseteq S \\ |T| > 1}} PV_T$ defined as follows. An element $\underline{M}$ in $\prod_{\substack{T \subseteq S \\ |T| > 1}} PV_T$ lies in $\overline{\mathcal{M}}_S$ if and only if for every pair of sets $U \subseteq T \subseteq S$ we have $M_T \leq (\pi_U^T)^{-1} M_U$ where each $M_U$ is the $U$ component of $\underline{M}$.

In set theoretic notation, this is

$$\overline{\mathcal{M}}_S = \left\{ \underline{M} \in \prod_{\substack{T \subseteq S \\ |T| > 1}} PV_T \ \middle| \ M_T \leq (\pi_U^T)^{-1} M_U \text{ for all } U \subseteq T \subseteq S \right\}$$

EXAMPLE 2.0.6. For $n = 0, 1$ $\overline{\mathcal{M}}_S$ is empty. For $n = 2$, $\overline{\mathcal{M}}_S$ is a point and for $n = 3$ $\overline{\mathcal{M}}_S$ is a copy of $\mathbb{C}P^1$. For $n = 4$, $\overline{\mathcal{M}}_S$ is the blowup of $PV_S$ over the four special marked points $[0:0:0:1], [0:0:1:0], [0:1:0:0], [1:0:0:0]$

In order to understand the topology of the space $\overline{\mathcal{M}}_S$ it is natural to generalize the above construction to more general collections $\mathcal{L}$ of subsets of $S$. This we will do in the next definition. We shall then analyze these spaces to deduce information about our original space.

DEFINITION 2.0.7. For any collection $\mathcal{L}$ of subsets of $S$ with $|T| > 1$ for every $T \in \mathcal{L}$ we define,

$$\overline{\mathcal{M}}_\mathcal{L} = \left\{ \underline{M} \in \prod_{T \in \mathcal{L}} PV_T \ \middle| \ M_T \leq (\pi_U^T)^{-1} M_U \text{ for all } U, T \in \mathcal{L} \text{ with } U \subseteq T \right\}.$$

REMARK 2.0.8. It is clear that for such collections $\mathcal{L}$ of subsets of $S$, $\overline{\mathcal{M}}_\mathcal{L}$ is a closed projective variety. Throughout this thesis whenever we specify a collection $\mathcal{L}$ of subsets of $S$ we will implicitly assume that every element $T \in \mathcal{L}$ has size at least 2.

DEFINITION 2.0.9.

$$\begin{aligned} L_T &= \text{ tautological vector bundle over } PV_T \\ y_T &= e(L_T) \in H^2 PV_T \text{ the Euler class of } L_T \end{aligned}$$

LEMMA 2.0.10. *Let $V$ be a complex vector space of dimension d then $H^*(PV) = \mathbb{Z}[y]/y^d$ where $y = e(L)$ and $L$ is the tautological vector bundle over $PV$.* □

REMARK 2.0.11. Let $\mathcal{L}$ be a collection of subsets of $S$ with $|T| > 1$ for every $T \in \mathcal{L}$. We will sometimes describe the space $\overline{\mathcal{M}}_\mathcal{L}$ in a slightly different way using the partial maps $\rho_U^T : PV_T \dashrightarrow PV_U$ for any $U \subseteq T \subseteq S$ previously defined.

$$\overline{\mathcal{M}}_\mathcal{L} = \left\{ \underline{M} \in \prod_{T \in \mathcal{L}} PV_T \;\middle|\; \rho_U^T(M_T) = M_U \text{ whenever } \rho_U^T(M_T) \text{ is defined} \right\}$$

We have suppressed the notation $U, T \in \mathcal{L}$ with $U \subseteq T$. It is clear that both definitions are precisely the same.

DEFINITION 2.0.12. Define $\widetilde{U}_S = \{ f \in F(S, \mathbb{C}) \mid f \text{ is injective } \}$. This is naturally the configuration space $\mathbb{C}_0^{|S|}$ i.e. the complement of all diagonals in $\mathbb{C}^{|S|}$. We then define the space $U_S = p_S q_S(\widetilde{U}_S)$ where $q_S : F(S, \mathbb{C}) \to V_S$ and $p_S : V_S^\times \to PV_S$ are the quotient maps.

REMARK 2.0.13. The space $U_S \subseteq PV_S$ appears throughout this thesis. It is straightforward to check that $U_S$ is a Zariski open dense set. It is well known that, modulo the fact we have not defined them yet, the set of generic curves give a copy of $U_S$. We also note here that the restricted map $\pi_S : \pi_S^{-1}(U_S) \to U_S$ gives us an isomorphism where $\pi_S : \overline{\mathcal{M}}_S \to PV_S$. These facts will be proven later in chapter 6, section 2.

Here we state some conventions that we will be applying throughout this thesis. Let $R$ be a ring and $r \in R$ be a non zero element. In general whene we write $r = \sum_{i \in I} a_i r_i$ with $a_i \in \mathbb{Z}$ we suppose that $a_i r_i \neq 0$ and $r_i \neq r_j$ for all $i \neq j$. If $r = \prod_{i \in I} r_i^{n_i}$ we suppose $n_i > 0$. We may on occasion relax these conventions but this should be clear from the context.

# CHAPTER 3

# Combinatorial results

An important combinatorial invariant of an n-pointed stable curve of genus zero is that of the associated tree. Here we introduce a notion of trees, more generally forests that will appear throughout this paper in the structure of $\overline{\mathcal{M}}_S$ and its cohomology ring. In this chapter we gather some miscellaneous results about trees and other objects. Here we will develop some simple combinatorial properties of these objects. We will consider them in more detail later, as they are required. In particular we will show that the relationship between the classical associated combinatorial tree and our trees are equivalent. The sense in which they are equivalent will be explained later.

## 3.1. Forests and the length of a tree

DEFINITION 3.1.1. A *forest* on $S$ is a collection $\mathcal{F}$ of subsets of $S$ such that,

$$\text{for all} \quad T \in \mathcal{F} \quad |T| > 1,$$
$$\text{if} \quad U, T \in \mathcal{F} \text{ then either } U \cap T \text{ is empty or } U \subseteq T \text{ or } T \subseteq U.$$

An $S-tree$ is a forest $\mathcal{T}$ on $S$ with $S \in \mathcal{T}$. We say that a forest $\mathcal{F}$ is *proper* if it is not a tree.

REMARK 3.1.2. If a collection of subsets $\mathcal{F}$ of $S$ is a forest, then so is any subset of $\mathcal{F}$.

DEFINITION 3.1.3. Let $\mathcal{F}$ be a forest and $T$ any element of $\mathcal{F}$ then we define the T-tree $\mathcal{F}|_T$ by $\mathcal{F}|_T = \{\, U \in \mathcal{F} \mid U \subseteq T \,\}$ and call this set the *restriction* of $\mathcal{F}$ to $T$.

DEFINITION 3.1.4. Let $\mathcal{L}$ be a collection of subsets of $S$ then we define the sets $\mathbb{F}_\mathcal{L} = \{\, \mathcal{F} \subseteq \mathcal{L} \mid \mathcal{F} \text{ is a forest } \,\}$ and $\mathbb{T}_\mathcal{L} = \{\, \mathcal{T} \subseteq \mathcal{L} \mid \mathcal{T} \text{ is an S-tree } \,\}$, that is the set of forests respectively trees of $\mathcal{L}$.

We next introduce the concept of the depth of an S-tree $\mathcal{T}$. This will be analogous to the concept of the length of a rooted tree, which is defined as the length of longest path from the root to a terminal vertex without repetitions. It will turn out that this concept will provide us with a nice inductive way of studying the space $\overline{\mathcal{M}}_\mathcal{T}$, this is the subject of the next chapter.

DEFINITION 3.1.5. Let $\mathcal{T}$ be an S-tree. A *chain* of $\mathcal{T}$ is a subset $\mathcal{C}$ of $\mathcal{T}$ such that for all $U, V \in \mathcal{C}$ either $U \subseteq V$ or $V \subseteq U$. We write $\mathcal{C}_\mathcal{T}$ for the set of all chains on $\mathcal{T}$ and define the *depth* of $\mathcal{T}$ denoted $\mathrm{d}(\mathcal{T})$ by $\mathrm{d}(\mathcal{T}) = \max\{\,|\mathcal{C}|\,|\,\mathcal{C} \in \mathcal{C}_\mathcal{T}\,\}$.

LEMMA 3.1.6. *Let $\mathcal{T}$ be an S-tree and $T \subset S$ an element of $\mathcal{T}$ then $\mathcal{T}|_T$ is a T-tree and $\mathrm{d}(\mathcal{T}|_T) < \mathrm{d}(\mathcal{T})$.*

PROOF. The proof of this is clear.

□

## 3.2. Partitions of forests

The power set $P(S)$ of $S$ has its usual natural order given by inclusion. For any collection $\mathcal{L}$ of subsets of $S$ we can use the order induced on it to give a natural decomposition of $\mathcal{L}$ by sets $\mathcal{L}_T$, one for each $T \in \mathcal{L}$ given by, for any $U \in \mathcal{L}$, $U \in \mathcal{L}_T$ if and only if $U \subset T$ and is maximal in $T$. This decomposition, when applied to trees will actually turn out to give us a natural partition of it. We shall later see that the combinatorics of such a partition will turn up a lot. We next explain this notion more precisely.

DEFINITION 3.2.1. For any collection $\mathcal{L}$ of subsets of $S$ and for every $V \in \mathcal{L}$, we say $V$ is *maximal* in $\mathcal{L}$ if there is no $T \in \mathcal{L}$ such that $V \subset T$ and write $M(\mathcal{L})$ for the set of all maximal elements, that is $M(\mathcal{L}) = \{\,V \in \mathcal{L}\,|\,V \text{ is maximal in } \mathcal{L}\,\}$. We call this the set of *maximal elements* of $\mathcal{L}$. We say $V$ is *minimal* in $\mathcal{L}$ if there is no $T \in \mathcal{L}$ with $T \subset V$.

Let $U, T \in \mathcal{L}$, we say $U$ is *maximal in $\mathcal{L}$ under $T$* if $U \subset T$ and there is no $V \in \mathcal{L}$ with $U \subset V \subset T$. Let $T \in \mathcal{L}$ we write $M(\mathcal{L}, T) = \{\,U \in \mathcal{L}\,|\,U \text{ is maximal in } \mathcal{L} \text{ under } T\,\}$ and call this the set of *maximal objects of $\mathcal{L}$ under $T$*.

LEMMA 3.2.2. *For any collection of subsets $\mathcal{L}$ of $S$ we have $\mathcal{L} = M(\mathcal{L}) \cup \bigcup_{T \in \mathcal{L}} M(\mathcal{L}, T)$.*

PROOF. Clearly we have $M(\mathcal{L}) \cup \bigcup_{T \in \mathcal{L}} M(\mathcal{L}, T) \subseteq \mathcal{L}$. Let $U \in \mathcal{L}$, then either $U \in M(\mathcal{L})$ or $U \notin M(\mathcal{L})$. In the second case we need to show that $U \in M(\mathcal{L}, T)$ for some $T \in \mathcal{L}$. Suppose for a contradiction that $U \notin M(\mathcal{L}, T)$ for any $T \in \mathcal{L}$. Let $T \in \mathcal{L}$ be a set of minimal size containing $U$, there must exist such a set since $U \notin M(\mathcal{L})$. Then we know that $U \notin M(\mathcal{L}, T)$ and $U \subset T$ so there exists $V \in \mathcal{L}$ such that $U \subset V \subset T$, this contradicts the minimality of $T$. □

LEMMA 3.2.3. *Let $\mathcal{F}$ be a collection of subsets of $S$, then $\mathcal{F}$ is a forest if and only if for every $T \in \mathcal{F}$ the set $M(\mathcal{F}, T)$ is a collection of pairwise disjoint elements and $M(\mathcal{F})$ is also a collection of pairwise disjoint elements.*

PROOF. $\Rightarrow$ Suppose that for some $T \in \mathcal{F}$ there are distinct elements $U, V \in M(\mathcal{F}, T)$ with $U \cap V$ non-empty. Then since $\mathcal{F}$ is a forest, we have that either $U \subset V \subset T$ or $V \subset U \subset T$. This contradicts maximality in $\mathcal{F}$ under $T$.

If $S \in \mathcal{F}$ then $M(\mathcal{F}) = \{S\}$, so we may suppose $S \notin \mathcal{F}$. Then, we are required to show that $M(\mathcal{F})$ is disjoint. Suppose there are distinct elements $U, V \in M(\mathcal{F})$ with $U \cap V$ non empty. Then either $U \subset V \subset S$ or $V \subset U \subset S$, contradicting maximality in $\mathcal{T}$

$\Leftarrow$ We prove this inductively. Define

$$\begin{aligned}\mathcal{F}_1 &= M(\mathcal{F}) \\ \mathcal{F}_{i+1} &= \mathcal{F}_i \cup \bigcup_{T \in \mathcal{F}_i} M(\mathcal{F}, T) \\ &= \mathcal{F}_i \amalg \coprod_{T \in \mathcal{U}_i} M(\mathcal{F}, T)\end{aligned}$$

where $\mathcal{U}_i = \mathcal{F}_i \setminus \mathcal{F}_{i-1}$. Then it is clear $\mathcal{F}_1$ is a forest and given $\mathcal{F}_i$ is a forest then $\mathcal{F}_{i+1}$ is a forest. Clearly there is some $n$ such that $\mathcal{F} = \mathcal{F}_m$ for all $m \geq n$ and we are done. □

LEMMA 3.2.4. *For every forest $\mathcal{F}$ and for every distinct $U, V \in \mathcal{F}$ $M(\mathcal{F}, U) \cap M(\mathcal{F}, V)$ is empty and for every $T \in \mathcal{F}$ $M(\mathcal{F}) \cap M(\mathcal{F}, T)$ is empty, that is all the sets are disjoint.*

PROOF. Suppose for a contradiction there are elements $U, V \in \mathcal{F}$ with $M(\mathcal{F}, U) \cap M(\mathcal{F}, V)$ not empty. Let $W \in M(\mathcal{F}, U) \cap M(\mathcal{F}, V)$, then $W \subset U$ and $W \subset V$ therefore $U \cap V$ is not empty. Since $\mathcal{F}$ is a forest we must have either $U \subset V$ or $V \subset U$. If $U \subset V$ then $W \subset U \subset V$. If $V \subset U$ then $W \subset V \subset U$. In both cases we have a contradiction. A similar result holds for the second case.

□

COROLLARY 3.2.5. *For any forest $\mathcal{F}$ of $S$ we have* $\quad \mathcal{F} = M(\mathcal{F}) \amalg \coprod_{T \in \mathcal{F}} M(\mathcal{F}, T)$　　□

COROLLARY 3.2.6. *For any forest $\mathcal{F}$ of $S$ we have,*

$$\sum_{T \in \mathcal{F}} |M(\mathcal{F}, T)| = |\mathcal{F}| - |M(\mathcal{F})|$$

$$\sum_{T \in \mathcal{F}} \sum_{U \in M(\mathcal{F}, T)} |U| = \sum_{V \in \mathcal{F} \setminus M(\mathcal{F})} |V|$$

□

## 3.3. Numbers associated to forests

We next define some numbers associated to each forest $\mathcal{F}$ and develop some simple combinatorial properties of them. These number will appear throughout this document in both the algebraic and topological aspects of our project.

DEFINITION 3.3.1. For each forest $\mathcal{F}$ and for any $T \in \mathcal{F}$ we define

$$m(\mathcal{F}, T) = (|T| - 1) - \sum_{U \in M(\mathcal{F}, T)} (|U| - 1)$$

$$n(\mathcal{F}, T) = (|T| - 1) - \sum_{U \in M(\mathcal{F}, T)} (|U| - 2)$$

LEMMA 3.3.2. *For each forest $\mathcal{F}$ and for any $T \in \mathcal{F}$ we have,*

$$\sum_{T \in \mathcal{F}} m(\mathcal{F}, T) = \left( \sum_{T \in M(\mathcal{F})} |T| \right) - |M(\mathcal{F})|$$

$$\sum_{T \in \mathcal{F}} n(\mathcal{F}, T) = \left( \sum_{T \in M(\mathcal{F})} |T| \right) + |\mathcal{F}| - 2|M(\mathcal{F})|$$

PROOF.

$$\begin{aligned}
\sum_{T \in \mathcal{F}} m(\mathcal{F}, T) &= \sum_{T \in \mathcal{F}} (|T| - 1) - \sum_{T \in \mathcal{F}} \sum_{U \in M(\mathcal{F}, T)} (|U| - 1) \\
&= \left( \sum_{T \in \mathcal{F}} |T| \right) - |\mathcal{F}| - \sum_{T \in \mathcal{F}} \sum_{U \in M(\mathcal{F}, T)} |U| + \sum_{T \in \mathcal{F}} |M(\mathcal{F}, T)| \\
&= \left( \sum_{T \in M(\mathcal{F})} |T| \right) - |M(\mathcal{F})|
\end{aligned}$$

To prove the first case we used corollary 3.2.6. To prove the second case we note that $n(\mathcal{F}, T) = m(\mathcal{F}, T) + |M(\mathcal{F}, T)|$ and again use corollary 3.2.6 together with the first case. □

CONSTRUCTION 3.3.3. Let $T$ be a subset of $S$ and define $S/T$ to be the set of equivalence classes of $S$ under the equivalence relation $\sim$ given by $u \sim v \iff u = v$ or $u, v \in T$. We then define $q_T : S \to S/T$ to be the evident quotient map. Let $\mathcal{L}_+$ be a collection of subsets of $S$ and $T$ an element of $\mathcal{L}_+$ of minimal size. Then we write $\mathcal{L}$ to be the set $\mathcal{L}_+$ with $T$ removed and put $\overline{\mathcal{L}} = q_T(\mathcal{L})$, note that by the minimality of $T$ every element of $\overline{\mathcal{L}}$ has size at least two. It should also be noted that we do not display the dependance of $\mathcal{L}$ on $T$ as this avoids extra notation and should cause no confusion.

DEFINITION 3.3.4. Let $\mathcal{L}$ be a collection of subsets of $S$ and $T$ a subset of $S$ then we define $\text{Forests}(\mathcal{L}, T) = \{ \mathcal{F} \in \mathbb{F}_{\mathcal{L}} \mid T \in \mathcal{F} \}$. We also write $\text{Forests}(\mathcal{L}) = \mathbb{F}_{\mathcal{L}}$

LEMMA 3.3.5. *Let $S, T, \mathcal{L}_+, \mathcal{L}$ be as in the previous construction then there is a bijective correspondence $b : \text{Forests}(\mathcal{L}_+, T) \leftrightarrow \text{Forests}(\overline{\mathcal{L}})$ such that*

$$m(\mathcal{F}, U) = \begin{cases} m(b(\mathcal{F}), \overline{U}) & if \quad U \neq T \\ |T| - 1 & if \quad U = T \end{cases}$$

$$n(\mathcal{F}, U) = \begin{cases} n(b(\mathcal{F}), \overline{U}) & if \quad U \neq T, U \neq W \\ n(b(\mathcal{F}), \overline{U}) + 1 & if \quad U = W \\ |T| - 1 & if \quad U = T \end{cases}$$

where $W$ is the minimal element of $\mathcal{F}$ with $W \supset T$ should it exist.

PROOF. Define the function $b : \text{Forests}(\mathcal{L}, T) \to \text{Forests}(\overline{\mathcal{L}})$ by $b(\mathcal{T}) = q_T(\mathcal{T} \setminus \{T\})$ and $c : \text{Forests}(\overline{\mathcal{L}}) \to \text{Forests}(\mathcal{L}, T)$ by $c(\mathcal{U}) = q_T^{-1}(\mathcal{U}) \amalg \{T\}$. Then $b$ is a bijection with inverse $c$. To prove the final part of the claim we only need to check elements $U \in \mathcal{T}$ with $T \subseteq U$. First suppose $T \subset U$ and consider $M(\mathcal{T}, U)$, since $T \in \mathcal{T}$ and the size of $T$ is minimal there must be an element $V \in M(\mathcal{T}, U)$ with $T \subseteq V$. If $T \subset V$ then $\overline{V} \in \overline{\mathcal{L}}$ and for any other $W \in M(\mathcal{T}, U)$ $|W| = |\overline{W}|$ thus $(|\overline{U}| - 1) - (|\overline{V}| - 1) = (|U| - |T|) - (|V| - |T|) = (|U| - 1) - (|V| - 1)$ and $m(b(\mathcal{T}), \overline{U}) = m(\mathcal{T}, U)$. If $V = T$ then $|\overline{U}| - 1 = |U| - |T| = (|U| - 1) - (|T| - 1)$ and again $m(b(\mathcal{T}), \overline{U}) = m(\mathcal{T}, U)$. The case when $U = T$ is clear. The proof of the second statement is similar. □

DEFINITION 3.3.6. Let $\mathcal{F}$ be a forest and $\mathcal{U} \subseteq \mathcal{F}$ then we define $P_\mathcal{U} = \prod_{T \in \mathcal{U}} (m(\mathcal{U}, T) - 1)$, $P_\emptyset = 1$ and $m(\mathcal{F}) = \sum_{\mathcal{U} \subseteq \mathcal{F}} P_\mathcal{U}$. We also define $n(\mathcal{F}) = \prod_{T \in \mathcal{F}} n(\mathcal{F}, T)$.

LEMMA 3.3.7. Let $\mathcal{T}_+$ be a tree and $T \in \mathcal{T}_+$ be an element of minimal size. Put $\mathcal{T}$ to be the set $\mathcal{T}_+$ with $T$ removed and $\overline{\mathcal{T}} = q_T(\mathcal{T})$ where $q_T : S \to S/T$ is the collapsing map. Let $b : \text{Forests}(\mathcal{T}_+, T) \to \text{Forests}(\overline{\mathcal{T}})$ be the bijection of lemma 3.3.5 then for any forest $\mathcal{F} \in \text{Forests}(\mathcal{T}_+, T)$ we have $P_\mathcal{F} = (|T| - 2) P_{b(\mathcal{F})}$.

PROOF. Using the construction in lemma 3.3.5 we see that $b(\mathcal{F}) = q_T(\mathcal{U})$ where $\mathcal{U} = \mathcal{F} \setminus \{T\}$. Again by lemma 3.3.5 we see that for every $U \in \mathcal{U}$ that $m(\mathcal{F}, U) = m(b(\mathcal{F}), \overline{U})$. Since the evident induced map $q_T : \mathcal{U} \to b(\mathcal{F})$ is a bijection we see that

$$P_{b(F)} = \prod_{\overline{U} \in b(\mathcal{F})} (m(b(\mathcal{F}), \overline{U}) - 1)$$

$$= \prod_{U \in \mathcal{U}} (m(b(\mathcal{F}), \overline{U}) - 1)$$

$$= \prod_{U \in \mathcal{U}} (m(\mathcal{F}, U) - 1)$$

and so $P_{\mathcal{F}} = (|T| - 2)P_{b(\mathcal{F})}$.

$\square$

Next we prove that $m(\mathcal{T}) = n(\mathcal{T})$ for any tree $\mathcal{T}$, more generally it is true for any forest $\mathcal{F}$. We require this fact but prove it more efficiently later.

LEMMA 3.3.8. *Let $\mathcal{T}$ be an S-tree then $m(\mathcal{T}) = n(\mathcal{T})$.*

PROOF. To prove the claim we will show that both sides satisfy a certain recurrence relation with the same initial conditions. The particular relation is given by $a_{\mathcal{T}_+} = a_{\mathcal{T}} + (|T| - 2)a_{\overline{\mathcal{T}}}$, $a_{\emptyset} = 1$. First Let $W \in \mathcal{T}$ be the set of minimal size containing $T$ then

$$\begin{aligned}
a_{\mathcal{T}_+} &= a_{\mathcal{T}} + (|T| - 2)a_{\overline{\mathcal{T}}} \\
&= \prod_{U \in \mathcal{T}} n(\mathcal{T}, U) + (|T| - 2) \prod_{V \in \overline{\mathcal{T}}} n(\overline{\mathcal{T}}, V) \\
&= \Big[ \prod_{\substack{U \in \mathcal{T}_+ \\ U \neq W\ U \neq T}} n(\mathcal{T}_+, U)) \Big] n(\mathcal{T}, W) + (|T| - 2) \Big[ \prod_{\substack{U \in \mathcal{T}_+ \\ U \neq W\ U \neq T}} n(\mathcal{T}_+, U) \Big] (n(\mathcal{T}_+, W) - 1) \\
&= \Big[ \prod_{\substack{U \in \mathcal{T}_+ \\ U \neq W\ U \neq T}} n(\mathcal{T}_+, U) \Big] \Big[ n(\mathcal{T}_+, W) + (|T| - 2) + (|T| - 2)(n(\mathcal{T}_+, W) - 1) \Big] \\
&= \Big[ \prod_{\substack{U \in \mathcal{T}_+ \\ U \neq W\ U \neq T}} n(\mathcal{T}_+, U)) \Big] (|T| - 1)(n(\mathcal{T}_+, W)) \\
&= \prod_{U \in \mathcal{T}_+} n(\mathcal{T}_+, U)
\end{aligned}$$

for the second numbers we have

$$
\begin{aligned}
a_{\mathcal{T}_+} &= a_{\mathcal{T}} + (|T| - 2)a_{\overline{\mathcal{T}}} \\
&= \sum_{\mathcal{F} \subseteq \mathcal{T}} p_{\mathcal{F}} + (|T| - 2) \sum_{\overline{\mathcal{F}} \subseteq \overline{\mathcal{T}}} p_{\overline{\mathcal{F}}} \\
&= \sum_{\mathcal{F} \subseteq \mathcal{T}} p_{\mathcal{F}} + \sum_{\overline{\mathcal{F}} \subseteq \overline{\mathcal{T}}} (|T| - 2) p_{\overline{\mathcal{F}}} \\
&= \sum_{\substack{\mathcal{F} \subseteq \mathcal{T}_+ \\ T \notin \mathcal{F}}} p_{\mathcal{F}} + \sum_{\substack{\mathcal{F} \subseteq \mathcal{T}_+ \\ T \in \mathcal{F}}} p_{\mathcal{F}} \text{ by lemma 3.3.7} \\
&= \sum_{\mathcal{F} \subseteq \mathcal{T}_+} p_{\mathcal{F}} \\
&= m(\mathcal{T}_+)
\end{aligned}
$$

clearly these formula agree on the initial terms and we have proven our result. We will use this result to later see that the equality is also true for forests $\mathcal{F}$.

□

DEFINITION 3.3.9. For any finite set $S$ we define a function from the power set $P^2(S)$ of $P(S)$ to $P(S)$ called the *support* $\text{supp} : P^2(S) \to P(S)$ by $\text{supp}(\mathcal{L}) = \bigcup_{U \in \mathcal{L}} U$.

LEMMA 3.3.10. $\text{supp}(\mathcal{U} \cup \mathcal{V}) = \text{supp}(\mathcal{U}) \cup \text{supp}(\mathcal{V})$

PROOF. The proof is immediate from the definitions.

□

# CHAPTER 4

# The topology of $\overline{\mathcal{M}}_\mathcal{L}$ for trees

In this chapter we consider the topology of the space $\overline{\mathcal{M}}_\mathcal{L}$ for the case when $\mathcal{L}$ is an $S$-tree $\mathcal{T}$, this case will be the easiest to consider. We will show that $\overline{\mathcal{M}}_\mathcal{T}$ can be seen as the total space of an iterated projective bundle, we make this notion precise in what follows. This will enable us to show that $\overline{\mathcal{M}}_\mathcal{T}$ is a smooth irreducible projective variety and thus also a smooth complex manifold, its complex dimension is $|S| - 2$. This description will then make it easy to compute its cohomology ring. We consider this in the next chapter. We will use the embedding $i : \overline{\mathcal{M}}_\mathcal{T} \to \prod PV_T$ to identify the tangent space of $\overline{\mathcal{M}}_\mathcal{T}$ as a sub-bundle of the tangent space of $\prod PV_T$. In the last section of this chapter we will show that the space $\overline{\mathcal{M}}_\mathcal{T}$ can also be seen as an iterated blowup of $PV_S$ over certain linear subspaces. This result will then be used to prove a more general result that we consider later. The understanding of the spaces $\overline{\mathcal{M}}_\mathcal{T}$ will turn out to be both convenient and neat for the understanding of the spaces $\overline{\mathcal{M}}_\mathcal{L}$ for more general collections $\mathcal{L}$, enabling us to deduce the respective analogous statements about them. This will be explained later in the relevant chapters.

## 4.1. The projective bundle description for $\overline{\mathcal{M}}_\mathcal{T}$

Here we prove the projective bundle description for $\overline{\mathcal{M}}_\mathcal{T}$. Specifically if we remove $S$ from $\mathcal{T}$, we are left with a forest $\mathcal{F}$, which can be written as a disjoint union of trees, $\mathcal{F} = \coprod_{T \in M(\mathcal{F})} \mathcal{T}|_T$. It is clear from the definitions that $\overline{\mathcal{M}}_\mathcal{F} = \prod_{T \in M(\mathcal{F})} \overline{\mathcal{M}}_{\mathcal{T}|_T}$. There is an evident projection map $\pi : \overline{\mathcal{M}}_\mathcal{T} \to \overline{\mathcal{M}}_\mathcal{F}$ and we will show that this lifts to an isomorphism $\overline{\mathcal{M}}_\mathcal{T} = PW_\mathcal{T}$ for some vector bundle $W_\mathcal{T}$ over $\overline{\mathcal{M}}_\mathcal{F}$. This fact will then enable us to deduce a number of results about $\overline{\mathcal{M}}_\mathcal{T}$, inductively, on the depth of the tree $\mathcal{T}$, which we defined in the previous chapter. Later these will be used to obtain general result about $\overline{\mathcal{M}}_\mathcal{L}$. We will exhibit the vector bundle $W_\mathcal{T}$ over $\overline{\mathcal{M}}_\mathcal{F}$ directly as the pullback over an algebraic map of the tautological bundle of a certain grassmann space.

LEMMA 4.1.1. *Suppose $\mathcal{F}$ is a forest then $\overline{\mathcal{M}}_{\mathcal{F}} = \prod_{T \in M(\mathcal{F})} \overline{\mathcal{M}}_{\mathcal{F}|_T}$ and each $\mathcal{F}|_T$ is a tree.*

PROOF. It is clear that $\mathcal{F} = \coprod_{T \in M(\mathcal{F})} \mathcal{F}|_T$ and each $\mathcal{F}|_T$ is a $T$-tree. The proof is then clear from the definitions. □

REMARK 4.1.2. The last lemma tells us that to study $\overline{\mathcal{M}}_{\mathcal{F}}$ for forests $\mathcal{F}$ it is enough to understand the cases for trees $\mathcal{T}$. We next define a map of vector spaces whose importance will be seen in the following lemma. We will then use this map to prove our claim that $\overline{\mathcal{M}}_{\mathcal{T}}$ is the total space of a projective bundle over $\overline{\mathcal{M}}_{\mathcal{F}}$.

DEFINITION 4.1.3. Let $\mathcal{F}$ be a forest and for any $U, T \in \mathcal{F}$ with $U \subseteq T$ let $\pi_U^T : V_T \to V_U$ be the usual restriction map. Suppose $M(\mathcal{F}, T)$ is non-empty then define the map $\pi_{\mathcal{F},T} : V_T \to \bigoplus_{U \in M(\mathcal{F},T)} V_U$ by $\pi_{\mathcal{F},T} = \prod_{U \in M(\mathcal{F},T)} \pi_U^T$.

LEMMA 4.1.4. *Suppose $\mathcal{T}$ is an $S$-tree and $\mathrm{d}(\mathcal{T}) > 1$. Put $\mathcal{F}$ to be $\mathcal{T}$ with $S$ removed then,*

$$\overline{\mathcal{M}}_{\mathcal{T}} = \left\{ \underline{M} \in \overline{\mathcal{M}}_{\mathcal{F}} \times PV_S \,\Big|\, M_S \leq \pi_{\mathcal{T},S}^{-1}(\oplus M_U) \right\}$$

PROOF.
$$\begin{aligned}
\overline{\mathcal{M}}_{\mathcal{T}} &= \left\{ \underline{M} \in \prod_{T \in \mathcal{T}} PV_T \,\Big|\, M_T \leq (\pi_U^T)^{-1} M_U \text{ for all } U, T \in \mathcal{T} \text{ with } U \subseteq T \right\} \\
&= \left\{ \underline{M} \in \left( \prod_{T \in M(\mathcal{T},S)} \overline{\mathcal{M}}_{\mathcal{T}|_T} \right) \times PV_S \,\Big|\, M_S \leq \bigcap_{T \in M(\mathcal{T},S)} (\pi_T^S)^{-1}(M_T) \right\} \\
&= \left\{ \underline{M} \in \overline{\mathcal{M}}_{\mathcal{F}} \times PV_S \,\Big|\, M_S \leq \pi_{\mathcal{T},S}^{-1}(\oplus M_U) \right\}
\end{aligned}$$

The second step requires some comment. First the inclusion $\subseteq$ is clear. We now prove the reverse inclusion. Let $\underline{M}$ be an element of the second set. Then to prove $\underline{M} \in \overline{\mathcal{M}}_{\mathcal{T}}$ it is enough to prove that for every $U \in \mathcal{T}$ we have $M_S \leq (\pi_U^S)^{-1} M_U$. Let $U \in \mathcal{T}$, the claim is clear if $U = S$ so we suppose otherwise. Now we can factor $U \subset S$ with some $T \in M(\mathcal{T}, S)$. That is $U \subseteq T \subset S$ and $U, T \in \mathcal{T}|_T$. Thus we have $M_S \leq (\pi_T^S)^{-1} M_T$ and $M_T \leq (\pi_U^T)^{-1} M_U$ thus $M_S \leq (\pi_T^S)^{-1}(\pi_U^T)^{-1} M_U = (\pi_U^S)^{-1} M_U$. This proves the claim. □

Next, given $\underline{N} \in \overline{\mathcal{M}}_{\mathcal{F}}$ we put $W_{\mathcal{T},\underline{N}} = \pi_{\mathcal{T},S}^{-1}(\oplus N_U)$. If we can assemble these vector spaces into a vector bundle $W_{\mathcal{T}}$ over $\overline{\mathcal{M}}_{\mathcal{F}}$, the above lemma will show that $\overline{\mathcal{M}}_{\mathcal{T}} = PW_{\mathcal{T}}$ and thus $\overline{\mathcal{M}}_{\mathcal{T}}$ is a projective sub-bundle of $\overline{\mathcal{M}}_{\mathcal{F}} \times PV_S$. In what follows we will explain how to do this assembly.

LEMMA 4.1.5. *Let $\mathcal{F}$ be a forest and $T \in \mathcal{F}$ such that $M(\mathcal{F}, T)$ is non-empty. Then the map $\pi_{\mathcal{F},T} : V_T \to \bigoplus_{U \in M(\mathcal{F},T)} V_U$ is surjective, the dimension of the kernel is $m(\mathcal{F}, T)$ and $\dim(\pi_{\mathcal{F},T}^{-1}(\oplus M_U)) = n(\mathcal{F}, T)$. These numbers are defined in definition 3.3.1*

PROOF. For any $\overline{y} \in \bigoplus_{U \in M(\mathcal{F},T)} V_U$ and any $U \in M(\mathcal{F}, T)$ let $\overline{y}_U \in V_U$ be the $U$ component and $y_U : U \to \mathbb{C}$ be a representative for $\overline{y}_U$. Next put $W = T \setminus \bigcup_{U \in M(\mathcal{F},T)} U$ and define $x : T \to \mathbb{C}$ by

$$x(t) = \begin{cases} y_U(t) & \text{if } t \in U \\ 0 & \text{if } t \in W \end{cases}$$

This construction is well defined since the sets $U \in M(\mathcal{F}, T)$ are disjoint. Then $\overline{x} \in V_T$ and $\pi(\overline{x}) = \overline{y}$ thus $\pi_{\mathcal{F},T}$ is surjective. The dimension of the kernel is now clear using the rank nullity formula for linear transformations.

Next put $W_{\mathcal{F},\underline{M}} = \pi_{\mathcal{F},T}^{-1}(\oplus M_U)$ then the restricted map $\pi_{\mathcal{F},T} : W_{\mathcal{F},\underline{M}} \to \bigoplus_{U \in M(\mathcal{F},T)} M_U$ is surjective and clearly has the same kernel as the unrestricted map. We then apply the rank nullity formula again, the previous calculation tells us the dimension of the kernel. □

DEFINITION 4.1.6. Let $\mathcal{T}$ be an S-tree with $d(\mathcal{T}) > 1$. Then we have the surjective map $\pi_{\mathcal{T},S} : V_S \to \bigoplus_{T \in M(\mathcal{T},S)} V_T$. We define the induced map $\overline{\pi}_{\mathcal{T},S} : \prod PV_T \to \text{Grass}(V_S, k)$ by $\overline{\pi}_{\mathcal{T},S}(\underline{M}) = \pi_{\mathcal{T},S}^{-1}(\underline{M})$, where $k = n(\mathcal{T}, S)$. It is clear that this map is a morphism of varieties.

PROPOSITION 4.1.7. *Suppose $\mathcal{T}$ is an S-tree with $d(\mathcal{T}) > 1$ then there is a smooth algebraic vector sub-bundle $W_{\mathcal{T}} \subseteq \overline{\mathcal{M}}_{\mathcal{F}} \times V_S$ of dimension $n(\mathcal{T}, S)$ over $\overline{\mathcal{M}}_{\mathcal{F}}$ such that $\overline{\mathcal{M}}_{\mathcal{T}} = PW_{\mathcal{T}}$, the projective bundle of $W_{\mathcal{T}}$. If $d(\mathcal{T}) = 1$ then $\overline{\mathcal{M}}_{\mathcal{T}} = PV_S$.*

PROOF. Let $L$ be the $k$ dimensional tautological vector bundle over $\text{Grass}(V_S, k)$ where $k = n(\mathcal{T}, S)$. We have by the previous definition the map $\overline{\pi}_{\mathcal{T},S} : \prod PV_T \to \text{Grass}(V_S, k)$ and define the smooth vector bundle $V_{\mathcal{T}}$ over $\prod PV_T$ to be the pullback of $L$ over $\overline{\pi}_{\mathcal{T},S}$. Define $\pi_{\mathcal{T}} : \overline{\mathcal{M}}_{\mathcal{F}} \to \prod PV_T$ to be the projection where the product is taken over $M(\mathcal{F})$. Next define the vector bundle $W_{\mathcal{T}}$ over $\overline{\mathcal{M}}_{\mathcal{F}}$ to be the pullback of $V_{\mathcal{T}}$ over $\pi_{\mathcal{T}} : \overline{\mathcal{M}}_{\mathcal{F}} \to \prod PV_T$. Then it is clear that this is the required vector bundle. To see that $W_{\mathcal{T}}$ is smooth we use induction on the depth of the tree $\mathcal{T}$. For any S-tree $\mathcal{T}$ with $d(\mathcal{T}) = 2$ we have that $W_{\mathcal{T}} = V_{\mathcal{T}}$ which we all ready know is smooth. Suppose the result is true for any S-tree $\mathcal{T}$ with $d(\mathcal{T}) \leq n-1$. Let $\mathcal{T}$ be any S-tree with $d(\mathcal{T}) = n$ then we know $W_{\mathcal{T}}$ is the pullback of the smooth vector bundle $V_{\mathcal{T}}$ over the algebraic map $\pi_{\mathcal{T}} : \overline{\mathcal{M}}_{\mathcal{F}} \to \prod PV_T$. For each $T \in M(\mathcal{F})$ we have $\mathcal{T}|_T$ is a $T$-tree and $d(\mathcal{T}|_T) < d(\mathcal{T})$. If $d(\mathcal{T}|_T) > 1$ then by induction each $\overline{\mathcal{M}}_{\mathcal{T}|_T}$ is the projectivization of a smooth algebraic bundle and so is smooth. If $d(\mathcal{T}|_T) = 1$ then $\overline{\mathcal{M}}_{\mathcal{T}|_T} = PV_T$ which is smooth. Thus $\overline{\mathcal{M}}_{\mathcal{F}} = \prod \overline{\mathcal{M}}_{\mathcal{T}|_T}$ is smooth. Therefore we have that $W_{\mathcal{T}}$ is the pullback of an algebraic vector bundle with smooth base over an algebraic map and so is a smooth algebraic vector bundle.

□

DEFINITION 4.1.8. Let $\mathcal{F}$ be a forest and $\mathcal{U} \subseteq \mathcal{F}$ then we say that $\mathcal{U}$ is *closed downwards* if for every $U \in \mathcal{U}$ and $V \in \mathcal{F}$ with $V \subseteq U$ we have $V \in \mathcal{U}$.

COROLLARY 4.1.9. *Let $\mathcal{F}$ be a forest and $\mathcal{U} \subseteq \mathcal{F}$ be closed downwards. Then the projection $\pi : \overline{\mathcal{M}}_{\mathcal{F}} \to \overline{\mathcal{M}}_{\mathcal{U}}$ is a fibre bundle map.* □

DEFINITION 4.1.10. Let $\mathcal{T}$ be an $S$-tree and for each $T \in \mathcal{T}$ let $L_T$ be the tautological vector bundle over $PV_T$ and $\pi_T : \overline{\mathcal{M}}_{\mathcal{F}} \to PV_T$ be the projection map. Then we define the vector bundle $N_T$ to be the pullback of $L_T$ over $\pi_T$. Suppose $d(\mathcal{T}) > 1$ then let $V_{\mathcal{T}}$ be the $m(\mathcal{T}, S)$-dimensional trivial vector bundle over $\overline{\mathcal{M}}_{\mathcal{F}}$ whose fibre is $\ker(\pi_{\mathcal{T},S})$ where $\pi_{\mathcal{T},S} : V_S \to \bigoplus_{T \in M(\mathcal{T},S)} V_T$ is the usual map.

LEMMA 4.1.11. *Let $\mathcal{T}$ be an $S$-tree with $d(\mathcal{T}) > 1$ then we have the smooth isomorphism of vector bundles $W_{\mathcal{T}} \cong V_{\mathcal{T}} \oplus X_{\mathcal{T}}$ where $X_{\mathcal{T}} = \bigoplus_{T \in M(\mathcal{T},S)} N_T$.*

PROOF. For each $\underline{N} \in \overline{\mathcal{M}}_{\mathcal{F}}$ we have the short exact sequence on the fibres of the vector bundles

$$
\begin{array}{ccccc}
\ker(\pi_{\mathcal{T},S}) & \longrightarrow & \pi_{\mathcal{T},S}^{-1}(\oplus N_T) & \xrightarrow{\pi_{\mathcal{T},S}} & \oplus N_T \\
\downarrow & & \downarrow & & \downarrow \\
\ker(\pi_{\mathcal{T},S}) & \longrightarrow & V_S & \xrightarrow{\pi_{\mathcal{T},S}} & \oplus V_T
\end{array}
$$

This induces a short exact sequence of vector bundles $V_{\mathcal{T}} \rightarrowtail W_{\mathcal{T}} \twoheadrightarrow X_{\mathcal{T}}$ and a choice of smooth inner product on the vector space $V_S$ gives us the splitting $W_{\mathcal{T}} = V_{\mathcal{T}} \oplus V_{\mathcal{T}}^{\perp}$. We then use the map $\pi_{\mathcal{T},S}$ to give us the isomorphism of the perpendicular bundle $V_{\mathcal{T}}^{\perp}$ with $X_{\mathcal{T}}$. Thus we obtain the splitting $W_{\mathcal{T}} \cong V_{\mathcal{T}} \oplus X_{\mathcal{T}}$. □

## 4.2. The irreducibility and dimension of $\overline{\mathcal{M}}_{\mathcal{T}}$

Here we will prove that $\overline{\mathcal{M}}_{\mathcal{T}}$ is a smooth irreducible algebraic manifold of dimension $|S| - 2$. In particular this will mean that any non-empty Zariski open set is dense in the Zariski topology and moreover in the classical topology.

LEMMA 4.2.1. *Let $\mathcal{T}$ be an S-tree then $\overline{\mathcal{M}}_{\mathcal{T}}$ is a smooth irreducible projective variety of dimension $|S| - 2$.*

PROOF. Since $\overline{\mathcal{M}}_{\mathcal{T}} = P(W_{\mathcal{T}})$ and $W_{\mathcal{T}}$ is a smooth algebraic vector bundle it is clear that $\overline{\mathcal{M}}_{\mathcal{T}}$ is a smooth projective variety. To show $\overline{\mathcal{M}}_{\mathcal{T}}$ is irreducible and compute its dimension we use induction on the depth of the tree $\mathcal{T}$. The cases when $d(\mathcal{T}) = 1$ are well known. Assume the cases to be true for $d(\mathcal{T}) \leq n-1$, then for $d(\mathcal{T}) = n$ and for every $T \in M(\mathcal{F})$ each $\mathcal{T}|_T$ is a $T$-tree and $d(\mathcal{T}|_T) < d(\mathcal{T})$. We first show that $\overline{\mathcal{M}}_{\mathcal{T}}$ is irreducible. By induction each $\overline{\mathcal{M}}_{\mathcal{T}|_T}$ is irreducible, thus $\overline{\mathcal{M}}_{\mathcal{F}} = \prod_{T \in M(\mathcal{F})} \overline{\mathcal{M}}_{\mathcal{T}|_T}$ is irreducible, therefore $\overline{\mathcal{M}}_{\mathcal{T}}$ is the projectivization of a smooth algebraic vector bundle whose base $\overline{\mathcal{M}}_{\mathcal{F}}$ is irreducible. Thus $\overline{\mathcal{M}}_{\mathcal{T}}$ is irreducible. We next compute the dimension. By the above description of $\overline{\mathcal{M}}_{\mathcal{F}}$ we see inductively that the dimension of $\overline{\mathcal{M}}_{\mathcal{F}}$ is $\sum_{T \in M(\mathcal{F})} |T| - 2$. Now the dimension of the projective bundle follows from proposition 4.1.7 and is $n(\mathcal{T}, S) - 1$, thus we have

$$\begin{aligned}
\dim(\overline{\mathcal{M}}_\mathcal{T}) &= n(\mathcal{T},S) - 1 + \sum_{T \in M(\mathcal{F})} (|T| - 2) \\
&= |S| - 2 - \sum_{T \in M(\mathcal{T},S)} (|T| - 2) + \sum_{T \in M(\mathcal{F})} (|T| - 2) \\
&= |S| - 2 \quad \text{since } M(\mathcal{T},S) = M(\mathcal{F})
\end{aligned}$$

□

## 4.3. The tangent space to $\overline{\mathcal{M}}_\mathcal{T}$

We next compute the tangent bundle of the space $\overline{\mathcal{M}}_\mathcal{T}$. Using the embedding of $\overline{\mathcal{M}}_\mathcal{T}$ in $\prod PV_T$ we identify it as a sub-bundle of the tangent bundle $T \prod PV_T = \prod TPV_T$. This result will then enable us to calculate the tangent bundle for more general collections $\mathcal{L}$ of subsets of $S$, this we consider later. In particular we will be able to compute the tangent space of $\overline{\mathcal{M}}_S$. Before we can begin our computations we will need to introduce some definitions. These definitions will be generalized in later sections and are important throughout this thesis. Here we also prove special cases of more general results as they are particularly easy and serve as a warmup for later results.

DEFINITION 4.3.1. Let $\mathcal{T}$ be an $S$-tree and $\underline{M} \in \overline{\mathcal{M}}_\mathcal{T}$. We define the *type* of $\underline{M}$ denoted type($\underline{M}$) as,

$$\text{type}(\underline{M}) = \{\, U \in \mathcal{T} \mid \text{ for all } T \supset U \text{ with } T \in \mathcal{T} \implies \pi_U^T(M_T) = 0 \,\}$$

REMARK 4.3.2. It is clear that type($\underline{M}$) is an $S$-tree.

DEFINITION 4.3.3. Let $\mathcal{T}$ and $\mathcal{U} \subseteq \mathcal{T}$ be $S$-trees then we define a function root : $\mathcal{T} \to \mathcal{U}$ by root($U$) = $T$ where $T \in \mathcal{U}$ is the element of minimal size containing $U$. We say $T$ is the *root* of $U$ in $\mathcal{U}$.

REMARK 4.3.4. Let $\mathcal{T}$ and $\mathcal{U} \subseteq \mathcal{T}$ be $S$-trees. Let $U \in \mathcal{T}$ then root($U$) = $U$ in $\mathcal{U}$ if and only if $U \in \mathcal{U}$. Further if $U \subseteq V$ then root($U$) $\subseteq$ root($V$). On occasion we may extend this definition to elements $T \subseteq S$ other than those in $\mathcal{T}$.

LEMMA 4.3.5. *Let $\mathcal{T}$ be an S-tree and $\underline{M} \in \overline{\mathcal{M}}_\mathcal{T}$ an element of tree type $\mathcal{U}$. Let $U \in \mathcal{T}$ and put $T = \mathrm{root}(U)$ in $\mathcal{U}$ then $\pi_U^T M_T = M_U$.*

PROOF. Suppose for a contradiction $\pi_U^T M_T = 0$. Then clearly $U \notin \mathcal{U}$ for otherwise $T = \mathrm{root}(U) = U$ and $\pi_U^T M_T = M_T$. Let $V \in \mathcal{T}$ with $V \supset U$ then since $T$ is minimal $V \supseteq T \supset U$ and $\pi_U^V M_V = \pi_U^T \pi_T^V M_V \leq \pi_U^T M_T = 0$ thus $U \in \mathcal{U}$ a contradiction. $\square$

LEMMA 4.3.6. *Let $V$ be a finite dimensional complex vector space, $L$ the tautological vector bundle over $PV$ and $M \in PV$. Define $\mathrm{inc}(M, V)$ to be the subset of $\hom(M, V)$ consisting of the non zero functions. We then define the algebraic map $h : \mathrm{inc}(M, V) \to PV$ by $h(\alpha) = \alpha(M)$ and $q_M : V \to V/M$ to be the quotient map then*

(1) *The map $h$ is locally a product projection in the Zariski topology.*
(2) *The derivative $\mathrm{d}h : \hom(M, V) \to T_N PV$ is surjective at every point.*
(3) *The kernel of the derivative $\mathrm{d}h$ of $h$ is $\hom(M, N)$.*
(4) *$TPV \cong \hom(L, V/L)$*
(5) *The derivative $\mathrm{d}h : \hom(M, V) \to \hom(M, V/M)$ is $\mathrm{d}h(f) = q_M \circ f$*

PROOF. We first prove that the map $h$ is locally a product projection. The fibre of the map is $h^{-1}(N) = \mathrm{inc}(M, N)$. It is clear that we can form the smooth algebraic vector bundle $W$ over $PV$ with fibre $W_N = \hom(M, N)$. Write $\pi : W \to PV$ for the projection morphism and let $U$ be a Zariski open set containing $N$ so that $\pi^{-1}(U) \simeq U \times \mathbb{C}$ is trivial. Let $W^\times = \{(N, v) \in W \mid v \neq 0\}$ and write $\pi^\times : W^\times \to PV$ for the evident morphism. Restricting $\mathbb{C}$ to $\mathbb{C}^\times$ we obtain the smooth isomorphism $\pi^{\times-1}(U) \simeq U \times \mathbb{C}^\times$ with fibre at $N$ being $h^{-1}(N)$, this show us that $h$ is locally trivial. Since $h : \mathrm{inc}(M, V) \to PV$ is a product projection near $N \in PV$ with $h(g) = N$ where $g : M \to V$ we see that the differential $\mathrm{d}h$ at $g : M \to V$ is surjective. Thus we then have the surjective map $\mathrm{d}h : \hom(M, V) \to T_N PV$. We have $W_N = \hom(M, N)$ is a subspace of $\hom(M, V)$ and it is clear $h$ sends $W_N^\times = \mathrm{inc}(M, N)$ to a point and $g \in W_M^\times$ thus taking the derivative we see that the kernel of the map contains $W_N$, but the kernel of the map is one dimensional so $W_N = \ker(\mathrm{d}h)$ at $g$. When $N = M$ we have $\overline{\mathrm{d}h} : \hom(M, V)/W_M \to T_M PV$ is an isomorphism. The former space may be identified with $\hom(M, V/M)$, we can do this locally for each $M \in PV$. Let $U$ be an open set containing $M$. If we put $\ker(h_U)$ for the vector bundle corresponding to the kernel of the derivative we obtain the local isomorphisms $\overline{dh}_U : \hom(M, V)/\ker(h_U) \to T|_U(PV)$, we may then glue these to obtain the tangent bundle. Then the tangent space at each point $M \in PV$ is represented by $\hom(M, V/M)$. Now the spaces $\hom(M, V/M)$ fit together over $PV$ to form a smooth vector bundle and the evident quotient map into the tangent space is an isomorphism. One then sees that we may then take the derivative map $dh : \hom(M, V) \to \hom(M, V/M)$ by $dh(f) = q_M \circ f$. □

DEFINITION 4.3.7. Let $\mathcal{T}$ be an S-tree and $\underline{M}$ a point in $\overline{\mathcal{M}}_{\mathcal{T}}$ then for any $U$ and $T$ in $\mathcal{T}$ with $U \subseteq T$ we have the usual map $\pi_U^T : V_T \to V_U$. Then $\pi_U^T(M_T) \leq M_U$, so we have induced maps $\pi_U^T : M_T \to M_U$ and $\overline{\pi}_U^T : V_T/M_T \to V_U/M_U$. We then define a complex vector space $\sigma_{\underline{M}}$ by

$$\sigma_{\mathcal{T},\underline{M}} = \left\{ \underline{\alpha} \in \prod_{T \in \mathcal{T}} \hom(M_T, V_T/M_T) \,\bigg|\, \overline{\pi}_U^T \alpha_T = \alpha_U \pi_U^T \text{ for all } U \subseteq T \text{ with } U, T \in \mathcal{T} \right\}$$

We then define $\sigma_{\mathcal{T}} = \coprod_{\underline{M}} \sigma_{\mathcal{T},\underline{M}}$. Note that the condition in the braces is just the requirement that the following diagram is commutative for all $U, T \in \mathcal{T}$ with $U \subseteq T$.

$$\begin{array}{ccc} M_T & \xrightarrow{\alpha_T} & V_T/M_T \\ \pi_U^T \downarrow & & \downarrow \overline{\pi}_U^T \\ M_U & \xrightarrow{\alpha_U} & V_U/M_U \end{array}$$

DEFINITION 4.3.8. Let $\underline{M} \in \overline{\mathcal{M}}_{\mathcal{T}}$ with tree type $\mathcal{U}$ and $T \in \mathcal{U}$. Suppose $M(\mathcal{U}, T)$ is non-empty then let $\pi_{\mathcal{U},T} : V_T \to \bigoplus_{U \in M(\mathcal{U},T)} V_U$ be the usual map. Then we define the vector space $W_{\underline{M},T} = \pi_{\mathcal{U},T}^{-1}(\oplus M_U)$. If $M(\mathcal{U}, T)$ is empty we define $W_{\underline{M},T} = V_T$.

CONSTRUCTION 4.3.9. Let $\mathcal{T}$ be an $S$-tree and $\underline{M} \in \overline{\mathcal{M}}_{\mathcal{T}}$. Put $\mathcal{U} = \text{type}(\underline{M})$. Let $U \in \mathcal{T}$ and put $T = \text{root}(U)$ in $\mathcal{U}$ then we define the surjective restriction maps $\theta_U^T : \hom(M_T, V_T/M_T) \to \hom(M_U, V_U/M_U)$ as follows. Since $T = \text{root}(U)$ we have $\pi_U^T M_T = M_U$ and thus the restricted map $\pi_U^T : M_T \to M_U$ is an isomorphism. We also have the induced quotient maps $\overline{\pi}_U^T : V_T/M_T \to V_U/M_U$. Then given an element $\alpha_T \in \hom(M_T, V_T/M_T)$ we may define the element $\theta_U^T(\alpha_T) \in \hom(M_U, V_U/M_U)$ by $\theta_U^T(\alpha_T) = \overline{\pi}_U^T \alpha_T (\pi_U^T)^{-1}$. We then define a map $\theta : \prod_{T \in \mathcal{U}} \hom(M_T, V_T/M_T) \to \prod_{U \in \mathcal{T}} \hom(M_U, V_U/M_U)$ by $\theta = \prod_{U \in \mathcal{T}} \theta_U^{\text{root}(U)}$.

REMARK 4.3.10. We observe here that $\theta_U^T = \mathrm{d}\rho_U^T$ where $\rho_U^T : PV_T \setminus \ker(\pi_U^T) \to PV_U$ and $\pi_U^T : V_T \to V_U$ is the usual map. That is over the part where $\rho_U^T : PV_T \dashrightarrow PV_U$ is regular.

LEMMA 4.3.11. Let $\mathcal{T}$ be an $S$-tree and $\underline{M} \in \overline{\mathcal{M}}_{\mathcal{T}}$. Put $\mathcal{U} = \mathrm{type}(\underline{M})$ and write $\underline{N}$ for the image of $\underline{M}$ in $\overline{\mathcal{M}}_{\mathcal{U}}$ then

(1) The projection $\pi : \sigma_{\mathcal{T},\underline{M}} \to \sigma_{\mathcal{U},\underline{N}}$ is an isomorphism of vector spaces.
(2) $\sigma_{\mathcal{U},\underline{N}} = \prod_{U \in \mathcal{U}} \hom(M_U, W_{\underline{M},U}/M_U)$
(3) $\dim(\sigma_{\mathcal{U},\underline{N}}) = |S| - 2$
(4) The map $\theta$ restricts to $\theta : \sigma_{\mathcal{U},\underline{N}} \to \sigma_{\mathcal{T},\underline{M}}$ which is inverse to $\pi$

The proof of this is deferred to lemma 6.3.4 where we give the statement in greater generality.

LEMMA 4.3.12. Let $U, W \subseteq S$ with $U \subseteq W$ and put $\mathcal{T}$ to be the tree $\{U, W\}$. Then the tangent bundle to $\overline{\mathcal{M}}_{\mathcal{T}}$ is $\sigma_{\mathcal{T}}$.

PROOF. Let $\underline{M} \in \overline{\mathcal{M}}_{\mathcal{T}}$ and for each $T \in \mathcal{T}$ let $\mathrm{inc}(M_T, V_T)$ be the subset of $\hom(M_T, V_T)$ consisting of the non zero functions and define a map $h_T : \mathrm{inc}(M_T, V_T) \to PV_T$ by $h_T(\alpha_T) = \alpha_T(M_T)$. Define $Y_{\mathcal{T}}$ to be the variety $\hom(M_W, M_U) \times \prod \mathrm{inc}(M_T, V_T)$ and $X_{\mathcal{T}}$ to be the closed subvariety of $Y_{\mathcal{T}}$ given by $X_{\mathcal{T}} = \{\, (\gamma_U^W, \alpha_W, \alpha_U) \in Y_{\mathcal{T}} \mid \pi_U^W \alpha_W = \alpha_U \gamma_U^W \,\}$ that is the following diagram commutes.

$$\begin{array}{ccc} M_W & \xrightarrow{\alpha_W} & V_W \\ {\scriptstyle \gamma_U^W} \downarrow & & \downarrow {\scriptstyle \pi_U^W} \\ M_U & \xrightarrow{\alpha_U} & V_U \end{array}$$

Next consider the following commutative diagram where $q : X_{\mathcal{T}} \to \overline{\mathcal{M}}_{\mathcal{T}}$ is given by $q(\gamma, \underline{\alpha}) = (h_W(\alpha_W), h_U(\alpha_U))$ and $r : \text{inc}(M_W, V_W) \times \text{inc}(M_U, V_U) \to PV_W \times PV_U$ is given by $r = h_W \times h_U$, $p$ is the projection and $j$ the inclusion.

$$\begin{array}{ccc} X_{\mathcal{T}} & \xrightarrow{p} & \text{inc}(M_W, V_W) \times \text{inc}(M_U, V_U) \\ {\scriptstyle q}\downarrow & & \downarrow{\scriptstyle r} \\ \overline{\mathcal{M}}_{\mathcal{T}} & \xrightarrow{j} & PV_W \times PV_U \end{array}$$

I first claim that this algebraically is a pullback. Let $Z_{\mathcal{T}}$ be the pullback and $p : X_{\mathcal{T}} \to Z_{\mathcal{T}}$ be the comparison map. We produce an inverse $s : Z_{\mathcal{T}} \to X_{\mathcal{T}}$ as follows. Let $(\underline{N}, \alpha_W, \alpha_U)$ be a point of $Z_{\mathcal{T}}$ then by definition $\underline{N} = (\alpha_W(M_W), \alpha_U(M_U))$, so that the restricted maps with the same names onto their images give us the isomorphism of vector spaces $\alpha_W : M_W \to N_W$ and $\alpha_U : M_U \to N_U$. We also have $\pi_U^W N_W \leq N_U$. Now we need to construct a linear map $\gamma_U^W : M_W \to M_U$ so that $\pi_U^W \alpha_W = \alpha_U \gamma_U^W$. Then it is clear that $\gamma_U^W$ must be the unique map defined by $\gamma_U^W = \alpha_U^{-1} \pi_U^W \alpha_W$, so we define $s : Z_{\mathcal{T}} \to X_{\mathcal{T}}$ by $s(\underline{N}, \alpha_W, \alpha_U) = (\alpha_U^{-1} \pi_U^W \alpha_W, \alpha_W, \alpha_U)$. This gives us the desired inverse morphism. Next by lemma 4.3.6 $r$ is locally a product projection. Thus we see that the derivative $\mathrm{d}r$ is surjective at every point and since $\overline{\mathcal{M}}_{\mathcal{T}}$ is smooth we see that the pullback $Z_{\mathcal{T}}$ is a smooth algebraic variety thus $X_{\mathcal{T}}$ is a smooth algebraic subvariety of $Y_{\mathcal{T}}$.

Next as $r$ is a product projection and the diagram is a pullback we see that $q$ is a product projection. Thus the differential of $q$ is surjective at every point. In particular we have $(\pi, \text{inc}, \text{inc}) \in X_{\mathcal{T}}$ where $\pi = \pi_U^W$ and the differential is surjective at this point. To prove our claim we are required to show that $\mathrm{d}j(T_M \overline{\mathcal{M}}_{\mathcal{T}}) \subseteq \sigma_M$ or equivalently that $\mathrm{d}(jq)(T_{(\pi,\text{inc},\text{inc})} X_{\mathcal{T}}) \subseteq \sigma_{\underline{M}}$. Let $(x, y_W, y_U) \in T_{(\pi,\text{inc},\text{inc})} X_{\mathcal{T}}$ then we can choose a smooth path $p : \mathbb{R} \to X_{\mathcal{T}}$ given by $p(t) = (\gamma(t), \alpha_W(t), \alpha_U(t))$ such that $p(0) = (\pi, \text{inc}, \text{inc})$ and $\mathrm{d}p(0) = (x, y_W, y_U)$ where $\mathrm{d}p(t) = (\gamma'(t), \alpha'_W(t), \alpha'_U(t))$. Since $p(t) \in X_{\mathcal{T}}$ then we have $\pi_U^W \circ \alpha_W(t) = \alpha_U(t) \circ \gamma(t)$. We then take derivatives with respect to $t$ using a Leibniz argument and evaluate at $t = 0$ to see that $\pi_U^W \circ \alpha'_W(0) = \alpha'_U(0) \circ \gamma(0) + \alpha_U(0) \circ \gamma'(0)$ therefore $\pi_U^W \circ y_W = y_U \circ \pi_U^W + \text{inc} \circ x$. Now the differential of the map $h_T : \text{Inc}(M_T, V_T) \to PV_T$ at the inclusion $\text{inc} : M_T \to V_T$ is the map $\mathrm{d}h_T : \hom(M_T, V_T) \to \hom(M_T, V_T/M_T)$ induced from the quotient map $q_T : V_T \to V_T/M_T$ and we have $q_U \pi_U^W = \overline{\pi}_U^W q_W$ and $q_U \circ \text{inc} = 0$ thus composing with $q_U$ we have $q_U \pi_U^W \circ y_W = q_U y_U \circ \pi_U^W$ that is $\overline{\pi}_U^W \circ (q_W y_W) = (q_U y_U) \circ \pi_U^W$ ie we have the following commutative diagram

$$\begin{array}{ccccc} M_W & \xrightarrow{y_W} & V_W & \xrightarrow{q_W} & V_W/M_W \\ {\scriptstyle \pi_U^W} \downarrow & & & & \downarrow {\scriptstyle \overline{\pi}_U^W} \\ M_U & \xrightarrow{y_U} & V_U & \xrightarrow{q_U} & V_U/M_U \end{array}$$

if we put $\alpha_W = q_W y_W$ and $\alpha_U = q_U y_U$ then we have $\underline{\alpha} \in T_{\underline{M}} \overline{\mathcal{M}}_{\mathcal{T}}$ and $\overline{\pi}_U^W \circ \alpha_T = \alpha_U \circ \pi_U^W$ therefore $T_{\underline{M}} \overline{\mathcal{M}}_{\mathcal{T}} \leq \sigma_{\mathcal{T},\underline{M}}$ and a dimension check gives us an equality.

$\square$

COROLLARY 4.3.13. *Let $\mathcal{T}$ be an S-tree then the tangent bundle of $\overline{\mathcal{M}}_{\mathcal{T}}$ is $\sigma_{\mathcal{T}}$.*

PROOF. Let $\underline{M} \in \overline{\mathcal{M}}_{\mathcal{T}}$ and $U, T \in \mathcal{T}$ with $U \subseteq T$. Put $\mathcal{V} = \{U, T\}$ and consider the following diagram.

$$\begin{array}{ccc} \overline{\mathcal{M}}_{\mathcal{T}} & \xrightarrow{i} & \prod_{T \in \mathcal{T}} PV_T \\ {\scriptstyle \pi} \downarrow & & \downarrow {\scriptstyle \pi} \\ \overline{\mathcal{M}}_{\mathcal{V}} & \xrightarrow{j} & \prod_{T \in \mathcal{V}} PV_T \end{array}$$

where $i$ and $j$ are the inclusions. Let $\underline{\alpha} \in \mathrm{d}i(T_{\underline{M}} \overline{\mathcal{M}}_{\mathcal{T}})$ then by taking the differential of this diagram and using the previous lemma on the space $\overline{\mathcal{M}}_{\mathcal{V}}$ we have that $\underline{\alpha}|_{\mathcal{V}} \in \sigma_{\mathcal{V},\underline{N}}$ where $\underline{N} = \pi(\underline{M})$ for each such pair $U$ and $T$. Thus $\underline{\alpha} \in \sigma_{\mathcal{T},\underline{M}}$. Now the dimension of $\sigma_{\mathcal{T},\underline{M}}$ is $|S| - 2$ by lemma 4.3.11 thus we may take the tangent space $T_{\underline{M}} \overline{\mathcal{M}}_{\mathcal{T}}$ to be $\sigma_{\mathcal{T},\underline{M}}$. $\square$

## 4.4. The blowup description for $\overline{\mathcal{M}}_{\mathcal{T}}$

In this section we will show that for any S-tree $\mathcal{T}$ the space $\overline{\mathcal{M}}_{\mathcal{T}}$ is an iterated blowup of the projective space $PV_S$ taken over certain linear subspaces. We will use this to later prove a more general result about blowups over linear subspaces of $PV_S$ without much further effort. Let $\mathrm{inc} : \overline{\mathcal{M}}_{\mathcal{T}} \to \overline{\mathcal{M}}_{\mathcal{F}} \times PV_S$ be the inclusion, $\pi_1 : \overline{\mathcal{M}}_{\mathcal{F}} \times PV_S \to \overline{\mathcal{M}}_{\mathcal{F}}$ be the projection onto the first factor and $\pi_2 : \overline{\mathcal{M}}_\mathcal{F} \times PV_S \to PV_S$ the projection onto the second. Then section 4.1 told us the composition $\pi_1 \circ \mathrm{inc} : \overline{\mathcal{M}}_\mathcal{T} \to \overline{\mathcal{M}}_\mathcal{F}$ is a projective bundle. Here we prove that the composition $\pi_2 \circ \mathrm{inc} : \overline{\mathcal{M}}_\mathcal{T} \to PV_S$ is an iterated blowup. We reduce this calculation to the most basic case of a tree of depth 2. This case is easily computed.

Let $X$ be a variety, $Y$ a subvariety and write $\pi : \mathrm{Bl}_Y X \to X$ for the blowup of $X$ over $Y$. In particular for a blowup we have (1) the restricted map $\pi : \pi^{-1}(X \setminus Y) \to X \setminus Y$ is an isomorphism and (2) the fibre of the map over $Y$ is a projective space of dimension $\mathrm{codim}(X, Y) - 1$. In the case when $X$ and $Y$ are smooth $\pi^{-1}(Y)$ may be identified with the projective normal bundle of $Y$ in $X$. It is stated in [4, page 604 remark 4] that conditions (1) and (2) characterize the blowup in the smooth case. This would greatly simplify our situation. However no proof of this is supplied and I cannot prove it so we proceed differently. In our situation everything is smooth, this simplifies our analysis. In the category of smooth complex manifolds there is a local approach to the process of blowing up. We give a reference to [4] for more details on this approach.

LEMMA 4.4.1. *Let $i : W \to V$ be an injective map of vector spaces and $\bar{i} : PW \to PV$ be the induced embedding. Let $N \in PW$ , $M = \bar{i}(N)$ and write $i : N \to M$ for restricted map $i$ which is then an isomorphism , $j : W/N \to V/M$ for the induced quotient map. Then the derivative $\mathrm{d}\bar{i} : \hom(N, W/N) \to \hom(M, V/M)$ is given by $\mathrm{d}\bar{i}(f) = jfi^{-1}$ and the image of $\mathrm{d}\bar{i}$ is $\hom(M, \overline{W}/M)$ where $\overline{W} = i(W)$. Further let $L_W$ be the tautological vector bundle over $PW$ then $N(\bar{i}) = \hom(L_W, V/\overline{W})$ where $N(\bar{i})$ is the normal bundle.*

PROOF. Consider the following commutative diagram,

$$\begin{array}{ccc} PW & \xrightarrow{\bar{i}} & PV \\ h \uparrow & & \uparrow h \\ \mathrm{inc}(N, W) & \xrightarrow{k} & \mathrm{inc}(M, V) \end{array}$$

where $k : \mathrm{inc}(N, W) \to \mathrm{inc}(M, V)$ is given by $k(f) = ifi^{-1}$. Then $k(\mathrm{inc}(N, N)) = \mathrm{inc}(M, M)$ and taking derivatives we obtain the desired result. The results on the image and normal bundle are then clear. □

In particular given $T \subseteq S$ we define an equivalence relation on $S$ by $u \sim v$ if and only if $u = v$ or $u, v \in T$ and write $S/T$ for the equivalence classes. We will then define a natural embedding of $PV_{S/T}$ into $PV_S$. This embedding, which we define next, will appear throughout the thesis.

CONSTRUCTION 4.4.2. Let $T \subseteq S$ and $q_T : S \to S/T$ be the usual quotient map. For each $U \subseteq S$ we write $\overline{U} = q_T(U)$ and $q_U : U \to \overline{U}$ for the restricted quotient map. Then we define maps a $r_U : PV_{\overline{U}} \to PV_U$ as follows. We first define a map $r_U : F(\overline{U}, \mathbb{C}) \to F(U, \mathbb{C})$. For any $f_{\overline{U}} \in F(\overline{U}, \mathbb{C})$ we define $f_U \in F(U, \mathbb{C})$ by $f_U = f_{\overline{U}} \circ q_U$. Now $r_U$ sends constants to constants so we can define the induced injective map $r_U : V_{\overline{U}} \to V_U$ and in turn an injective map $r_U : PV_{\overline{U}} \to PV_U$, this is the first of our desired maps. One readily verifies that we have a short exact sequence $V_{\overline{U}} \to V_U \to V_{T \cap U}$ and we define $\overline{V}_U = \ker(\pi^U_{T \cap U})$. We then define the isomorphisms $s_U : \overline{V}_U \to V_{\overline{U}}$ by the short exact sequence and $s_U : P\overline{V}_U \to PV_{\overline{U}}$ to be the induced map. This is the second of our desired maps.

LEMMA 4.4.3. Let $T \subseteq S$ and $q_T : S \to S/T$ the quotient map. Then with the notation as above we have maps $r_U : PV_{\overline{U}} \to PV_U$ and $s_U : P\overline{V}_U \to PV_{\overline{U}}$. We also write $r_U : PV_{\overline{U}} \to P\overline{V}_U$ to be the restricted map of $r_U$ onto its image. Then $r_U$ and $s_U$ are inverses for each other. Let $U \subseteq W \subseteq S$ then we have the following commutative diagrams,

$$\begin{array}{ccc} V_{\overline{W}} \xrightarrow{r_W} \overline{V}_W & & \overline{V}_W \xrightarrow{s_W} V_{\overline{W}} \\ \pi^{\overline{W}}_{\overline{U}} \downarrow \quad \downarrow \pi^W_U & & \pi^W_U \downarrow \quad \downarrow \pi^{\overline{W}}_{\overline{U}} \\ V_{\overline{U}} \xrightarrow{r_U} \overline{V}_U & & \overline{V}_U \xrightarrow{s_U} V_{\overline{U}} \end{array}$$

In particular if $|\overline{U}| = 1$ then we obtain $\pi^W_U r_W = 0$

DEFINITION 4.4.4. We say that $(\mathcal{T}, S, T)$ is *compatible* if $\mathcal{T}$ is an $S$-tree, $T \subseteq S$ and for every $U \in \mathcal{T}$ with $T \cap U$ non-empty we have $T \subset U$. Then $\mathcal{T} \amalg \{T\}$ is an $S$-tree. We then write $\mathcal{T}_+$ for the tree $\mathcal{T} \amalg \{T\}$ and $T$ is a minimal element of $\mathcal{T}_+$.

CONSTRUCTION 4.4.5. Let $(\mathcal{T}, S, T)$ compatible and write $\overline{\mathcal{T}} = q_T(\mathcal{T})$ where $q_T : S \to S/T$ is the collapsing map. Then we define a map $j_{\mathcal{T}} : \overline{\mathcal{M}}_{\overline{\mathcal{T}}} \to \overline{\mathcal{M}}_{\mathcal{T}}$ as follows. Given $\underline{N} \in \overline{\mathcal{M}}_{\overline{\mathcal{T}}}$ and $U \in \mathcal{T}$ we define $M_U = r_U(N_{\overline{U}})$. This clearly gives an injective map $j_{\mathcal{T}} : \overline{\mathcal{M}}_{\overline{\mathcal{T}}} \to \prod_{U \in \mathcal{T}} PV_U$. We must show that its image lies in $\overline{\mathcal{M}}_{\mathcal{T}}$. Given $U, W \in \mathcal{T}$ with $U \subseteq W$ we have $\overline{U} \subseteq \overline{W}$ and it is then clear from the first commutative diagram of the last lemma that we have a map $j_{\mathcal{T}} : \overline{\mathcal{M}}_{\overline{\mathcal{T}}} \to \overline{\mathcal{M}}_{\mathcal{T}}$ and we have the following lemma.

LEMMA 4.4.6. *Let $(\mathcal{T}, S, T)$ be compatible then the map $j_{\mathcal{T}} : \overline{\mathcal{M}}_{\overline{\mathcal{T}}} \to \overline{\mathcal{M}}_{\mathcal{T}}$ is injective and its image is the closed variety*

$$Z(\mathcal{T}, T) = \{\, \underline{M} \in \overline{\mathcal{M}}_{\mathcal{T}} \mid \text{ for all } U \in \mathcal{T} \text{ with } U \supset T \implies \pi_T^U M_U = 0 \,\}$$

PROOF. The inclusion image($j_{\mathcal{T}}$) $\subseteq Z(\mathcal{T}, T)$ is immediate since for every $U \supset T$ we have the short exact sequence $V_{U/T} \to V_U \to V_T$. Let $\underline{M} \in Z(\mathcal{T}, T)$ then by definition $\pi_T^U M_U = 0$ for every $U \supset T$ with $U \in \mathcal{T}$ thus $M_U \in P\overline{V}_U$. We have the map $s_U : P\overline{V}_U \to PV_{U/T}$ and write $N_{\overline{U}}$ for the image of $M_U$ then by the second commutative diagram of the last lemma we immediately see that $\underline{N} \in \overline{\mathcal{M}}_{\overline{\mathcal{T}}}$ and $j_{\mathcal{T}}(\underline{N}) = \underline{M}$ thus the image is as claimed. It is clearly Zariski closed from the definition. □

CONSTRUCTION 4.4.7. Let $(\mathcal{T}, S, T)$ be compatible. We then define a map $i_{\mathcal{T}} : PV_T \times \overline{\mathcal{M}}_{\overline{\mathcal{T}}} \to \overline{\mathcal{M}}_{\mathcal{T}_+}$ by $i_{\mathcal{T}}(M_T, \underline{M}) = (M_T, j_{\mathcal{T}}(\underline{M}))$. Because $T$ has minimal size we see by the last lemma that this does indeed lie in $\overline{\mathcal{M}}_{\mathcal{T}_+}$

LEMMA 4.4.8. *Let $(\mathcal{T}, S, T)$ be compatible then the following diagram is a pullback.*

$$\begin{array}{ccc} \overline{\mathcal{M}}_{\overline{\mathcal{T}}} \times PV_T & \xrightarrow{i_{\mathcal{T}}} & \overline{\mathcal{M}}_{\mathcal{T}_+} \\ {\scriptstyle p}\downarrow & & \downarrow{\scriptstyle \pi} \\ \overline{\mathcal{M}}_{\overline{\mathcal{T}}} & \xrightarrow{j_{\mathcal{T}}} & \overline{\mathcal{M}}_{\mathcal{T}} \end{array}$$

□

LEMMA 4.4.9. *Let $(\mathcal{T}, S, T)$ be compatible and $U = \text{root}(T)$ in $\mathcal{T}$. Let $\pi_U : \overline{\mathcal{M}}_\mathcal{T} \to PV_U$ be the projection map and $r_U : PV_{U/T} \to PV_U$ be the usual embedding then $\pi_U \pitchfork r_U$.*

PROOF. We begin by working out the images of the respective tangent maps. Then we proceed to work out their sum. In this proof we will be making extensive use of lemma 4.1.5 to work out the dimensions of the various vector spaces we consider. This should be clear from the context and will not be explicitly stated when used.

Let $\underline{M} \in \overline{\mathcal{M}}_\mathcal{T}$ and $N \in PV_{U/T}$ such that $\pi_U(\underline{M}) = r_U(N)$. Then it is clear by lemma 4.4.1 that the image of the tangent space of $k$ is $\hom(M_U, \mathbf{W}/M_U)$ where $\mathbf{W} = \ker(\pi_T^U)$ and $\pi_T^U : V_U \to V_T$ is the usual restriction map. Next put $\mathcal{U} = \text{type}(\underline{M})$ and $W = \text{root}(U)$ in $\mathcal{U}$. To compute the image of the tangent space of $\pi_U$ at $\underline{M}$ we need to compute the image of the projection map $\pi : \sigma_{\mathcal{T},\underline{M}} \to \hom(M_U, V_U/M_U)$. By lemma 4.3.11 we have the isomorphism $\theta : \sigma_{\mathcal{U},\underline{N}} \to \sigma_{\mathcal{T},\underline{M}}$ where $\underline{N}$ is the image of $\underline{M}$ in $\overline{\mathcal{M}}_\mathcal{U}$ and $\sigma_{\mathcal{U},\underline{N}} = \prod_{U \in \mathcal{U}} \hom(M_U, W_{\underline{M},U}/M_U)$. Then we see that we must compute the image of $\hom(M_W, W_{\underline{M},W}/M_W)$ under the map $\theta_U^W : \hom(M_W, V_W/M_W) \to \hom(M_U, V_U/M_U)$. Equivalently let $\pi_U^W : V_W \to V_U$ be the usual projection map then the image of $\hom(M_W, W_{\underline{M},W}/M_W)$ under $\theta_U^W$ is the space $\hom(M_U, \pi_U^W(W_{\underline{M},W})/M_U)$. Define $A(U) = \{ V \in M(\mathcal{U}, W) \mid V \subseteq U \}$. If $A(U)$ is empty the whole claim is trivial since $\pi_U^W(W_{\underline{M},W}) = V_U$ so we suppose otherwise. Then we have the following commutative diagram where all the maps are surjective.

$$\begin{array}{ccc} V_W & \xrightarrow{\pi_U^W} & V_U \\ \pi_{\mathcal{U},W} \downarrow & & \downarrow t_U \\ \bigoplus_{X \in M(\mathcal{U},W)} V_X & \xrightarrow{p} & \bigoplus_{Y \in A(U)} V_Y \end{array}$$

Put $\mathbf{Z} = t_U^{-1}(\oplus M_Y)$ then we prove that the image of the tangent space of $\pi_U$ at $\underline{M}$ is $\hom(M_U, \mathbf{Z}/M_U)$. The diagram induces the map $\gamma_U^W : \pi_{\mathcal{U},W}^{-1}(\oplus M_X) \to \pi_U^{-1}(\oplus M_Y)$ where $W_{\underline{M},W} = \pi_{\mathcal{U},W}^{-1}(\oplus M_X)$, $\gamma_U^W$ is the restricted map of $\pi_U^W$ and we are using the same limits as above. We need to prove that this is surjective. Put $B(U)$ to be the complement of the set $A(U)$ in $M(\mathcal{U}, W)$ and write $q : V_W \to V_U \oplus \bigoplus_{X \in B(U)} V_X$ for the evident map. This map is surjective as $Y \cap U$ is empty for each $Y \in B(U)$. One then sees that the kernel of $\gamma_U^W$ is $q^{-1}(0_U \oplus \bigoplus M_Y)$. Given this we have the following dimensions,

$$\dim(\ker(\gamma_U^W)) = |B(U)| + |W| - 1 - \sum_{Y \in B(U)} (|Y| - 1) - (|U| - 1)$$

$$\dim(\mathbf{Z}) = |U| - 1 - \sum_{X \in A(U)} (|X| - 2)$$

$$\dim(\pi_{\mathcal{U},W}^{-1}(\oplus M_X)) = |W| - 1 - \sum_{X \in M(\mathcal{U},W)} (|X| - 2)$$

One then readily checks using the observation $A(U) \amalg B(U) = M(\mathcal{U}, W)$ that we have

$$\dim(\pi_{\mathcal{U},W}^{-1}(\oplus M_X)) = \dim(\ker(\gamma_U^W)) + \dim(\mathbf{Z})$$

This then shows us that the image of the projection map $\pi : \sigma_{T,\underline{M}} \to \hom(M_U, V_U/M_U)$ is $\hom(M_U, \mathbf{Z}/M_U)$.

For a transverse intersection we must show that $\hom(M_U, \mathbf{Z}/M_U) + \hom(M_U, \mathbf{W}/M_U) = \hom(M_U, V_T/M_U)$ or equivalently $\mathbf{Z} + \mathbf{W} = V_U$. Now let $p : V_U \to V_T \oplus \bigoplus_{Y \in A(U)} V_Y$ be the evident map then this map is surjective because $Y \cap T$ is empty for every $Y \in A(U)$ since $U$ is the minimal element containing $T$ and we know the dimension of the kernel. Moreover the space $\mathbf{Z} \cap \mathbf{W}$ is $p^{-1}(0_T \oplus \bigoplus M_X)$ thus we have the following dimensions,

$$\dim(\mathbf{Z}) = |U| - 1 - \sum_{X \in A(U)} |X| - 2$$

$$\dim(\mathbf{W}) = |U| - |T|$$

$$\dim(\mathbf{Z} \cap \mathbf{W}) = |A(U)| + |U| - 1 - \sum_{X \in A(U)} (|X| - 1) - (|T| - 1)$$

Thus we obtain

$$\dim(\mathbf{Z} + \mathbf{W}) = \dim(\mathbf{Z}) + \dim(\mathbf{W}) - \dim(\mathbf{Z} \cap \mathbf{W}) = |U| - 1$$

and we are done. □

LEMMA 4.4.10. *Let $X$ be a smooth variety and $Y, Z$ smooth subvarieties. Suppose $Y$ and $Z$ intersect transversely and put $W = Y \cap Z$. Write $\pi : \mathrm{Bl}_Y X \to X$ for the blowup of $X$ along $Y$ then the restricted map $\pi : \pi^{-1}(Z) \to Z$ is the blowup of $Z$ along $W$.*

For a proof see for example [**13**, page 74]

LEMMA 4.4.11. *Let $\pi : \mathrm{Bl}_Y X \to X$ be the blowup of $X$ along $Y$ and $Z$ be any variety then $\pi \times \mathrm{id} : \mathrm{Bl}_Y X \times Z \to X \times Z$ is the blowup of $X \times Z$ along $Y \times Z$.* □

LEMMA 4.4.12. *Let $V$ be a complex vector space, $W$ a subspace of $V$ and $\pi : V \to V/W$ the quotient map then the blowup of $PV$ along $PW$ is the projection $p : B(V, W) \to PV$ where*

$$B(V, W) = \{\, (L, M) \in PV \times P(V/W) \mid L \leq \pi^{-1}(M) \,\}$$

PROOF. We have $PW = \{\, M \in PV \mid \pi(M) = 0 \,\}$ thus

$$\begin{aligned}\mathrm{Bl}_{PW} PV &= \mathrm{cl}\{(M, N) \in PV \times P(V/W) \mid M \in PV \setminus PW, \quad N = \pi(M)\} \\ &\leq B(V, W)\end{aligned}$$

One easily sees this is the closure as $B(V, W)$ is irreducible and we are done. □

COROLLARY 4.4.13. *Let $S$ be a finite set and $T \subseteq S$. Let $r_S : PV_{S/T} \to PV_S$ be the embedding of construction 4.4.2 then the blowup of $PV_S$ over $PV_{S/T}$ is $\pi : \overline{\mathcal{M}}_{\{S,T\}} \to PV_S$.*

PROOF. Consider the short exact sequence $V_{S/T} \to V_S \to V_T$. The last lemma tells us that the blowup of $PV_S$ along $PV_{S/T}$ is $B(V_S, V_{S/T})$. Write $\pi : V_S \to V_S/V_{S/T}$ for the quotient map and $\overline{\pi}_T^S : V_S/V_{S/T} \to V_T$ for the induced isomorphism where $\pi_T^S : V_S \to V_T$ it the standard map. Then we have $\pi_T^S = \overline{\pi}_T^S \pi$. Consider the isomorphism $\theta : PV_S \times P(V_S/V_{S/T}) \to PV_S \times PV_T$ induced by $\overline{\pi}_T^S$ which is the identity on the first factor. Then $\theta$ restricts to an isomorphism $\theta : B(V_S, V_{S/T}) \to \overline{\mathcal{M}}_{\{S,T\}}$ that commutes with projections to $PV_S$. Thus the claim is proven. □

LEMMA 4.4.14. *Let $(\mathcal{T}, S, T)$ be compatible then the blowup of $\overline{\mathcal{M}}_\mathcal{T}$ along $\overline{\mathcal{M}}_{\overline{\mathcal{T}}}$ is given by $\pi : \overline{\mathcal{M}}_{\mathcal{T}_+} \to \overline{\mathcal{M}}_\mathcal{T}$.*

PROOF. Let $U = \mathrm{root}(T)$ in $\mathcal{T}$ and put $\mathcal{V}$ to be the set $\mathcal{T}$ with $U$ removed. Consider the following commutative diagram where $j_\mathcal{T}$ is the map of construction 4.4.5 and $k \circ i$ is the obvious factoring.

$$\begin{array}{ccc}
\overline{\mathcal{M}}_{\mathcal{T}_+} & \xrightarrow{\mathrm{inc}} & \overline{\mathcal{M}}_{\{U,T\}} \times \overline{\mathcal{M}}_\mathcal{V} \\
\pi \downarrow & & \downarrow l \\
\overline{\mathcal{M}}_\mathcal{T} & \xrightarrow{\mathrm{inc}} & PV_U \times \overline{\mathcal{M}}_\mathcal{V} \\
j_\mathcal{T} \uparrow & & \uparrow k \\
\overline{\mathcal{M}}_{\overline{\mathcal{T}}} & \xrightarrow{i} & PV_{U/T} \times \overline{\mathcal{M}}_\mathcal{V}
\end{array}$$

I first claim that $\mathrm{image}(k) \cap \overline{\mathcal{M}}_\mathcal{T} = \mathrm{image}(j_\mathcal{T})$. Let $\underline{M} \in \mathrm{image}(k) \cap \overline{\mathcal{M}}_\mathcal{T}$ then $\pi_T^U(M_U) = 0$. Now let $V$ be any element containing $T$ then since $\mathcal{T}$ is a tree we have $V \supseteq U$ and $\pi_U^V(M_V) \leq M_U$ and therefore $\pi_T^V M_V = 0$ thus by lemma 4.4.6 we have $\underline{M} = j(\underline{N})$ for some unique $\underline{N} \in \mathcal{M}_{\overline{\mathcal{T}}}$. The other inclusion is automatically true. Then using lemmas 4.4.13 and 4.4.11 we see that the blowup of $PV_U \times \overline{\mathcal{M}}_\mathcal{V}$ along $PV_{U/T} \times \overline{\mathcal{M}}_\mathcal{V}$ is given by the morphism $l$. It is easy to see by lemma 4.4.9 that $k \pitchfork \mathrm{inc}$ then using lemma 4.4.10 we see that the blowup of $\overline{\mathcal{M}}_\mathcal{T}$ along $\overline{\mathcal{M}}_{\overline{\mathcal{T}}}$ is the restricted morphism $l : l^{-1}(\overline{\mathcal{M}}_\mathcal{T}) \to \overline{\mathcal{M}}_\mathcal{T}$. We are required to prove that $\overline{\mathcal{M}}_{\mathcal{T}_+} = l^{-1}(\overline{\mathcal{M}}_\mathcal{T})$. The inclusion $\subseteq$ is clear by the commutative diagram. To see the reverse inclusion let $\underline{M} \in l^{-1}(\overline{\mathcal{M}}_\mathcal{T})$ then for any $V \supseteq T$ with $V \in \mathcal{T}$ it will be enough to prove that $\pi_T^V M_V \leq M_T$. Because $U$ is the smallest element of $\mathcal{T}$ containing $T$ we have $T \subset U \subseteq V$. By construction we have $\pi_T^U M_U \leq M_T$ and $\pi_U^V M_V \leq M_U$ thus $\pi_T^V M_V = \pi_T^U \pi_U^V M_V \leq \pi_T^U M_U \leq M_T$. This proves $\underline{M} \in \overline{\mathcal{M}}_{\mathcal{T}_+}$ and we have the equality. Then clearly we have $\pi : \overline{\mathcal{M}}_{\mathcal{T}_+} \to \overline{\mathcal{M}}_\mathcal{T}$ is the desired blowup, as claimed. □

COROLLARY 4.4.15. *Let $(\mathcal{T}, S, T)$ be compatible and $j_\mathcal{T} : \overline{\mathcal{M}}_{\overline{\mathcal{T}}} \to \overline{\mathcal{M}}_\mathcal{T}$ be the usual embedding. Then the normal bundle is $\hom(N_{U/T}, V_T)$ where $U = \mathrm{root}(T)$ in $\mathcal{T}$. The projectivization of the normal bundle is the trivial bundle $\overline{\mathcal{M}}_{\overline{\mathcal{T}}} \times PV_T$.*

PROOF. Consider the following commutative diagram where $\mathcal{V}$ is $\mathcal{T}$ with $U$ removed and all the maps are embeddings

$$\begin{array}{ccc} \overline{\mathcal{M}}_\mathcal{T} & \longrightarrow & PV_U \times \overline{\mathcal{M}}_\mathcal{V} \\ j_\mathcal{T} \uparrow & & \uparrow k \\ \overline{\mathcal{M}}_{\overline{\mathcal{T}}} & \longrightarrow & PV_{U/T} \times \overline{\mathcal{M}}_\mathcal{V} \end{array}$$

Then we have an induced map $N(j_\mathcal{T}) \to N(k)$ and a transverse intersection tells us that the map is surjective on fibres. The dimension of the fibres are both $|T|-1$ thus the map is an isomorphism on fibres. Clearly $N(k) \cong N(r_U)$ where $r_U : PV_{U/T} \to PV_U$ and we have the short exact sequence $V_{U/T} \to V_U \to V_T$. Thus by lemma 4.4.1 $N(k) \cong \hom(L_{U/T}, V_T)$. Let $N_{U/T}$ be the pullback of the tautological bundle $L_{U/T}$ over the projection $\pi : \overline{\mathcal{M}}_{\overline{\mathcal{T}}} \to PV_{U/T}$. Then we see that the normal bundle is $\hom(N_{U/T}, V_T)$ where we write $V_T$ for the trivial bundle with the same fibre. This can also be written as $N^*_{U/T} \otimes V_T$ where the star denotes the dual bundle. This gives the projectivization as $\overline{\mathcal{M}}_{\overline{\mathcal{T}}} \times PV_T$ the expected answer. □

DEFINITION 4.4.16. Let $U_S$ be the usual open set of $PV_S$ and $\pi : \overline{\mathcal{M}}_\mathcal{T} \to PV_S$ be the projection map. Then we have proven that $\pi$ is a composition of blowups and is therefore surjective. We then define the non-empty Zariski open set $\mathcal{M}_\mathcal{T} = \pi^{-1}(U_S)$. Then $\mathcal{M}_\mathcal{T} \subseteq \overline{\mathcal{M}}_\mathcal{T}$ is dense in the classical topology since $\overline{\mathcal{M}}_\mathcal{T}$ is smooth and irreducible.

CONSTRUCTION 4.4.17. Given an $S$-tree $\mathcal{T}$ we fix an order on it as follows. We define $T_1 = S$ and $\mathcal{T}_1 = \{S\}$. Suppose $T_1, ..., T_r$ and $\mathcal{T}_1, ..., \mathcal{T}_r$ are defined then we define $T_{r+1} \in \mathcal{T} \setminus \mathcal{T}_r$ to be an element of maximal size and $\mathcal{T}_{r+1} = \mathcal{T}_r \amalg \{T_{r+1}\}$. Then $\mathcal{T}_i = \{ T_j \mid j \leq i \}$ and define $S_1^U = P(\ker(\pi_U^S))$, $B_1 = PV_S$ and $B_{i+1} = \mathrm{Bl}_{S_i^T} B_i$ where $T = T_{i+1}$ and write $\pi_{i+1} : \overline{\mathcal{M}}_{\mathcal{T}_{i+1}} \to \overline{\mathcal{M}}_{\mathcal{T}_i}$ for the projection map. Then for any $U \in \mathcal{T}$ we define $S_{i+1}^U = $ strict transform of $S_i^U$ in $B_i$.

LEMMA 4.4.18. *Let $\mathcal{T}$ be an S-tree and $\mathcal{T}_i$ be an order as above then $B_i = \overline{\mathcal{M}}_{\mathcal{T}_i}$ and given $U \in \mathcal{T}$ with $U = T_k$ we have $S_i^U = j_U(\overline{\mathcal{M}}_{\overline{\mathcal{T}}_i})$ for every $i < k$ and is zero otherwise where we quotient in $q_U : S \to S/U$ and write $j_U = j_{\overline{\mathcal{L}}_i}$ in particular $B_m = \overline{\mathcal{M}}_{\mathcal{T}}$ where $m = |\mathcal{T}|$.*

PROOF. We prove the claim by an induction on $i$. The case when $i = 1$ is clear from the definitions. Put $T = T_{i+1}$ then inductively $B_i = \overline{\mathcal{M}}_{\mathcal{T}_i}$ and $S_i^T = j_T(\overline{\mathcal{M}}_{\overline{\mathcal{T}}_i})$. Then by lemma 4.4.14 we see that the blowup of $B_i$ along $S_i^T$ is $B_{i+1} = \overline{\mathcal{M}}_{\mathcal{T}_{i+1}}$. We are left to prove the claim for the strict transforms. Choose $i + 1 < k$ and put $U = T_k$ then define $W_i^U = S_i^U \setminus S_i^T$ where inductively $S_i^U = j_U(\overline{\mathcal{M}}_{\overline{\mathcal{T}}_i})$ then we need to calculate $\mathrm{cl}(\pi_{i+1}^{-1}(W_i^U))$. We first prove that

$$j_U(\mathcal{M}_{\overline{\mathcal{T}}_{i+1}}) \subseteq \pi_{i+1}^{-1}(W_i^U) \subseteq j_U(\overline{\mathcal{M}}_{\overline{\mathcal{T}}_{i+1}})$$

We prove the left hand inclusion first. Let $\underline{N} \in \mathcal{M}_{\overline{\mathcal{T}}_{i+1}}$ and put $\underline{M} = j_U(\underline{N})$. Then by lemma 4.4.6 we see that for every $V \in \mathcal{T}_{i+1}$ with $V \supset U$ we have $\pi_U^V M_V = 0$. Put $\underline{L}$ to be the image of $\underline{M}$ under $\pi_{i+1}$. Then again by lemma 4.4.6 we see that $\underline{L} \in S_i^U$. Because $T \not\subseteq U$ by choice of order we see $\pi_T^S M_S \neq 0$ so that $\underline{L} \notin S_i^T$ thus $\underline{M} \in \pi_{i+1}^{-1}(W_i^U)$ and the left hand inclusion is true.

Next let $\underline{M} \in \pi_{i+1}^{-1}(W_i^U)$ and put $\underline{L}$ to be its image under $\pi_{i+1}$. Consider $W_i^U \subseteq \overline{\mathcal{M}}_{\mathcal{T}_i}$, because $\underline{L} \notin S_i^T$ by lemma 4.4.6 we can find a $V \in \mathcal{T}_i$ containing $T$ such that $\pi_T^V M_V = M_T$. Because $\underline{L} \in S_i^U$ we deduce that $\pi_U^W M_W = 0$ for all $W \in \mathcal{T}_i$ with $W \supset U$. In particular this then tells us that $\pi_U^T M_T = 0$ if $T \supset U$. This proves the right hand inclusion.

Then taking the closure of both sides and recalling that $j_U$ is a proper map and $\mathcal{M}_{\overline{\mathcal{T}}_{i+1}}$ is dense in $\overline{\mathcal{M}}_{\overline{\mathcal{T}}_{i+1}}$ we see that $S_{i+1}^U = \mathrm{cl}(\pi_{i+1}^{-1}(W^T)) = j_U(\overline{\mathcal{M}}_{\overline{\mathcal{T}}_{i+1}})$. This completes the induction.

□

# CHAPTER 5

# The cohomology of $\overline{\mathcal{M}}_{\mathcal{L}}$ for forests

In this chapter we compute the cohomology of the space $\overline{\mathcal{M}}_{\mathcal{F}}$ for any forest $\mathcal{F}$. We begin by defining a ring $R_{\mathcal{F}}$ and showing that it has a natural decomposition into a tensor product of rings $R_{\mathcal{F}} \cong \bigotimes_{T \in M(\mathcal{F})} R_{\mathcal{F}|_T}$ where each $\mathcal{F}|_T$ is a $T$-tree. For the case of an $S$-tree $\mathcal{T}$ the result of proposition 4.1.7 represents $\overline{\mathcal{M}}_{\mathcal{T}}$ as the total space of a tower of projective bundles. Effectively, this will enable us to compute its cohomology ring and prove it is $R_{\mathcal{T}}$ for the tree case. We then prove more generally that the cohomology of $\overline{\mathcal{M}}_{\mathcal{F}}$ is $R_{\mathcal{F}}$. This is just the Künneth theorem. This will help us to understand the cohomology ring of $\overline{\mathcal{M}}_S$, which is the main subject of chapter 9. In particular we will show that the cohomology ring of $\overline{\mathcal{M}}_{\mathcal{F}}$ is finitely generated by its elements of degree two and free as an abelian group, there is a natural choice for these generators, and the rank of $R_{\mathcal{F}}^2$ is the number of elements in $\mathcal{F}$ whose sizes $|T|$ are at least three. We also compute the rank of this ring. This will enable us to produce two different descriptions of a basis for it. The first description will be the natural one associated to the cohomology of a projective bundle. The second basis will provide some insight into the structure of the cohomology ring of $\overline{\mathcal{M}}_S$ and will be formulated in the language of trees. Later we will see that the second description extends naturally to rings of greater generality. We will also specify a necessary and sufficient condition for a monomial $x$ in $R_{\mathcal{F}}$ to be zero, this result will be useful, and developed further when we later examine the zeros of more general rings and in particular the cohomology ring of $\overline{\mathcal{M}}_S$. We recall here that we will be using the conventions in chapter 2 regarding rings.

## 5.1. The cohomology ring $R_{\mathcal{F}}$

In this section we will compute the ordinary integral cohomology ring of the space $\overline{\mathcal{M}}_{\mathcal{F}}$, the proof presented here relies on the previous chapter for the description of $\overline{\mathcal{M}}_{\mathcal{T}}$ as the total space of an iterated projective bundle. We then proceed inductively (on the depth of the tree $\mathcal{T}$) to compute its cohomology ring using a general result that relates the cohomology of a projective bundle to that of it base space and the Chern polynomial of the underlying vector bundle, a notion we define next as this is non standard terminology.

DEFINITION 5.1.1. Let $\pi : V \to X$ be an $n$ dimensional complex vector bundle over $X$ and $c_m(V) \in H^{2m}(X)$ be the $m^{\text{th}}$ Chern class of $V$ then we define the *Chern polynomial* $f_V(t) \in H^*(X)[t]$ by $f_V(t) = \sum_{i=0}^{n} c_i(V) t^{n-i}$

We next remind the reader of some standard results about Chern polynomials, these follow immediately from the properties of Chern classes. We then state a result that relates the cohomology of $X$ with the cohomology of $P(V)$, the projective bundle of $V$ over $X$ together with the Chern polynomial of $V$. For the proofs see any standard texts on characteristics classes [**11**] for example.

PROPOSITION 5.1.2. *Let $V$ and $W$ be vector bundles over a base space $X$ then the Chern polynomial of the vector bundle $U = V \oplus W$ is $f_{V \oplus W} = f_V f_W$.*

□

PROPOSITION 5.1.3. *Let $V$ be a trivial complex vector bundle over $X$ of dimension $d$ then $f_V(t) = t^d$.*

□

PROPOSITION 5.1.4. *Let $V$ be a complex line bundle over $X$ then the Chern polynomial of $V$ is $f_V(t) = t - e(V)$ where $e(V)$ is the euler class of $V$.*

□

THEOREM 5.1.5. *let $V$ be an $n$ dimensional vector bundle over $X$ and $f_V(t) \in H^*(X)[t]$ be the associated Chern polynomial then $H^*(PV) = H^*(X)[t]/f_V(t)$ where the identification is natural. Moreover $\{\, 1, t, ..., t^{n-1} \,\}$ form a basis making $H^*(PV)$ free of rank $n$ over $H^*(X)$*

□

CONSTRUCTION 5.1.6. Let $\mathcal{T}$ be an S-tree with $\mathrm{d}(\mathcal{T}) > 1$ and for every $T \in \mathcal{T}$ define $L_T$ to be the tautological line bundle over $PV_T$ and $y_T = e(L_T)$, the Euler class of $L_T$. Let $\mathcal{F}$ be the set $\mathcal{T}$ with $S$ removed. For every $T \in M(\mathcal{F})$ let $\pi_T : \overline{\mathcal{M}}_\mathcal{F} \to PV_T$ be the projection map and put $x_T = \pi_T^*(y_T)$. We define the polynomial $f_\mathcal{T} \in H^*(\overline{\mathcal{M}}_\mathcal{F})[t]$ by

$$f_\mathcal{T}(t) = t^{m(\mathcal{T},S)} \prod_{T \in M(\mathcal{T},S)} (t - x_T).$$

LEMMA 5.1.7. *Let $\mathcal{T}$ be an $S$-tree with $\mathrm{d}(\mathcal{T}) > 1$ then $f_\mathcal{T}$ is the Chern polynomial of the vector bundle $W_\mathcal{T}$ where $W_\mathcal{T}$ is the vector bundle of Proposition 4.1.7.*

PROOF. Lemma 4.1.11 gives us the splitting $W_\mathcal{T} \cong V_\mathcal{T} \oplus X_\mathcal{T}$ where $X_\mathcal{T} = \bigoplus_{T \in M(\mathcal{T},S)} N_T$. Then using the previously discussed results on Chern polynomials we have

$$\begin{aligned} f_W &= f_{V_\mathcal{T}} \prod_{T \in M(\mathcal{T},S)} f_{N_T} \\ &= t^{m(\mathcal{T},S)} \prod_{T \in M(\mathcal{T},S)} (t - e(N_T)) \\ &= t^{m(\mathcal{T},S)} \prod_{T \in M(\mathcal{T},S)} (t - x_T) \end{aligned}$$

□

DEFINITION 5.1.8. For any forest $\mathcal{F}$ of $S$ let $I_\mathcal{F}$ be the ideal in the polynomial ring $\mathbb{Z}_\mathcal{F} = \mathbb{Z}[\, y_T \mid T \in \mathcal{F}\, ]$ generated by elements of the form $y_T^{m(\mathcal{F},T)} \prod_{U \in M(\mathcal{F},T)} (y_T - y_U)$ where $T$ is an element of $\mathcal{F}$. Note that when $M(\mathcal{F},T)$ is empty i.e. $T$ is a minimal in $\mathcal{F}$ we take by convention the relation to be $y_T^{m(\mathcal{F},T)} = y_T^{|T|-2}$. We then define $S_\mathcal{F} = \mathbb{Z}_\mathcal{F}/I_\mathcal{F}$.

LEMMA 5.1.9. *Let $\mathcal{F}$ be a forest then $S_\mathcal{F} \cong \bigotimes_{T \in M(\mathcal{F})} S_{\mathcal{F}|T}$* □

PROPOSITION 5.1.10. *Let $\mathcal{F}$ be a forest then $S_\mathcal{F}$ is the cohomology ring of $\overline{\mathcal{M}}_\mathcal{F}$ and is free as an abelian group with finite rank. The identification sends $\overline{y}_T$ to $x_T = \pi_T^*(z_T)$ where $\pi_T : \overline{\mathcal{M}}_\mathcal{F} \to PV_T$ and $z_T$ is the standard generator of $H^*(PV_T)$.*

PROOF. We prove this for the case of an $S$-tree first. We proceed by an induction argument on the depth of the tree $\mathcal{T}$. If $\mathrm{d}(\mathcal{T}) = 1$ then we have that $\overline{\mathcal{M}}_\mathcal{T} = PV_S$ thus the claim is clear in this case. Assume the result is true for a tree of depth $\mathrm{d}(\mathcal{T}) \leq n-1$ then consider any tree of depth $\mathrm{d}(\mathcal{T}) = n$. We have by proposition 4.1.7 that $\pi : \overline{\mathcal{M}}_\mathcal{T} \to \overline{\mathcal{M}}_\mathcal{F}$ is a projective bundle where $\mathcal{F}$ is the set $\mathcal{T}$ with $S$ removed and we now apply theorem

5.1.5 to compute its cohomology. Then,

$$
\begin{aligned}
H^*(\overline{\mathcal{M}}_{\mathcal{T}}) &\cong \frac{H^*(\prod \overline{\mathcal{M}}_{\mathcal{T}|_T})[t]}{f_{\mathcal{T}}(t)} \text{ by proposition 5.1.5} \\
&\cong \frac{\left[\bigotimes H^*(\overline{\mathcal{M}}_{\mathcal{T}|_T})\right][t]}{f_{\mathcal{T}}(t)} \text{ inductively by the Künneth theorem} \\
&\cong \frac{\left[\bigotimes S_{\mathcal{T}|_T}\right][t]}{f_{\mathcal{T}}(t)} \text{ by induction} \\
&\cong \frac{S_{\mathcal{F}}[t]}{f_{\mathcal{T}}(t)} \text{ by lemma 5.1.9} \\
&\cong \frac{\mathbb{Z}_{\mathcal{T}}}{I_{\mathcal{T}}} \\
&= S_{\mathcal{T}}
\end{aligned}
$$

Thus we see $S_{\mathcal{T}}$ is finitely generated and free as an abelian group and $S_{\mathcal{T}}$ is the cohomology ring of $\overline{\mathcal{M}}_{\mathcal{T}}$. Note that by a slight abuse of notation we have not identified the various $f_{\mathcal{T}}(t)$ under each isomorphism to avoid extra notation. This should not lead to any confusion. The arbitrary case of $\overline{\mathcal{M}}_{\mathcal{F}}$ then follows from the description given in lemma 4.1.1 the Künneth theorem and the previous tree case. □

COROLLARY 5.1.11. *Let $\mathcal{T}$ be an S-tree with $\mathrm{d}(\mathcal{T}) > 1$ then $S_{\mathcal{T}} \cong \frac{[\bigotimes S_{\mathcal{T}|_T}][t]}{f_{\mathcal{T}}(t)}$* □

Having proven $S_{\mathcal{F}}$ is the cohomology ring of $\overline{\mathcal{M}}_{\mathcal{F}}$ we now proceed to describe a different ring $R_{\mathcal{F}}$ that is naturally isomorphic to $S_{\mathcal{F}}$, in fact they are equal but this is not immediate from the definitions. This ring will be in a more natural form to enable us to compare it with the cohomology ring of $\overline{\mathcal{M}}_S$ which we consider in chapter 9.

DEFINITION 5.1.12. For any forest $\mathcal{F}$ of $S$ let $J_{\mathcal{F}}$ be the ideal in the polynomial ring $\mathbb{Z}_{\mathcal{F}} = \mathbb{Z}[\, y_T \mid T \in \mathcal{F}\,]$ generated by elements of the form $y_T^{m(\mathcal{T},T)} \prod_{U \in M(\mathcal{T},T)} (y_T - y_U)$ where $T \in \mathcal{F}$ and $\mathcal{T} \subseteq \mathcal{F}$ is a $T$-tree of depth 2 or less. Note that in the case that $\mathcal{T}$ has depth one we have the special relation $y_T^{m(\mathcal{T},T)} = y_T^{|T|-2}$ corresponding to the fact that $M(\mathcal{T},T)$ is empty. We define $R_{\mathcal{F}} = \mathbb{Z}_{\mathcal{F}}/J_{\mathcal{F}}$ and write $B_{\mathcal{F}}$ for the standard monomial basis of $\mathbb{Z}_{\mathcal{F}}$.

REMARK 5.1.13. More explicitly let $T, U_1, ..., U_r \in \mathcal{F}$ with the $U_i \subseteq T$ and disjoint then we have the relation $y_T^m \prod_i (y_T - y_{U_i}) = 0$ where $m = |T| - 1 - \sum_i (|U_i| - 1)$. The ideal $I_\mathcal{F}$ relates to the case when $\{U_1, ..., U_r\} = M(\mathcal{F}, T)$.

REMARK 5.1.14. It is immediate that we have $I_\mathcal{F} \subseteq J_\mathcal{F}$ and that for every $\mathcal{U} \subseteq \mathcal{F}$ we have the following commutative diagram where $\pi : \overline{\mathcal{M}}_\mathcal{F} \to \overline{\mathcal{M}}_\mathcal{U}$ and $s$ is induced from the inclusion.

$$\begin{array}{ccc} H^*(\overline{\mathcal{M}}_\mathcal{U}) & \xrightarrow{\pi^*} & H^*(\overline{\mathcal{M}}_\mathcal{F}) \\ \downarrow & & \downarrow \\ S_\mathcal{U} & \xrightarrow{s} & S_\mathcal{F} \end{array}$$

We will next show that this inclusion of ideals is actually an equality so that $R_\mathcal{F} = S_\mathcal{F}$. Topologically from the above diagram this will be easy to prove, since we may put $\mathcal{U}$ to be the $T$-tree associated to the relation we wish to prove. Then the equivalent relation will hold in $H^*(\overline{\mathcal{M}}_\mathcal{U})$, we obtain the result by applying $\pi^*$. However it would be better if we could understand this ring from a purely algebraic viewpoint.

LEMMA 5.1.15. *Let $\mathcal{F}$ be a forest then $R_\mathcal{F} = S_\mathcal{F}$*

PROOF. We prove this for $S$-trees first. We use the description of $S_\mathcal{T}$ as $\dfrac{[\bigotimes S_{\mathcal{T}|_T}][t]}{f_\mathcal{T}(t)}$ where the product is taken over $M(\mathcal{T}, S)$ to prove inductively on the depth of the tree $\mathcal{T}$ that $S_\mathcal{T}$ and $R_\mathcal{T}$ are equal. When $d(\mathcal{T}) = 1$ this case is clear. Assume the result to be true for any $S$-tree with $d(\mathcal{T}) \leq n - 1$ then we have $S_\mathcal{T} \cong \dfrac{[\bigotimes R_{\mathcal{T}|_T}][t]}{f_\mathcal{T}(t)}$ and so for any $T$-tree $\mathcal{U}$ of depth 2 or less where $T \neq S$ we obtain the desired relation. We also have the relation $y_S^{m(\mathcal{T}, S)} \prod_{T \in M(\mathcal{T}, S)} (y_S - y_T)$. Now let $\mathcal{V}$ be an $S$-tree contained in $M(\mathcal{T}, S) \cup \{S\}$ then we use a downward induction on the size of $\mathcal{V}$ to prove that we have the relation $y_S^{m(\mathcal{V}, S)} \prod_{T \in M(\mathcal{V}, S)} (y_S - y_T)$. Put $m = |M(\mathcal{T}, S)| + 1$ then when $|\mathcal{V}| = m$ the claim is immediately true. Suppose the claim is true for some $\mathcal{V}$ of size $l + 1 < m$ and let $\mathcal{W}$ be an $S$-tree contained in $M(\mathcal{T}, S) \cup \{S\}$ of size $l$. Then extend $\mathcal{W}$ to a set $\mathcal{V}$ by adding one disjoint element $W$ of $M(\mathcal{T}, S)$. Then by induction we have the relation.

$$
\begin{aligned}
x_S^{m(\mathcal{V},S)} \prod_{T \in M(\mathcal{V},S)} (x_S - x_T) &= 0 \\
x_S^{m(\mathcal{V},S)+1} \prod_{T \in M(\mathcal{W},S)} (x_S - x_T) &= x_W x_S^{m(\mathcal{V},S)} \prod_{T \in M(\mathcal{W},S)} (x_S - x_T) \text{ expanding } (x_S - x_W) \\
x_S^{m(\mathcal{V},S)+|W|-1} \prod_{T \in M(\mathcal{W},S)} (x_S - x_T) &= x_W^{|W|-1} x_S^{m(\mathcal{V},S)} \prod_{T \in M(\mathcal{W},S)} (x_S - x_T) \text{ iterating} \\
x_S^{m(\mathcal{W},S)} \prod_{T \in M(\mathcal{W},S)} (x_S - x_T) &= 0 \text{ as inductively } x_W^{|W|-1} = 0
\end{aligned}
$$

Thus by induction we obtain the result. Next given any other $S$-tree $\mathcal{U}$ of depth 2 we prove that the equivalent relation holds in our ring. For each $T \in M(\mathcal{T}, S)$ put $M(T) = \{\, U \in \mathcal{U} \mid U \subseteq T \,\}$ and $\mathcal{U}_T = M(T) \cup \{T\}$. First observe that the depth of $\mathcal{U}_T$ is less than or equal to 2. Next put

$$
\begin{aligned}
\mathcal{W} &= \{\, T \in M(\mathcal{T}, S) \mid M(T) \text{ is non-empty} \,\} \\
\mathcal{W}' &= \{T \in \mathcal{W} \mid |\mathcal{U}_T| = 2\} \\
\mathcal{W}'' &= \{T \in \mathcal{W} \mid |\mathcal{U}_T| \neq 2\}
\end{aligned}
$$

then $\mathcal{W} = \mathcal{W}' \amalg \mathcal{W}''$. Put $\mathcal{V} = \mathcal{W}'' \cup \{S\}$ so that $\mathcal{W}'' = M(\mathcal{V}, S)$ and we have

$$
M(\mathcal{U}, S) = \coprod_{T \in \mathcal{W}'} M(\mathcal{U}_T, T) \amalg M(\mathcal{V}, S)
$$

For each $T \in \mathcal{W}'$ we have by induction the relation $x_T^{m(\mathcal{U}_T,T)} \prod_{U \in M(\mathcal{U}_T,T)} (x_T - x_U) = 0$. Put $R_{\mathcal{T},T}$ to be the quotient of the ring $R_{\mathcal{T}}$ by the ideal $(x_S - x_T)$ then we obtain $x_S^{m(\mathcal{U}_T,T)} \prod_{U \in M(\mathcal{U}_T,T)} (x_S - x_U) = 0$ in $R_{\mathcal{T},T}$ thus $x_S^{m(\mathcal{U}_T,T)} \prod_{U \in M(\mathcal{U}_T,T)} (x_S - x_U) = f_T(x_S - x_T)$.
Next put $\mathcal{W}_S = \mathcal{W} \cup \{S\}$ and multiply the above equations for each $T \in \mathcal{W}'$ together with $(x_S - x_T)$ for each $T \in M(\mathcal{V}, S) = \mathcal{W}''$ to obtain

$$
\begin{aligned}
x_S^m \prod_{U \in M(\mathcal{U},S)} (x_S - x_U) &= f \prod_{U \in M(\mathcal{W}_S,S)} (x_S - x_U) \\
x_S^{m(\mathcal{U},S)} \prod_{U \in M(\mathcal{U},S)} (x_S - x_U) &= f x_S^{m(\mathcal{W}_S,S)} \prod_{U \in M(\mathcal{W}_S,S)} (x_S - x_U) = 0
\end{aligned}
$$

where $m = \sum_{T \in \mathcal{W}'} m(\mathcal{U}_T, T)$ and we need to show that $m + m(\mathcal{W}_S, S) = m(\mathcal{U}, S)$

$$\begin{aligned} m(\mathcal{W}_S, S) + m &= m(\mathcal{W}_S, S) + \sum_{T \in \mathcal{W}'} m(\mathcal{U}_T, T) \\ &= |S| - 1 - \sum_{T \in \mathcal{W}} (|T| - 1) + \sum_{T \in \mathcal{W}'} \left[ (|T| - 1) - \sum_{U \in M(\mathcal{U}_T, T)} (|U| - 1) \right] \\ &= |S| - 1 - \sum_{\mathcal{W}''} (|T| - 1) - \sum_{T \in \mathcal{W}'} \sum_{U \in M(\mathcal{U}_T, T)} (|U| - 1) \\ &= |S| - 1 - \sum_{T \in M(\mathcal{U}, S)} (|T| - 1) \\ &= m(\mathcal{U}, S) \end{aligned}$$

This completes the induction. Hence we have proven that $I_\mathcal{T} = J_\mathcal{T}$ thus $S_\mathcal{T} = R_\mathcal{T}$. Now consider an arbitrary forest $\mathcal{F}$ then we have each $\mathcal{F}|_T$ is a $T$-tree and

$$S_\mathcal{F} \cong \bigotimes_{T \in M(\mathcal{F})} S_{\mathcal{F}|_T} = \bigotimes_{T \in M(\mathcal{F})} R_{\mathcal{F}|_T} \cong R_\mathcal{F}$$

The composition of these maps $s : S_\mathcal{F} \to R_\mathcal{F}$ is induced from the identity map on $\mathbb{Z}_\mathcal{F}$. Then since this composition is an isomorphism we see that $I_\mathcal{F} = J_\mathcal{F}$ and $S_\mathcal{F} = R_\mathcal{F}$ □

REMARK 5.1.16. The next lemma which is now immediate from corollary 5.1.11 is the key to a purely algebraic understanding of the ring $R_\mathcal{T}$ we shall use this description to that end.

LEMMA 5.1.17. *Let $\mathcal{T}$ be a $S$-tree with $d(\mathcal{T}) > 1$ then $R_\mathcal{T} \cong \dfrac{[\bigotimes R_{\mathcal{T}|_T}][t]}{f_\mathcal{T}(t)}$ where the product is taken over $M(\mathcal{T}, S)$ and $R_\mathcal{F} \cong \bigotimes R_{\mathcal{T}|_T}$. Moreover $\{1, t, ..., t^{d-1}\}$ form a basis for $R_\mathcal{T}$ over $R_\mathcal{F}$ where $d = n(\mathcal{T}, S)$.* □

COROLLARY 5.1.18. *Let $\mathcal{F}$ be a forest then the rank of $R_\mathcal{F}$ is $n(\mathcal{F}) = \prod_{T \in \mathcal{F}} n(\mathcal{F}, T)$*

PROOF. We first prove this for the case of an $S$-tree $\mathcal{T}$. We proceed by an induction on the depth of the tree $\mathcal{T}$. For $d(\mathcal{T}) = 1$ the claim is clear. Assume the case for trees $\mathcal{T}$ with $d(\mathcal{T}) \leq n - 1$ then for $d(\mathcal{T}) = n$ we see by the previous lemma that $R_\mathcal{T} \cong [\bigotimes R_{\mathcal{T}|_T}]\{1, t, ..., t^{d-1}\}$ as free modules where $d = \deg(f_\mathcal{T}) = n(\mathcal{T}, S)$. Put $\mathcal{U}$ to be

the set $\mathcal{T}$ with $S$ removed then,

$$\begin{aligned}
\operatorname{rank}(R_\mathcal{T}) &= n(\mathcal{T},S) \prod_{T \in M(\mathcal{U})} \operatorname{rank}(R_{\mathcal{T}|_T}) \\
&= n(\mathcal{T},S) \prod_{T \in M(\mathcal{U})} \prod_{U \in \mathcal{T}|_T} n(\mathcal{T}|_T, U) \text{ by induction} \\
&= n(\mathcal{T},S) \prod_{T \in \mathcal{U}} n(\mathcal{T}, T) \\
&= \prod_{T \in \mathcal{T}} n(\mathcal{T}, T)
\end{aligned}$$

The third line is valid because for each $U \in \mathcal{T}|_T$ we have $M(\mathcal{T}|_T, U) = M(\mathcal{T}, U)$ so that $n(\mathcal{T}|_T, U) = n(\mathcal{T}, U)$ and we also have $\mathcal{U} = \coprod_{T \in M(\mathcal{U})} \mathcal{T}|_T$. This completes the induction.

Now let $\mathcal{F}$ be any forest. We have $R_\mathcal{F} \cong \bigotimes_{T \in M(\mathcal{F})} R_{\mathcal{F}|_T}$. Thus

$$\begin{aligned}
\operatorname{rank}(R_\mathcal{F}) &= \prod_{T \in M(\mathcal{F})} n(\mathcal{F}|_T) \\
&= \prod_{T \in M(\mathcal{F})} \prod_{U \in \mathcal{F}|_T} n(\mathcal{F}|_T, U) \\
&= \prod_{T \in \mathcal{F}} n(\mathcal{F}, T)
\end{aligned}$$

$\square$

## 5.2. A basis for $R_\mathcal{F}$

We complete our description of the ring $R_\mathcal{F}$ by specifying basis $\overline{A}[\mathcal{F}]$ and $\overline{B}[\mathcal{F}]$ for it. The basis $\overline{A}[\mathcal{F}]$ will naturally follow from the projective bundle description of $\overline{\mathcal{M}}_\mathcal{T}$ for trees $\mathcal{T}$ and we will see that the set $\overline{B}[\mathcal{F}]$ uses the combinatorics of forests. Later we will see how the set $\overline{B}[\mathcal{F}]$ can be used more generally to produce a basis for the cohomology ring of $\overline{\mathcal{M}}_S$.

DEFINITION 5.2.1. Let $B_\mathcal{F}$ be the monomial basis for $\mathbb{Z}_\mathcal{F}$. For any monomial $y = \prod_{T \in \mathcal{U}} y_T^{n_T} \in B_\mathcal{F}$ we define functions called the *shape*, shape : $B_\mathcal{F} \to P^2(S)$ and called the *support*, supp : $B_\mathcal{F} \to P(S)$ by shape$(y) = \mathcal{U}$ and supp$(y) = \bigcup_{T \in \mathcal{U}} T$.

DEFINITION 5.2.2. For any monomial $y$ with $y = \prod_{T \in \mathcal{U}} y_T^{n_T}$ so that shape$(y) = \mathcal{U}$ and for any $\mathcal{V} \subseteq \mathcal{U}$ we write $y|_\mathcal{V} = \prod_{V \in \mathcal{V}} y_V^{n_V}$ and call $y|_\mathcal{V}$ the *restriction* of $y$ to $\mathcal{V}$.

DEFINITION 5.2.3. Let $\mathcal{F}$ be a forest then for each $\mathcal{U} \subseteq \mathcal{F}$ we define,

$$A[\mathcal{F}] = \left\{ \prod_{T \in \mathcal{F}} y_T^{n_T} \mid 0 \leq n_T < n(\mathcal{F}, T) \text{ for every } T \in \mathcal{F} \right\}$$

$$\overline{A}[\mathcal{F}] = q_\mathcal{F}(A[\mathcal{F}]) \text{ where } q_\mathcal{F} : \mathbb{Z}_\mathcal{F} \to R_\mathcal{F} \text{ is the quotient map .}$$

and

$$B[\mathcal{F}][\mathcal{U}] = \left\{ \prod_{T \in \mathcal{U}} y_T^{m_T} \mid 1 \leq m_T < m(\mathcal{U}, T) \text{ for every } T \in \mathcal{U} \right\}$$

$$B[\mathcal{F}] = \coprod_{\mathcal{U} \subseteq \mathcal{F}} B[\mathcal{F}][\mathcal{U}]$$

$$\overline{B}[\mathcal{F}] = q_\mathcal{F}(B[\mathcal{F}])$$

The numbers $n(\mathcal{F}, T)$ and $m(\mathcal{F}, T)$ are defined in 3.3.1

LEMMA 5.2.4. *Let $\mathcal{T}$ be an S-tree with* d$(\mathcal{T}) > 1$, $\mathcal{F}$ *a forest and* $\mathcal{U} \subseteq \mathcal{F}$ *then in natural notation*

$$A[\mathcal{T}] = \left( \prod_{T \in M(\mathcal{T},S)} A[\mathcal{T}|_T] \right) \{1, ...., y_S^{d-1}\} \text{ where } d = n(\mathcal{T}, S)$$

$$A[\mathcal{F}] = \prod_{T \in M(\mathcal{F})} A[\mathcal{F}|_T]$$

$$B[\mathcal{F}] = \prod_{T \in M(\mathcal{F})} B[\mathcal{F}|_T]$$

$$B[\mathcal{U}] \subseteq B[\mathcal{F}]$$

LEMMA 5.2.5. *Let $\mathcal{F}$ be a forest $\mathcal{F}$ then $n(\mathcal{F}) = |A[\mathcal{F}]|, m(\mathcal{F}) = |B[\mathcal{F}]|$ and $n(\mathcal{F}) = m(\mathcal{F})$.*

PROOF. Let $\mathcal{F}$ be a forest. From the previous lemma we deduce that $n(\mathcal{F}) = |A[\mathcal{F}]|$ and it is clear from the definitions that $m(\mathcal{F}) = |B[\mathcal{F}]|$ where these numbers are defined in 3.3.6. For the case when $\mathcal{T}$ is an $S$-tree lemma 3.3.8 tells us that $m(\mathcal{T}) = n(\mathcal{T})$. We now use the previous lemma to see that

$$n(\mathcal{F}) = |A[\mathcal{F}]| = \prod_{T \in M(\mathcal{F})} |A[\mathcal{F}|_T]| = \prod_{T \in M(\mathcal{F})} |B[\mathcal{F}|_T]| = |B[\mathcal{F}]| = m(\mathcal{F})$$

$\square$

LEMMA 5.2.6. *For any forest $\mathcal{F}$ we have that $\overline{A}[\mathcal{F}]$ is a basis for $R_\mathcal{F}$.*

PROOF. We begin by proving it for an $S$-tree $\mathcal{T}$. We proceed by induction on the depth $\mathrm{d}(\mathcal{T})$ of the tree $\mathcal{T}$. If $\mathrm{d}(\mathcal{T}) = 1$ then we are reduced to the case of the cohomology of $PV_S$ which is well known. Suppose the case is true for $\mathrm{d}(\mathcal{T}) < n$ for some $n > 1$ then for $\mathrm{d}(\mathcal{T}) = n$ lemma 5.1.17 to tells us that $R_\mathcal{T} \cong [\bigotimes R_{\mathcal{T}|_T}]\{1, t, ..., t^{d-1}\}$ as free modules where $d = n(\mathcal{T}, S)$. Inductively we see that $\overline{A}[\mathcal{T}|_T]$ is a basis for $R_{\mathcal{T}|_T}$ and one can check by apply the above correspondence that in natural notation $\overline{A}[\mathcal{T}] = \prod \overline{A}[\mathcal{T}|_T]\{1, x_S, ..., x_S^{d-1}\}$ to deduce our result. To complete the proof let $\mathcal{F}$ be a forest then $R_\mathcal{F} \cong \bigotimes_{T \in M(\mathcal{F})} R_{\mathcal{F}|_T}$ and we use the tree case to deduce our result.

$\square$

LEMMA 5.2.7. *Let $\mathcal{F}$ be a forest and $y = \prod_{T \in \mathcal{U}} y_T^{m_T} \in B[\mathcal{F}]$ then*

(1) $\sum_{T \in \mathcal{U}} m_T \leq \sum_{T \in M(\mathcal{F})} (|T| - 2)$ *with equality* $\iff \mathcal{U} = M(\mathcal{F})$ *and* $m_T = |T| - 2$ *for each* $T \in M(\mathcal{F})$.

(2) *For each* $U \in \mathcal{U}$ *we have* $\sum_{T \in \mathcal{U}|_U} m_T \leq |U| - 2$ *with equality* $\iff \mathcal{U}|_U = \{U\}$ *and* $m_U = |U| - 2$

(3) *In particular if $\mathcal{F}$ is an $S$-tree* $\sum_{T \in \mathcal{U}} m_T \leq |S| - 2$ *with equality* $\iff \mathcal{U} = \{S\}$ *and* $m_S = |S| - 2$.

We first prove the special case when $\mathcal{F}$ is an $S$-tree.

PROOF.

$$\begin{aligned}
\sum_{T \in \mathcal{U}} m_T &\leq \sum_{T \in \mathcal{U}} (m(\mathcal{U}, T) - 1) \\
&= \sum_{T \in \mathcal{U}} m(\mathcal{U}, T) - |\mathcal{U}| \\
&= \sum_{T \in M(\mathcal{U})} |T| - |M(\mathcal{U})| - |\mathcal{U}| \text{ by lemma 3.3.2} \\
&= |\coprod_{T \in M(\mathcal{U})} T| - |M(\mathcal{U})| - |\mathcal{U}| \\
&\leq |S| - 2 \text{ as } |M(\mathcal{U})|, |\mathcal{U}| \geq 1
\end{aligned}$$

It is now clear that we have equality if and only if $\mathcal{U} = \{S\}$ and $m_S = |S| - 2$. To prove the general case we observe that we have the partition $\mathcal{F} = \coprod_{T \in M(\mathcal{F})} \mathcal{F}|_T$ where each $\mathcal{F}|_T$ is a $T$-tree and apply the special case. $\square$

## 5.3. A filtration for $R_\mathcal{F}$

We next define a weight function on $R_\mathcal{F}$ which we use to give us a filtration. We will use this to show that the set $\overline{B}[\mathcal{F}]$ spans $R_\mathcal{F}$. Since we have shown that the size of $B[\mathcal{F}]$ is the rank of the ring $R_\mathcal{F}$, which is a finite free module, this then shows us that $\overline{B}[\mathcal{F}]$ is a basis for $R_\mathcal{F}$. We do not provide many of the proofs in this section as they are particulary simple and unilluminating.

DEFINITION 5.3.1. Let $B_\mathcal{F}$ be the monomial basis for $\mathbb{Z}_\mathcal{F}$ and $y \in B_\mathcal{F}$ so that $y = \prod_{T \in \mathcal{F}} y_T^{n_T}$ say. We define a function $\mathrm{wt} : B_\mathcal{F} \to \mathbb{N}$ called the *weight* by $\mathrm{wt}(y) = \sum_{T \in \mathcal{F}} n_T |T|$. In particular this gives the monomial $y_T$ weight $|T|$. We also put $\deg : R_\mathcal{F} \to \mathbb{N}$ to be the cohomological degree function, so that $\deg(y) = 2 \sum_{T \in \mathcal{F}} n_T$

LEMMA 5.3.2. $\mathrm{wt}(xy) = \mathrm{wt}(x) + \mathrm{wt}(y)$ *and* $\deg(x) \leq \mathrm{wt}(x) \leq \frac{n}{2} \deg(x)$ *where* $n = |\mathrm{supp}(\mathcal{F})|$. $\square$

DEFINITION 5.3.3. We define filtrations on $\mathbb{Z}_\mathcal{F}$ by

$$F_k \mathbb{Z}_\mathcal{F} = \text{span}\{\, y \in B_\mathcal{F} \mid \text{wt}(y) \geq k \,\}$$
$$G_k \mathbb{Z}_\mathcal{F} = \text{span}\{\, y \in B_\mathcal{F} \mid \deg(y) \geq k \,\}$$

LEMMA 5.3.4. *Put $H = F$ or $G$ then $H_k \mathbb{Z}_\mathcal{T}$ is a convergent decreasing filtration of $\mathbb{Z}_\mathcal{F}$, that is,*

$$\mathbb{Z}_\mathcal{F} = H_0 \mathbb{Z}_\mathcal{F} \supseteq H_1 \mathbb{Z}_\mathcal{F} \supseteq H_2 \mathbb{Z}_\mathcal{F} \supseteq \ldots \text{ and } \bigcap_{k \geq 0} H_k \mathbb{Z}_\mathcal{F} = \{0\}$$

$\square$

DEFINITION 5.3.5. We define a function $\text{wt} : \mathbb{Z}_\mathcal{F} \to \mathbb{N} \cup \{\infty\}$ called the *weight* as follows. For every non-zero $y \in \mathbb{Z}_\mathcal{F}$ we know by the previous lemma that there exists a largest k such that $y \in F_k \mathbb{Z}_\mathcal{F}$ but $y \notin F_{k+1} \mathbb{Z}_\mathcal{F}$. We define $\text{wt}(y) = k$ and $\text{wt}(0) = \infty$

LEMMA 5.3.6. *The weight function has the following properties*

(1) *If $y = \sum_{i \in I} a_i y_i$ with $y_i \in B_\mathcal{F}$ and $a_i \neq 0$ then $\text{wt}(y) = \min\{\, \text{wt}(y_i) \mid i \in I \,\}$*
(2) $\text{wt}(xy) = \text{wt}(x) + \text{wt}(y)$.

$\square$

DEFINITION 5.3.7. We define filtrations on $R_\mathcal{F}$ by $F_k R_\mathcal{F} = q_\mathcal{F}(F_k \mathbb{Z}_\mathcal{F})$ and $G_k R_\mathcal{F} = q_\mathcal{F}(G_k \mathbb{Z}_\mathcal{F})$ where $q_\mathcal{F} : \mathbb{Z}_\mathcal{F} \to R_\mathcal{F}$ is the quotient map.

LEMMA 5.3.8. *Put $H = F$ or $G$ then $H_k R_\mathcal{F}$ is a convergent decreasing filtration of $R_\mathcal{F}$, that is,*

$$R_\mathcal{F} = H_0 R_\mathcal{F} \supseteq H_1 R_\mathcal{F} \supseteq H_2 R_\mathcal{F} \supseteq \ldots \quad \text{and} \bigcap_{k \geq 0} H_k R_\mathcal{F} = \{0\}$$

PROOF. The first part is clear. By lemma 5.3.2 it is clear that $F_k\mathbb{Z}_\mathcal{F} \subseteq G_l\mathbb{Z}_\mathcal{F}$ where $l = [2k/n]$ and $[a]$ is the integer part of $a$ thus $F_kR_\mathcal{F} = q_\mathcal{F}(F_k\mathbb{Z}_\mathcal{F}) \subseteq q_\mathcal{F}(G_l\mathbb{Z}_\mathcal{F}) = G_lR_\mathcal{F}$. Therefore $\bigcap_k F_kR_\mathcal{F} \subseteq \bigcap_l G_lR_\mathcal{F} = \{0\}$ as $R_\mathcal{F}$ is graded by degree.

□

DEFINITION 5.3.9. We next define a function wt : $R_\mathcal{F} \to \mathbb{N} \cup \{\infty\}$ called the *weight* as follows. For every non-zero $x \in R_\mathcal{F}$ we know by the previous lemma that there exists a largest $k \in \mathbb{N}$ such that $x \in F_kR_\mathcal{F}$, but $x \notin F_{k+1}R_\mathcal{F}$. We define $\mathrm{wt}(x) = k$ and $\mathrm{wt}(0) = \infty$.

LEMMA 5.3.10. *The function* wt *has the following properties.*

(1) *If* $x = \sum_{i \in I} a_i x_i$ *with* $x_i \in R_\mathcal{F}$ *then* $\mathrm{wt}(x) \geq \min\{\mathrm{wt}(x_i) \mid i \in I\}$
(2) *There is an* $m \in \mathbb{N}$ *such that* $\mathrm{wt}(x) \leq m$ *for all non-zero* $x$
(3) $\mathrm{wt}(xy) \geq \mathrm{wt}(x) + \mathrm{wt}(y)$
(4) $\mathrm{wt}(q_\mathcal{F}(x)) \geq \mathrm{wt}(x)$

□

DEFINITION 5.3.11. We define $F_kB[\mathcal{F}] = \{x \in B[\mathcal{F}] \mid \mathrm{wt}(x) \geq k\}$ and the induced set $F_k\overline{B}[\mathcal{F}] = q_\mathcal{F}(F_kB[\mathcal{F}])$

DEFINITION 5.3.12. Let $y \in \mathbb{Z}_\mathcal{F}$ be a monomial. We say $y$ is *admissible* if $y \in B[\mathcal{F}]$ and *inadmissible* if $y \notin B[\mathcal{F}]$. We also say $y \in \mathbb{Z}_\mathcal{F}$ is *minimally inadmissible* if it has one of the following forms

(1) $y_T^{m(\mathcal{T},T)} \prod_{U \in M(\mathcal{T},T)} y_U$ with $\mathcal{T}$ a $T$-tree of depth 2 and $m(\mathcal{T},T)$ is as defined in 3.3.1
(2) $y_T^{|T|-1}$ for some $T \in \mathcal{F}$.

LEMMA 5.3.13. *Let* $z \in \mathbb{Z}_\mathcal{F}$ *be a non-zero monomial, then* $z$ *is inadmissible if and only if* $z = xy$ *with* $x$ *minimally inadmissible.* □

LEMMA 5.3.14. *If* $y \in \mathbb{Z}_\mathcal{F}$ *is minimally inadmissible then* $\mathrm{wt}(q_\mathcal{F}(y)) > \mathrm{wt}(y)$.

PROOF. We have by a relation that $x_T^{m(\mathcal{T},T)} \prod_{U \in M(\mathcal{T},T)} (x_U - x_T) = 0$. Upon expansion of this expression we have

$$x = x_T^m \prod_{U \in M(\mathcal{T},T)} x_U = \sum_{i=1}^{n} \epsilon(i) x_T^{m+i} x_i$$

where $\epsilon : I \to \{\pm 1\}$ and $m = m(\mathcal{T},T)$. Put $\mathcal{T}_i = \text{shape}(x_i)$ then $\mathcal{T}_i \subset M(\mathcal{T},T)$, $|\mathcal{T}_i| = |\mathcal{T}| - 1 - i$ and for every $U \in \mathcal{T}_i$ we have $|U| < |T|$. Consider an element $x_T^{m+i} x_i$ of this sum then

$$\begin{aligned}
\text{wt}(x_T^{m+i} x_i) &\geq m|T| + i|T| + \sum_{U \in \mathcal{T}_i} |U| \\
&> m|T| + \sum_{U \in M(\mathcal{T},T)} |U| \\
&= \text{wt}(y)
\end{aligned}$$

thus we have $\text{wt}(q_{\mathcal{F}}(y)) \geq \min\{\text{wt}(x_T^{m+i} x_i) \mid 1 \leq i \leq n\} > \text{wt}(y)$

$\square$

COROLLARY 5.3.15. *Let $z \in \mathbb{Z}_{\mathcal{F}}$ be a monomial. Then if $z$ is inadmissible, that is $z \notin B[\mathcal{F}]$ we have $\text{wt}(q_{\mathcal{F}}(z)) > \text{wt}(z)$.*

PROOF. Let $z \in \mathbb{Z}_{\mathcal{F}}$ be an inadmissible element then we may write $z = xy$ with $x$ minimally admissible thus we have

$$\begin{aligned}
\text{wt}(q_{\mathcal{F}}(z)) &= \text{wt}(q_{\mathcal{F}}(xy)) \\
&= \text{wt}(q_{\mathcal{F}}(x) q_{\mathcal{F}}(y)) \\
&\geq \text{wt}(q_{\mathcal{F}}(x)) + \text{wt}(q_{\mathcal{F}}(y)) \\
&> \text{wt}(x) + \text{wt}(y) \text{ by lemma 5.3.10 part 4 and 5.3.14} \\
&= \text{wt}(xy) \\
&= \text{wt}(z)
\end{aligned}$$

$\square$

LEMMA 5.3.16. $R_\mathcal{F} = \mathrm{span} B[\mathcal{F}]$

PROOF. The proof will follow from a downward induction on the weight $w$ of the statement that for every $k \in \mathbb{N}$ $F_k R_\mathcal{F} = \mathrm{span} F_k \overline{B}[\mathcal{F}]$. For $k \gg 0$ we know by lemmas 5.3.10 part 2 and 5.2.7 that $F_k R_\mathcal{F} = 0 = \mathrm{span} F_k \overline{B}[\mathcal{F}]$. Suppose it is true for all $k > w$. Let $x \in F_w R_\mathcal{F} = q_\mathcal{F}(F_w \mathbb{Z}_\mathcal{F})$, then $x = q_\mathcal{F}(y)$ where $y = \sum_{i \in I} a_i y_i \in F_w \mathbb{Z}_\mathcal{F}$ and $y_i \in B_\mathcal{F}$ thus $\mathrm{wt}(y_i) \geq w$. So it is enough to show that for every $i \in I$, $q_\mathcal{F}(y_i) \in \mathrm{span} F_w \overline{B}[\mathcal{F}]$. If $y_i \in B[\mathcal{F}]$ then this holds so we may assume that $y_i \notin B[\mathcal{F}]$. In this case we know that $\mathrm{wt}(q_\mathcal{F}(y_i)) > \mathrm{wt}(y_i) \geq w$ and so by induction $q_\mathcal{F}(y_i) \in F_{k+1} R_\mathcal{F} = \mathrm{span} F_{k+1} \overline{B}[\mathcal{F}] \subseteq \mathrm{span} F_k \overline{B}[\mathcal{F}]$. Therefore $R_\mathcal{F} = F_0 R_\mathcal{F} = \mathrm{span} F_0 \overline{B}[\mathcal{F}] = \mathrm{span} \overline{B}[\mathcal{F}]$ and we are done. □

LEMMA 5.3.17. *For each forest $\mathcal{F}$ the set $\overline{B}[\mathcal{F}]$ forms a basis for $R_\mathcal{F}$. Put $x_\mathcal{F} = \prod_{T \in M(\mathcal{F})} x_T^{|T|-2}$ and $n = \deg(x_\mathcal{F})$ then we have $R_\mathcal{F}^{2n} = \mathbb{Z}[x_\mathcal{F}]$ and $R_\mathcal{F}^i = 0$ for $i > 2n$. In particular if $\mathcal{F}$ is an $S$-tree then $n = |S| - 2$ and $R_\mathcal{F}^{2(|S|-2)} = \mathbb{Z}[x_S^{|S|-2}]$.*

PROOF. The previous lemma shows us that $R_\mathcal{F} = \mathrm{span} \overline{B}[\mathcal{F}]$ and lemma 5.2.5 shows that $|B[\mathcal{F}]| = \mathrm{rank}(R_\mathcal{F})$. Thus $|B[\mathcal{F}]| = \mathrm{rank}(R_\mathcal{F}) \leq |\overline{B}[\mathcal{F}]| \leq |B[\mathcal{F}]|$ so that $|\overline{B}[\mathcal{F}]| = \mathrm{rank}(R_\mathcal{F})$. Since $x_\mathcal{F} \in B[\mathcal{F}]$ we see that $x_\mathcal{F}$ is non-zero and by lemma 5.2.7 part 1 the degree of $x_\mathcal{F}$ is maximal and is the only element of maximal degree. This completes the proof. □

LEMMA 5.3.18. *Let $\mathcal{F}$ be a forest and $\mathcal{U} \subseteq \mathcal{F}$. Let $i : \mathbb{Z}_\mathcal{U} \to \mathbb{Z}_\mathcal{F}$ be the inclusion map. We have the inclusion of ideal $J_\mathcal{U} \subseteq J_\mathcal{F}$ and an induced map $\phi : R_\mathcal{U} \to R_\mathcal{F}$. Then the map $\phi$ is injective.*

PROOF. The proof is trivial since we have an inclusion of basis $B[\mathcal{U}] \subseteq B[\mathcal{F}]$ □

COROLLARY 5.3.19. *Let $y \in \mathbb{Z}_\mathcal{F}$ with $\mathrm{shape}(y) = \mathcal{V}$ say and put $x = q_\mathcal{F}(y)$ where $q_\mathcal{F} : \mathbb{Z}_\mathcal{F} \to R_\mathcal{F}$ is the usual quotient map. Suppose there exists a subset $\mathcal{U}$ of $\mathcal{V}$ such that $\deg(y|_\mathcal{U}) \geq 2(|\mathrm{supp}(\mathcal{U})| - 1)$ then $x = 0$.*

PROOF. Put $T = \text{supp}(\mathcal{U})$ and put $\mathcal{W}$ to be the $T$-tree $\mathcal{U} \cup \{T\}$. by considering $x|_\mathcal{U}$ as an element in $R_\mathcal{W}$ we see by the previous result that $x|_\mathcal{U} = 0$ in $R_\mathcal{W}$. Let $\phi : R_\mathcal{W} \to R_\mathcal{F}$ be the ring map induced by the inclusion $i : \mathbb{Z}_\mathcal{W} \to \mathbb{Z}_\mathcal{F}$ then since $x_\mathcal{U}$ is zero in $R_\mathcal{W}$ we see that $x = x|_\mathcal{U} \cdot x'$ is zero in $R_\mathcal{F}$. □

## 5.4. The zero condition for monomials of $R_\mathcal{F}$

In this section we prove precisely which monomials of the ring $R_\mathcal{F}$ are zero. Later we will see how to use this to deduce which monomials are zero in more general rings. In particular for the cohomology ring of $\overline{\mathcal{M}}_S$. The last corollary gives us sufficient conditions for monomials of $R_\mathcal{F}$ to be zero. We next prove that they are also necessary. In this section we will be using the ordinary degree of a monomial.

DEFINITION 5.4.1. Let $\mathcal{F}$ be a forest. We define $N[\mathcal{F}]$ to be the set of monomials of the form $y = \prod_{T \in \mathcal{F}} y_T^{n_T}$ such that for every $T \in \mathcal{F}$ we have $\sum_{\substack{U \in \mathcal{F} \\ U \subseteq T}} n_U \leq |T| - 2$. We then define $\overline{N}[\mathcal{F}] = q_\mathcal{F}(N[\mathcal{F}])$ where $q_\mathcal{F} : \mathbb{Z}_\mathcal{F} \to R_\mathcal{F}$ is the quotient map.

LEMMA 5.4.2. Let $\mathcal{T}$ be an $S$-tree. Let $y = \prod_{T \in \mathcal{U}} y_T^{n_T}$ be an element of $N[\mathcal{T}]$ and $\mathcal{V}$ be the set $\mathcal{U}$ with $S$ removed. Suppose that $M(\mathcal{V})$ is non-empty then $m + \deg(y) - n_S \leq |S| - 2$ where $m = |S| - 1 - \sum_{T \in M(\mathcal{V})} (|T| - 1)$

PROOF. We have,

$$\begin{aligned}
m + \deg(y) - n_S &= m + \sum_{T \in \mathcal{V}} n_T \\
&= m + \sum_{T \in M(\mathcal{V})} \sum_{U \in \mathcal{V}|_T} n_U \\
&\leq m + \sum_{T \in M(\mathcal{V})} (|T| - 2) \text{ by definition 5.4.1} \\
&= |S| - 1 - |M(\mathcal{V})| \\
&\leq |S| - 2 \text{ as } |M(\mathcal{V})| \geq 1
\end{aligned}$$

□

LEMMA 5.4.3. *Let $x$ be an element of $R_\mathcal{F}$ and $y \in \mathbb{Z}_\mathcal{F}$ with $x = q_\mathcal{F}(y)$ then*

(a) $q_\mathcal{F}(y) = 0$ *if and only if* $y \notin N[\mathcal{F}]$.

(b) *If $y \in N[\mathcal{F}]$ and* shape$(y)$ *is a $T$-tree of maximal degree $|T| - 2$ then $x = x_T^{|T|-2}$.*

PROOF. We first prove the special case for $S$-trees $\mathcal{T}$. We have already proven in lemma 5.3.19 that if $y \notin N[\mathcal{T}]$ then $x = q_\mathcal{T}(y) = 0$ so it remains to prove the other direction. We prove that for any $x \in \overline{N}[\mathcal{T}]$ we have $x_S^{|S|-2-\deg(x)} x = x_S^{|S|-2}$. From this it immediately follows that $x \neq 0$. Let $y = \prod_{U \in \mathcal{U}} y_U^{n_U}$ be an element of $N[\mathcal{T}]$, we may suppose by lemma 5.3.17 that $\mathcal{U} \neq \{S\}$. Put $\mathcal{V}$ to be the set $\mathcal{U}$ with $S$ removed and $m = |S| - 1 - \sum_{V \in M(\mathcal{V})} (|V| - 1)$ then we have the relation

$$x_S^m \prod_{V \in M(\mathcal{V})} (x_S - x_V) = 0$$

$$x_S^m \prod_{V \in \mathcal{V}} (x_S - x_V) = 0$$

$$x_S^m \prod_{V \in \mathcal{V}} x_V^{n_V - 1} \prod_{V \in \mathcal{V}} (x_V - x_S) = 0$$

$$x_S^{|S|-2-\deg(x)+n_S} \prod_{V \in \mathcal{V}} (x_V^{n_V} - x_S x_V^{n_V - 1}) = 0$$

$$x_S^{|S|-2-\deg(x)} x_S^{n_S} \prod_{V \in \mathcal{V}} (x_V^{n_V} - x_S x_V^{n_V - 1}) = 0$$

For the fourth step we use lemma 5.4.2 to see that $n = |S| - 2 + n_S - m - \deg(x) \geq 0$ and multiply by $x_S^n$. Then expanding the last equation we see that,

$$x_S^{|S|-2-\deg(x)} x = \sum_{i \in I} \epsilon(i) x_S^{|S|-2-\deg(x_i)} x_i$$

where the $x_i$ properly divide $x$. Therefore the $x_i$ lie in $\overline{N}[\mathcal{T}]$ and $\deg(x_i) < \deg(x)$. We also have that $\epsilon : I \to \{\pm 1\}$ and $\sum_{i \in I} \epsilon(i) = 1$. To prove the claim we now apply an induction argument on the degree of $x$. For $\deg(x) = 0$ the claim is clear, then assume the claim is true for any monomials $x$ with $\deg(x) < n$ then for $\deg(x) = n$ we have

$$\begin{aligned}
x_S^{|S|-2-\deg(x)}x &= \sum_{i\in I}\epsilon(i)x_S^{|S|-2-\deg(x_i)}x_i \\
&= \sum_{i\in I}\epsilon(i)x_S^{|S|-2} \\
&= x_S^{|S|-2}\sum_{i\in I}\epsilon(i) \\
&= x_S^{|S|-2}
\end{aligned}$$

this completes the induction. Since we know by lemma 5.3.17 that $x_S^{|S|-2}$ is non-zero we have proven the first part of the claim. The second claim is now immediate. We now prove the claim for an arbitrary forest $\mathcal{F}$. One first checks that in the natural notation $N(\mathcal{F}) = \prod_{T\in M(\mathcal{F})} N(\mathcal{F}|_T)$ and then apply lemma 5.1.9 to deduce our result. □

## 5.5. Summary of results

SUMMARY 5.5.1. Here we summarize the results of this chapter. Let $\mathcal{F}$ be a forest, for neatness we suppose that every $T \in \mathcal{F}$ has $|T| > 2$. Clearly if $|T| = 2$ we have $x_T = 0$.

(0) $S_\mathcal{F} = R_\mathcal{F}$ and the specified relations of $I_\mathcal{F}$ are minimal with $|\mathcal{F}|$ generators.

(1) The cohomology ring of $\overline{\mathcal{M}}_\mathcal{F}$ is $R_\mathcal{F}$.

(2) $R_\mathcal{F}$ is a finite free module generated by its elements of degree 2.

(3) The rank of $R_\mathcal{F}^2$ is $|\mathcal{F}|$.

(4) $R_\mathcal{F} \cong \bigotimes_{T\in M(\mathcal{F})} R_{\mathcal{F}|_T}$ and each $\mathcal{F}|_T$ is a $T$-tree.

(5) The rank of $R_\mathcal{F}$ is $n(\mathcal{F}) = \prod_{T\in\mathcal{F}} n(\mathcal{F},T)$ and $n(\mathcal{F}) = m(\mathcal{F})$.

(6) $R_\mathcal{T} \cong \dfrac{[\bigotimes R_{\mathcal{T}|_T}][t]}{f_\mathcal{T}(t)}$ where $f_\mathcal{T}$ is defined in construction 5.1.6.

(7) $\overline{A}[\mathcal{F}]$ and $\overline{B}[\mathcal{F}]$ form basis for $R_\mathcal{F}$.

(8) For any $\mathcal{U} \subseteq \mathcal{F}$ the evident map $\phi : R_\mathcal{U} \to R_\mathcal{F}$ is injective.

(9) $\overline{N}[\mathcal{F}]$ are the non zero elements of $R_\mathcal{F}$.

(10) Put $x_\mathcal{F} = \prod_{T\in M(\mathcal{F})} x_T^{|T|-2}$ and $n = \deg(x_\mathcal{F})$ then $R_\mathcal{F}^{2n} = \mathbb{Z}[x_\mathcal{F}]$ and $R_\mathcal{F}^i = 0$ for $i > 2n$.

# CHAPTER 6

# The topology of $\overline{\mathcal{M}}_\mathcal{L}$ for non tree sets

We shall next consider the space $\overline{\mathcal{M}}_\mathcal{L}$ for certain collections $\mathcal{L}$ of subsets of $S$ other than forests. The restrictions we will impose on these sets will be mild for our purposes but necessary for our approach. For such sets $\mathcal{L}$ this will enable us to introduce the notion of the type of an element $\underline{M} \in \overline{\mathcal{M}}_\mathcal{L}$. This will be an $S$-tree whose elements are members of $\mathcal{L}$. It is this construction that will enable us to compare the spaces $\overline{\mathcal{M}}_\mathcal{L}$ with $\overline{\mathcal{M}}_\mathcal{T}$. The latter space has already been studied in some detail in chapter 4 and the map between them is just the evident projection map. By the end of this chapter we will prove that $\overline{\mathcal{M}}_\mathcal{L}$ is a smooth irreducible projective variety of dimension $|S| - 2$ and we use the embedding $i : \overline{\mathcal{M}}_\mathcal{L} \to \prod PV_T$ to compute the tangent bundle of $\overline{\mathcal{M}}_\mathcal{L}$. Given a result from chapter 4 this reduces to computing the dimension of a certain vector space. We also gather in this chapter a number of results concerning the cohomology of $\overline{\mathcal{M}}_\mathcal{L}$. These will be required in chapter 9 where we give a presentation for the cohomology ring of $\overline{\mathcal{M}}_S$. We will prove that for each permitted $\mathcal{L}$ its cohomology ring is a finite free module generated by certain characteristic classes in degree two and the rank of $H^2(\overline{\mathcal{M}}_\mathcal{L})$ is the number of elements in $\mathcal{L}$ with size $|T|$ at least three. The same techniques will allow us to analyze the Chow ring. In particular this all works for $\overline{\mathcal{M}}_S$. For the final section of this chapter we prove using a blowup argument that $\overline{\mathcal{M}}_S$ and $\overline{\mathcal{X}}_S$ are in fact isomorphic, however this only implies the existence of such a morphism between them. The approach relies on a paper by Kapranov [**8**] and also provides some motivation for the construction of an isomorphism.

## 6.1. The associated tree to elements of $\overline{\mathcal{M}}_\mathcal{L}$

We next define the notion of the type of an element $\underline{M} \in \overline{\mathcal{M}}_\mathcal{L}$ which for a large collection of sets $\mathcal{L}$ will be a tree. The point of this definition is that it encodes the minimal amount of information that one requires to reproduce elements of $\overline{\mathcal{M}}_\mathcal{L}$ with the specified tree type and allows us to compare the space $\overline{\mathcal{M}}_\mathcal{L}$ with the various tree spaces $\overline{\mathcal{M}}_\mathcal{T}$ under the projection maps. For the space $\overline{\mathcal{M}}_S$ we will later see that this notion of the type of an element is equivalent to the ordinary notion of the tree attached to an S-curve. We will explain this equivalence precisely in the relevant chapter. First we define the notion of a thicket $\mathcal{L}$ which is a collection of subsets of $S$. We will see in a moment that the notion of a thicket is precisely the condition required to ensure that the type of every element in $\overline{\mathcal{M}}_\mathcal{L}$ is a tree.

DEFINITION 6.1.1. We say that a collection of subsets $\mathcal{L}$ of $S$ is a *thicket* if $S \in \mathcal{L}$ and for all $T \in \mathcal{L}$ $|T| > 1$ and for any $U, V \in \mathcal{L}$ with $U \cap V$ non-empty we have $W = U \cup V \in \mathcal{L}$.

LEMMA 6.1.2. *Let $U, W \subseteq S$ and put $T = U \cup W$ then the map $\pi : V_T \to V_U \oplus V_W$ is injective if and only if $U \cap W$ is non-empty.* □

DEFINITION 6.1.3. Let $\mathcal{L}$ be a thicket and $\underline{M} \in \overline{\mathcal{M}}_\mathcal{L}$ then we define the *type* of $\underline{M}$ denoted type($\underline{M}$) as,

$$\text{type}(\underline{M}) = \{\, U \in \mathcal{L} \mid \text{ for all } T \supset U \text{ with } T \in \mathcal{L} \implies \pi_U^T(M_T) = 0 \,\}$$

LEMMA 6.1.4. *Let $\mathcal{L}$ be a thicket and $\underline{M} \in \overline{\mathcal{M}}_\mathcal{L}$ then* type($\underline{M}$) *is an $S$-tree.*

PROOF. Suppose for a contradiction that type($\underline{M}$) is not a forest. Then there are elements $U, V \in$ type($\underline{M}$) such that $U \not\subseteq V$, $V \not\subseteq U$ and $U \cap V$ is non-empty. Put $W = U \cup V$, as $\mathcal{L}$ is a thicket $W \in \mathcal{L}$ and $W \supset U, V$ therefore $\pi_U^W M_W = 0$ and $\pi_V^W M_W = 0$. Then by lemma 6.1.2 we see that $M_W = 0$ a contradiction since $M_W \in PV_W$, thus type($\underline{M}$) is a forest. It is also clear that $S \in$ type($\underline{M}$) therefore type($\underline{M}$) is an $S$-tree.
□

LEMMA 6.1.5. *Let $\mathcal{L}$ be a thicket and $\underline{M} \in \overline{\mathcal{M}}_\mathcal{L}$ with $\mathcal{T} = $ type($\underline{M}$). Suppose*

(1) $T \in \mathcal{T}$
(2) $U \in \mathcal{L}$ and $U \subseteq T$
(3) For all $V \in \mathcal{L}$ with $U \subset V \subseteq T$ we have $\pi_U^V M_V = 0$

*Then $U \in \mathcal{T}$.*

PROOF. Choose any $V \in \mathcal{L}$ with $V \supset U$, we must show that $\pi_U^V M_U = 0$. This is clear if either $V \subseteq T$ or $T \subseteq V$ so suppose neither $V \not\subseteq T$ nor $T \not\subseteq V$. Put $W = V \cup T$ then $V \cap T \supseteq U$ and so $W \in \mathcal{L}$ since $\mathcal{L}$ is a thicket also $W \supset T, V$ thus $\pi_T^W M_W = 0$ as $T \in \mathcal{T}$. Then by lemma 6.1.2 we must have $\pi_V^W M_W = M_V$ as $T$ and $V$ intersect non trivially and $M_W \in PV_W$. Thus $\pi_U^V M_V = \pi_U^V \pi_V^W M_W = \pi_U^W M_W = \pi_U^T \pi_T^W M_W = \pi_U^T(0_T) = 0_U$ and $U \in \mathcal{T}$. $\square$

DEFINITION 6.1.6. For any thicket $\mathcal{L}$ of $S$ we define $\mathbb{T}_\mathcal{L} = \{\, \mathcal{T} \subseteq \mathcal{L} \mid \mathcal{T} \text{ is an } S-\text{tree}\,\}$

DEFINITION 6.1.7. Let $\mathcal{L}$ be a thicket and $\mathcal{T} \in \mathbb{T}_\mathcal{L}$ then we define a function root $: \mathcal{L} \to \mathcal{T}$ by $\text{root}(U) = T$ where $T \in \mathcal{T}$ is the element of minimal size containing $U$. We say $T$ is the *root* of $U$ in $\mathcal{T}$. This always exists because $S \in \mathcal{T}$.

LEMMA 6.1.8. *Let $\mathcal{L}$ be a thicket, $\underline{M} \in \prod_{U \in \mathcal{L}} PV_U$ and $\mathcal{T} \in \mathbb{T}_\mathcal{L}$. Further if we define $\pi_\mathcal{T} : \prod_{U \in \mathcal{L}} PV_U \to \prod_{U \in \mathcal{T}} PV_U$ to be the projection map and suppose that $\pi_\mathcal{T}(\underline{M}) \in \overline{\mathcal{M}}_\mathcal{T}$ then the following are equivalent*

(1) $\underline{M} \in \overline{\mathcal{M}}_\mathcal{L}$ *and* $\text{type}(\underline{M}) \subseteq \mathcal{T}$
(2) *For every $T \in \mathcal{T}$ and $U \in \mathcal{L}$ with $T = \text{root}(U)$ we have $\pi_U^T M_T = M_U$.*

*moreover $\underline{M}$ is the unique element with respect to the second property.*

PROOF. $\implies$ Suppose for a contradiction that there is a $U \in \mathcal{L}$ with $T = \text{root}(U)$ such that $\pi_U^T M_T = 0$. Then clearly $U \notin \mathcal{T}$, define $\mathcal{L}_U^T = \{V \in \mathcal{L} \mid U \subset V \subseteq T\,\}$. Since $U \notin \mathcal{T}$ then by the previous lemma there exists an element $V \in \mathcal{L}_U^T$ with $\pi_U^V M_V = M_U$ and clearly $V \neq T$. Choose such an element $V$ whose size is maximal with respect to this property, then for all $W \in \mathcal{L}_V^T$ we have $\pi_V^W M_W = 0$ and again by the previous lemma $V \in \mathcal{T}$ and $\text{root}(U) \subseteq V \subset T$, a contradiction thus $\pi_U^T M_T = M_U$.

For the converse let $U$ and $V$ be elements of $\mathcal{L}$ with $U \subseteq V$. Put $T = \text{root}(U)$ and $W = \text{root}(V)$ then $T \subseteq W$ and by condition 2 we have $M_V = \pi_V^W M_W$, $M_U = \pi_U^T M_T$ and $\pi_T^W M_W \leq M_T$ by the assumption that $\pi_{\mathcal{T}}(\underline{M}) \in \overline{\mathcal{M}}_{\mathcal{T}}$. Therefore

$$\begin{aligned}
\pi_U^V M_V &= \pi_U^V \pi_V^W M_W \\
&= \pi_U^W M_W \\
&= \pi_U^T \pi_T^W M_W \\
&\leq \pi_U^T M_T \\
&= M_U
\end{aligned}$$

thus $\underline{M} \in \overline{\mathcal{M}}_{\mathcal{L}}$. Now given $U \in \mathcal{L} \setminus \mathcal{T}$ put $T = \text{root}(U)$ then $T \supset U$ and $\pi_U^T M_T = M_U$ thus $U \notin \text{type}(\underline{M})$ therefore $\text{type}(\underline{M}) \subseteq \mathcal{T}$. The uniqueness issue is clear, this completes the proof. $\square$

LEMMA 6.1.9. *The map* $\text{type} : \overline{\mathcal{M}}_{\mathcal{L}} \to \mathbb{T}_{\mathcal{L}}$ *is surjective.*

PROOF. To prove this it suffices to prove the special case that $\text{type} : \overline{\mathcal{M}}_S \to \mathbb{T}_S$ is surjective. Then given any $\mathcal{T} \in \mathbb{T}_{\mathcal{L}}$ we can choose an element $\underline{M} \in \overline{\mathcal{M}}_S$ with $\text{type}(\underline{M}) = \mathcal{T}$. Let $\pi_{\mathcal{L}} : \overline{\mathcal{M}}_S \to \overline{\mathcal{M}}_{\mathcal{L}}$ be the projection map and put $\underline{N} = \pi_{\mathcal{L}}(M)$. Then I claim that $\text{type}(\underline{N}) = \mathcal{T}$. It is clear from the definitions that $\text{type}(\underline{N}) \supseteq \mathcal{T}$. Suppose for a contradiction that we can find $U \in \mathcal{L} \setminus \mathcal{T}$ with $U \in \text{type}(\underline{N})$ that is to say $\pi_U^V M_V = 0$ for every $V \in \mathcal{L}$ strictly containing $U$. Because $U \notin \text{type}(\underline{M})$ we can find $W \subseteq S$ with $W \notin \mathcal{L}$ strictly containing $U$ such that $\pi_U^W M_W \neq 0$. Then put $T = \text{root}(W) \in \mathcal{T} \subseteq \mathcal{L}$ so that by the previous lemma $\pi_W^T M_T = M_W$. Thus $\pi_U^T M_T = \pi_U^W \pi_W^T M_T = \pi_U^W M_W \neq 0$ a contradiction. Thus $\text{type}(\underline{N}) = \mathcal{T}$ and the map is surjective. We will prove the special case in lemma 7.2.3 $\square$

## 6.2. The smoothness and irreducibility of $\overline{\mathcal{M}}_{\mathcal{L}}$

Let $\mathcal{L}$ be a thicket on $S$. In this section we show that $\overline{\mathcal{M}}_{\mathcal{L}}$ is a smooth irreducible projective variety of dimension $|S| - 2$. In particular this will mean that every non-empty Zariski open set is dense in the classical topology. Lemma 6.1.8 tells us that the data encoded by $\mathcal{T} = \text{type}(\underline{M})$ enables us to uniquely determine a point $\underline{M} \in \overline{\mathcal{M}}_{\mathcal{L}}$ from its image $\underline{N} \in \overline{\mathcal{M}}_{\mathcal{T}}$ under the projection map $\pi_{\mathcal{T}} : \overline{\mathcal{M}}_{\mathcal{L}} \to \overline{\mathcal{M}}_{\mathcal{T}}$ and is the minimal such data. This suggests that the following definitions will aid in the understanding of these spaces.

DEFINITION 6.2.1. Let $\mathcal{L}$ be a thicket and $\mathcal{T} \in \mathbb{T}_\mathcal{L}$ then we define,

$$\begin{aligned}\mathcal{N}_\mathcal{L}(\mathcal{T}) &= \{\, \underline{M} \in \overline{\mathcal{M}}_\mathcal{L} \mid \text{type}(\underline{M}) \subseteq \mathcal{T} \,\} \\ \mathcal{N}'_\mathcal{L}(\mathcal{T}) &= \{\, \underline{N} \in \overline{\mathcal{M}}_\mathcal{T} \mid \text{ for all } U \in \mathcal{L} \text{ we have } \pi_U^{\text{root}(U)} N_{\text{root}(U)} \neq 0 \,\}\end{aligned}$$

LEMMA 6.2.2. *Let $\mathcal{L}$ be a thicket, $\mathcal{T} \in \mathbb{T}_\mathcal{L}$ and $\pi_\mathcal{T} : \overline{\mathcal{M}}_\mathcal{L} \to \overline{\mathcal{M}}_\mathcal{T}$ be the usual projection map then*

(1) $\mathcal{N}'_\mathcal{L}(\mathcal{T}) = \pi_\mathcal{T}(\mathcal{N}_\mathcal{L}(\mathcal{T}))$ *and* $\mathcal{N}_\mathcal{L}(\mathcal{T}) = \pi_\mathcal{T}^{-1}(\mathcal{N}'_\mathcal{L}(\mathcal{T}))$

(2) $\mathcal{N}_\mathcal{L}(\mathcal{T})$ *and* $\mathcal{N}'_\mathcal{L}(\mathcal{T})$ *are non-empty Zariski open*

(3) $\mathcal{N}_\mathcal{L}(\mathcal{T}) \cap \mathcal{N}_\mathcal{L}(\mathcal{U}) = \mathcal{N}_\mathcal{L}(\mathcal{T} \cap \mathcal{U})$

(4) $\bigcup_{\mathcal{T} \in \mathbb{T}_\mathcal{L}} \mathcal{N}_\mathcal{L}(\mathcal{T}) = \overline{\mathcal{M}}_\mathcal{L}$ *and* $\bigcap_{\mathcal{T} \in \mathbb{T}_\mathcal{L}} \mathcal{N}_\mathcal{L}(\mathcal{T}) = \mathcal{N}_\mathcal{L}(S)$

PROOF. This result is essentially lemma 6.1.8. Let $\underline{N} \in \mathcal{N}'_\mathcal{L}(\mathcal{T})$ and define an element $\underline{M} \in \prod_{U \in \mathcal{L}} PV_U$ by $M_U = \pi_U^T N_T$ where $T = \text{root}(U)$ for each $U \in \mathcal{L}$. This makes sense since by construction each $M_U$ is non-zero. Furthermore $M_T = N_T$ for each $T \in \mathcal{T}$. Now applying lemma 6.1.8 we see that $\underline{M} \in \overline{\mathcal{M}}_\mathcal{L}$ and $\text{type}(\underline{M}) \subseteq \mathcal{T}$ thus $\underline{M} \in \mathcal{N}_\mathcal{L}(\mathcal{T})$. Then since $\pi_\mathcal{T}(\underline{M}) = \underline{N}$ we see that $\underline{N} \in \pi_\mathcal{T}(\mathcal{N}_\mathcal{L}(\mathcal{T}))$. To see the reverse inclusion let $\underline{N} \in \pi_\mathcal{T}(\mathcal{N}_\mathcal{L}(\mathcal{T}))$ then $\underline{N} = \pi_\mathcal{T}(\underline{M})$ for some $\underline{M} \in \mathcal{N}_\mathcal{L}(\mathcal{T})$. That is to say $\text{type}(\underline{M}) \subseteq \mathcal{T}$. Thus applying lemma 6.1.8 again, we see immediately that $\underline{N} \in \mathcal{N}'_\mathcal{L}(\mathcal{T})$ and we have an equality. The second part of the first statement is now immediate from the last part of lemma 6.1.8. To see the second claim for each $U \in \mathcal{L}$ define $C_U = \{\underline{N} \in \overline{\mathcal{M}}_\mathcal{T} \mid \pi_U^T N_T = 0 \text{ where } T = \text{root}(U)\,\}$ then $C_U$ is clearly Zariski closed and

$$\begin{aligned}\mathcal{N}'_\mathcal{L}(\mathcal{T}) &= \{\, \underline{N} \in \overline{\mathcal{M}}_\mathcal{T} \mid \text{ for all } U \in \mathcal{L} \text{ we have } \pi_U^{\text{root}(U)} N_{\text{root}(U)} \neq 0 \,\} \\ &= \overline{\mathcal{M}}_\mathcal{T} \setminus \coprod_{U \in \mathcal{L}} C_U\end{aligned}$$

thus $\mathcal{N}'_\mathcal{L}(\mathcal{T})$ is Zariski open. Now since $\pi_\mathcal{T} : \overline{\mathcal{M}}_\mathcal{L} \to \overline{\mathcal{M}}_\mathcal{T}$ is a morphism of varieties and $\mathcal{N}_\mathcal{L}(\mathcal{T}) = \pi_\mathcal{T}^{-1}(\mathcal{N}'_\mathcal{L}(\mathcal{T}))$ we see that $\mathcal{N}_\mathcal{L}(\mathcal{T})$ is also Zariski open. They are non-empty by lemma 6.1.9. The third and fourth parts of the claim are clear from the definitions. □

LEMMA 6.2.3. *Let $\mathcal{L}$ be a thicket, $\mathcal{T} \in \mathbb{T}_\mathcal{L}$ and $\pi_\mathcal{T} : \overline{\mathcal{M}}_\mathcal{L} \to \overline{\mathcal{M}}_\mathcal{T}$ be the usual projection map then the restricted map $\pi_\mathcal{T} : \mathcal{N}_\mathcal{L}(\mathcal{T}) \to \mathcal{N}'_\mathcal{L}(\mathcal{T})$ is a isomorphism of Zariski open varieties.*

PROOF. We first construct an injective map $\theta_\mathcal{T} : \mathcal{N}'_\mathcal{L}(\mathcal{T}) \rightarrowtail \overline{\mathcal{M}}_\mathcal{L}$ which is inverse to $\pi_\mathcal{T}$ as sets. We will then proceed further in the argument to prove it is a morphism of varieties. Given $\underline{N} \in \mathcal{N}'_\mathcal{L}(\mathcal{T})$ and any $U \in \mathcal{L}$ put $T = \text{root}(U)$ in $\mathcal{T}$ then we may define $M_U \in PV_U$ by $M_U = \rho^T_U(M_T)$ where $\rho^T_U : PV_T \dashrightarrow PV_U$ is the usual partial map. Now put $\underline{M} = \prod M_T$, then by the lemma 6.1.8 we have that $\underline{M} \in \mathcal{N}_\mathcal{L}(\mathcal{T})$ and $\pi_\mathcal{T} \circ \theta_\mathcal{T} = \text{inc}$. I claim that $\theta_\mathcal{T}$ is the inverse for the restricted map of $\pi_\mathcal{T}$. For brevity we put $\mathbb{U}_\mathcal{L} = \mathcal{N}_\mathcal{L}(\mathcal{T})$ and $\mathbb{V}_\mathcal{L} = \mathcal{N}'_\mathcal{L}(\mathcal{T})$ then it suffices to prove that $\pi_\mathcal{T}^{-1}(\mathbb{V}_\mathcal{L}) = \theta_\mathcal{T}(\mathbb{V}_\mathcal{L})$, as then $\mathcal{N}_\mathcal{L}(\mathcal{T}) \subseteq \pi_\mathcal{T}^{-1}(\mathbb{V}_\mathcal{L}) = \theta_\mathcal{T}(\mathbb{V}_\mathcal{L}) \subseteq \mathcal{N}_\mathcal{L}(\mathcal{T})$ thus $\theta_\mathcal{T}(\mathbb{V}_\mathcal{L}) = \mathcal{N}_\mathcal{L}(\mathcal{T})$. Let $\underline{M} \in \pi_\mathcal{T}^{-1}(\mathbb{V}_\mathcal{L})$ and put $\underline{N} = \pi_\mathcal{T}(\underline{M}) \in \mathbb{V}_\mathcal{L}$ then by lemma 6.1.8 there is only one element $\underline{M} \in \overline{\mathcal{M}}_\mathcal{L}$ with this property, namely $\underline{M} = \theta_\mathcal{T}(\underline{N})$ thus $\underline{M} \in \theta_\mathcal{T}(\mathbb{V}_\mathcal{L})$, the other inclusion is automatic. We must check that $\theta_\mathcal{T} : \mathbb{V}_\mathcal{L} \to \mathbb{U}_\mathcal{L}$ is a morphism of varieties. For any $T \in \mathcal{T}$ let $j^T : \mathbb{V}_\mathcal{L} \to PV_T$ be the projection map and put $X^T = j^T(\mathbb{V}_\mathcal{L})$. Then for any $U \in \mathcal{L}$ with $\text{root}(U) = T$ define $j^T_U : \mathbb{V}_\mathcal{L} \to PV_U$ by $j^T_U = \rho^T_U j^T$ where $\rho^T_U : PV_T \dashrightarrow PV_U$ is the usual partial map, then $\theta_\mathcal{T} = \prod_{U \in \mathcal{L}} j^{\text{root}(U)}_U$ restricted to $\mathbb{V}_\mathcal{L}$. For each $U \in \mathcal{L}$ with $\text{root}(U) = T$ must show that $X^T$ is contained in the regular part of $\rho^T_U$ that is $\rho^T_U : PV_T \backslash \ker(\pi^T_U) \to PV_U$ where $\pi^T_U : V_T \to V_U$ is the usual map. Clearly by lemma 6.1.8 we have $X^T \subseteq PV_T \backslash \ker(\pi^T_U)$ and $j^T_U = \rho^T_U j^T$ thus $\theta_\mathcal{T}$ is a morphism of varieties as required. □

DEFINITION 6.2.4. Let $\mathcal{L}$ be a thicket then we define $\mathcal{M}_\mathcal{L} = \pi_S^{-1}(U_S)$ and $U_\mathcal{L} = \pi_S(\mathcal{N}_\mathcal{L}(S))$ where $\pi_S : \overline{\mathcal{M}}_\mathcal{L} \to PV_S$ is the usual map.

REMARK 6.2.5. For every thicket $\mathcal{L}$ we see that $U_S \subseteq U_\mathcal{L}$ and by the last lemma we know that $\pi_S : \mathcal{N}_\mathcal{L}(S) \to U_\mathcal{L}$ is an isomorphism with $\mathcal{N}_\mathcal{L}(S) = \pi_S^{-1}(U_\mathcal{L})$ so that $\mathcal{M}_\mathcal{L} \subseteq \mathcal{N}_\mathcal{L}(S)$ and $\pi_S : \mathcal{M}_\mathcal{L} \to U_S$ is an isomorphism. In the case when $\mathcal{L}$ is all subsets $W \subseteq S$ with $|W| > 1$ we see that $\mathcal{M}_S = \mathcal{N}_\mathcal{L}(S)$ which consists of the generic curves and $U_S = U_\mathcal{L}$.

EXAMPLE 6.2.6. One can easily verify that we have a homotopy equivalence $\mathcal{M}_S \simeq \mathbb{C}_0^{|S|}/S^1$ where we take $S^1$ to be the set of complex number with length 1 and take $S^1$ to act by multiplication. Again one checks that we then get $\mathcal{M}_S \times S^1 \simeq \mathbb{C}_0^{|S|}$.

LEMMA 6.2.7. *Let $X$ be a topological space and suppose we have an open cover of $X$ so that $X = \bigcup_{i \in I} U_i$, each $U_i$ is irreducible and $\bigcap_{i \in I} U_i$ is non-empty, then $X$ is irreducible.*

PROOF. Suppose for a contradiction that $X$ is not irreducible. Then we may write $X = A \cup B$ for some proper non trivial closed subsets $A$ and $B$. Now consider any element of our cover, $U$ say, then since $U$ is also irreducible we must have $U$ contained in $A$ or $U$ contained in $B$. Without loss of generality we may suppose $U \subseteq A$. Now consider any other element $V$ of our cover then I claim $V$ is also contained in $A$. Since $U \cap V$ is non-empty and open in the irreducible set $U$ we have $\mathrm{cl}(U \cap V) = \mathrm{cl}(U) \subseteq A$. But also $\mathrm{cl}(U \cap V) = \mathrm{cl}(V)$ thus $V \subseteq \mathrm{cl}(V) = \mathrm{cl}(U) \subseteq A$ and the claim is true. Then $X = \bigcup_{i \in I} U_i \subseteq A$. As $A$ is proper this provides a contradiction, thus our claim is true.

□

COROLLARY 6.2.8. *Let $\mathcal{L}$ be a thicket on $S$ then $\overline{\mathcal{M}}_{\mathcal{L}}$ is a smooth irreducible projective variety of dimension $|S| - 2$ and $\mathcal{M}_{\mathcal{L}}$ is a Zariski open subset that is dense in the classical topology. In particular this is true for $\overline{\mathcal{M}}_S$.*

PROOF. For any thicket $\mathcal{L}$ and S-tree $\mathcal{T}$ of $\mathcal{L}$ we have seen in lemma 4.2.1 that $\overline{\mathcal{M}}_{\mathcal{T}}$ is a smooth irreducible projective variety of dimension $|S| - 2$ and $\mathcal{N}'_{\mathcal{L}}(\mathcal{T})$ is a Zariski open subset and thus irreducible. Now by the previous lemma we have the morphism of varieties $\pi_{\mathcal{T}} : \overline{\mathcal{M}}_{\mathcal{L}} \to \overline{\mathcal{M}}_{\mathcal{T}}$ which restricts to an isomorphism $\pi_{\mathcal{T}} : \mathcal{N}_{\mathcal{L}}(\mathcal{T}) \to \mathcal{N}'_{\mathcal{L}}(\mathcal{T})$ thus $\mathcal{N}_{\mathcal{L}}(\mathcal{T})$ is a Zariski open irreducible subset of $\overline{\mathcal{M}}_{\mathcal{L}}$. Now, the $\mathcal{N}_{\mathcal{L}}(\mathcal{T})$ form an open cover of $\overline{\mathcal{M}}_{\mathcal{L}}$ and $\bigcap_{\mathcal{T} \in \mathbb{T}_{\mathcal{L}}} \mathcal{N}_{\mathcal{L}}(\mathcal{T}) = \mathcal{N}_{\mathcal{L}}(S)$ therefore by the last lemma $\overline{\mathcal{M}}_{\mathcal{L}}$ is irreducible. It is smooth of dimension $|S| - 2$ as this is true for the open cover under the isomorphisms. For the final part of the claim as $\overline{\mathcal{M}}_{\mathcal{L}}$ is irreducible we see that $\mathcal{M}_{\mathcal{L}}$ is dense in the Zariski topology but since $\overline{\mathcal{M}}_{\mathcal{L}}$ is smooth we see that it is dense in the classical topology.

□

## 6.3. The tangent bundle to $\overline{\mathcal{M}}_{\mathcal{L}}$

We next proceed to compute the tangent bundle $\sigma_{\mathcal{L}}$ of $\overline{\mathcal{M}}_{\mathcal{L}}$ for when $\mathcal{L}$ is a thicket. We will use the embedding $i : \overline{\mathcal{M}}_{\mathcal{L}} \to \prod_{T \in \mathcal{L}} PV_T$ to identify it as a sub-bundle of $T \prod PV_T$ where $T_{\underline{M}} \prod PV_T = \prod \hom(M_T, V_T/M_T)$. Given our calculation in section 4.3 this reduces to showing that the dimension of a certain vector space which is a natural extension to that defined in 4.3.7 is $|S| - 2$. This space we define below.

DEFINITION 6.3.1. Let $\underline{M} \in \overline{\mathcal{M}}_\mathcal{L}$ with tree type $\mathcal{T}$ and $T \in \mathcal{T}$. Suppose $M(\mathcal{T}, T)$ is non-empty then let $\pi_{\mathcal{T},T} : V_T \to \bigoplus_{U \in M(\mathcal{T},T)} V_U$ be the usual map. Then we define the vector space by $W_{\underline{M},T} = \pi_{\mathcal{T},T}^{-1}(\oplus M_U)$. If $M(\mathcal{T}, T)$ is empty we define $W_{\underline{M},T} = V_T$.

DEFINITION 6.3.2. Let $\underline{M}$ be a point in $\overline{\mathcal{M}}_\mathcal{L}$ and $\pi_U^T : V_T \to V_U$ be the usual map. Then $\pi_U^T(M_T) \leq M_U$, so we have induced maps $\pi_U^T : M_T \to M_U$ and $\overline{\pi}_U^T : V_T/M_T \to V_U/M_U$ and define $\sigma_{\mathcal{L},\underline{M}}$ by

$$\sigma_{\mathcal{L},\underline{M}} = \left\{ \underline{\alpha} \in \prod_{T \in \mathcal{L}} \hom(M_T, V_T/M_T) \;\middle|\; \overline{\pi}_U^T \alpha_T = \alpha_U \pi_U^T \text{ for all } U \subseteq T \text{ with } U, T \in \mathcal{L} \right\}$$

We then define $\sigma_\mathcal{L} = \coprod_{\underline{M}} \sigma_{\mathcal{L},\underline{M}}$. Note that the condition in the braces is just the requirement that the following diagram is commutative. For all $U \subseteq T$ with $U, T \in \mathcal{L}$

$$\begin{array}{ccc} M_T & \xrightarrow{\alpha_T} & V_T/M_T \\ \pi_U^T \downarrow & & \downarrow \overline{\pi}_U^T \\ M_U & \xrightarrow{\alpha_U} & V_U/M_U \end{array}$$

CONSTRUCTION 6.3.3. Let $\mathcal{L}$ be a thicket on $S$ and $\underline{M} \in \overline{\mathcal{M}}_\mathcal{L}$. Put $\mathcal{T} = \text{type}(\underline{M})$. Let $U \in \mathcal{L}$ and put $T = \text{root}(U)$ in $\mathcal{T}$ then we define the surjective restriction maps $\theta_U^T : \hom(M_T, V_T/M_T) \to \hom(M_U, V_U/M_U)$ as follows. Since $T = \text{root}(U)$ we have $\pi_U^T M_T = M_U$ and thus the restricted map $\pi_U^T : M_T \to M_U$ is an isomorphism. We also have the induced quotient maps $\overline{\pi}_U^T : V_T/M_T \to V_U/M_U$. Then given an element $\alpha_T \in \hom(M_T, V_T/M_T)$ we may define the element $\theta_U^T(\alpha_T) \in \hom(M_U, V_U/M_U)$ by $\theta_U^T(\alpha_T) = \overline{\pi}_U^T \alpha_T (\pi_U^T)^{-1}$. We then define a map $\theta : \prod_{T \in \mathcal{T}} \hom(M_T, V_T/M_T) \to \prod_{U \in \mathcal{L}} \hom(M_U, V_U/M_U)$ by $\theta = \prod_{U \in \mathcal{L}} \theta_U^{\text{root}(U)}$.

LEMMA 6.3.4. *Let $\mathcal{L}$ be a thicket and $\underline{M} \in \overline{\mathcal{M}}_\mathcal{L}$. Put $\mathcal{T} = \text{type}(\underline{M})$ and write $\underline{N}$ to be the image of $\underline{M}$ in $\overline{\mathcal{M}}_\mathcal{T}$ then*

(1) *The projection $\pi : \sigma_{\mathcal{L},\underline{M}} \to \sigma_{\mathcal{T},\underline{N}}$ is an isomorphism of vector spaces.*
(2) $\sigma_{\mathcal{T},\underline{N}} = \prod_{T \in \mathcal{T}} \hom(M_U, W_{\underline{M},T}/M_U)$
(3) $\dim(\sigma_{\mathcal{T},\underline{N}}) = |S| - 2$
(4) *The map $\theta$ restricts to $\theta : \sigma_{\mathcal{T},\underline{N}} \to \sigma_{\mathcal{L},\underline{M}}$ which is inverse to $\pi$*

PROOF. It is clear that $\sigma_{\mathcal{L},\underline{M}}$ is a vector subspace of $\prod_{T \in \mathcal{L}} \hom(M_T, V_T/M_T)$. We next examine how the commutative diagrams restricts the functions $\alpha_U, \alpha_T$ for $U \subseteq T$ that we are interested in.

Let $\underline{M} \in \overline{\mathcal{M}}_\mathcal{L}$ and $U, T \in \mathcal{L}$ with $U \subseteq T \subseteq S$. Then by construction either $\pi_U^T(M_T) = M_U$ or $\pi_U^T(M_T) = 0$. In the first case we have that $\pi_U^T : M_T \to M_U$ is an isomorphism and so given $\alpha_T \in \hom(M_T, V_T/M_T)$ there exists a unique $\alpha_U \in \hom(M_U, V_U/M_U)$ with the required property, namely $\alpha_U = \overline{\pi}_U^T \alpha_T (\pi_U^T)^{-1} = \theta_U^T(\alpha_T)$. In the second case for any $\alpha_U \in \hom(M_U, V_U/M_U)$ we have $\alpha_U \pi_U^T = 0$ and so by the commutative diagram we must have $\overline{\pi}_U^T \alpha_T = 0$ that is we require $\text{im}(\alpha_T) \le \ker(\overline{\pi}_U^T)$. It is clear that $\ker(\overline{\pi}_U^T) = K_U^T/M_T$ where $K_U^T = (\pi_U^T)^{-1}(M_U)$ and so for any $\alpha_U \in \hom(M_U, V_U/M_U)$ and $\alpha_T \in \hom(M_T, V_T/M_T)$ the required diagram commutes if and only if $\alpha_T \in \hom(M_T, K_U^T/M_T)$ where we are using the obvious notation.

Let $\underline{M} \in \overline{\mathcal{M}}_\mathcal{L}$ with $\mathcal{T} = \text{type}(\underline{M})$, then given any $T \in \mathcal{T}$ and $U \in \mathcal{L}$ with $U \subseteq T$ it is clear that $\pi_U^T(M_T) = 0$ if and only if $U \subseteq W$ for some $W \in M(\mathcal{T}, T)$. Put $n_T = |M(\mathcal{T}, T)|$ and suppose $n_T > 0$, then given $\alpha_T \in \hom(M_T, V_T/M_T)$ we must have that $\alpha_T \in \hom(M_T, (\pi_U^T)^{-1}M_U/M_T)$ for all $U \in M(\mathcal{T}, T)$. That is to say $\alpha_T \in \hom(M_T, W_{\underline{M},T}/M_T)$ where $W_{\underline{M},T} = \bigcap_{U \in M(\mathcal{T},T)} (\pi_U^T)^{-1}M_U \le V_T$. If $n_T = 0$ then we take $W_{\underline{M},T} = V_T$. Next put $\sigma'_{\mathcal{T},\underline{N}} = \prod_{T \in \mathcal{T}} \hom(M_T, W_{\underline{M},T}/M_T)$, this proves we have a projection map $\pi : \sigma_{\mathcal{L},\underline{M}} \to \sigma'_{\mathcal{T},\underline{N}}$. We next construct the inverse. For any $T \in \mathcal{T}$ let $\alpha_T \in \hom(M_T, W_{\underline{M},T}/M_T)$. Let $U \in \mathcal{L}$ with $T = \text{root}(U)$ then by construction $\pi_U^T(M_T) = M_U$ and we define $\alpha_U$ to be the only choice possible from this commutative diagram, that is $\alpha_U = \overline{\pi}_U^T \alpha_T (\pi_U^T)^{-1} = \theta_U^T(\alpha_T)$. We first need to check that the $\underline{\alpha}$ constructed in this way does indeed lie in $\sigma_{\mathcal{L},\underline{M}}$. Let $U \subseteq V \subseteq S$ and put $T = \text{root}(U), W = \text{root}(V)$, then $T \subseteq W$ and we have $\pi_U^T(M_T) = M_U, \pi_V^W M_W = M_V$. First suppose that $T = W$ then given a linear map $\alpha_T$ we are forced to take $\alpha_U$ and $\alpha_V$ to be the unique function satisfying the conditions that $\alpha_U \pi_U^T = \overline{\pi}_U^T \alpha_T$ and $\alpha_V \pi_V^T = \overline{\pi}_V^T \alpha_T$. We are required to prove that, $\alpha_U \pi_U^V = \overline{\pi}_U^V \alpha_V$. We have that $\pi_V^T M_T = M_V$, $\pi_U^V \pi_V^T = \pi_U^T$ and the map $\pi_V^T$ is invertible so $\pi_U^V = \pi_U^T (\pi_V^T)^{-1}$ also $\overline{\pi}_U^T = \overline{\pi}_U^V \overline{\pi}_V^T$. We then have,

$$\begin{aligned}
\alpha_U \pi_U^V &= \alpha_U \pi_U^T (\pi_V^T)^{-1} \\
&= \overline{\pi}_U^T \alpha_T (\pi_V^T)^{-1} \\
&= \overline{\pi}_U^V \overline{\pi}_V^T \alpha_T (\pi_V^T)^{-1} \\
&= \overline{\pi}_U^V \alpha_V
\end{aligned}$$

Now suppose $T \subset W$ then we can find an $X \in M(\mathcal{T}, W)$ with $T \subseteq X$. From this we see that $\pi_T^W M_W = 0$ and therefore $\overline{\pi}_T^W \alpha_W = 0$. Now $\alpha_V$ is the unique function satisfying $\alpha_V \pi_V^W = \overline{\pi}_V^W \alpha_W$ and since $\pi_U^V M_V = 0$ we are required to prove that $\overline{\pi}_U^V \alpha_V = 0$

$$\begin{aligned}
\overline{\pi}_U^V \alpha_V &= \overline{\pi}_U^V \overline{\pi}_V^W \alpha_W (\pi_V^W)^{-1} \\
&= \overline{\pi}_U^W \alpha_W (\pi_V^W)^{-1} \\
&= \overline{\pi}_U^T (\overline{\pi}_T^W \alpha_W)(\pi_V^W)^{-1} \\
&= 0
\end{aligned}$$

This proves that $\underline{\alpha}$ is a point in $\sigma_{\mathcal{L},\underline{M}}$, thus we have the injective map $\theta : \sigma'_{\mathcal{T},\underline{N}} \to \sigma_{\mathcal{L},\underline{M}}$. This is clearly the inverse map. Next put $H_T = \hom(M_T, W_{\underline{M},T}/M_T)$ and $d_T = \dim_{\mathbb{C}}(H_T)$. We are required to prove that $\dim_{\mathbb{C}}(\sigma_{\mathcal{L},\underline{M}}) = \sum_{T \in \mathcal{T}} d_T = |S| - 2$. For each $T \in \mathcal{T}$, $d_T = \dim_{\mathbb{C}} W_{\underline{M},T} - 1$. If $n_T > 0$ let $\pi_{\mathcal{T},T} : V_T \to \bigoplus_{U \in M(\mathcal{T},T)} V_U$ be the usual map, then we have $W_{\underline{M},T} = \pi_{\mathcal{T},T}^{-1}(\oplus M_T)$ and by corollary 4.1.5 we have $\dim_{\mathbb{C}} W_{\underline{M},T} = n(\mathcal{T},T)$ and so $d_T = n(\mathcal{T},T) - 1$. If $n_T = 0$ we see that $W_{\underline{M},T} = V_T$ and so in all cases $d_T = n(\mathcal{T},T) - 1$. Therefore,

$$\begin{aligned}
\dim_{\mathbb{C}}(\sigma_{\mathcal{L},\underline{M}}) &= \sum_{T \in \mathcal{T}}(n(\mathcal{T},T) - 1) \\
&= |S| - 2 \text{ by lemma 3.3.2}
\end{aligned}$$

Next we have the projection map $\pi : \sigma_{\mathcal{L},\underline{M}} \to \sigma_{\mathcal{T},\underline{N}}$ and the map $\theta : \sigma'_{\mathcal{T},\underline{N}} \to \sigma_{\mathcal{L},\underline{M}}$ thus composing we see that $\sigma'_{\mathcal{T},\underline{N}} \subseteq \sigma_{\mathcal{T},\underline{N}}$. It is also clear from our calculations that we have the projection $p : \sigma_{\mathcal{T},\underline{N}} \to \sigma'_{\mathcal{T},\underline{N}}$. This shows $\sigma_{\mathcal{T},\underline{N}} \subseteq \sigma'_{\mathcal{T},\underline{N}}$ thus $\sigma_{\mathcal{T},\underline{N}} = \sigma'_{\mathcal{T},\underline{N}}$ and our claims are proven. □

LEMMA 6.3.5. *Let $\mathcal{L}$ be a thicket then the tangent bundle of $\overline{\mathcal{M}}_{\mathcal{L}}$ is $\sigma_{\mathcal{L}}$*

PROOF. For $|\mathcal{L}| \leq 2$ we see that $\mathcal{L}$ is an $S$-tree and we have already proven this claim in lemma 4.3.12 so we may suppose $|\mathcal{L}| > 2$. Let $\underline{M} \in \overline{\mathcal{M}}_{\mathcal{L}}$ and $U, T \in \mathcal{L}$ with $U \subset T$. Put $\mathcal{U} = \{U, T\}$ and consider the following commutative diagram,

$$\begin{array}{ccc} \overline{\mathcal{M}}_{\mathcal{L}} & \xrightarrow{i} & \prod_{T \in \mathcal{L}} PV_T \\ \pi \downarrow & & \downarrow p \\ \overline{\mathcal{M}}_{\mathcal{U}} & \xrightarrow{j} & \prod_{T \in \mathcal{U}} PV_T \end{array}$$

Put $\underline{N}$ to be the image of $\underline{M}$ in $\overline{\mathcal{M}}_{\mathcal{U}}$. Then the tangent space $T_{\underline{N}}\overline{\mathcal{M}}_{\mathcal{U}} = \sigma_{\mathcal{U},\underline{N}}$ thus we may identify the tangent space $T_{\underline{M}}\overline{\mathcal{M}}_{\mathcal{L}} \leq \sigma_{\mathcal{L},\underline{M}}$ but $\dim(T_{\underline{M}}\overline{\mathcal{M}}_{\mathcal{L}}) = |S| - 2 = \dim(\sigma_{\mathcal{L},\underline{M}})$ by the previous lemma thus $T\overline{\mathcal{M}}_{\mathcal{L}} = \sigma_{\mathcal{L}}$.

□

COROLLARY 6.3.6. *For any finite set $S$, $\overline{\mathcal{M}}_S$ is a smooth irreducible projective variety of dimension $|S| - 2$ with tangent bundle $\sigma_S$.*

## 6.4. Results on the cohomology ring of $\overline{\mathcal{M}}_{\mathcal{L}}$

In this section we will produce a certain pullback diagram that will give rise to a Mayer Vietoris type sequence. This diagram will be one useful way of deducing many pleasant results about the cohomology rings of $\overline{\mathcal{M}}_{\mathcal{L}}$ and in particular $\overline{\mathcal{M}}_S$. In the following section we will use the same approach to analyze the Chow ring. This diagram will also turn out to have other nice properties that we discuss in the last section.

DEFINITION 6.4.1. Let $\mathcal{L}$ be a thicket on $S$. Then for each $T \in \mathcal{L}$ we define the elements $x_T \in H^2(\overline{\mathcal{M}}_{\mathcal{L}})$ by $x_T = \pi_T^*(y_T)$ where $\pi_T : \overline{\mathcal{M}}_{\mathcal{L}} \to PV_T$ is the projection map and $y_T$ is the standard generator of $H^2(PV_T)$, that is $y_T = e(L_T)$ where $L_T$ is the tautological line bundle over $PV_T$.

DEFINITION 6.4.2. We say a triple $(\mathcal{L}, S, T)$ is *admissible* if $\mathcal{L} \subseteq P(S)$, $\mathcal{L}$ is a thicket, $T \subseteq S$ and there is no $U \in \mathcal{L}$ with $U \subseteq T$. Further if $\mathcal{L}_+ = \mathcal{L} \amalg \{T\}$ is a thicket we call the admissible triple an *admissible thicket*.

CONSTRUCTION 6.4.3. Let $T \subseteq S$ and $\sim$ be an equivalence relation on $T$. Write $\overline{T}$ for the set of equivalence classes of $T$ and let $q_T : T \to \overline{T}$ be the quotient map then we define a map $r_T : PV_{\overline{T}} \to PV_T$ as follows. We first define a map $r_T : F(\overline{T}, \mathbb{C}) \to F(T, \mathbb{C})$. For any $f_{\overline{T}} \in F(\overline{T}, \mathbb{C})$ we define $f_T \in F(T, \mathbb{C})$ by $f_T = f_{\overline{T}} \circ q_T$. Then $r_T$ send constants to constants so we can define the induced injective map $r_T : V_{\overline{T}} \to V_T$ and in turn an injective map $r_T : PV_{\overline{T}} \to PV_T$, this is the first of our desired maps. One readily checks we have a short exact sequence $V_{\overline{T}} \to V_T \to \bigoplus_{U \in \overline{T}} V_U$ and we define $\overline{V}_T = \ker(\pi_{\overline{T}})$. Thus we have isomorphisms $s_T : \overline{V}_T \to V_{\overline{T}}$ induced from this sequence and $s_T : P\overline{V}_T \to PV_{\overline{T}}$. The latter is the second of our desired maps.

LEMMA 6.4.4. *Let $T \subseteq S$ and $\sim$ be a equivalence relation on $T$. Then with the notation as above we have maps $r_T : PV_{\overline{T}} \to PV_T$ and $s_T : P\overline{V}_T \to PV_{\overline{T}}$. We also write $r_T : PV_{\overline{T}} \to P\overline{V}_T$ to be the restricted map of $r_T$ onto its image. Then $r_T$ and $s_T$ are inverses for each other. Let $U$ be a subset of $T$ and give $U$ the induced equivalence relation then we have the following commutative diagrams,*

$$\begin{array}{ccc} V_{\overline{T}} \xrightarrow{r_T} \overline{V}_T & & \overline{V}_T \xrightarrow{s_T} V_{\overline{T}} \\ \pi_{\overline{U}}^{\overline{T}} \downarrow \quad \downarrow \pi_U^T & & \pi_U^T \downarrow \quad \downarrow \pi_{\overline{U}}^{\overline{T}} \\ V_{\overline{U}} \xrightarrow{r_U} \overline{V}_U & & \overline{V}_U \xrightarrow{s_U} V_{\overline{U}} \end{array}$$

*In particular if $|\overline{U}| = 1$ then we obtain $\pi_U^T r_T = 0$*

PROOF. The proofs of these claims are clear.

$\square$

DEFINITION 6.4.5. Let $T \subseteq S$. Then in this section we make extensive use of the equivalence relation on $S$ by $T$ defined by $u \sim v$ if and only if $u = v$ or $u, v \in T$. We write $S/T$ for the equivalence classes and $q_S : S \to S/T$ for the quotient map. Let $U \subseteq S$ and put $\overline{U} = q_S(U)$ then we write $q_U : U \to \overline{U}$ for the induced quotient map on $U$.

CONSTRUCTION 6.4.6. Let $(\mathcal{L}, S, T)$ be admissible then we define a map $j_\mathcal{L} : \overline{\mathcal{M}}_{\overline{\mathcal{L}}} \hookrightarrow \overline{\mathcal{M}}_\mathcal{L}$ as follows. Given any set $U \in \mathcal{L}$ we write $\overline{U} = q_T(U)$ where $q_T : S \to S/T$ is the quotient map then $|\overline{U}| > 1$ because $U \not\subseteq T$. We define $M_U = r_U(N_{\overline{U}})$ and $j_\mathcal{L}(\underline{N}) = \prod r_U(N_{\overline{U}})$. To see that the image of $j_\mathcal{L}$ lies in $\overline{\mathcal{M}}_\mathcal{L}$ we observe that by the last lemma we have the following commutative diagram. Let $U, W \in \mathcal{L}$ with $U \subseteq W$ and $\overline{U}, \overline{W}$ be their images in $\overline{\mathcal{L}}$ then,

$$\begin{array}{ccc} V_{\overline{W}} & \xrightarrow{r_W} & V_W \\ \pi^{\overline{W}}_{\overline{U}} \downarrow & & \downarrow \pi^W_U \\ V_{\overline{U}} & \xrightarrow{r_U} & V_U \end{array}$$

LEMMA 6.4.7. Let $(\mathcal{L}, S, T)$ be admissible then

$$j_\mathcal{L}(\overline{\mathcal{M}}_{\overline{\mathcal{L}}}) = \{\, \underline{M} \in \overline{\mathcal{M}}_\mathcal{L} \mid \text{ for all } U \in \mathcal{L} \text{ with } U \cap T \text{ non-empty} \implies \pi^U_{U \cap T} M_U = 0 \,\}$$

PROOF. Let $\underline{M} \in \overline{\mathcal{M}}_\mathcal{L}$ and $U \in \mathcal{L}$ with $U \cap T$ non-empty. Then because $\pi^U_{U \cap T} M_U = 0$ we see that $M_U \in P(\overline{V}_U)$. Then given $U \in \mathcal{L}$ we may define an element $N_{\overline{U}} \in PV_{\overline{U}}$ by $N_{\overline{U}} = s_U(M_U)$. We must check this is well defined. Let $W \in \mathcal{L}$ be another element with $q_T(W) = q_T(U)$. Then by the first part we readily deduce that $\text{root}(U) = \text{root}(W)$ in $\text{type}(\underline{M})$ so that we are well defined. Put $\underline{N} = \prod N_{\overline{U}}$ then it is clear from the commuting diagram below of 6.4.4 that $\underline{N} \in \overline{\mathcal{M}}_{\overline{\mathcal{L}}}$ and $\underline{M} = j_\mathcal{L}(\underline{N})$

$$\begin{array}{ccc} \overline{V}_T & \xrightarrow{s_T} & V_{\overline{T}} \\ \pi^T_U \downarrow & & \downarrow \pi^{\overline{T}}_{\overline{U}} \\ \overline{V}_U & \xrightarrow{s_U} & V_{\overline{U}} \end{array}$$

The reverse inclusion is clear as we have the short exact sequence $V_{U/T} \to V_U \to V_{U \cap T}$. $\square$

CONSTRUCTION 6.4.8. Let $(\mathcal{L}, S, T)$ be admissible. We then define a map $i_\mathcal{L} : PV_T \times \overline{\mathcal{M}}_\mathcal{L} \to \overline{\mathcal{M}}_{\mathcal{L}_+}$ by $i_\mathcal{L}(M_T, \underline{M}) = (M_T, j_\mathcal{L}(\underline{M}))$. Because $T$ has minimal size we see by the last lemma that this does indeed lie in $\overline{\mathcal{M}}_{\mathcal{L}_+}$

The following lemma is a refinement of the previous one that works for admissible thickets $(\mathcal{L}, S, T)$. This will be important for several results we will require later.

LEMMA 6.4.9. *Let $(\mathcal{L}, S, T)$ be an admissible thicket then*

$$j_\mathcal{L}(\overline{\mathcal{M}}_{\overline{\mathcal{L}}}) = \{\, \underline{M} \in \overline{\mathcal{M}}_\mathcal{L} \mid \text{ for all } U \in \mathcal{L} \text{ with } U \supseteq T \implies \pi_T^U M_U = 0 \,\}$$

PROOF. Suppose we have $\underline{M} \in \overline{\mathcal{M}}_\mathcal{L}$ and for all $U \in \mathcal{L}$ with $U \supset T$ we have $\pi_T^U M_U = 0$. Then for each $U \in \mathcal{L}$ with $U \cap T$ non-empty. Put $V = U \cup T \supset T$ and $W = U \cap T$ then as $\mathcal{L}_+$ is a thicket we have $V \in \mathcal{L}$ and $\pi_T^V M_V = 0$. Then by lemma 6.1.2 we must have $\pi_U^V M_V = M_U$ as $U$ and $T$ have non-trivial intersection. Then $\pi_W^U M_U = \pi_W^U \pi_U^V M_V = \pi_W^V M_V = 0$ because $W \subseteq T$ and $\pi_T^V M_V = 0$. Let $q_T : S \to S/T$ be the usual quotient map. For each $U \in \mathcal{L}$ put $\overline{U} = q_T(U)$ then because $M_U \in P\overline{V}_U$ we may define $N_{\overline{U}} = s_U(M_U)$ and $\underline{N} = \prod N_{\overline{U}}$, to see this is well defined let $W \in \mathcal{L}$ with $q_T(U) = q_T(W)$ then by the first part we deduce that $\text{root}(U) = \text{root}(W)$ and we are well defined. Next I claim $\underline{N} \in \overline{\mathcal{M}}_{\overline{\mathcal{L}}}$, this is immediate from the previous commutative diagram of lemma 6.4.4

Then it is clear $j_\mathcal{L}(\underline{N}) = \underline{M}$ therefore $\underline{M} \in j_\mathcal{L}(\overline{\mathcal{M}}_{\overline{\mathcal{L}}})$. The reverse inclusion of this result is clear because of the short exact sequence $V_{U/T} \to V_U \to V_T$.

$\square$

PROPOSITION 6.4.10. *Let $(\mathcal{L}, S, T)$ be an admissible thicket. Then we have the following pullback diagram,*

$$\begin{array}{ccc} \overline{\mathcal{M}}_{\overline{\mathcal{L}}} \times PV_T & \xrightarrow{i_\mathcal{L}} & \overline{\mathcal{M}}_{\mathcal{L}_+} \\ {\scriptstyle p} \downarrow & & \downarrow {\scriptstyle \pi} \\ \overline{\mathcal{M}}_{\overline{\mathcal{L}}} & \xrightarrow{j_\mathcal{L}} & \overline{\mathcal{M}}_\mathcal{L} \end{array}$$

$i_{\mathcal{L}}$ and $j_{\mathcal{L}}$ are injective closed maps and $\pi : (\overline{\mathcal{M}}_{\mathcal{L}_+}, A_{\mathcal{L}}) \to (\overline{\mathcal{M}}_{\mathcal{L}}, B_{\mathcal{L}})$ is a relative isomorphism where $A_{\mathcal{L}} = i_{\mathcal{L}}(\overline{\mathcal{M}}_{\overline{\mathcal{L}}} \times PV_T)$ and $B_{\mathcal{L}} = j_{\mathcal{L}}(\overline{\mathcal{M}}_{\overline{\mathcal{L}}})$ are Zariski closed moreover the map $\pi$ is surjective and the square is a pushout.

PROOF. The commutativity of this diagram is clear from the definitions of the various maps. We next prove that the map $j_{\mathcal{L}} : \overline{\mathcal{M}}_{\overline{\mathcal{L}}} \hookrightarrow \overline{\mathcal{M}}_{\mathcal{L}}$ is injective. Suppose we are given $\underline{N}, \underline{M} \in \overline{\mathcal{M}}_{\overline{\mathcal{L}}}$ with $j_{\mathcal{L}}(\underline{N}) = \underline{L} = j_{\mathcal{L}}(\underline{M})$, then for any $U \in \mathcal{L}$ we have $L_U = r_U(N_{\overline{U}}) = r_U(M_{\overline{U}})$. Since $r_U$ is injective we must have $M_{\overline{U}} = N_{\overline{U}}$ and therefore $\underline{N} = \underline{M}$. From this it follows that the map $i_{\mathcal{L}}$ is also injective.

We next prove that the map of pairs $\pi : (\overline{\mathcal{M}}_{\mathcal{L}_+}, A_{\mathcal{L}}) \to (\overline{\mathcal{M}}_{\mathcal{L}}, B_{\mathcal{L}})$ is a relative isomorphism. For this the following result will be useful. We have that,

$$|\pi^{-1}(\underline{M})| = 1 \iff \text{there exists } U \in \mathcal{L} \text{ with } U \supset T \text{ such that } \pi_T^U M_U \neq 0$$

To see the first direction suppose there is no such $U$ then by lemma 6.4.9 we see $\underline{M} \in \text{image}(j_{\mathcal{L}})$, and by the commutativity of the diagram $\pi^{-1}(\underline{M})$ contains a copy of $PV_T$ therefore $|\pi^{-1}(\underline{M})| = \infty$. To see the converse choose a set $U \supset T$ of minimal size with $\pi_T^U M_U \neq 0$. Write $M_T = \pi_T^U M_U$. Let $V \in \mathcal{L}$ be any element with $V \supset T$ then we need to show that $\pi_T^V M_V \leq M_T$. Put $W = U \cup V$ then $U \cap V \supseteq T$ therefore $W \in \mathcal{L}$ as $\mathcal{L}$ is a thicket. Clearly we can't have both $\pi_U^W M_W = 0$ and $\pi_V^W M_W = 0$. First suppose that $\pi_U^W M_W \neq 0$, then $\pi_U^W M_W = M_U$. Because $\pi_T^U M_U = M_T$ we have that $\pi_T^W M_W = M_T$. Therefore $\pi_V^W M_W = M_V$ as $V \supset T$ and $\pi_T^V M_V = \pi_T^V \pi_V^W M_W = \pi_T^W M_W = M_T$. Next suppose $\pi_U^W M_W = 0$ then $\pi_T^W M_W = 0$ as $T \subset U$ and $\pi_V^W M_W = M_V$ and therefore $\pi_T^V M_V = \pi_T^V \pi_V^W M_W = \pi_T^W M_W = 0$. Then in all cases we have $\pi_T^V M_V \leq M_T$.

By the previous lemma we have for any $\underline{M} \in \overline{\mathcal{M}}_{\mathcal{L}}$ then $\underline{M} \in j_{\mathcal{L}}(\overline{\mathcal{M}}_{\overline{\mathcal{L}}})$ if and only if for every $U \in \mathcal{L}$ with $U \supset T$ we have $\pi_T^U M_U = 0$. Thus this proves our map $\pi : (\overline{\mathcal{M}}_{\mathcal{L}_+}, A_{\mathcal{L}}) \to (\overline{\mathcal{M}}_{\mathcal{L}}, B_{\mathcal{L}})$ is a bijection. We now construct the inverse over our relative spaces for the function $\pi$. Write $A_{\mathcal{L}}^c$ for the complement of $A_{\mathcal{L}}$ in $\overline{\mathcal{M}}_{\mathcal{L}}$ and $B_{\mathcal{L}}^c$ for the complement of $B_{\mathcal{L}}$ in $\overline{\mathcal{M}}_{\mathcal{L}_+}$. Then as $\overline{\mathcal{M}}_{\mathcal{L}}$ is a complete variety for every $\mathcal{L}$ we see that $A_{\mathcal{L}}$ and $B_{\mathcal{L}}$ are Zariski closed. Next put $\mathcal{U}_T = \{\, U \in \mathcal{L} \mid U \supset T \,\}$ and define a function $\mu_T : B_{\mathcal{L}}^c \to PV_T$ by $\mu_T(\underline{L}) = \bigcup_{U \in \mathcal{U}_T} \pi_T^U(M_U)$. We then define $\theta : B_{\mathcal{L}}^c \to A_{\mathcal{L}}^c$ by $\theta(\underline{N}) = (\underline{N}, \mu_T(\underline{N}))$. This is clearly the inverse map, we must show that $\mu_T$ is a morphism of varieties. For each tree $\mathcal{T} \in \mathbb{T}_\mathcal{L}$ put $\mathcal{O}_\mathcal{L}(\mathcal{T}) = \mathcal{N}_\mathcal{L}(\mathcal{T}) \cap B_\mathcal{L}^c$. Then each $\mathcal{O}_\mathcal{L}(\mathcal{T})$ is Zariski open and $B_\mathcal{L}^c = \bigcup_{\mathcal{T} \in \mathbb{T}_\mathcal{L}} \mathcal{O}_\mathcal{L}(\mathcal{T})$. For each $\mathcal{T} \in \mathbb{T}_\mathcal{L}$ put $U = \text{root}(T)$ in $\mathcal{T}$. Then I claim the partial map $\rho_T^U : PV_U \dashrightarrow PV_T$ is defined at $M_U$. Note that this is not a consequence of lemma 6.1.8 as $T \notin \mathcal{L}$ and we are extending the definition of the root function for this case (in the evident natural way) to the element $T$. For suppose otherwise that $\pi_T^U M_U = 0$. Let $V \in \mathcal{L}$ with $V \supset T$ then $\text{root}(V) \supseteq \text{root}(T) = U$. Put $W = \text{root}(V)$ so that $\pi_V^W M_W = M_V$ thus $\pi_T^V M_V = \pi_T^V \pi_V^W M_W = \pi_T^W M_W = \pi_T^U \pi_U^W M_W \leq \pi_T^U M_U = 0$ contrary to the assumption that $\underline{M} \in B_\mathcal{L}^c$. Thus the restricted map of $\mu_T$ to $\mathcal{O}_\mathcal{L}(\mathcal{T})$ is given by $\mu_T(\underline{M}) = \rho_T^U(M_U)$ where $\rho_T^U : PV_U \setminus P(\ker \pi_T^U) \to PV_T$ and is therefore a regular morphism of varieties.

Finally we prove that the diagram is a pushout. Let $W_\mathcal{L}$ be the pushout of this diagram and $p : W_\mathcal{L} \to \overline{\mathcal{M}}_\mathcal{L}$ be the evident map. Then since $W_\mathcal{L}$ is compact and $\overline{\mathcal{M}}_\mathcal{L}$ is Hausdorph it will be enough to prove that $p$ is a bijection. This is then clear since $j_\mathcal{L}$ is injective $\pi$ is surjective and $\pi : (\overline{\mathcal{M}}_{\mathcal{L}_+}, A_\mathcal{L}) \to (\overline{\mathcal{M}}_\mathcal{L}, B_\mathcal{L})$ is a relative bijection.

$\square$

PROPOSITION 6.4.11. *Let $(\mathcal{L}, S, T)$ be an admissible thicket, then we have the long exact sequence,*

$$H^*(\overline{\mathcal{M}}_\mathcal{L}) \xrightarrow{f_\mathcal{L}^*} H^*(\overline{\mathcal{M}}_{\overline{\mathcal{L}}}) \oplus H^*(\overline{\mathcal{M}}_{\mathcal{L}_+}) \xrightarrow{g_\mathcal{L}^*} H^*(\overline{\mathcal{M}}_{\overline{\mathcal{L}}} \times PV_T) \xrightarrow{\delta_\mathcal{L}} H^{*+1}(\overline{\mathcal{M}}_\mathcal{L}).$$

*where $f_\mathcal{L}^* = (j_\mathcal{L}^*, \pi^*)$ and $g_\mathcal{L}^* = i_\mathcal{L}^* - p^*$*

PROOF. To prove this claim it will be enough to show that the diagram of Proposition 6.4.10 is a homotopy pushout. Since $i_\mathcal{L}$ is an embedding of smooth compact manifolds we see it is a cofibration and by lemma 6.4.10 the diagram is a pushout thus the diagram is a homotopy push out and we obtain the specified Mayer-Vietoris type sequence.

$\square$

REMARK 6.4.12. We will see later that the odd cohomology is zero and all the modules are free. Therefore the exact sequence decouples and we will obtain a series of split exact sequences in even degrees.

LEMMA 6.4.13. *Let $i_\mathcal{L} : \overline{\mathcal{M}}_{\overline{\mathcal{L}}} \times PV_T \to \overline{\mathcal{M}}_{\mathcal{L}_+}$ and $j_\mathcal{L} : \overline{\mathcal{M}}_{\overline{\mathcal{L}}} \to \overline{\mathcal{M}}_\mathcal{L}$ be the maps of proposition* 6.4.10, *then for any $U \in \mathcal{L}$ we have $i_\mathcal{L}^*(x_U) = x_{\overline{U}} \otimes 1$, $j_\mathcal{L}^*(x_U) = x_{\overline{U}}$ and $i_\mathcal{L}^*(x_T) = 1 \otimes x_T$. Note the slight abuse of notation, this should cause no problems.*

PROOF. We prove this for $j_\mathcal{L}$ the other case being similar. It is clear from the definition of $j_\mathcal{L}$ that we have the following commutative diagram.

$$\begin{array}{ccc} \overline{\mathcal{M}}_{\overline{\mathcal{L}}} & \xrightarrow{j_\mathcal{L}} & \overline{\mathcal{M}}_\mathcal{L} \\ \downarrow & & \downarrow \\ PV_{\overline{U}} & \xrightarrow{r_U} & PV_U \end{array}$$

Now it is easy to check that $y_{\overline{U}} = r_U^*(y_U)$ thus by definition we see that $x_{\overline{U}} = j_\mathcal{L}^*(x_U)$. □

REMARK 6.4.14. Before we prove a list of results concerning the cohomology ring of $\overline{\mathcal{M}}_\mathcal{L}$ we need one definition and a lemma so that we may compute the Poincaré series. We also recall a result from chapter 3 that we need. Let $\mathcal{L}$ be a thicket and $T \subseteq S$ with $q_T : S \to S/T$ the usual collapsing map. Recall in chapter 3 we defined a map $b : \text{Forests}(\mathcal{L}, T) \to \text{Forests}(\overline{\mathcal{L}})$ which we proved in lemma 3.3.5 was a bijection with the property of preserving the numbers $m(\mathcal{F}, T)$. We will use this to give a (rather crude) combinatorial description of the Poincaré series of our cohomology ring. This description will be in the right form to later allow us to compare the rings and work out a basis for each cohomology ring. This is considered in the final chapter. We will write $PS(H^*X)$ for the Poincaré series of the cohomology ring of $X$ and in the case of $X = \overline{\mathcal{M}}_\mathcal{L}$ we may abbreviate this to $PS_\mathcal{L}$.

DEFINITION 6.4.15. For any thicket $\mathcal{L}$ of $S$ and for any forest $\mathcal{F} \in \mathbb{F}_\mathcal{L}$ let $U \in \mathcal{F}$ and write $p_{\mathcal{F},U} = \sum_{i=1}^{m(\mathcal{F},U)-1} t^i$ and $p_\mathcal{F} = \prod_{U \in \mathcal{F}} p_{\mathcal{F},U}$. We then define $P_\mathcal{L} = \sum_{\mathcal{F} \in \mathbb{F}_\mathcal{L}} p_\mathcal{F}$

REMARK 6.4.16. It is possible that $p_{\mathcal{F},U}$ and so $p_\mathcal{F}$ can be zero. This can happen if and only if there is an element $U \in \mathcal{F}$ of size 2 and corresponds to the fact that $x_U = 0$ or $M(\mathcal{F}, U) = \{W\}$ and $|W| = |U| - 1$ this will corresponds to the relation $x_U(x_U - x_W)$ that we have in this case.

LEMMA 6.4.17. *Let $(\mathcal{L}, S, T)$ be an admissible thicket and $b : \text{Forests}(\mathcal{L}_+, T) \to \text{Forests}(\overline{\mathcal{L}})$ be the bijection of lemma 3.3.5 where $\overline{\mathcal{L}} = q_T(\mathcal{L})$ then for any forest $\mathcal{F} \in \text{Forests}(\mathcal{L}_+, T)$ we have $p_{\mathcal{F}} = p_{b(\mathcal{F})} p_T$ where $p_T = \sum_{i=1}^{|T|-2} t^i$.*

PROOF. Using the construction in lemma 3.3.5 we see that $b(\mathcal{F}) = q_T(\mathcal{U})$ where $\mathcal{U} = \mathcal{F} \setminus \{T\}$. Again by lemma 3.3.5 we see that for every $U \in \mathcal{U}$ that $m(\mathcal{F}, U) = m(b(\mathcal{F}), \overline{U})$ thus $p_{\mathcal{F}, U} = p_{b(\mathcal{F}), \overline{U}}$. Since the evident induced map $q_T : \mathcal{U} \to b(\mathcal{F})$ is a bijection we see that

$$\begin{aligned} p_{b(F)} &= \prod_{\overline{U} \in b(\mathcal{F})} p_{b(\mathcal{F}), \overline{U}} \\ &= \prod_{U \in \mathcal{U}} p_{b(\mathcal{F}), \overline{U}} \\ &= \prod_{U \in \mathcal{U}} p_{\mathcal{F}, U} \end{aligned}$$

and so $p_{\mathcal{F}} = p_{b(\mathcal{F})} p_T$. □

PROPOSITION 6.4.18. *For any finite set $S$ and thicket $\mathcal{L}$ of $S$ we have the following.*

$$\begin{aligned} H^{2m+1} \overline{\mathcal{M}}_{\mathcal{L}} &= 0 \\ H^m \overline{\mathcal{M}}_{\mathcal{L}} &\quad \text{is a finite free module} \\ H^2 \overline{\mathcal{M}}_{\mathcal{L}} &= \mathbb{Z}\{\, x_T \mid T \in \mathcal{L},\, |T| > 2 \,\} \\ \text{rank}(H^2 \overline{\mathcal{M}}_{\mathcal{L}}) &= |\{\, T \in \mathcal{L} \mid |T| > 2 \,\}| \\ H^* \overline{\mathcal{M}}_{\mathcal{L}} &\quad \text{is generated as a ring by } H^2 \overline{\mathcal{M}}_{\mathcal{L}} \\ PS(H^* \overline{\mathcal{M}}_{\mathcal{L}}) &= P_{\mathcal{L}} \\ H^{2(|S|-2)} \overline{\mathcal{M}}_{\mathcal{L}} &= \mathbb{Z}[x_S^{|S|-2}] \\ H^m \overline{\mathcal{M}}_{\mathcal{L}} &= 0 \text{ for } m > 2(|S|-2) \end{aligned}$$

PROOF. We use an induction argument to prove our claims. For $|S| = 3$ the claims are trivial to check for any thicket $\mathcal{L}$ of $S$. Suppose the claims are true for $|S| = n-1$ for any thicket $\mathcal{L}$ of $S$. Then for sets $S$ of size $|S| = n$ let $\mathcal{L}$ be any thicket. We use induction on the size $|\mathcal{L}|$. For $|\mathcal{L}| = 1$ the results are again clear. Suppose the result is true for some thicket $\mathcal{L}$ of $S$ of size $|\mathcal{L}| = m - 1$. Then for any thicket $\mathcal{L}_+$ of $S$ with $|\mathcal{L}_+| = m$, $m > 1$ put $\mathcal{L}$ to be the set $\mathcal{L}_+$ with $T$ removed where $T$ is an element in $\mathcal{L}_+$ of minimal size. Then it is clear that $\mathcal{L}$ is a thicket and $|\mathcal{L}| = m - 1$. We first suppose $|T| = 2$ in this case it is clear that $\overline{\mathcal{M}}_{\mathcal{L}_+} \simeq \overline{\mathcal{M}}_\mathcal{L}$ and $x_T = 0$ thus by induction all of our claims hold. Next we suppose $|T| > 2$ then in this case $\overline{\mathcal{L}}$ is a thicket of $S/T$ and $|S/T| < |S|$.

We first prove that $H^*\overline{\mathcal{M}}_\mathcal{L}$ is generated in degree 2 by $\{x_T | T \in \mathcal{L}\}$. Let $w \in H^*\overline{\mathcal{M}}_{\mathcal{L}_+}$ then inductively $i_\mathcal{L}^*(w)$ is a polynomial in $H^2(\mathcal{M}_{\overline{\mathcal{L}}}) \oplus H^2 PV_T$ which is generated by $\{x_{\overline{U}}, x_T | \overline{U} \in \overline{\mathcal{L}}\}$. Choose a section $s_T : \overline{\mathcal{L}} \to \mathcal{L}$ ie $q_T s_T = $ id and put $w'$ to be the obvious same polynomial in $\{x_{s(\overline{U})}, ..., x_T\} \subseteq H^*\overline{\mathcal{M}}_{\mathcal{L}_+}$ so that $i_\mathcal{L}^*(w) = i_\mathcal{L}^*(w')$. Now consider the long exact sequence from proposition 6.4.11. Then $(0, w - w') \in \ker(g_\mathcal{L}^*) = \text{im}(f_\mathcal{L}^*)$ and so there is some $u \in H^*\overline{\mathcal{M}}_\mathcal{L}$ with $w - w' = \pi^*(u)$ but $H^*\overline{\mathcal{M}}_\mathcal{L}$ is generated by degree 2 and thus $w = w' + \pi^*(u)$ where $u$ is a polynomial in degree 2 elements. From this it immediately follows that there is no odd cohomology.

We next prove that $H^*\overline{\mathcal{M}}_\mathcal{L}$ is a finite free module. Since the odd cohomology is zero the exact sequence of proposition 6.4.11 decouples and inductively both $H^{2q}(\overline{\mathcal{M}}_\mathcal{L})$ and $H^{2q}(\overline{\mathcal{M}}_{\overline{\mathcal{L}}} \times PV_T)$ are finite free modules thus the short exact sequence splits and we obtain that $H^*(\overline{\mathcal{M}}_{\mathcal{L}_+})$ is a finite free module. By considering the $H^2$ term of our sequence we also obtain that

$$\begin{aligned}
\text{rank}(H^2\overline{\mathcal{M}}_{\mathcal{L}_+}) &= 1 + \text{rank}(H^2\overline{\mathcal{M}}_\mathcal{L}) \\
&= 1 + |\{U \in \mathcal{L} \mid |U| > 2\}| \\
&= |\{U \in \mathcal{L}_+ \mid |U| > 2\}| \text{ because } |T| > 2
\end{aligned}$$

Next we prove that the inclusion $k : \mathbb{Z}[x_T | T \in \mathcal{L}, |T| > 2] \to H^2\overline{\mathcal{M}}_\mathcal{L}$ is an isomorphism. The map $k$ is surjective by the first part and by the previous part we deduced that the ranks are equal so $k$ is an isomorphism.

For the Poincaré series we have by the split exact sequence that $PS(\overline{\mathcal{M}}_{\mathcal{L}_+}) + PS(\overline{\mathcal{M}}_{\overline{\mathcal{L}}}) = PS(\overline{\mathcal{M}}_\mathcal{L}) + PS(\overline{\mathcal{M}}_{\overline{\mathcal{L}}})PS(PV_T)$ and $PS(PV_T) = \sum_{i=0}^{|T|-2} t^i$ thus rearranging we find that $PS_{\mathcal{L}_+} = PS_\mathcal{L} + PS_{\overline{\mathcal{L}}} P_T$ thus inductively we find that

$$\begin{aligned} PS_{\mathcal{L}_+} &= PS_\mathcal{L} + PS_{\overline{\mathcal{L}}} P_\mathcal{T} \\ &= P_\mathcal{L} + \sum_{\mathcal{U} \subseteq \overline{\mathcal{L}}} P_\mathcal{U} P_\mathcal{T} \\ &= \sum_{\mathcal{T} \subseteq \mathcal{L}} P_\mathcal{T} + \sum_{T \in \mathcal{T} \subseteq \mathcal{L}_+} P_\mathcal{T} \text{ by lemma 6.4.17} \\ &= \sum_{T \notin \mathcal{T} \subseteq \mathcal{L}_+} P_\mathcal{T} + \sum_{T \in \mathcal{T} \subseteq \mathcal{L}_+} P_\mathcal{T} \\ &= \sum_{\mathcal{T} \subseteq \mathcal{L}_+} P_\mathcal{T} \\ &= P_{\mathcal{L}_+} \end{aligned}$$

Here we prove the final two parts of the claim. The spaces $\overline{\mathcal{M}}_\mathcal{L}$ and $PV_S$ are compact and have the same dimension. Consider the map $\pi_S : \overline{\mathcal{M}}_\mathcal{L} \to PV_S$ which is an isomorphism $\pi_S : \mathcal{M}_\mathcal{L} \to U_\mathcal{L}$, $\mathcal{M}_\mathcal{L} = \pi_S^{-1}(U_\mathcal{L})$. Put $d = \dim(\overline{\mathcal{M}}_\mathcal{L}) = 2(|S| - 2)$ then Poincaré duality tells us that $H^d(\overline{\mathcal{M}}_\mathcal{L}) = \mathbb{Z}u$ and $H^d(PV_S) = \mathbb{Z}v$. Therefore $\pi_S^*(v) = ku$ where $k = \deg(\pi_S)$. Since $\pi : \pi^{-1}U_\mathcal{L} \to U_\mathcal{L}$ is an isomorphism we see that $\deg(\pi) = \pm 1$. Because $\pi$ is analytic we see $\deg(\pi) = 1$. Thus $\pi_S^d : H^d(PV_S) \to H^d(\overline{\mathcal{M}}_\mathcal{L})$ is an isomorphism. We know $u = y_S^{n-2}$ and $mu \neq 0$ for every non-zero $m \in \mathbb{Z}$ so we deduce $mx_S^{n-2} = \pi_S^*(mu) \neq 0$. The last part is clear. □

REMARK 6.4.19. In particular let $i : \overline{\mathcal{M}}_\mathcal{L} \to \prod PV_T$ be the inclusion map then we have shown that the map $i^* : H^*(\prod PV_T) \to H^*(\overline{\mathcal{M}}_\mathcal{L})$ is surjective and in degree 2 an isomorphism. Then $H^*(\overline{\mathcal{M}}_\mathcal{L})$ is a quotient of $\mathbb{Z}[\,x_T \,|\, T \in \mathcal{L} \text{ and } |T| > 2\,]$ and the kernel contains the ideal generated by the $x_T^{|T|-2}$.

## 6.5. Comparisons with the Chow ring

In this section we consider the natural map $\text{cl} : A^*(\overline{\mathcal{M}}_\mathcal{L}) \to H^{2*}(\overline{\mathcal{M}}_\mathcal{L})$ [**3**, Chapter 19] and prove it is an isomorphism for thickets $\mathcal{L}$. In particular we will show that $\overline{\mathcal{M}}_S$ is a homology isomorphism. This will be easy given our current results. We then given an explicit description for this map in chapter 9 for the case of $\overline{\mathcal{M}}_S$.

LEMMA 6.5.1. *Let $\mathcal{L}$ be a thicket then the natural map* cl : $A_*(\overline{\mathcal{M}}_\mathcal{L}) \to H_{2*}(\overline{\mathcal{M}}_\mathcal{L})$ *is an isomorphism. In particular every analogous result to proposition* 6.4.18 *holds for* $A_*(\overline{\mathcal{M}}_S)$.

PROOF. Consider the following commutative diagram where the bottom row is short exact and the top row is right exact by [**3**, Example 1.8.1]

$$
\begin{array}{ccccc}
A_*(\overline{\mathcal{M}}_{\overline{\mathcal{L}}} \times PV_T) & \longrightarrow & A_*(\overline{\mathcal{M}}_{\overline{\mathcal{L}}}) \oplus A_*(\overline{\mathcal{M}}_{\mathcal{L}_+}) & \longrightarrow & A_*(\overline{\mathcal{M}}_\mathcal{L}) \\
\downarrow \mathrm{cl} & & \downarrow \mathrm{cl} & & \downarrow \mathrm{cl} \\
H_{2*}(\overline{\mathcal{M}}_{\overline{\mathcal{L}}} \times PV_T) & \longrightarrow & H_{2*}(\overline{\mathcal{M}}_{\overline{\mathcal{L}}}) \oplus H_{2*}(\overline{\mathcal{M}}_{\mathcal{L}_+}) & \longrightarrow & H_{2*}(\overline{\mathcal{M}}_\mathcal{L})
\end{array}
$$

Then it is well known that cl : $A_*(PV) \to H_{2*}(PV)$ is an isomorphism because for example $PV$ has a cellular decomposition. Thus we may use an induction argument as in proposition 6.4.18 to suppose that the outside maps are isomorphisms and then use a simple diagram chase to deduce that the middle map is an isomorphism and thus cl : $A_*(\overline{\mathcal{M}}_\mathcal{L}) \to H_{2*}(\overline{\mathcal{M}}_\mathcal{L})$ is an isomorphism. We leave the details to the interested reader as previous calculations are similar. □

LEMMA 6.5.2. *Let $\mathcal{L}$ be a thicket then the natural map* cl : $A^*(\overline{\mathcal{M}}_\mathcal{L}) \to H^{2*}(\overline{\mathcal{M}}_\mathcal{L})$ *is an isomorphism.*

PROOF. Consider the following commutative diagram

$$
\begin{array}{ccc}
A^*(\prod PV_T) & \longrightarrow & A^*(\overline{\mathcal{M}}_\mathcal{L}) \\
\downarrow \mathrm{cl} & & \downarrow \mathrm{cl} \\
H^{2*}(\prod PV_T) & \longrightarrow & H^{2*}(\overline{\mathcal{M}}_\mathcal{L})
\end{array}
$$

then the left hand vertical map is an isomorphism and the bottom horizontal map is surjective thus the map cl : $A^*(\overline{\mathcal{M}}_\mathcal{L}) \to H^{2*}(\overline{\mathcal{M}}_\mathcal{L})$ is surjective. We deduce using the last lemma that both are finite free modules and have the same rank. Thus the map is an isomorphism. □

## 6.6. The blowup description for $\overline{\mathcal{M}}_\mathcal{L}$

In the final section of this chapter we extend the result of chapter 4 section 3 by proving that the commutative diagram of section 4 is actually a blowup diagram and that the space $\overline{\mathcal{M}}_S$ may be seen as an iterated blowup of projective space. This should be compared with Kapranov's approach [8]. It will also become clear that the order in which the blowup is obtained is not of great importance although some care is required. Although this section answers the question of whether the two spaces $\overline{\mathcal{X}}_S$ and $\overline{\mathcal{M}}_S$ are the same in the affirmative it does not provide us with an actual isomorphism between the two so as such is not a satisfactory answer to this question. We will provide a more careful analysis in later chapters to produce an explicit isomorphism.

We first introduce an order on the set $P^+(S)$, this is the subset of $P(S)$ whose elements $T$ have $|T| > 1$. This is the order introduced by Kapranov in [8]. We will use this to prove that $\overline{\mathcal{X}}_S$ the moduli space of stable $n+1$ pointed curves of genus zero and $\overline{\mathcal{M}}_S$ are isomorphic by blowing up the same initial space along the same subspaces in the same order. Here whenever $S$ is a finite set with $|S| = n$ we take $S = \{1, 2, .., n\}$.

DEFINITION 6.6.1. Let $S$ be a finite set with $|S| = n$ and put $S_m = \{\,1, 2, ..., n-m\,\}$ for $0 \leq m < n$ then we define

$$P^S_{m,i} = \{\,T \subseteq S_m \mid |T| = |S_m| - i,\ n-m \in T\,\}.$$

DEFINITION 6.6.2. Let $S$ be a finite set of size $n$ with the evident order then we define the binary ordering $\leq$ on $P^+(S)$ by $U \leq V \iff \sum_{u \in U} 2^u \leq \sum_{v \in V} 2^v$. For any $T \subseteq S$ we have the map $q_T : S \to S/T$ and give $S/T$ the induced order under the injection $q_T : S \setminus T \to S/T$ where we use $j = \max T$ for the equivalence class $T$.

DEFINITION 6.6.3. We define a total ordering $\preccurlyeq$ on $P^+(S)$ by

$$
\begin{aligned}
U \preccurlyeq V \iff \max V &\leq \max U \text{ or} \\
\max V &= \max U \text{ and } |V| \leq |U| \text{ or} \\
\max V &= \max U \text{ and } |V| = |U| \text{ and } U \leq V
\end{aligned}
$$

We write $T_1, T_2, \ldots\ldots, T_m$ for this order where $m = 2^n - 1 - n$.

REMARK 6.6.4. Let $U, V \in P^+(S)$ then $U \in P_{l,i}$ and $V \in P_{m,j}$ then the above ordering is equivalent to the following,

$$
\begin{aligned}
U \preccurlyeq V \iff m &< l \quad \text{or} \\
m &= l \quad \text{and } j < i \text{ or} \\
m &= l \quad \text{and } j = i \text{ and } U \leq V
\end{aligned}
$$

LEMMA 6.6.5. *If $U \subseteq T$ then $T \preccurlyeq U$.* □

REMARK 6.6.6. The last result tells us that if for each $U \subseteq S$ we put $\mathcal{L}_U = \{T \subseteq S | T \preccurlyeq U\}$ then $\mathcal{L}_U$ is a thicket.

REMARK 6.6.7. Here we state the following fact, the proof is unilluminating so we do not offer one. Define the following recurrence relation

$$
\begin{aligned}
V_{n-1}^n &= n - 1 \\
V_{k-1}^n &= V_k^n + \sum_{i=1}^{k-2} (k-1-i)\binom{k}{i} V_{i+1}^{n+i-k}
\end{aligned}
$$

then $V_1^n = \dim(H^*(\overline{\mathcal{M}}_S))$ where $n = |S|$

LEMMA 6.6.8. *Let $(\mathcal{L}, S, T)$ be an admissible thicket and put $\mathcal{T}_\mathcal{L}(T) = \{\, \mathcal{T} \in \mathbb{T}_\mathcal{L} \mid \mathcal{T} \amalg \{T\}$ is a tree $\}$ then for any tree $\mathcal{T} \in \mathcal{T}(T)$ we have the following commutative diagram*

$$\begin{array}{ccc} \overline{\mathcal{M}}_\mathcal{L} & \xrightarrow{\pi_1} & \overline{\mathcal{M}}_\mathcal{T} \\ \uparrow{j_\mathcal{L}} & & \uparrow{j_\mathcal{T}} \\ \overline{\mathcal{M}}_{\overline{\mathcal{L}}} & \xrightarrow{\pi_2} & \overline{\mathcal{M}}_{\overline{\mathcal{T}}} \end{array}$$

$$\pi_1^{-1} j_\mathcal{T}(\overline{\mathcal{M}}_{\overline{\mathcal{T}}}) \cap \mathcal{N}_\mathcal{L}(\mathcal{T}) = j_\mathcal{L}(\overline{\mathcal{M}}_{\overline{\mathcal{L}}}) \cap \mathcal{N}_\mathcal{L}(\mathcal{T})$$

$$\bigcup_{\mathcal{T} \in \mathcal{T}_\mathcal{L}(T)} \mathcal{N}_\mathcal{L}(\mathcal{T}) \supseteq j_\mathcal{L}(\overline{\mathcal{M}}_{\overline{\mathcal{L}}})$$

PROOF. Let $\underline{M} \in \pi_1^{-1} j_\mathcal{T}(\overline{\mathcal{M}}_{\overline{\mathcal{T}}}) \cap \mathcal{N}_\mathcal{L}(\mathcal{T})$ then $\underline{M} \in \mathcal{N}_\mathcal{L}(\mathcal{T})$ and $\underline{N} = \pi_1(x) \in j_\mathcal{T}(\overline{\mathcal{M}}_{\overline{\mathcal{T}}}) \cap \mathcal{N}'_\mathcal{L}(\mathcal{T})$ since $\mathcal{N}_\mathcal{L}(\mathcal{T})$ is a saturated set. Because $j_\mathcal{T}$ is injective there is a unique $\underline{L} \in \overline{\mathcal{M}}_{\overline{\mathcal{T}}}$ with $j_\mathcal{T}(\underline{L}) = \underline{N}$. We prove that $\underline{L} \in \mathcal{N}'_{\overline{\mathcal{L}}}(\overline{\mathcal{T}})$. We first observe that the map $q_\mathcal{T} : \mathcal{T} \to \overline{\mathcal{T}}$ is surjective. Now let $U \in \mathcal{L}$ and put $T = \text{root}(U)$ in $\mathcal{T}$. Then $\overline{T} = \text{root}(\overline{U})$ in $\overline{\mathcal{T}}$. As $\underline{N} \in \mathcal{N}'_\mathcal{L}(\mathcal{T})$ then by definition $\pi_U^T N_T \neq 0$. By the definition of $j_\mathcal{T}$ we see that $N_T = r_T(L_{\overline{T}})$ and $N_U = r_U(L_{\overline{U}})$. Thus we may apply lemma 6.4.4 to deduce that $\pi_{\overline{U}}^{\overline{T}} N_{\overline{T}} \neq 0$ and so $\underline{L} \in \mathcal{N}'_{\overline{\mathcal{L}}}(\overline{\mathcal{T}})$. Now $\pi_1^{-1}(\underline{N}) = \underline{M}$ and there is a unique $\underline{K} \in \mathcal{N}_{\overline{\mathcal{L}}}(\overline{\mathcal{T}})$ with $\pi_2(\underline{K}) = \underline{L}$ and therefore $\underline{M} = j_\mathcal{L}(\underline{K})$. Thus $\pi_1^{-1} j_\mathcal{T}(\overline{\mathcal{M}}_{\overline{\mathcal{T}}}) \cap \mathcal{N}_\mathcal{L}(\mathcal{T}) \subseteq j_\mathcal{L}(\overline{\mathcal{M}}_{\overline{\mathcal{L}}}) \cap \mathcal{N}_\mathcal{L}(\mathcal{T})$ and the reverse inclusion is automatically true.

Next given $\underline{M} \in j_\mathcal{L}(\overline{\mathcal{M}}_{\overline{\mathcal{L}}})$ put $\mathcal{T} = \text{type}(\underline{M})$. Let $U \in \mathcal{L}$ with $U \cap T$ non-empty and $U \not\supseteq T$. By the minimality of $T$ we have $U \not\subseteq T$. It follows that $W = U \cup T \supset T$ and $W \in \mathcal{L}$ as $\mathcal{L}_+$ is a thicket. Then since by lemma 6.4.9 we have $\pi_T^W(M_W) = 0$ we must have $\pi_U^W(M_W) = M_U$ thus $U \notin \mathcal{T}$. Therefore $\mathcal{T} \amalg \{T\}$ is an $S$-tree.

$\square$

COROLLARY 6.6.9. *Let $(\mathcal{L}, S, T)$ be an admissible thicket then the following diagram is a pullback and $\pi : \overline{\mathcal{M}}_{\mathcal{L}_+} \to \overline{\mathcal{M}}_\mathcal{L}$ is the blowup of $\overline{\mathcal{M}}_\mathcal{L}$ along $\overline{\mathcal{M}}_{\overline{\mathcal{L}}}$.*

$$\begin{array}{ccc} \overline{\mathcal{M}}_{\overline{\mathcal{L}}} \times PV_T & \xrightarrow{i_{\mathcal{L}}} & \overline{\mathcal{M}}_{\mathcal{L}_+} \\ {\scriptstyle p} \downarrow & & \downarrow {\scriptstyle \pi} \\ \overline{\mathcal{M}}_{\overline{\mathcal{L}}} & \xrightarrow{j_{\mathcal{L}}} & \overline{\mathcal{M}}_{\mathcal{L}} \end{array}$$

PROOF. Let $b : \mathrm{Bl}_Y X \to X$ to be the blowup of $X$ along $Y$ where $X = \overline{\mathcal{M}}_{\mathcal{L}}$ and $Y = \overline{\mathcal{M}}_{\overline{\mathcal{L}}}$. Let $\pi : \overline{\mathcal{M}}_{\mathcal{L}_+} \to \overline{\mathcal{M}}_{\mathcal{L}}$ be the standard projection. Then to prove $\pi$ is the blowup of $X$ along $Y$ it will be enough to cover $\overline{\mathcal{M}}_{\mathcal{L}}$ by Zariski open sets $\{U_i\}_{i \in I}$ so that $\pi^{-1}(U_i)$ and $b^{-1}(U_i)$ are isomorphic by a morphism (necessarily unique) that commutes with projection onto $U_i$. Let $\mathcal{T}$ be an element of $\mathcal{T}_{\mathcal{L}}(T)$ and consider the following commutative diagram

$$\begin{array}{ccc} \overline{\mathcal{M}}_{\mathcal{L}_+} & \xrightarrow{\pi_3} & \overline{\mathcal{M}}_{\mathcal{T}_+} \\ {\scriptstyle \pi} \downarrow & & \downarrow {\scriptstyle \pi_4} \\ \overline{\mathcal{M}}_{\mathcal{L}} & \xrightarrow{\pi_2} & \overline{\mathcal{M}}_{\mathcal{T}} \\ {\scriptstyle j_{\mathcal{L}}} \uparrow & & \uparrow {\scriptstyle j_{\mathcal{T}}} \\ \overline{\mathcal{M}}_{\overline{\mathcal{L}}} & \xrightarrow{\pi_1} & \overline{\mathcal{M}}_{\overline{\mathcal{T}}} \end{array}$$

By lemma 6.2.3 $\pi_2$ maps $\mathcal{N}_{\mathcal{L}}(\mathcal{T})$ isomorphically onto $\mathcal{N}'_{\mathcal{L}}(\mathcal{T})$. By the previous lemma $\pi_2$ maps $j_{\mathcal{L}}(\overline{\mathcal{M}}_{\overline{\mathcal{L}}}) \cap \mathcal{N}_{\mathcal{L}}(\mathcal{T})$ isomorphically onto $j_{\mathcal{T}}(\overline{\mathcal{M}}_{\overline{\mathcal{T}}}) \cap \mathcal{N}'_{\mathcal{L}}(\mathcal{T})$ as saturated sets and by lemma 4.4.14 $\pi_4 : \overline{\mathcal{M}}_{\mathcal{T}_+} \to \overline{\mathcal{M}}_{\mathcal{T}}$ is the blowup of $\overline{\mathcal{M}}_{\mathcal{T}}$ along $\overline{\mathcal{M}}_{\overline{\mathcal{T}}}$. Thus we have $b^{-1}(\mathcal{N}_{\mathcal{L}}(\mathcal{T}))$ is isomorphic to $\pi_4^{-1}(\mathcal{N}'_{\mathcal{L}}(\mathcal{T}))$ commuting through the restricted map of $\pi_2$ onto their respective bases. Now put $U_{\mathcal{T}} = \pi^{-1}(\mathcal{N}_{\mathcal{L}}(\mathcal{T}))$ and $V_{\mathcal{T}} = \pi_4^{-1}(\mathcal{N}'(\mathcal{L}))$ then it is enough to show that the map $\pi_3 : U_{\mathcal{T}} \to V_{\mathcal{T}}$ is an isomorphism. This is clear since $T$ is an element of minimal size. Thus we obtain $\pi^{-1}(\mathcal{N}_{\mathcal{L}}(\mathcal{T}))$ is isomorphic to $b^{-1}(\mathcal{N}_{\mathcal{L}}(\mathcal{T}))$ with projections commuting through $\mathcal{N}_{\mathcal{L}}(\mathcal{T}) \subseteq \overline{\mathcal{M}}_{\mathcal{L}}$. By the previous lemma we have that

$$\bigcup_{\mathcal{T} \in \mathcal{T}_{\mathcal{L}}(T)} \mathcal{N}_{\mathcal{L}}(\mathcal{T}) \supseteq j_{\mathcal{L}}(\overline{\mathcal{M}}_{\overline{\mathcal{L}}})$$

and on each of these open sets $\mathcal{N}_{\mathcal{L}}(\mathcal{T})$ we have shown that $\pi^{-1}(\mathcal{N}_{\mathcal{L}}(\mathcal{T}))$ is isomorphic to $b^{-1}(\mathcal{N}_{\mathcal{L}}(\mathcal{T}))$ with projections commuting through $\mathcal{N}_{\mathcal{L}}(\mathcal{T})$. Now by lemma 6.4.10 we see that $\pi : \overline{\mathcal{M}}_{\mathcal{L}_+} \to \overline{\mathcal{M}}_{\mathcal{L}}$ is an isomorphism away from $j_{\mathcal{L}}(\overline{\mathcal{M}}_{\overline{\mathcal{L}}})$ thus we see that $\pi$ is the blowup of $\overline{\mathcal{M}}_{\mathcal{L}}$ along $\overline{\mathcal{M}}_{\mathcal{L}_+}$.

□

COROLLARY 6.6.10. *Let $(\mathcal{L}, S, T)$ be an admissible thicket and $j_{\mathcal{L}} : \overline{\mathcal{M}}_{\overline{\mathcal{L}}} \to \overline{\mathcal{M}}_{\mathcal{L}}$ be the embedding. Then the normal bundle is $\hom(W_{\mathcal{L}}, V_T)$ and the projectivization of the normal bundle is the trivial bundle $\overline{\mathcal{M}}_{\overline{\mathcal{L}}} \times PV_T$. We define $W_{\mathcal{L}}$ below.*

PROOF. Consider the following commutative diagram where $\mathcal{T} \in \mathcal{T}_{\mathcal{L}}(T)$

$$\begin{array}{ccc} \overline{\mathcal{M}}_{\mathcal{L}} & \xrightarrow{\pi_1} & \overline{\mathcal{M}}_{\mathcal{T}} \\ j_{\mathcal{L}} \uparrow & & \uparrow j_{\mathcal{T}} \\ \overline{\mathcal{M}}_{\overline{\mathcal{L}}} & \xrightarrow{\pi_2} & \overline{\mathcal{M}}_{\overline{\mathcal{T}}} \end{array}$$

Lemma 6.6.8 essentially tells us that over $\mathcal{N}_{\overline{\mathcal{L}}}(\overline{\mathcal{T}})$ the normal bundle $N(j_{\mathcal{L}})$ is $N(j_{\mathcal{T}})$ and this is $\hom(N_{U/T}, V_T)$ where $U = \text{root}(T)$ in $\mathcal{T}$. Doing this for each tree in $\overline{\mathcal{L}}$ we obtain local description for the normal bundle. One can check that the $N_{U/T}$ are compatible over pairwise intersection in our cover and glue to obtain the vector bundle $W_{\mathcal{L}}$, thus we may take $N(j_{\mathcal{L}}) = \hom(W_{\mathcal{L}}, V_T)$. This can also be written as $W_{\mathcal{L}}^* \otimes V_T$ where the star denotes the dual bundle. This gives the projectivization as $\overline{\mathcal{M}}_{\overline{\mathcal{L}}} \times PV_T$ the expected answer.

□

CONSTRUCTION 6.6.11. Given a thicket $\mathcal{L}$ we define an order on it as follows. We define $T_1 = S$ and $\mathcal{L}_1 = \{S\}$. Suppose $T_1, ..., T_r$ and $\mathcal{L}_1, ..., \mathcal{L}_r$ are defined. Then we define $T_{r+1} \in \mathcal{L} \setminus \mathcal{L}_r$ to be an element of maximal size and $\mathcal{L}_{r+1} = \mathcal{L}_r \amalg \{T_{r+1}\}$. Then it is clear that the $(\mathcal{L}_i, S, T_{i+1})$ form a sequence of admissible thickets and $\mathcal{L}_i = \{T_j \mid j \leq i\}$. We then define a sequence of spaces $B_i, S_i^U$ by $B_1 = \overline{\mathcal{M}}_{\mathcal{L}_1} = PV_S$, $S_1^U = P(\ker(\pi_U^S))$ and given $B_i$ and $S_i^U$ are defined. We define $B_{i+1} = Bl_{S_i^T} B_i$ where $T = T_{i+1}$ and write $\pi_{i+1} : \overline{\mathcal{M}}_{\mathcal{L}_{i+1}} \to \overline{\mathcal{M}}_{\mathcal{L}_i}$ for the projection map. Then for any $U \in \mathcal{L}$ we define $S_{i+1}^U =$ strick transform of $S_i^U$ in $B_i$.

LEMMA 6.6.12. *Let $\mathcal{L}$ be a thicket and $\mathcal{L}_i$ an associated sequence of thickets then we have the following sequence of commutative diagrams*

$$\begin{array}{ccc} \overline{\mathcal{M}}_{\overline{\mathcal{L}}_n} \times PV_T & \xrightarrow{i} & \overline{\mathcal{M}}_{\mathcal{L}_{n+1}} \\ p \downarrow & & \downarrow \pi \\ \overline{\mathcal{M}}_{\overline{\mathcal{L}}_n} & \xrightarrow{j} & \overline{\mathcal{M}}_{\mathcal{L}_n} \end{array}$$

*Then $B_i = \overline{\mathcal{M}}_{\mathcal{L}_i}$ and given $U \in \mathcal{L}$ with $U = T_k$ we have $S_i^U = j_U(\overline{\mathcal{M}}_{\overline{\mathcal{L}}_i})$ for every $i < k$ and is zero otherwise where we quotient in $q_U : S \to S/U$ and we write $j_U = j_{\overline{\mathcal{L}}_i}$.*

PROOF. We prove the claim by an induction on $i$. The case when $i = 1$ is clear from the definitions. Put $T = T_{i+1}$ then inductively $B_i = \overline{\mathcal{M}}_{\mathcal{L}_i}$ and $S_i^T = j_T(\overline{\mathcal{M}}_{\overline{\mathcal{L}}_i})$. Then by lemma 6.6.9 we see that the blowup of $B_i$ along $S_i^T$ is $B_{i+1} = \overline{\mathcal{M}}_{\mathcal{L}_{i+1}}$. We are left to prove the claim for the strict transforms. Choose $i + 1 < k$ and put $U = T_k$ then define $W_i^U = S_i^U \setminus S_i^T$ where inductively $S_i^U = j_U(\overline{\mathcal{M}}_{\overline{\mathcal{L}}_i})$ then we need to calculate $\mathrm{cl}(\pi_{i+1}^{-1}(W_i^U))$. We first prove that

$$j_U(\mathcal{M}_{\overline{\mathcal{L}}_{i+1}}) \subseteq \pi_{i+1}^{-1}(W_i^U) \subseteq j_U(\overline{\mathcal{M}}_{\overline{\mathcal{L}}_{i+1}})$$

We prove the left hand inclusion first. Let $\underline{N} \in \mathcal{M}_{\overline{\mathcal{L}}_{i+1}}$ and put $\underline{M} = j_U(\underline{N})$. Then by lemma 6.4.7 we see that for every $V \in \mathcal{L}_{i+1}$ with $V \cap U$ non-empty we have $\pi_{V \cap U}^V M_V = 0$. Put $\underline{L}$ to be the image of $\underline{M}$ under $\pi_{i+1}$. Then again by lemma 6.4.7 we see that $\underline{L} \in S_i^U$. Because $T \not\subseteq U$ by choice of order we see $\pi_T^S M_S \neq 0$ so that $\underline{L} \notin S_i^T$ thus $\underline{M} \in \pi_{i+1}^{-1}(W_i^U)$ and the left hand inclusion is true.

Next let $\underline{M} \in \pi_{i+1}^{-1}(W_i^U)$ and put $\underline{L}$ to be its image under $\pi_{i+1}$. Consider $W_i^U \subseteq \overline{\mathcal{M}}_{\mathcal{L}_i}$. Now $\mathcal{L}_i$ and $\mathcal{L}_{i+1}$ are both thickets and $\underline{L} \notin S_i^T$. Thus by lemma 6.4.9 we can find a $V \in \mathcal{L}_i$ containing $T$ such that $\pi_T^V M_V = M_T$. Because $\underline{L} \in S_i^U$ we deduce that $\pi_{W \cap U}^W M_W = 0$ for all $W \in \mathcal{L}_i$ such that $W \cap U$ is non-empty. In particular this then tells us that $\pi_{T \cap U}^T M_T = 0$ if $T \cap U$ is non-empty. This proves the right hand inclusion.

Then taking the closure of both sides and recalling that $j_U$ is a proper map and $\mathcal{M}_{\overline{\mathcal{L}}_{i+1}}$ is dense in $\overline{\mathcal{M}}_{\overline{\mathcal{L}}_{i+1}}$ we see that $S_{i+1}^U = \mathrm{cl}(\pi_{i+1}^{-1}(W^T)) = j_U(\overline{\mathcal{M}}_{\overline{\mathcal{L}}_{i+1}})$. This completes the induction. $\square$

COROLLARY 6.6.13. *The spaces $\overline{\mathcal{X}}_S$ and $\overline{\mathcal{M}}_S$ are isomorphic as projective varieties .*

PROOF. Here we take $\mathcal{L} = P^+(S)$ and use the order from definition 6.6.3. We then apply the previous result. $\square$

As a result of the previous corollary and given the uniqueness of blowups we have the following commutative diagram where $\theta_S : \overline{\mathcal{X}}_S \to \overline{\mathcal{M}}_S$ (necessarily unique) is an isomorphism, $\sigma_S : \overline{\mathcal{X}}_S \to PV_S$ is Kapranov's regular map defined in [**7**] and $\pi_S : \overline{\mathcal{M}}_S \to PV_S$ is the standard projection map.

$$\begin{array}{ccc} \overline{\mathcal{X}}_S & \xrightarrow{\theta_S} & \overline{\mathcal{M}}_S \\ & \searrow^{\sigma_S} \swarrow_{\pi_S} & \\ & PV_S & \end{array}$$

clearly the restricted map $\sigma_S : \mathcal{X}_S \to U_S$ is an isomorphism. Suppose for each $T \subseteq S$ we could show that the following diagram is commutative where $\pi_T^S : \overline{\mathcal{X}}_S \to \overline{\mathcal{X}}_T$ is a natural map discussed in [**9**]

$$\begin{array}{ccc} \overline{\mathcal{X}}_S & \xrightarrow{\theta_S} & \overline{\mathcal{M}}_S \\ \pi_T^S \downarrow & & \downarrow \pi_T^S \\ \overline{\mathcal{X}}_T & \xrightarrow{\theta_T} & \overline{\mathcal{M}}_T \end{array}$$

then it is easy to see that the map $\theta_S$ is given by $\prod \sigma_T \pi_T^S$. This is not an unreasonable expectation given that to prove this diagram commutes we would only have to restrict ourself to the isomorphisms $\sigma_V : \mathcal{X}_V \to U_V$ and intuitively it is clear what form this map must have. However we do not explicitly analyze this map over the spaces $\mathcal{X}_V$. In a later chapter we will describe an equivalent map ( ie differing by an isomorphism ) $\sigma_S : \overline{\mathcal{X}}_S \to PV_S$ that is a natural extension for the usual identification of $\mathcal{X}_S$ with $U_S$. This will have the desired properties and will give us an explicit isomorphism as described above.

# CHAPTER 7

# The topology of $\overline{\mathcal{M}}_S$

In this section we turn our attention back to the main object of interest, the space $\overline{\mathcal{M}}_S$. We will introduce various notions of trees and compare them by constructing natural bijections between them. We will then consider the different representations equivalent by these correspondences. We also prove a number of analogous results from the theory of $\overline{\mathcal{X}}_S$ associated to subspaces $\mathcal{M}_S(\mathcal{T})$ of $\overline{\mathcal{M}}_S$ that consists of elements of tree type $\mathcal{T}$. These spaces are of particular importance for the study of $\overline{\mathcal{M}}_S$ and we have already used similar constructions in chapter 6. The results we prove will highlight the tree structures of our space. We will need these later when we compare in more detail the spaces $\overline{\mathcal{M}}_S$ with $\overline{\mathcal{X}}_S$. A detailed proof of some of the results in this chapter would be unilluminating so we do not supply them. Instead we give (carefully chosen) examples that should highlight the main ideas. These results are self contained and will not be mentioned elsewhere.

## 7.1. The combinatorial structure of $\overline{\mathcal{M}}_S$

In this section we examine in more detail the combinatorial structure of $\overline{\mathcal{M}}_S$. In particular we will define a space $\mathbb{P}_S$ that is constructed in the same spirit as $\overline{\mathcal{M}}_S$ modulo some natural combinatorial structure. This should be thought of as a universal way of assigning to elements of $\overline{\mathcal{M}}_S$ the structure of a tree. We then consider the combinatorial notion of a tree and our original definition in 3.1.1 and explain the way that we consider them equivalent, that is produce natural bijections between them.

DEFINITION 7.1.1. For every $T \subseteq S$ with $|T| > 1$ we define,

$$Q_T = \{ \text{ partitions of } T \text{ into at least 2 blocks } \}$$

and $\theta_T : PV_T \to Q_T$ by $i \sim_{\theta_T(M_T)} j$ in $T \iff m(i) = m(j)$ for all $m \in q_T^{-1}(M_T)$ where $q_T : F(T, \mathbb{C}) \to V_T$ is the quotient map.

REMARK 7.1.2. Equivalently we may define the previous map as follows. Given an element $M_T \in PV_T$ and $U \subseteq T$ we define the image $\theta_T(M_T)$ by $U \in \theta_T(M_T) \iff \pi_U^T(M_T) = 0$ and for any other set $V \supset U$ we have $\pi_V^T M_T \neq 0$, that is $U$ is maximal with respect to this property.

LEMMA 7.1.3. *For every $U \subseteq T \subseteq S$ the restriction of partitions gives partial maps $\rho_U^T : Q_T \dashrightarrow Q_U$ defined by $\rho_U^T(V) = V \cap U$ compatible with $\theta_T$ that is, the following diagram commutes whenever it is defined. This happens if and only if $\rho_U^T(M_T)$ is defined or equivalently if $\rho_U^T(\omega)$ is defined where $\omega = \theta_T(M_T)$. The latter just means $U$ is not contained in a single block of the partition $\omega$ of $T$.*

$$\begin{array}{ccc} PV_T & \xrightarrow{\theta_T} & Q_T \\ \rho_U^T \downarrow & & \downarrow \rho_U^T \\ PV_U & \xrightarrow{\theta_U} & Q_U \end{array}$$

PROOF. This is clear given the last remark. □

DEFINITION 7.1.4.

$$\mathbb{P}_S = \left\{ \underline{\omega} \in \prod_{\substack{T \subseteq S \\ |T|>1}} Q_T : \rho_U^T(\omega_T) = \omega_U \text{ for all } U \subseteq T \subseteq S \text{ whenever } \rho_U^T(\omega_T) \text{ is defined} \right\}$$

$$\mathbb{T}_S = \{ \text{ S-trees } \mathcal{T} \}$$

$$\mathbb{U}_S = \left\{ \begin{array}{c} \text{isomorphisms classes of rooted trees, all vertices of valence at least 3,} \\ \text{leaves labelled bijectively with S} \end{array} \right\}$$

In the last definition we assume that $0 \notin S$ and call the root 0.

CONSTRUCTION 7.1.5. Here we construct a map $\overline{\theta}_S : \overline{\mathcal{M}}_S \to \mathbb{P}_S$ that we will later see is surjective. Given $\underline{M} \in \overline{\mathcal{M}}_S$ we define $\omega_T \in Q_T$ by $\omega_T = \theta_T(M_T)$. Then by the partially commuting diagram of lemma 7.1.3 and remark 2.0.11 it is easily seen that $\underline{\omega} = \prod \omega_T \in \mathbb{P}_S$.

DEFINITION 7.1.6. Let $t \in \mathbb{U}_S$ then for each internal vertex $v$ of $t$ we put $T_v$ to be the set of leaves of $t$ lying below $v$ away from the root 0.

CONSTRUCTION 7.1.7. Here we produce natural maps $T_1 : \mathbb{U}_S \to \mathbb{P}_S$, $T_2 : \mathbb{P}_S \to \mathbb{T}_S$ and $T_3 : \mathbb{T}_S \to \mathbb{U}_S$. We will then explain why these are bijections.

We first construct the map $T_1 : \mathbb{U}_S \to \mathbb{P}_S$. Given $t \in \mathbb{U}_S$ and $U \subseteq S$ with $|U| > 1$ let $t_U$ be the subtree of $t$ that spans $U$, that is the (rooted) subtree whose leaves are on $U$, $r_U$ the root of $t_U$ and $\omega(t)_U$ the partition of $T$ by components of $t_T \setminus \{r_T\}$. Equivalently let $t$ be a tree of $\mathbb{U}_S$ and write $V_I(t)$ for the set of internal vertices of $t$. Then for each internal vertex $v$ in $V_I(t)$ let $T_v$ be as above. Then $T_v$ is a subset of $S$ and it is clear that $\mathcal{T} = \{ T_v \mid v \text{ is an internal vertex of } t \}$ is an $S$-tree. Now given $U \subseteq S$ with $|U| > 1$ put $T = \text{root}(U)$ in $\mathcal{T}$ then we define $\omega(t)_U$ to be the induced partition on $U$ from the partition of $T$ defined by $u \sim v$ if and only if $u = v$ or $u, v \in W$ for some $W \in M(\mathcal{T}, T)$. Put $\underline{\omega} = \underline{\omega}(t)$. Then to see $\underline{\omega}$ lies in $\mathbb{P}_S$ we observe that given $U \subseteq V \subseteq S$ then $\rho_U^V(\omega_V)$ is defined if an only if $\text{root}(U) = \text{root}(V)$ in $\mathcal{T}$ and when this is the case we see $\rho_U^V(\omega_V) = \omega_U$.

EXAMPLE 7.1.8.

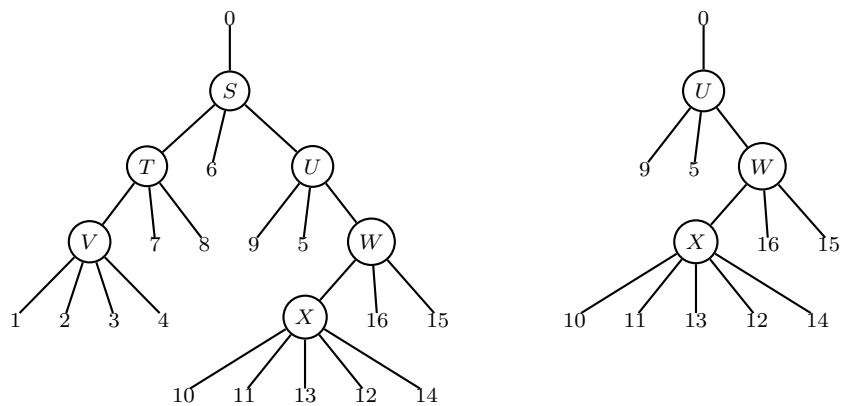

Eg if $t$ is the first tree and $Y = \{5, 9, 10, 13, 15\}$ then $\text{root}(Y) = U$ , $U$ is partitioned as $\{\{5\}, \{9\}, \{10, 11, 12, 13, 14, 15, 16\}\}$ then $\omega(t)_Y = \{\{9\}, \{5\}, \{10, 13, 15\}\}$

We next construct a map $T_2 : \mathbb{P}_S \to \mathbb{T}_S$. Given $\underline{\omega} \in \mathbb{P}_S$ define,

$$T_2(\underline{\omega}) = \{\, U \subseteq S \mid |U| > 1 \text{ and for every } U \subset T \subseteq S \implies \rho_U^T(\omega_T) \text{ is undefined } \}$$

I claim this gives us an element in $\mathbb{T}_S$. Suppose for a contradiction that $\mathcal{T} = T_2(\underline{\omega})$ is not a forest. Then we may choose elements $U, V \in \mathcal{T}$ with $U \nsubseteq V, V \nsubseteq U$ and $U \cap V$ non-empty. Put $T = U \cup V$ then $\omega_T \in Q_T$ is a partition of $T$ into at least 2 blocks. Since $\rho_U^T(\omega_T)$ and $\rho_V^T(\omega_T)$ are undefined the induced partitions of $\omega_T$ restricted to $U$ and $V$ are just $\{U\}$ and $\{V\}$ but then we can choose elements $X, Y \in \omega_T$ with $U \subseteq X$ and $V \subseteq Y$ thus $X \cap Y$ contains $U \cap V$ and so is non-empty, as $\omega_T$ is a partition we have $X = Y$ and therefore $\omega_T = \{T\}$. This contradicts the fact that $\omega_T$ has at least 2 blocks. Thus $\mathcal{T}$ is a forest and it is clear that $S \in \mathcal{T}$ so that $\mathcal{T} \in \mathbb{T}_S$.

We now construct a map $T_3 : \mathbb{T}_S \to \mathbb{U}_S$. Given an $S$-tree $\mathcal{T} \in \mathbb{T}_S$ we define an element $t \in \mathbb{U}_S$ as follows. We define the set of internal vertices of $t$ by $V_I(t) = \mathcal{T}$ and the set of external vertices by $V_E(t) = S \amalg \{0\}$. For every $U, V \in V_I(t)$ we then connect $U$ to $V$ by an edge if and only if $U \in M(\mathcal{T}, V)$ or $V \in M(\mathcal{T}, U)$. For every $T \in \mathcal{T}$ we write $U_T = T \setminus \coprod_{U \in M(\mathcal{T},T)} U$ then observe $\coprod_{T \in \mathcal{T}} U_T = S$. We connect every external vertex $s \in U_T$ to the internal vertex $T$. The vertex $0$ is the root of the tree and is connected to the internal vertex $S$. We then take the isomorphism class of this tree.

Here we produce the inverse for the bijection $T_3 : \mathbb{T}_S \to \mathbb{U}_S$ as we will need it explicitly later. We first construct a map $U_3 : \mathbb{U}_S \to \mathbb{T}_S$. Let $t$ be a tree of $\mathbb{U}_S$ and write $V_I(t)$ for the set of internal vertices of $t$. Then for each internal vertex $v$ in $V_I(t)$ let $T_v$ be the set of leaves of $t$ lying below $v$ away from the root $0$. Then $T_v$ is a subset of $S$ and it is clear that $\mathcal{T} = \{\, T_v \mid v \text{ is an internal vertex of } t\,\}$ is an $S$-tree. It is easy to see that this map is an inverse for $T_3$, moreover $U_3 = T_2 \circ T_1$.

LEMMA 7.1.9. *The maps $T_1, T_2$ and $T_3$ are natural bijections $\mathbb{P}_S \simeq \mathbb{U}_S \simeq \mathbb{T}_S$ such that we have the following commutative diagram.*

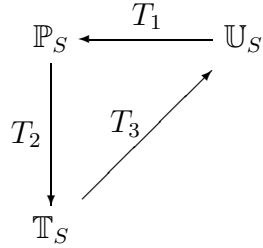

☐

## 7.2. The associated tree to elements of $\overline{\mathcal{M}}_S$

In this section we define the tree type attached to each element $\underline{M}$ in $\overline{\mathcal{M}}_S$ and consider in more detail the combinatorial notion of trees in $\mathbb{U}_S$. We will then introduce a partial ordering on the sets $\mathbb{T}_S$ and $\mathbb{U}_S$. The ordering on $\mathbb{U}_S$ will be the standard one which is used in the study of $\overline{\mathcal{X}}_S$. We will show that the bijection $U_3 : \mathbb{U}_S \to \mathbb{T}_S$ will respect several natural properties including the partial orderings and in particular we will construct in a later chapter an isomorphism of $\overline{\mathcal{X}}_S$ with $\overline{\mathcal{M}}_S$ such that the induced maps onto trees are preserved by our bijection.

DEFINITION 7.2.1. Recall from chapter 6 that we constructed the map type : $\overline{\mathcal{M}}_S \to \mathbb{T}_S$. We use this map to define the *combinatorial tree type* of an element $\underline{M}$ of $\overline{\mathcal{M}}_S$.

REMARK 7.2.2. It is useful be kept in mind that any of the previous constructions of trees are equivalent to the above definition in the obvious sense.

LEMMA 7.2.3. *The map* type : $\overline{\mathcal{M}}_S \to \mathbb{T}_S$ *is surjective.*

PROOF. The proof of this may be seen from proposition 7.3.6 part 2.

☐

DEFINITION 7.2.4. For every $T \subseteq S$ we define a map $\pi_T^S : \mathbb{T}_S \to \mathbb{T}_T$ by,

$$\pi_T^S(\mathcal{T}) = \{\, V \cap T \mid V \in \mathcal{T} \text{ and } |V \cap T| > 1 \,\}$$

We also define a map $\pi_T^S : \mathbb{U}_S \to \mathbb{U}_T$ as follows. Given a rooted tree $t \in \mathbb{U}_S$ we define the tree $t|_T = \pi_T^S(t) \in \mathbb{U}_T$ to be the restriction of the leaves $S$ of $t$ to the set $T$ subject to the constraint that if some vertex has valency 2 after restricting we remove that vertex and make that edge rigid. We call this process stably forgetting the set $S \setminus T$. See below for an example.

EXAMPLE 7.2.5.

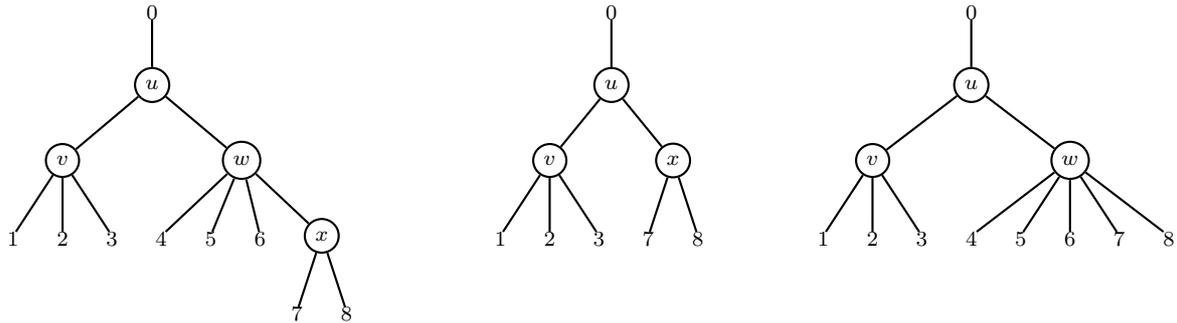

The second tree is the restriction of the first to the set $\{1, 2, 3, 7, 8\}$ and the third tree is the contraction of vertex $x$ to $w$

LEMMA 7.2.6. Let $U \subseteq T \subseteq S$ then for each of the maps $\pi_T^S$ above we have the composition rule $\pi_U^S = \pi_U^T \pi_T^S$ and we have the following commutative diagram

$$\begin{array}{ccc} \mathbb{U}_S & \xrightarrow{U_3} & \mathbb{T}_S \\ \pi_T^S \downarrow & & \downarrow \pi_T^S \\ \mathbb{U}_T & \xrightarrow{U_3} & \mathbb{T}_T \end{array}$$

$\square$

DEFINITION 7.2.7. Here we define orderings on $\mathbb{T}_S$ and $\mathbb{U}_S$. For $\mathbb{T}_S$ we define $\mathcal{T} \leq \mathcal{U} \iff \mathcal{T} \subseteq \mathcal{U}$. For $\mathbb{U}_S$ we define $t \leq u \iff t$ is obtained from $u$ by contracting some internal edges. See 7.2.5 for an example.

LEMMA 7.2.8. Let $t, u \in \mathbb{U}_S$ and put $\mathcal{T} = U_3(t)$, $\mathcal{U} = U_3(u)$ then $t \leq u \iff \mathcal{T} \leq \mathcal{U}$ that is the orders are equivalent under this bijection. See 7.2.9 for an example.   $\square$

EXAMPLE 7.2.9.

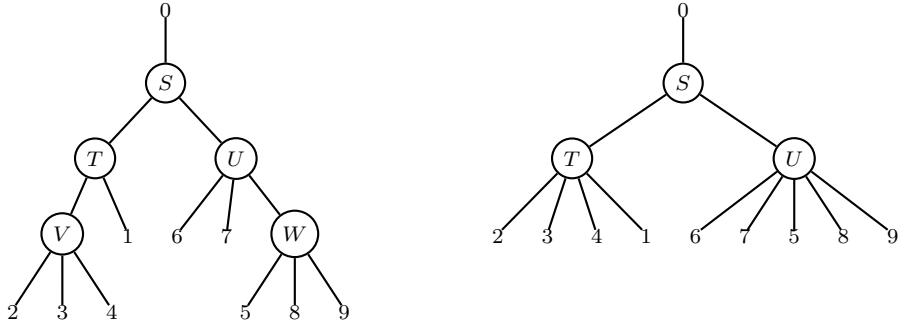

The first tree in set form is $\mathcal{U} = \{S, T, U, V, W\}$ the second tree is $\mathcal{T} = \{S, T, U\}$ and $\mathcal{T} \subseteq \mathcal{U}$ so that $\mathcal{T} \leq \mathcal{U}$.

REMARK 7.2.10. From now on we will use whichever definition of trees that is the most convenient and the choice should be clear from the context.

DEFINITION 7.2.11. For every $S$-tree $\mathcal{T} \in \mathbb{T}_S$ we define,

$$\begin{aligned}
\mathcal{M}_S(\mathcal{T}) &= \mathrm{type}^{-1}(\mathcal{T}) \\
\overline{\mathcal{M}}_S(\mathcal{T}) &= \coprod_{\mathcal{U} \supseteq \mathcal{T}} \mathcal{M}_S(\mathcal{U}) \\
\mathcal{M}_S &= \mathrm{type}^{-1}(\{S\})
\end{aligned}$$

REMARK 7.2.12. The space $\mathcal{M}_S(\mathcal{T})$ is non-empty because the map $\mathrm{type} : \overline{\mathcal{M}}_S \to \mathbb{T}_S$ is surjective. Recall from chapter 6 that we defined the Zariski open subvariety $U_S$ of $PV_S$ and proved that the map $\pi : \overline{\mathcal{M}}_S \to PV_S$ restricted to $\mathcal{M}_S$ gives us an isomorphism $\pi : \mathcal{M}_S \to U_S$ that classifies generic curves.

## 7.3. Tree isomorphisms

In this section we produce analogous results to those already known about $\overline{\mathcal{X}}_S$ regarding the structure of the spaces $\mathcal{M}_S(\mathcal{T})$ and $\overline{\mathcal{M}}_S(\mathcal{T})$. We will require these results in the next chapter when we analyze the space $\overline{\mathcal{X}}_S$ in more detail.

DEFINITION 7.3.1. For any S-tree $\mathcal{T}$ and given an element $T \in \mathcal{T}$ we define $C_T$ to be the set of equivalence classes of $T$ under the relation $u \sim v \iff u = v$ or $u, v \in W \in M(\mathcal{T}, T)$. We call $C_T$ the *children* of the vertex $T$ on the tree $\mathcal{T}$. This should be compared with its combinatorial counterpart. Observe the map $q_T : T \setminus \coprod_{W \in M(\mathcal{T}, T)} W \to C_T \setminus M(\mathcal{T}, T)$ is a bijection where $q_T : T \to C_T$ is the quotient map.

EXAMPLE 7.3.2.

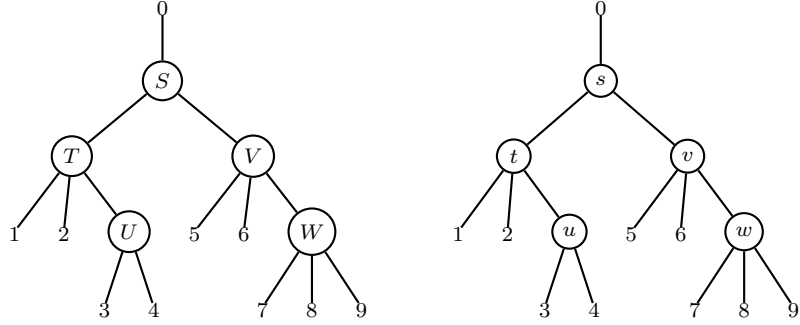

For the first tree $C_V = \{\{5\}, \{6\}, W\}$ and for second tree $C_v = \{5, 6, w\}$

REMARK 7.3.3. For each tree $t \in \mathbb{U}_S$ and corresponding S-tree $\mathcal{T}$ let $v$ be an internal vertex with $T = T_v$ then there is a natural bijection $b : C_v \to C_T$ sending $i \in S \cap C_v$ to $\{i\}$ and $u$ to $T_u$.

CONSTRUCTION 7.3.4. Let $\mathcal{T}$ be an $S$-tree and $q_T : T \to C_T$ be the usual quotient map. We construct a map $r : \prod_{T \in \mathcal{T}} \mathbb{T}_{C_T} \to \mathbb{T}_S$ by $r(\prod \mathcal{T}_T) = \coprod q_T^{-1}(\mathcal{T}_T)$. In more detail a point of $\prod \mathbb{T}_{C_T}$ consists of a system of trees $\mathcal{T}_T \subseteq P(C_T)$, one for each $T \in \mathcal{T}$. We let $q_T : T \to C_T$ be the usual quotient map and define $q_T^{-1}(\mathcal{T}_T) = \{ q_T^{-1}(U) \mid U \in \mathcal{T}_T \}$, so each $q_T^{-1}(\mathcal{T}_T) \subseteq P(T)$ is a $T$-tree. We then put $r(\prod \mathcal{T}_T) = \coprod q_T^{-1}(\mathcal{T}_T)$.

LEMMA 7.3.5. *The map* $r : \prod_{T \in \mathcal{T}} \mathbb{T}_{C_T} \to \mathbb{T}_S$ *is injective, and the image is the set of $S$-trees $\mathcal{U}$ containing $\mathcal{T}$.*

PROOF. It is clear that the map is injective and $r(\prod \mathcal{T}_T)$ is an $S$-tree. For every $T \in \mathcal{T}$ since $C_T \in \mathcal{T}_T$ we see that $T \in q_T^{-1}(\mathcal{T}_T)$ and so $\mathcal{T} \subseteq r(\prod \mathcal{T}_T)$. Given any $S-$tree $\mathcal{U}$ containing $\mathcal{T}$ let $T \in \mathcal{T}$, put $\mathcal{U}_T = \{ U \in \mathcal{U} \mid \text{root}(U) = T \text{ in } \mathcal{T}\}$ then we have an induced tree $\mathcal{T}_T = q_T(\mathcal{U}_T)$ and $r(\prod \mathcal{T}_T) = \mathcal{U}$. □

PROPOSITION 7.3.6. *Let $\mathcal{T}$ be an S-tree then we have that,*

$$\begin{aligned}
\overline{\mathcal{M}}_S &= \coprod_{\mathcal{T} \in \mathbb{T}_S} \mathcal{M}_S(\mathcal{T}) \\
\mathcal{M}_S(\mathcal{T}) &\cong \prod_{T \in \mathcal{T}} \mathcal{M}_{C_T} \\
\overline{\mathcal{M}}_S(\mathcal{T}) &\cong \prod_{T \in \mathcal{T}} \overline{\mathcal{M}}_{C_T} \\
\operatorname{cl}(\mathcal{M}_S(\mathcal{T})) &= \overline{\mathcal{M}}_S(\mathcal{T})
\end{aligned}$$

*The second isomorphism is an extension of the first and* $\operatorname{type}(\phi(\prod \underline{M}_{C_T})) = r(\prod \operatorname{type}(\underline{M}_{C_T}))$ *where $\phi$ is the inverse morphism to the second map.*

The proofs will be given at the end of this section. We observe here that in principle the closures could be different in the Zariski topology and the classical topology. However as our varieties are smooth and irreducible this result is valid in both topologies.

DEFINITION 7.3.7. For any finite set $S$ we define $P^+(S) = \{\, U \subseteq S \mid |U| > 1 \,\}$ to be the *reduced* power set of $S$.

LEMMA 7.3.8. *The space $\overline{\mathcal{M}}_S(\mathcal{T}) = \coprod_{\mathcal{U} \supseteq \mathcal{T}} \mathcal{M}_S(\mathcal{U})$ is a Zariski closed subvariety of $\overline{\mathcal{M}}_S$*

PROOF. The proof of this is easy since

$$\begin{aligned}
\overline{\mathcal{M}}_S(\mathcal{T}) &= \coprod_{\mathcal{U} \supseteq \mathcal{T}} \mathcal{M}_S(\mathcal{U}) \\
&= \overline{\mathcal{M}}_S \setminus \bigcup_{\mathcal{V} \not\supseteq \mathcal{T}} \mathcal{N}_S(\mathcal{V})
\end{aligned}$$

and each $\mathcal{N}_S(\mathcal{V})$ is Zariski open by lemma 6.2.2 part 2.

□

CONSTRUCTION 7.3.9. Let $\mathcal{T}$ an S-tree and $U \subseteq S$ any subset with $|U| > 1$. For each $T \in \mathcal{T}$ we have the quotient map $q_T : T \to C_T$. Suppose $U \subseteq T$ then we define $C_T^U = q_T(U)$, we then have $C_T^U \subseteq C_T$. If $T = \text{root}(U)$ in $\mathcal{T}$, that is the smallest element of $\mathcal{T}$ containing $U$ then by construction $|C_T^U| > 1$. This construction then gives us a map $s : P^+(S) \to \coprod_{T \in \mathcal{T}} P^+(C_T)$ defined by $s(U) = q_{\text{root}(U)}(U) \in P^+(C_{\text{root}(U)})$. It is easily seen that this map is surjective and has the following properties.

(1) $s(U) = s(V) \iff \text{root}(U) = T = \text{root}(V)$ and $q_T(U) = q_T(V)$.
(2) Let $T \in \mathcal{T}$ then for all $U \subseteq T$ we have $\text{root}(U) = T \iff |q_T(U)| > 1$.
(3) $s(U) \in P^+(C_T) \iff \text{root}(U) = T$.

For each $T \in \mathcal{T}$ and $U \subseteq T$ with $|U| > 1$ and associated set $C_T^U$ we define the induced quotient map $q_T^U : U \to C_T^U$ by $q_T^U = q_T|_U$ where $q_T$ is the quotient map $q_T : T \to C_T$

LEMMA 7.3.10. *For any $U \subseteq W$ with $T = \text{root}(W)$ then we have the following commutative diagrams where the vertical maps are the usual projections. The image of $r_W$ is $\overline{V}_W$ and the maps $r_W : V_{C_T^W} \to \overline{V}_W$ onto its image and $s_W : \overline{V}_W \to V_{C_T^W}$ are inverses.*

$$\begin{array}{ccc} V_{C_T^W} & \xrightarrow{r_W} & V_W \\ \pi_U^W \downarrow & & \downarrow \pi_U^W \\ V_{C_T^U} & \xrightarrow{r_U} & V_U \end{array} \qquad \begin{array}{ccc} \overline{V}_W & \xrightarrow{s_T} & V_{C_T^W} \\ \pi_U^W \downarrow & & \downarrow \pi_U^W \\ \overline{V}_U & \xrightarrow{s_U} & V_{C_T^U} \end{array}$$

*In particular if $\text{root}(U) \neq T$ then $|C_W^U| = 1$ and we have $\pi_U^W r_W = 0$*

PROOF. This is a particular case of lemma 6.4.4

□

CONSTRUCTION 7.3.11. For any S-tree $\mathcal{T}$ we define an injective map $\phi : \prod_{T \in \mathcal{T}} \overline{\mathcal{M}}_{C_T} \hookrightarrow \overline{\mathcal{M}}_S$ whose image is contained in $\overline{\mathcal{M}}_S(\mathcal{T})$ and $\text{type}(\phi(\prod \underline{M}_T)) = r(\prod \text{type}(\underline{M}_T))$.

Let $U \subseteq S$, $|U| > 1$ and $T = \text{root}(U)$ we define $M_U = r_U(M_{C_T^U})$. We then define $\underline{M} = \prod_U M_U$. I claim this gives us an element of $\overline{\mathcal{M}}_S$. For any $T \in \mathcal{T}$ define

$$\mathcal{L}_{S,T} = \{\, U \subseteq S \mid |U| > 1 \text{ and } \text{root}(U) = T \text{ in } \mathcal{T}\,\}$$

then

$$P^+(S) = \coprod_{T \in \mathcal{T}} \mathcal{L}_{S,T}$$

Now for every $U \subseteq V \subseteq S$ with $|U| > 1$ we have $U \in \mathcal{L}_{S,T}$ and $V \in \mathcal{L}_{S,W}$ where $T = \text{root}(U)$, $W = \text{root}(V)$ and $T \subseteq W$. Suppose first that $T = W$ then lemma 7.3.10 tells us that $M_V \le (\pi_U^V)^{-1} M_U$. If $T \subset W$ then again using lemma 7.3.10 we see that $\pi_U^V(M_V) = 0$ because in this case $|C_W^U| = 1$. Thus we have $\underline{M} \in \overline{\mathcal{M}}_S$. The construction of this map can be viewed more clearly by the following commutative diagram, the top maps are projections.

$$\begin{array}{ccc} \prod \overline{\mathcal{M}}_{C_T} \longrightarrow \overline{\mathcal{M}}_{C_{\text{root}(U)}} \longrightarrow PV_{s(U)} \\ \Big\downarrow \qquad\qquad\qquad\qquad\qquad \Big\downarrow r_U \\ \overline{\mathcal{M}}_S \longrightarrow PV_U \end{array}$$

We next prove that $\text{type}(\phi(\prod \underline{M}_{C_T})) = r(\prod \text{type}(\underline{M}_{C_T}))$. For each $T \in \mathcal{T}$ put $\mathcal{U}_T = \text{type}(\underline{M}_{C_T})$ and use the map $r : \prod_{T \in \mathcal{T}} \mathbb{T}_{C_T} \to \mathbb{T}_S$ to define an $S$-tree $\mathcal{U}$. Put $\underline{M} = \phi(\prod \underline{M}_{C_T})$, fix $T \in \mathcal{T}$ and choose $U \in q_T^{-1}(\mathcal{U}_T) \subseteq \mathcal{U}$. Then $\text{root}(U) = T$ in $\mathcal{T}$ and to prove $U \in \text{type}(\underline{M})$ we need to show that $\pi_U^V M_V = 0$ for all $U \subset V \subseteq S$. Since $U \subset V$ we have that $V \in \mathcal{L}_{S,W}$ where $W = \text{root}(V)$ in $\mathcal{T}$ and $T \subseteq W$. First suppose $T \subset W$ then as in the first paragraph $\pi_U^V M_V = 0$. Next suppose that $T = W$, since $U \subset V$ we see that $C_T^U \subseteq C_T^V$. Because $U \in q_T^{-1}(\mathcal{U}_T)$ we see $C_T^U \in \mathcal{U}_T$ and that $U$ is the maximal set with $q_T(U) = C_T^U$. Thus we see that $C_T^U \subset C_T^V$. Now put $X = C_T^U$ and $Y = C_T^V$ then as $X \in \mathcal{U}_T$ we see $\pi_X^Y M_Y = 0$. Thus by the first diagram of lemma 7.3.10 we have $\pi_U^V M_V = 0$. We have proven the condition for each $V \supset U$ thus $U \in \text{type}(\underline{M})$ and $\mathcal{U} \subseteq \text{type}(\underline{M})$.

Next choose a set $U \subseteq S$ so that $U$ is not in any $q_W^{-1}(\mathcal{U}_W)$ and put $T = \text{root}(U)$ in $\mathcal{T}$ so that $U \in \mathcal{L}_{S,T}$. If $C_T^U \notin \mathcal{U}_T$ then we can find $U \subset V \subseteq T$ with $C_T^U \subset C_T^V$, such that if we put $X = C_T^U$ and $Y = C_T^V$ then $\pi_X^Y M_Y = M_X$. Thus by lemma 7.3.10 we see that $\pi_U^V M_V = M_U$ and so $U \notin \text{type}(\underline{M})$. If $C_T^U \in \mathcal{U}_T$ put $V = q_T^{-1}(C_T^U)$ then as $V$ is the maximal set with $q_T(V) = C_T^U$ and $q_T(V) = q_T(U)$ we see that $V \supseteq U$. As $U \notin q_T^{-1}(\mathcal{T}_T)$ we see that $V \supsetneq U$ thus applying lemma 7.3.10 using the fact that $C_T^U = C_T^V$ we find that $\pi_U^V M_V = M_U$ and $U \notin \text{type}(\underline{M})$ therefore $r(\prod \mathcal{U}_T) = \text{type}(\underline{M})$ as claimed. The injectivity of this map is clear. It is now immediate from lemma 7.3.5 that $\text{type}(\underline{M}) \supseteq \mathcal{T}$. Thus the image of the map is contained in $\overline{\mathcal{M}}_S(\mathcal{T})$ as required.

CONSTRUCTION 7.3.12. For any S-tree $\mathcal{T}$ we construct a map $\theta : \overline{\mathcal{M}}_S(\mathcal{T}) \hookrightarrow \prod_{T \in \mathcal{T}} \overline{\mathcal{M}}_{C_T}$.

Consider a point $\underline{M} \in \overline{\mathcal{M}}_S(\mathcal{T})$. For each $T \in \mathcal{T}$, we must define an element $\underline{M}_{C_T} \in \overline{\mathcal{M}}_{C_T}$. With the previous notation we have $P^+(S) = \coprod_{T \in \mathcal{T}} \mathcal{L}_{S,T}$. Fix $T \in \mathcal{T}$, I first claim that for every $W \in \mathcal{L}_{S,T}$ and $\underline{M} \in \overline{\mathcal{M}}_S(\mathcal{T})$ we have $\pi_{W \cap V}^W(M_W) = 0$ for each $V \in M(\mathcal{T}, T)$ with $|W \cap V| > 1$ and therefore $M_W \in P\overline{V}_W$. To see this put $\mathcal{U} = \text{type}(\underline{M})$ and $X = \text{root}(W)$ in $\mathcal{U}$. Then as $\mathcal{U} \supseteq \mathcal{T}$ we have $V \subset X \subseteq T$, $\pi_W^X M_X = M_W$ and $\pi_V^X M_X = 0$. One then deduces the result. Then by construction we may use the second diagram of lemma 7.3.10 to define the following map. Let $\underline{M} \in \overline{\mathcal{M}}_S(\mathcal{T})$ and fix $T \in \mathcal{T}$ then for each $U \in \mathcal{L}_{S,T}$ we have $M_U \in P\overline{V}_U$ and define $M_{C_T^U} = s_T(M_U)$. We need to check this is well defined. Put $W = \text{root}(U)$ in $\mathcal{U}$ and choose another $V \in \mathcal{L}_{S,T}$ with $q_T(V) = C_T^U$. Then as $C_T^U = C_T^V$ one readily deduces that $\text{root}(V) = W$. From this we see our construction is well defined. Then for every such $T$ this define a element $\underline{M}_{C_T} = \prod M_{C_T^U}$. We are required to prove that $\underline{M}_{C_T} \in \overline{\mathcal{M}}_{C_T}$. This is now clear from the second diagram in lemma 7.3.10. Repeating this for every $T \in \mathcal{T}$ then gives as an element $\prod \underline{M}_{C_T} \in \prod \overline{\mathcal{M}}_{C_T}$. This completes the construction. It is a simple fact to verify that this is injective. One also readily verifies that over $\mathcal{M}_S(\mathcal{T})$ the map reduces on each component to

$$\begin{array}{ccc} \mathcal{M}_S(\mathcal{T}) & \longrightarrow U_{C_T} & \xrightarrow{\sim} \mathcal{M}_{C_T} \\ \downarrow & \downarrow & \\ \mathcal{M}_S & \longrightarrow PV_T & \end{array}$$

We are now in a position to prove proposition 7.3.6.

PROOF. It is clear that $\overline{\mathcal{M}}_S = \coprod_{\mathcal{T} \in \mathbb{T}_S} \mathcal{M}_S(\mathcal{T})$. We next prove that $\overline{\mathcal{M}}_S(\mathcal{T}) \cong \prod_{T \in \mathcal{T}} \overline{\mathcal{M}}_{C_T}$. By the previous constructions we have the map $\phi : \prod_{T \in \mathcal{T}} \overline{\mathcal{M}}_{C_T} \to \overline{\mathcal{M}}_S$ and the map $\theta : \overline{\mathcal{M}}_S(\mathcal{T}) \to \prod_{T \in \mathcal{T}} \overline{\mathcal{M}}_{C_T}$ with the image of $\phi$ contained in $\overline{\mathcal{M}}_S(\mathcal{T})$. Now by lemma 7.3.10 it is clear that these maps are inverses to each other. Thus the third claim is true. To prove that $\mathcal{M}_S(\mathcal{T}) \cong \prod_{T \in \mathcal{T}} \mathcal{M}_{C_T}$ we observe that the restricted map $\phi : \prod_{T \in C_T} \mathcal{M}_{C_T} \to \overline{\mathcal{M}}_S$ has its image in $\mathcal{M}_S(\mathcal{T})$ and $\theta : \mathcal{M}_S(\mathcal{T}) \to \prod_{T \in C_T} \mathcal{M}_{C_T}$. Finally we prove the last part of the claim that $\text{cl}(\mathcal{M}_S(\mathcal{T})) = \overline{\mathcal{M}}_S(\mathcal{T})$. By the previous part the map $\phi : \prod_{T \in \mathcal{T}} \overline{\mathcal{M}}_{C_T} \to \overline{\mathcal{M}}_S(\mathcal{T})$ restricts to a map $\phi : \prod_{T \in C_T} \mathcal{M}_{C_T} \to \mathcal{M}_S(\mathcal{T})$. By lemma 6.2.8 $\mathcal{M}_{C_T}$ is an open dense subset of $\overline{\mathcal{M}}_{C_T}$ thus $\prod_{T \in C_T} \mathcal{M}_{C_T}$ is an open dense subset of $\prod_{T \in \mathcal{T}} \overline{\mathcal{M}}_{C_T}$ and by lemma 7.3.8 the space $\overline{\mathcal{M}}_S(\mathcal{T})$ is Zariski closed and by definition contains $\mathcal{M}_S(\mathcal{T})$ thus,

$$\begin{aligned} \overline{\mathcal{M}}_S(\mathcal{T}) \supseteq \text{cl}(\mathcal{M}_S(\mathcal{T})) &= \text{cl}(\phi(\prod_{T \in C_T} \mathcal{M}_{C_T})) \\ &\supseteq \phi(\text{cl}(\prod_{T \in C_T} \mathcal{M}_{C_T})) \\ &= \phi(\prod_{T \in \mathcal{T}} \overline{\mathcal{M}}_{C_T}) \\ &= \overline{\mathcal{M}}_S(\mathcal{T}) \end{aligned}$$

The last argument works equally well in the Zariski or classical topology. This completes our results.

□

## 7.4. The natural maps from $\overline{\mathcal{M}}_S$ to $\overline{\mathcal{M}}_T$

Here we define natural projection maps $\pi_T^S : \overline{\mathcal{M}}_S \to \overline{\mathcal{M}}_T$. These will be analogous to the maps $\pi_T^S : \overline{\mathcal{X}}_S \to \overline{\mathcal{X}}_T$ that we consider later and are important in the study of $\overline{\mathcal{X}}_S$.

DEFINITION 7.4.1. let $T \subseteq S$ then we define the natural projection map $\pi_T^S : \overline{\mathcal{M}}_S \to \overline{\mathcal{M}}_T$ as follows. For every $U \subseteq T$ we define $\pi_T^S(\underline{M})_U = M_U$

LEMMA 7.4.2. *For every $U \subseteq T \subseteq S$ we have that the following diagram commutes*

$$\begin{array}{ccc} \overline{\mathcal{M}}_S & \xrightarrow{\pi_T^S} & \overline{\mathcal{M}}_T \\ & \searrow{\pi_U^S} & \downarrow{\pi_U^T} \\ & & \overline{\mathcal{M}}_U \end{array}$$

PROOF. The proof of this is clear from the definitions. □

LEMMA 7.4.3. *The map $\pi_T^S : \mathbb{T}_S \to \mathbb{T}_T$ of definition 7.2.4 commutes in the following diagram.*

$$\begin{array}{ccc} \overline{\mathcal{M}}_S & \xrightarrow{\text{type}} & \mathbb{T}_S \\ \pi_T^S \downarrow & & \downarrow \pi_T^S \\ \overline{\mathcal{M}}_T & \xrightarrow{\text{type}} & \mathbb{T}_T \end{array}$$

PROOF. Let $\underline{M} \in \overline{\mathcal{M}}_S$ put $\underline{N} = \pi_T^S(\underline{M})$ and $\mathcal{U} = \pi_T^S(\mathcal{T})$ where $\mathcal{T} = \text{type}(\underline{M})$. Let $W \in \mathcal{U}$ and $V \in \mathcal{T}$ so that $W = V \cap T$. Then we must show that $\pi_W^X M_X = 0$ for every $W \subset X \subseteq T$. For this it suffices to prove that in $\mathcal{T}$ we have $Y = \text{root}(X) \supset V$. Then $\pi_X^Y M_Y = M_X$ and $\pi_V^Y M_Y = 0$ because $V \in \mathcal{T}$ thus we deduce that $\pi_W^X M_X = 0$ so that $W \in \text{type}(\underline{N})$. Because $W \subset X$ there is some $x \in X \setminus W$. As $X \subseteq T$ we see that $x \notin S \setminus T$. Thus $x \notin W \amalg U = V$ for some $U \subseteq S \setminus T$, ie $x \in X \setminus V$ so that $X \nsubseteq V$. Put $Y = \text{root}(X)$ in $\mathcal{T}$ then as $W$ is non-empty and $W = V \cap T \subset X \subseteq Y$ we see that $Y \cap V$ is non-empty. Because $\mathcal{T}$ is a tree either $Y \subseteq V$ or $V \subseteq Y$. In the former case we would have $X \subseteq Y \subseteq V$ so that $X \subseteq V$ a contradiction, thus we deduce $V \subset Y = \text{root}(X)$.

Let $W \subseteq T$ with $|W| > 1$ and $W \notin \mathcal{U}$. Put $Y = \text{root}(W)$ in $\mathcal{T}$ then $\pi_W^Y M_Y = M_W$ therefore if we put $Z = Y \cap T$ we see that $Z \supseteq W$ so that $Z \in \mathcal{U}$. Then as $W \notin \mathcal{U}$ we see that $Z \supset W$ therefore we have $M_Z = \pi_Z^Y M_Y$ and we deduce that $\pi_W^Z M_Z = M_W$ therefore $W \notin \text{type}(\underline{N})$ thus we have shown that $\text{type}(\underline{N}) = \mathcal{U}$ and the diagram commutes.

$\square$

DEFINITION 7.4.4. Let $\mathcal{T}$ be an $S$-tree then for each $T \in \mathcal{T}$ and $i \in S$ we define the following $S+$ trees where $S+ = S \amalg \{+\}$.

$$A(\mathcal{T}, T) = \{ U+ \mid U \in \mathcal{T} \text{ and } U \supseteq T \} \amalg \{ V \mid V \in \mathcal{T} \text{ and } V \not\supseteq T \}$$
$$B(\mathcal{T}, T) = A(\mathcal{T}, T) \amalg \{T\}$$
$$C(\mathcal{T}, i) = \{ U+ \mid U \in \mathcal{T} \text{ and } i \in U \} \amalg \{ V \mid V \in \mathcal{T} \text{ and } i \notin V \} \amalg \{\{i, +\}\}$$

We then define the set of trees

$$E(\mathcal{T}) = \coprod_{T \in \mathcal{T}} \{A(\mathcal{T}, T)\} \amalg \coprod_{T \in \mathcal{T}} \{B(\mathcal{T}, T)\} \amalg \coprod_{i \in S} \{C(\mathcal{T}, i)\}$$

EXAMPLE 7.4.5.

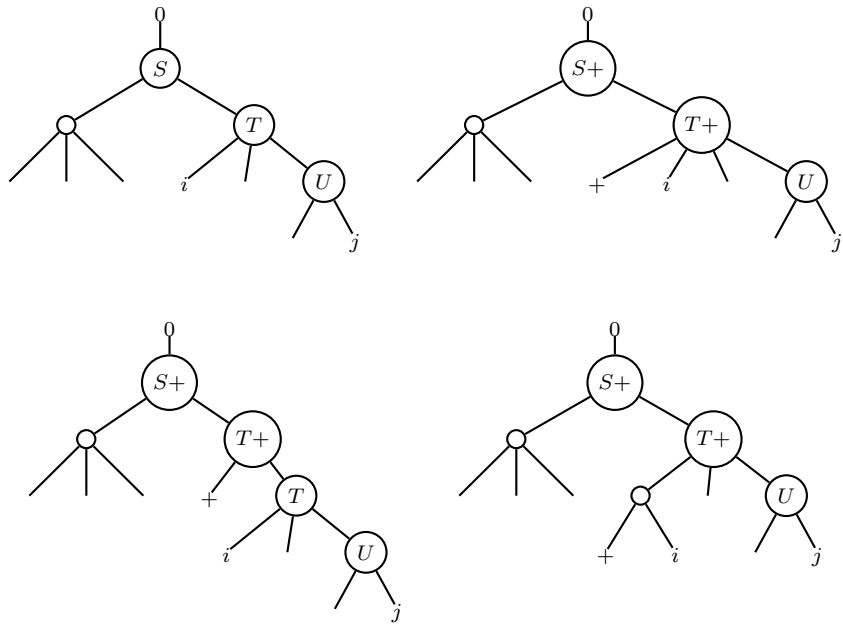

Top left is $\mathcal{T}$, top right $A(\mathcal{T}, T)$, bottom left $B(\mathcal{T}, T)$ and bottom right $C(\mathcal{T}, i)$

LEMMA 7.4.6. *For any finite set $S$ and any $\mathcal{T} \in \mathbb{T}_S$ let $\pi_S^{S+} : \mathbb{T}_{S+} \to \mathbb{T}_S$ be the usual map then $(\pi_S^{S+})^{-1}(\mathcal{T}) = E(\mathcal{T})$ and $|(\pi_S^{S+})^{-1}(\mathcal{T})| = |E(\mathcal{T})| = 2|\mathcal{T}| + |S|$*

PROOF. Suppose we can prove that $(\pi_S^{S+})^{-1}(\mathcal{T}) = E(\mathcal{T})$ then the second part of the claim is clear. To prove the first claim is not difficult but unilluminating. Instead we refer the reader to the diagram above for a picture of the process. □

DEFINITION 7.4.7. We introduce a grading on the set of trees $\mathbb{T}_S$ by

$$E_i^n = \{\, \mathcal{T} \in \mathbb{T}_S \mid |\mathcal{T}| = i \,\}$$

and define $h_n = \sum_{i>0} |E_i^n| t^i$, where $n = |S|$. Then $h_n(1) = |\mathbb{T}_S|$ and by the last lemma we have the following result.

LEMMA 7.4.8. $h_{n+1}(t) = t^2 h_n'(t) + t h_n'(t) + n t h_n(t)$ *and* $h_2 = t$ □

LEMMA 7.4.9. *For every $T \subseteq S$ and for any tree $\mathcal{T} \in \mathbb{T}_T$*

$$\begin{aligned}(\pi_T^S)^{-1}(\mathcal{M}_T(\mathcal{T})) &= \coprod_{\mathcal{U} \in (\pi_T^S)^{-1}(\mathcal{T})} \mathcal{M}_S(\mathcal{U}) \\ (\pi_T^S)^{-1}(\overline{\mathcal{M}}_T(\mathcal{T})) &= \bigcup_{\mathcal{U} \in (\pi_T^S)^{-1}(\mathcal{T})} \overline{\mathcal{M}}_S(\mathcal{U})\end{aligned}$$

PROOF. The first equality is immediate from lemma 7.4.3 and the second from the first using the description of $\overline{\mathcal{M}}_T(\mathcal{T})$. □

## 7.5. The projection from $\overline{\mathcal{M}}_{S+}$ to $\overline{\mathcal{M}}_S$

In this section we will investigate the local structure of the projection map $\pi : \overline{\mathcal{M}}_{S+} \to \overline{\mathcal{M}}_S$ over $\mathcal{M}_S$ and construct universal structure sections $\sigma : S_+ \times \overline{\mathcal{M}}_S \to \overline{\mathcal{M}}_{S+}$ where $S_+ = S \amalg \{0\}$. We will see later that the geometric fibre of this map $\pi^{-1}(\underline{M})$ together with the image of it structure sections at $\underline{M}$ is a stable $S$-curve and this will be an example of a family of $S$-curves over $\overline{\mathcal{M}}_S$. These notions are defined precisely in sections 3 and 6 of the next chapter. In particular it is the universal curve for our moduli space. $\overline{\mathcal{X}}_S$ will be the space of isomorphism classes of $S$-curves and we will have the induced map $\overline{\pi} : \overline{\mathcal{M}}_S \to \overline{\mathcal{X}}_S$. This will turn out to give us an isomorphism of varieties. We prove this further in the document.

LEMMA 7.5.1. *Let $\pi : \overline{\mathcal{M}}_{S+} \to \overline{\mathcal{M}}_S$ be the usual projection map then the restricted map $\pi : \pi^{-1}(\mathcal{M}_S) \to \mathcal{M}_S$ is a $\mathbb{C}P^1$ bundle. $\pi^{-1}(\mathcal{M}_S) = P(N_S \oplus \mathbb{C})$ where abusively we write $N_S$ for the restriction of the usual bundle $N_S$ over $\overline{\mathcal{M}}_S$ to the open set $\mathcal{M}_S$.*

PROOF. Consider the following commutative diagram

$$\begin{array}{ccc} \overline{\mathcal{M}}_{S+} & \xrightarrow{\pi} & \overline{\mathcal{M}}_S \\ {\scriptstyle p}\downarrow & & \downarrow{\scriptstyle p} \\ \overline{\mathcal{M}}_{\{S+,S\}} & \xrightarrow{\pi} & PV_S \end{array}$$

then we know $\pi : \overline{\mathcal{M}}_{\{S+,S\}} \to PV_S$ is the projectivization of a 2-dimensional algebraic vector bundle so it is a $\mathbb{C}P^1$ bundle and by lemma 6.2.3 the restricted map $p : \mathcal{M}_S \to U_S$ is an isomorphism with $\mathcal{M}_S = p^{-1}(U_S)$. Next put $X_S = \pi^{-1}(\mathcal{M}_S)$ and $W_S = \pi^{-1}(U_S)$ then $p^{-1}(W_S) = X_S$ and to prove the claim it will be sufficient to show that the restricted map $p : X_S \to W_S$ is an isomorphism. That is we need to construct an inverse for the restricted map $p : X_S \to W_S$. Let $N_S \times N_{S+} \in W_S$ and $T \subseteq S$ with $|T| > 2$. Then it is clear we may define $M_S = \rho_T^S N_S$ as $N_S \in U_S$ and $M_{T+} = \rho_{T+}^{S+} M_{S+}$. For any set $U \subseteq S+$ with $|U| = 2$ there is only one point in $PV_U$ so that $M_U$ is uniquely defined. We then define $\underline{M} = \prod_U M_U$. We see that this defines an injective map $r : W_S \to X_S$ and by construction its image $r(\underline{N})$ is the unique point in $\overline{\mathcal{M}}_S$ determined by $\underline{N}$ which is inverse to $p$. This completes the proof. $\square$

CONSTRUCTION 7.5.2. Let $\pi : \overline{\mathcal{M}}_{S+} \to \overline{\mathcal{M}}_S$ be the usual projection map then we define sections $\sigma : S_+ \times \overline{\mathcal{M}}_S \to \overline{\mathcal{M}}_{S+}$. For each $i \in S_+$ we define a section $\sigma_i : \overline{\mathcal{M}}_S \to \overline{\mathcal{M}}_{S+}$ as follows. For $i \in S$ consider the tree $\mathcal{T}_i = \{S+, \{i, +\}\}$. Proposition 7.3.6 tells us that $\overline{\mathcal{M}}_{S+}(\mathcal{T}_i) \cong \overline{\mathcal{M}}_{C_{S+}} \times \overline{\mathcal{M}}_{\{i,+\}} \cong \overline{\mathcal{M}}_S \times \text{pt} \cong \overline{\mathcal{M}}_S$ where the second isomorphism is induced from the bijection $b : S \to C_{S+}$ given by $b(j) = \overline{j}$. This gives us an isomorphism $\alpha_i : \overline{\mathcal{M}}_{S+}(\mathcal{T}_i) \to \overline{\mathcal{M}}_S$. One readily checks that the map implicit in this is just the inclusion $j_i : \overline{\mathcal{M}}_{S+}(\mathcal{T}_i) \to \overline{\mathcal{M}}_{S+}$ composed with the projection $\pi : \overline{\mathcal{M}}_{S+} \to \overline{\mathcal{M}}_S$. Thus we may define a section by $\sigma_i = j_i \circ \alpha_i^{-1}$. For $i = 0$ we define the tree $\mathcal{T}_0 = \{S+, S\}$. Then again by proposition 7.3.6 we obtain $\overline{\mathcal{M}}_{S+}(\mathcal{T}_0) \cong \overline{\mathcal{M}}_S \times \overline{\mathcal{M}}_{C_{S+}} \cong \overline{\mathcal{M}}_S \times \text{pt} \cong \overline{\mathcal{M}}_S$. This gives us the isomorphism $\alpha_0 : \overline{\mathcal{M}}_{S+}(\mathcal{T}_0) \to \overline{\mathcal{M}}_S$. Again one checks $\alpha_0 = \pi j_0$ thus we define a section $\sigma_0 = j_0 \circ \alpha_0^{-1}$. This completes our constructions and gives us a map $\sigma : S_+ \times \overline{\mathcal{M}}_S \to \overline{\mathcal{M}}_{S+}$.

REMARK 7.5.3. We will later see that the pair $\pi : \overline{\mathcal{M}}_{S+} \to \overline{\mathcal{M}}_S$ and $\sigma : S_+ \times \overline{\mathcal{M}}_S \to \overline{\mathcal{M}}_{S+}$ form the universal curve for our moduli space.

CHAPTER 8

# Definitions and properties of the space $\overline{\mathcal{X}}_S$

In this chapter we introduce the moduli space $\overline{\mathcal{X}}_S$ and using certain established results we consider the structure of this space. The most important of these will be Mumford's main theorem and the natural morphisms $\pi : \overline{\mathcal{X}}_{S+} \to \overline{\mathcal{X}}_S$ constructed by Deligne and Mumford. We begin by defining a generic point of $\overline{\mathcal{X}}_S$ and put $\mathcal{X}_S$ to be the Zariski open space consisting of all the generic points. We will then show that there is a natural isomorphism $\phi_S : \mathcal{X}_S \to U_S$. Next we proceed to introduce the non generic elements of our space. We then explain the construction of the tree associated to a stable $S$-curve and use this to stratify the space $\overline{\mathcal{X}}_S$. This is the approach we used on $\overline{\mathcal{M}}_S$ in chapter 7. After introducing the notion of a family of $S$-curves over a scheme $X$ we state Mumford's main theorem 8.6.3. We will then use this to construct a regular morphism of varieties $\theta_S : \overline{\mathcal{X}}_S \to PV_S$ that is a birational equivalence. In particular the restricted map $\theta_S : \mathcal{X}_S \to U_S$ will be an isomorphism of varieties that agrees with $\phi_S$. In the last section of this chapter we construct an isomorphism $\overline{\theta}_S : \overline{\mathcal{X}}_S \to \overline{\mathcal{M}}_S$ induced from the maps $\theta_T : \overline{\mathcal{X}}_T \to PV_T$. This will preserve certain natural properties such as the tree types as specified by the bijections under which we consider the various notions of trees equivalent. Parts of this chapter are closely related to unpublished notes by Professor N.P. Strickland. In this chapter we will write $\widetilde{V}_S$ for the vector space $F(S, \mathbb{C})$ for brevity, $S+ = S \amalg \{pt\}$ and $S_+ = S \amalg \{0\}$.

## 8.1. Introduction to moduli spaces

A moduli spaces typically consists of two pieces of data, a class of objects and a notion of an algebraic family of these objects over a scheme. Let $\mathcal{C}$ be a class of objects and for any base scheme $B$ let $S(B)$ be the set of all families over $B$. Let $\sim$ be an equivalence relation on $S(B)$ and consider the functor $F$ from the category of schemes to the category of sets given by $F(B) = S(B)/\sim$. $F$ is called the moduli functor of our moduli problem. Suppose $F$ is representable by a scheme $\mathcal{M}$ then we say the scheme is a fine moduli space for the moduli problem $F$.

## 8.2. Generic $S$-curves

DEFINITION 8.2.1. A *generic $S$-curve* is a pair $(C, x)$ where $C$ is an algebraic curve isomorphic to $\mathbb{C}P^1$ and $x : S_+ \to C$ is an injective map. We consider two generic $S$-curves $(C, x)$ and $(D, y)$ *isomorphic* if there is a algebraic isomorphism $s : C \to D$ such that $s \circ x = y$, that is the distinguished marked points of $C$ are sent to those of $D$ under the isomorphism $s$ with their order preserved. This is clearly an equivalence relation and we write $[C, x]$ for the isomorphism type of a generic $S$-curve $(C, x)$ under our equivalence. We define $\mathcal{X}_S$ to be the set of isomorphism classes of generic $S$-curves.

EXAMPLE 8.2.2. This is a generic curve with 5 marked points.

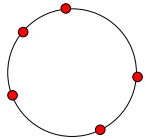

We next need to show how to give $\mathcal{X}_S$ the structure of a variety. We will do this by identifying it with a Zariski open subset of $PV_S$. This identification will be important for our understanding of $\overline{\mathcal{X}}_S$ which we consider later. Recall in 2.0.12 that we defined the Zariski open subset $U_S$ of $PV_S$ consisting of the 'injective functions'.

PROPOSITION 8.2.3. *There is a natural bijection $\phi_S : \mathcal{X}_S \to U_S$.*

PROOF. Let $[C, x] \in \mathcal{X}_S$ be a generic $S$-curve and $(D, u)$ be a representative for this equivalence class. Then we can choose an isomorphism $f : D \to \mathbb{C}P^1$ where we consider $\mathbb{C}P^1 = \mathbb{C} \cup \{\infty\}$ such that $f(u(0)) = \infty$. Since $u : S_+ \to D$ is injective and by construction $f(u(i)) \in \mathbb{C}$. Then we may define a function $z : S \to \mathbb{C}$ by $z(i) = f(u(i))$. We observe that this function is injective and so we have an induced element $\overline{z} \in U_S$. We must check that this construction is well defined. Let $(E, v)$ be another representative of $[C, x]$ and $g : E \to \mathbb{C}P^1$ be another isomorphism with $g(v(0)) = \infty$. Define an injective map $w : S \to \mathbb{C}$ in the same way, that is for each $i \in S$　　$w(i) = g(v(i))$. We can choose an isomorphism $\alpha : E \to D$ with $u = \alpha \circ v$ and put $h = f \alpha g^{-1}$, this is a map $h : \mathbb{C}P^1 \to \mathbb{C}P^1$ with $h(\infty) = \infty$ thus $h = az + b$ for some $(a, b) \in \mathbb{C}^\times \times \mathbb{C}$. Then for each $i \in S$ we have $z(i) = fu(i) = f\alpha v(i) = hgv(i) = hw(i) = aw(i) + b$ thus $\overline{z} = \overline{w}$ in $U_S$ and our map is well defined. There is an inverse for this map and so $\phi_S$ is a bijection as claimed. $\square$

Given the bijective map $\phi_S : \mathcal{X}_S \to U_S$ there is evidently a unique way in which to regard $\mathcal{X}_S$ as a variety for which $\phi_S$ is an isomorphism.

## 8.3. Stable $S$-curves

Here we define the compactification $\overline{\mathcal{X}}_S$ of $\mathcal{X}_S$. To do this we need a more general marked curve, this we define next.

DEFINITION 8.3.1. A *stable $S$-curve* is a pair $(C, x)$, where $C$ is a (possibly singular) algebraic curve over $\mathbb{C}$ and $x : S_+ \to C$ is an injective map, such that certain conditions are satisfied. To formulate these, we say that a point in $C$ is *marked* if it is in the image of $x$, and *special* if it is either singular or marked. The conditions are as follows.

(a) $C$ is reduced and connected, and any singular points are ordinary double points. Equivalently, the completion of the local ring at any point is isomorphic either to $\mathbb{C}[\![x]\!]$ (for a smooth point) or $\mathbb{C}[\![x,y]\!]/xy$ (for a singular point).
(b) All the marked points are nonsingular.
(c) Each irreducible component of $C$ is isomorphic to $\mathbb{C}P^1$, and contains at least three special points.
(d) $H^1(C; \mathcal{O}_C) = 0$.

We consider two $S$-curves $(C, x)$ and $(D, y)$ to be *isomorphic* if there is an algebraic isomorphism $s : C \to D$ such that $s \circ x = y$, that is $s$ send marked points to marked points and preserves the order and define $\overline{\mathcal{X}}_S$ to be the set of isomorphisms classes of stable $S$-curves, observe $\mathcal{X}_S \subset \overline{\mathcal{X}}_S$.

EXAMPLE 8.3.2. This is a non-generic curve with 12 marked points

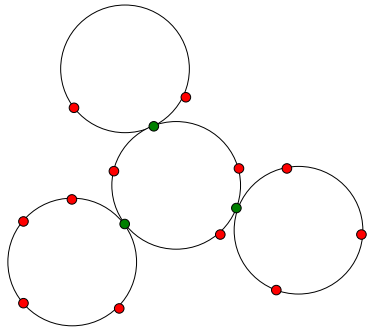

### 8.4. The tree associated to elements of $\overline{\mathcal{X}}_S$

There is a natural way to partition the set $\overline{\mathcal{X}}_S$ that reflects the combinatorics of an $S$-curve $(C, x)$. This is the notion of the tree type associated to a curve and is an important construction for the study of $\overline{\mathcal{X}}_S$. We next discus how to construct the associated tree, this will require several intermediate results.

LEMMA 8.4.1. *Let $(C, x)$ be a stable $S$-curve. Then every global regular function on $C$ is constant, so $H^0(C; \mathcal{O}_C) = \mathbb{C}$.*

PROOF. This is well-known for $\mathbb{C}P^1$, and every irreducible component of $C$ is a copy of $\mathbb{C}P^1$, so every regular function on $C$ is constant on irreducible components and thus takes only finitely many values. As $C$ is connected, the claim follows. □

LEMMA 8.4.2. *Let $s$ be the number of singular points in $C$. Then $H^1(C; \mathcal{O}_C) = 0$ if and only if there are precisely $s + 1$ irreducible components.*

PROOF. Let $C_1, \ldots, C_t$ be the irreducible components of $C$, so the claim is that $t = s + 1$ if and only if $H^1(C, \mathcal{O}_C) = 0$. Put $\widetilde{C} = \coprod_{i=1}^t C_i \simeq \coprod_{i=1}^t \mathbb{C}P^1$, and let $q : \widetilde{C} \to C$ be the obvious map from the disjoint union to the actual union. Now let $D \subset C$ be the finite set of singular points, and let $i : D \to C$ be the inclusion. Each point $d \in D$ has two preimages in $\widetilde{C}$; we choose one and call it $\sigma_0(d)$, and then we call the other one $\sigma_1(d)$. Given an open set $U \subseteq C$ and a regular function $f$ on $q^{-1}U$, we define $\delta(f) : D \cap U \to \mathbb{C}$ by $\delta(f)(d) = f(\sigma_0(d)) - f(\sigma_1(d))$. This construction gives a map $\delta : q_*\mathcal{O}_{\widetilde{C}} \to i_*\mathbb{C}$ of sheaves on $C$, which is easily seen to be an epimorphism with kernel $\mathcal{O}_C$. As $H^1(\widetilde{C}; \mathcal{O}_{\widetilde{C}}) = 0$ we have the following four term exact sequence,

$$H^0(C; \mathcal{O}_C) \rightarrowtail H^0(C; q_*\mathcal{O}_{\widetilde{C}}) \to H^0(C; i_*\mathbb{C}) \twoheadrightarrow H^1(C; \mathcal{O}_C),$$

or equivalently

$$H^0(C; \mathcal{O}_C) \rightarrowtail H^0(\widetilde{C}; \mathcal{O}_{\widetilde{C}}) \to H^0(D; \mathbb{C}) \twoheadrightarrow H^1(C; \mathcal{O}_C),$$

that is

$$\mathbb{C} \hookrightarrow \bigoplus_{i=1}^{t} \mathbb{C} \to \bigoplus_{d \in D} \mathbb{C} \twoheadrightarrow H^1(C; \mathcal{O}_C)$$

then counting dimensions we see $\dim H^1(C; \mathcal{O}_C) = 1 - t + s$. Thus $\dim H^1(C, \mathcal{O}_C) = 0 \iff t = s + 1$ as required.

$\square$

We next define a graph $G = G[C, x]$, more precisely an isomorphism class of a simplicial complex whose geometric realization is the associated graph to an isomorphism class of an $S$-curve $(C, x)$. We then proceed to prove that this graph is in fact a tree. Of-course our construction is equivalent to the ordinary graph associated to an $S$-curve that is well known to be a tree. Here we clarify this fact.

CONSTRUCTION 8.4.3. We define a graph $G = G(C, x)$ as follows. We let $V_0$ be the set of marked points (so $x$ gives a bijection $S_+ \simeq V_0$) and we let $V_1$ be the set of irreducible components of $C$. The vertex set of our graph is $V = V_0 \amalg V_1$. The vertices in $V_0$ are called *external*, and those in $V_1$ are called *internal*. For each marked point we have a corresponding external vertex, and also an internal vertex corresponding to the component of $C$ containing the marked point. The graph has an edge joining these two vertices; these are called *external edges*. Next, every singular point is a double point and so lies in the intersection of two irreducible components. We give the graph an edge joining the two corresponding internal vertices. (In principle, this could lead to multiple edges between the same two vertices.)

Note that we have one edge for each special point. If $C_0$ is an irreducible component, then the edges incident on the corresponding internal vertex biject with the special points in $C_0$, so the valence is at least three.

The combinatorial graph associated to the curve $[C, x] \in \overline{\mathcal{X}}_S$ is then defined to be the isomorphism class of the graph $G$. It is clear that this construction is well defined modulo isomorphism classes.

EXAMPLE 8.4.4. This example shows how we construct the graph.

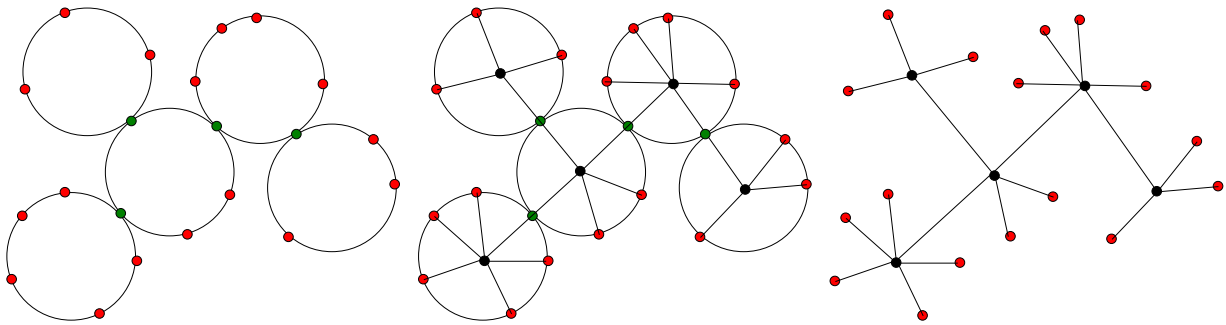

PROPOSITION 8.4.5. *The graph $G(C, x)$ is a tree (and thus has no multiple edges).*

PROOF. Let $|G|$ be the geometric realization of $G$, obtained by taking a copy of the unit interval for each edge and making the obvious gluings. We first claim that $|G|$ is connected. If not, choose a disconnection $|G| = X_0 \amalg X_1$ and let $E_i$ be the set of components of $C$ such that the corresponding internal edges lie in $X_i$. Now let $C_i$ be the union of the components in $E_i$, and observe that $C_i \neq \emptyset$ and $C = C_0 \amalg C_1$, contradicting the connectedness of $C$.

Now put $n = |S|$ and let $s$ be the number of singular points, so there are $s + 1$ irreducible components. There are then $n + 1$ external vertices, $n + 1$ external edges, $s + 1$ internal vertices, and $s$ internal edges. The Euler characteristic of $|G|$ is thus $\chi(|G|) = (n + 1 + s + 1) - (n + 1 + s) = 1$, which implies that $|G|$ is contractible and thus that $G$ is a tree. □

## 8.5. The regular map $\theta_S : \overline{\mathcal{X}}_S \to PV_S$

We have already shown that there is an isomorphism $\phi_S : \mathcal{X}_S \to U_S$. In order to continue our analysis of $\overline{\mathcal{X}}_S$ we wish to extend $\phi_S$ to a map $\theta_S : \overline{\mathcal{X}}_S \to PV_S$ of all of $\overline{\mathcal{X}}_S$. Clearly should such a map exists it must be unique as $PV_S$ is separated and $\mathcal{X}_S$ has the same dimension as $\overline{\mathcal{X}}_S$. Here we define $\theta_S$ as a map of sets. Later we will give a more careful construction that works for parameterized families of curves and thus gives a map of schemes. This map should be compared with Kapranov's construction in [**7**].

DEFINITION 8.5.1. Let $(C, x)$ be an $S$-curve. We define $M = M(C, x)$ to be the set of rational functions on $C$ with at worst a simple pole at $x(0)$, and no poles elsewhere.

REMARK 8.5.2. This is clearly a vector space containing the set $\mathbb{C}$ of constant functions. We will need two lemmas about $M$, the proofs are given later in this section.

LEMMA 8.5.3. For each $S$-curve $(C, x)$ we have $\dim_\mathbb{C}(M(C, x)) = 2$.

LEMMA 8.5.4. Let $(C, x)$ be a stable $S$-curve. Suppose that $f \in M(C, x)$ and $f$ vanishes at all marked points other than $x(0)$. Then $f = 0$.

CONSTRUCTION 8.5.5. Let $[C, x]$ be a stable $S$-curve and $(D, u)$ be a representative. Then we define a linear map $\sigma : M(D, u) \to F(S, \mathbb{C})$ as follows. For any $f \in M(D, u)$ and $i \in S$ define $\sigma(f)(i) = f(u(i))$, this makes sense as $u$ is injective and $f$ has a pole only at $x(0)$. By lemma 8.5.4 this map is injective and sends constants to constants, and thus induces an injective linear map $\overline{\sigma} : \overline{M}(C, x) \to V_S$ with image $L \in PV_S$ say. We need to check that this construction is well defined. Let $(E, v)$ be any other representative and given $g \in M(E, v)$ define $\sigma(g)(i) = g(v(i))$. Let $h : D \to E$ be an isomorphism with $h \circ u = v$ then clearly the map $h^* : M(E, v) \to M(D, u)$ defined by $h^*(g) = gh$ is an isomorphism of vector spaces. Then put $f = h^*(g)$ then $\sigma(f)(i) = ghu(i) = gv(i) = \sigma(g)(i)$. This shows that the images of the vector spaces $M(D, u)$ and $M(E, v)$ under $\sigma$ are the same thus our map is well defined. We now define a map $\theta_S : \overline{\mathcal{X}}_S \to PV_S$ by $\theta_S([C, x]) = L$.

PROPOSITION 8.5.6. If $[C, x]$ is a generic $S$-curve, then the above construction is the same as in proposition 8.2.3.

PROOF. Let $(C, x)$ be a representative for $[C, x]$ and choose an isomorphism $f : C \to \mathbb{C}P^1$ with $f(x(0)) = \infty$ as in Proposition 8.2.3. As $f$ is an isomorphism, we see that it is regular away from $x(0)$ and has a simple pole at $x(0)$ thus $\{1, f\}$ is a basis for $M = M(C, x)$. It follows that if we put $g = f \circ x$ and $\overline{g}$ to be the image in $PV_S$ then $\{\overline{g}\}$ is a basis for $L = \overline{\sigma}(M)$, and the claim follows directly. $\square$

We now prove the two deferred lemmas. We will make some detours along the way, to pick up results that will be useful later.

LEMMA 8.5.7. *Let $f$ be a rational function on $\mathbb{C}P^1 = \mathbb{C} \cup \{\infty\}$ with at worst a simple pole at $\infty$ and no poles elsewhere. Suppose also that $f$ vanishes at two distinct points in $\mathbb{C}$. Then $f = 0$.*

PROOF. As $f(z)$ is rational with no poles in $\mathbb{C}$, it must be polynomial. As it has at worst a simple pole at infinity, it must have degree at most one. The claim follows easily. □

PROOF OF LEMMA 8.5.4. Let $G$ be the tree defined in Construction 8.4.3, and let $v_0$ be the external vertex corresponding to $x(0)$. Let $C_0$ be the irreducible component in $C$ containing $x(0)$, and let $e_0$ be the corresponding external edge. Let $v_1$ be the internal vertex corresponding to $C_0$, so $e_0$ joins $v_0$ to $v_1$. Every edge $e \neq e_0$ corresponds to a special point in $C \setminus \{x(0)\}$, and we write $f(e)$ for the value of $f$ at that point. Every external vertex $v \neq v_0$ corresponds to a marked point different from $x(0)$, so $f(v) = 0$. Every internal vertex $v \neq v_1$ corresponds to a component on which $f$ is regular and thus constant, with value $f(v)$ say. If an edge $e$ and a vertex $v$ are incident, it is clear that $f(e) = f(v)$. Using the fact that $G$ is a tree, and working inwards from the external vertices, we see that $f(e) = 0$ for all $e \neq e_0$ and $f(v) = 0$ for all $v \notin \{v_0, v_1\}$. In particular, we see that $f$ vanishes at all special points other than $x(0)$. By assumption, there are at least two such points in the component $C_0$, and it follows from Lemma 8.5.7 that $f = 0$ on $C_0$ as well. □

Now let $D$ be the set of marked points other than $x(0)$ (so $x$ gives a bijection $S \to D$) and let $j: D \to C$ be the inclusion. Let $\mathcal{J}_0$ be the ideal sheaf of functions vanishing at $x(0)$, let $\mathcal{J}$ be the ideal sheaf of functions vanishing on $D$, and put $\mathcal{K} = \mathcal{J} \otimes \mathcal{J}_0^{-1}$. Thus $\mathcal{J}_0^{-1}$ is the sheaf of functions with at worst a simple pole at $x(0)$, and regular elsewhere, so $M = H^0(C; \mathcal{J}_0^{-1})$. Moreover, $\mathcal{K}$ is the subsheaf of such functions that vanish on $D$. Lemma 8.5.3 says that $\dim(M) = \dim(H^0(C; \mathcal{J}_0^{-1})) = 2$, so it is part of the following result

PROPOSITION 8.5.8. *If we write $h^i(\mathcal{F}) = \dim_{\mathbb{C}}(H^i(C;\mathcal{F}))$, then*

$$\begin{aligned} h^0(\mathcal{O}_C) &= 1 & h^1(\mathcal{O}_C) &= 0 \\ h^0(j_*\mathcal{O}_D) &= n & h^1(j_*\mathcal{O}_D) &= 0 \\ h^0(\mathcal{J}_0^{-1}) &= 2 & h^1(\mathcal{J}_0^{-1}) &= 0 \\ h^0(\mathcal{K}) &= 0 & h^1(\mathcal{K}) &= n-2. \end{aligned}$$

PROOF. Lemma 8.4.1 says that $h^0(\mathcal{O}_C) = 1$, and we are given that $h^1(\mathcal{O}_C) = 0$. As $j_*\mathcal{O}_D$ is a skyscraper sheaf, it is standard that $H^0(C; j_*\mathcal{O}_D) = \bigoplus_{d \in D} \mathbb{C} = \widetilde{V}_S$ and $H^1(C; j_*\mathcal{O}_D) = 0$, so $h^*(j_*\mathcal{O}_D)$ is as described. Next, put $\mathcal{N} = \mathcal{J}_0^{-1}/\mathcal{O}_C$. This is another skyscraper sheaf, with a one-dimensional stalk at $x(0)$, so $h^0(\mathcal{N}) = 1$ and $h^1(\mathcal{N}) = 0$. The short exact sequence $\mathcal{O}_C \to \mathcal{J}_0^{-1} \to \mathcal{N}$ gives a six-term sequence

$$\mathbb{C} \to M \to H^0(C;\mathcal{N}) \to 0 \to H^1(C;\mathcal{J}_0^{-1}) \to 0,$$

showing that $h^1(\mathcal{J}_0^{-1}) = 0$ and $\dim(M) = h^0(\mathcal{J}_0^{-1}) = h^0(\mathcal{O}_C) + h^0(\mathcal{N}) = 2$. Finally, Lemma 8.5.4 tells us that $h^0(\mathcal{K}) = 0$. There is a short exact sequence $\mathcal{K} \to \mathcal{J}_0^{-1} \to j_*\mathcal{O}_D$, giving a six-term sequence

$$0 \to M \to \widetilde{V}_S \to H^1(C;\mathcal{K}) \to 0 \to 0.$$

It follows that $h^1(\mathcal{K}) = \dim(\widetilde{V}_S) - \dim(M) = n - 2$. □

## 8.6. Families of stable $S$-curves

To make $\overline{\mathcal{X}}_S$ into a variety, the key point is to decide what we mean by an algebraically varying family of stable $S$-curves, parameterized by a scheme $X$. We should certainly have a stable curve $(C_a, x_a)$ for each point $a \in X$. Given these, we can form a set

$$C = \{(a, c) \mid a \in X \text{ and } c \in C_a\}.$$

We then have a map $\pi \colon C \to X$ given by $\pi(a, c) = a$, and a map $x \colon S_+ \times X \to C$ given by $x(i, a) = (a, x_a(i))$, so that $\pi(x(i, a)) = a$. It is certainly natural to require that $C$ should be a scheme, and that the functions $\pi$ and $x$ should be maps of schemes. One also needs some other technical conditions to make the theory work smoothly. The full definition is as follows.

DEFINITION 8.6.1. Let $X$ be a locally Noetherian scheme over $\mathbb{C}$. A *stable S-curve over* $X$ is a scheme $C$ equipped with maps $S_+ \times X \xrightarrow{x} C \xrightarrow{\pi} X$ of schemes such that

(a) $\pi$ is flat and proper
(b) $x$ is a closed inclusion
(c) $\pi \circ x$ is just the projection $S_+ \times X \to X$ (so for each $a \in X$ we have a fibre $C_a = \pi^{-1}\{a\}$ and a map $x_a : S_+ \to C_a$)
(d) each pair $(C_a, x_a)$ is a stable $S$-curve.

For any stable $S$-curve $C$ over $X$, we can define a function $\gamma_C : X \to \overline{\mathcal{X}}_S$ by

$$\gamma_C(a) = \text{ the isomorphism class of } (C_a, x_a) \in \overline{\mathcal{X}}_S.$$

A morphism of stable $S$-families $(C, X, \pi, x)$ and $(D, Y, \pi, y)$ is defined to be morphisms $g : C \to D$ and $h : X \to Y$ such that the following diagrams commute.

$$\begin{array}{ccc} C & \xrightarrow{g} & D \\ \pi \downarrow & & \downarrow \pi \\ X & \xrightarrow{h} & Y \end{array} \qquad \begin{array}{ccc} C & \xrightarrow{g} & D \\ x \uparrow & & \uparrow y \\ S_+ \times X & \xrightarrow{h} & S_+ \times Y \end{array}$$

A morphism is an isomorphism if $g$ and $h$ are. A morphism of families over the same base is a morphism of families with $X = Y$ and $h : X \to Y$ is the identity map.

The Noetherian hypothesis is not essential, but is included for technical convenience.

EXAMPLE 8.6.2. Take $S = \{1, 2, 3\}$ and $X = \mathbb{C} \setminus \{1\} = \text{spec}(\mathbb{C}[a][(a-1)^{-1}])$. Put

$$C = \{(a, [x : y : z]) \in X \times \mathbb{C}P^2 \mid xy = az^2\}$$

and $\pi(a, [x : y : z]) = a$. In other words, the fibre $C_a$ over a point $a \in X$ is the projective closure of the hyperbola $xy = a$. As $C$ is a closed subscheme of $X \times \mathbb{C}P^2$, the projection $\pi$ is certainly proper. We can rewrite the defining equation in terms of the variable $u = x - y$ as $y^2 + uy - az^2 = 0$, showing that the homogeneous coordinate ring of $C$ is a free module over $\mathcal{O}_X[u, z]$ with basis $\{1, y\}$; it follows easily that $\pi$ is also flat.

We then define $p : S_+ \times X \to S$ by

$$p(0, a) = (a, [1 : 0 : 0])$$
$$p(1, a) = (a, [0 : 1 : 0])$$
$$p(2, a) = (a, [a : 1 : 1])$$
$$p(3, a) = (a, [1 : a : 1]).$$

We removed the point $a = 1$ from $X$ to ensure that $p(2, a)$ is never equal to $p(3, a)$. Given this, one can check that $p$ is injective and that it gives an isomorphism of $S_+ \times X$ with the closed subscheme given by the equation $(x + y - (1 + a)z)z = 0$. Thus, $p$ is a closed inclusion. It is clear that $\pi \circ p : S_+ \times X \to X$ is just the projection. For $a \neq 0$ we can define an isomorphism $\mathbb{C}P^1 \to C_a$ by $[s : t] \mapsto [s^2 : at^2 : st]$, so $C_a$ is a generic $S$-curve. In the case $a = 0$ the curve $C_0$ is $\{[x : y : z] \in \mathbb{C}P^2 \mid xy = 0\}$. The irreducible components are $C_0' = \{[x : 0 : z] \mid [x : z] \in \mathbb{C}P^1\}$ and $C_0'' = \{[0 : y : z] \mid [y : z] \in \mathbb{C}P^1\}$, both of which are isomorphic to $\mathbb{C}P^1$. The components intersect only at the point $c = [0 : 0 : 1]$. The part of $C_0$ where $z = 1$ is an affine open neighborhood of $c$, with equation $xy = 0$ in $\mathbb{C}^2$. The completed local ring is thus $\mathbb{C}[\![x, y]\!]/xy$, so we have an ordinary double point. This shows that $C_0$ is a stable $S$-curve, so $C$ is a stable $S$-curve over $X$.

Results of Deligne and Mumford can be summarized as follows:

THEOREM 8.6.3. *One can make the set $\overline{\mathcal{X}}_S$ into a variety, and construct a universal stable $S$-curve $\mathcal{C}_S$ over $\overline{\mathcal{X}}_S$, such that*

(a) *The classifying map $\gamma_{\mathcal{C}_S} : \overline{\mathcal{X}}_S \to \overline{\mathcal{X}}_S$ is just the identity.*
(b) *For any locally Noetherian scheme $X$ over $\mathbb{C}$, and any stable $S$-curve $C$ over $X$, there is a unique map $\widetilde{\gamma}_C : X \to \overline{\mathcal{X}}_S$ of schemes such that $\widetilde{\gamma}_C^* \mathcal{C}_S$ is isomorphic to $C$. Moreover:*
   (i) *The induced map of complex points is just $\gamma_C$ as defined previously.*
   (ii) *The isomorphism $\widetilde{\gamma}_C^* \mathcal{C}_S \simeq C$ is unique.* □

As far as possible, we will use the universal property stated above, rather than any of the various constructions of $\overline{\mathcal{X}}_S$. The first exercise is to prove the following result:

PROPOSITION 8.6.4. *The inclusion $\mathcal{X}_S \to \overline{\mathcal{X}}_S$ is actually a morphism of varieties.*

PROOF. We have by lemma 4.1.7 that $\pi : \overline{\mathcal{M}}_{\{S+,S\}} \to PV_S$ is the projectivization of a Zariski locally trivial vector bundle of rank 2. We need to construct an $S$-curve over $U_S$. Put $C_S = \pi^{-1}(U_S)$ and write $\pi : C_S \to U_S$ for the restricted projection map. This is clearly flat and proper. We next define sections $x : S_+ \times U_S \to C_S$ as follows. Let $i \in S$ and $q_i : S+ \to S$ be the quotient map defined by $q_i(j) = j$ if $j \neq +$ and $q_i(+) = i$. Then the map $q_i$ induces a map $r_i : PV_S \to PV_{S+}$ in the usual way. We also define $r_0 : PV_S \to PV_{S+}$ by $r_0(N) = L$ where $L = \ker(\pi : V_{S+} \to V_S)$. This defines us an injective map $x : S_+ \times U_S \to C_S$ such that $\pi \circ x = \mathrm{id}$. It is clear that the fiber of $\pi : C_S \to U_S$ is a copy of $\mathbb{C}P^1$ thus $\pi : C_S \to U_S$ together with $x : S_+ \times U_S \to C_S$ is a stable family of $S$-curves over $U_S$. There is thus a morphism $\tilde{\gamma} : U_S \to \overline{\mathcal{X}}_S$ of schemes with $\tilde{\gamma}^* \mathcal{C}_S \simeq C_S$, whose effect on complex points is just the map $\gamma_{C_S}$. Each fibre of $C_S$ is clearly a generic $S$-curve, so $\gamma_{C_S}$ is a map from $U_S$ to $\mathcal{X}_S$. We claim that $\phi_S \gamma_{C_S} = \mathrm{id}$ as a map from $U_S$ to itself. In other words, we claim that for a point $N \in U_S$, the corresponding fibre $(C_N, x_N)$ of $C_S$ satisfies $\phi_S([C_N, x_N]) = N$. We can define a map $h : C_N \to \mathbb{C} \cup \{\infty\}$ by $h(N, L) = \frac{f(+)-f(1)}{f(2)-f(1)}$ where $f \in L$ and we have $h(x(0, N)) = \infty$. Next define $g : S \to \mathbb{C}$ by $g(i) = hx(i, N)$ we need to prove the image of $g$ in $PV_S$ is $N$ this is clear by construction, so $\phi_S([C_N, x_N]) = N$ as claimed. Since $\phi_S$ is a bijection we have $\gamma_{C_S} = \phi_S^{-1}$ and therefore the inclusion $i : \mathcal{X}_S \to \overline{\mathcal{X}}_S$ can be written as $\tilde{\gamma} \circ \phi_S$ so is a morphism of schemes. This should be compared with 7.5.1.

$\square$

## 8.7. A catalogue of results on $\overline{\mathcal{X}}_S$

In this section we discuss a number of results regarding the space $\overline{\mathcal{X}}_S$ that will be required in our analysis. We recall the results here and refer the reader to standard references for more details on these matters. We will discuss a map $\pi : \overline{\mathcal{X}}_{S+} \to \overline{\mathcal{X}}_S$ that together with certain structure sections $x : S_+ \times \overline{\mathcal{X}}_S \to \overline{\mathcal{X}}_{S+}$ form the universal curve for our moduli space. We will also discuss the structure of the tree partition $\overline{\mathcal{X}}_S = \coprod_t \mathcal{X}_S(t)$ and the structure of the spaces $\overline{\mathcal{X}}_S(t) = \mathrm{cl}\mathcal{X}_S(t)$ that are important in the study of $\overline{\mathcal{X}}_S$.

CONSTRUCTION 8.7.1. Here we describe a process which acts on the isomorphism classes of stable $S$-curve known as *stably forgetting*. Given any element $[C, x] \in \overline{\mathcal{X}}_{S+}$ we form the new curve of genus zero $[D, y] \in \overline{\mathcal{X}}_S$ as follows. We first forget the marked point $x(+)$. This new curve may be unstable. This can happen if and only if the component $C+$ of $C$ containing the marked point $x(+)$ contains only one other marked point $x(j)$ say for some $j \in S$. We then form the new $S$-curve $D$ by shrinking this component to a point and replacing the marked point $x(j)$ to be this new point, the corresponding singular point of $C$. This process then forms the stable $S$-curve $[D, y] \in \overline{\mathcal{X}}_S$ and the process is known as *stably forgetting* the marked point $+$.

EXAMPLE 8.7.2. This illustrates the last construction

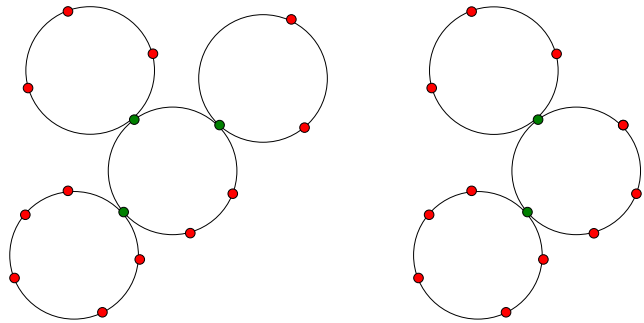

THEOREM 8.7.3. *There is a canonical morphism $\pi : \overline{\mathcal{X}}_{S+} \to \overline{\mathcal{X}}_S$ which acts on the isomorphism classes of stable $S+$-curves by stably forgetting the point $+$. This morphism is compatible with restriction to trees.*

$$\begin{array}{ccc} \overline{\mathcal{X}}_{S+} & \longrightarrow & \mathbb{U}_{S+} \\ \downarrow & & \downarrow \\ \overline{\mathcal{X}}_S & \longrightarrow & \mathbb{U}_S \end{array}$$

$\square$

CONSTRUCTION 8.7.4. Let $[D, x] \in \overline{\mathcal{X}}_S$ be a stable $S$-curve with combinatorial tree type $t$. Then we may disconnect $D$ into its irreducible components one for each internal vertex $v$ of $t$. Each of these is a copy of $\mathbb{C}P^1$. To each such component $D_v$ of $D$ we can consider it as a stable curve with its markings induced from $[D, x]$ together with the points where $D_v$ intersected the curve $D$. For each internal vertex $v$, together these marked points are $C_v$ the children of the vertex $v$ we write $[D_v, x_v]$ for each such curve. It is clear that each of these is a generic curve in the moduli space $\mathcal{X}_{C_v}$ and that the original curve $[D, x]$ can be reconstructed from the $[D_v, x_v]$ by gluing them along the intersection points, see 8.7.5 for an example. This construction then gives us a bijective map $\theta_t : \mathcal{X}_S(t) \to \prod_v \mathcal{X}_{C_v}$. Given a little more work it also turns out there is a bijection extending the previous map which we also (abusively) call $\theta_t$ given by $\theta_t : \overline{\mathcal{X}}_S(t) \to \prod_v \overline{\mathcal{X}}_{C_v}$, this construction is analogous in a reasonable sense to the first construction, where $\overline{\mathcal{X}}_S(t) = \coprod_{u \leq t} \mathcal{X}_S(u)$. We will not need the second map in our analysis.

EXAMPLE 8.7.5. An illustration of the previous construction

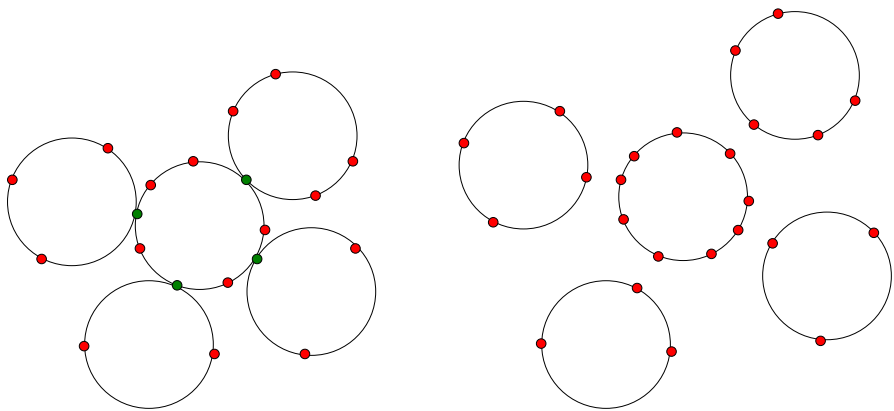

THEOREM 8.7.6. *Let $t$ be a tree then we have the following results*

(1) $\overline{\mathcal{X}}_S(t) = \text{cl}\mathcal{X}_S(t)$
(2) *The map $\theta_t : \mathcal{X}_S(t) \to \prod_v \mathcal{X}_{C_v}$ is an isomorphism.*
(3) *The map $\theta_t : \overline{\mathcal{X}}_S(t) \to \prod_v \overline{\mathcal{X}}_{C_v}$ is an isomorphism extending the previous map.*

CONSTRUCTION 8.7.7. Here we construct sections $x : S_+ \times \overline{\mathcal{X}}_S \to \overline{\mathcal{X}}_{S_+}$ that will endow the map $\pi : \overline{\mathcal{X}}_{S_+} \to \overline{\mathcal{X}}_S$ with the structure of a stable family. Let $[C, y] \in \overline{\mathcal{X}}_S$ be a stable curve and $t$ be its combinatorial tree type. For each $j \in S_+$ let $t_j$ be the tree $t$ with an extra edge attached to the middle of the edge whose external vertex is $j$ and label the vertex of this new edge $+$. Clearly $t_j$ is a tree on $S_+$ and we define $x(j, [C, y])$ to be the unique curve $[D, z]$ with $\pi([D, z]) = [C, y]$. Equivalently $D$ is the curve $C$ with an extra sphere connected to the marked point $y(j)$ and the markings $z$ of $D$ are induced from those of $y$ with the exception that the marked point for the labels $j$ and

+ are attached to the new sphere. We note that there is only one isomorphism type of a sphere with 3 marked points thus this construction uniquely define the isomorphism type of $[D, z]$. To be explicit let $t$ be the smallest tree, that is the tree associated to a generic curve and for each $i \in S$ let $t_i$ be the tree associated to $t$ as above and $u, v$ the internal vertices with $v$ such that $C_v = \{i, +\}$. Then as before consider the space $\overline{\mathcal{X}}_S(t_i)$. Then $\overline{\mathcal{X}}_S(t_i) \cong \overline{\mathcal{X}}_{C_u} \times \overline{\mathcal{X}}_{C_v} \cong \overline{\mathcal{X}}_{C_u} \times \text{pt} \cong \overline{\mathcal{X}}_{C_u} \cong \overline{\mathcal{X}}_S$. The last isomorphism is induced from the bijection $b : S \to C_u$ given by $b(j) = j$ is $j \neq i$ and $b(i) = v$. This gives us the isomorphism $\alpha_i : \overline{\mathcal{X}}_S(t_i) \to \overline{\mathcal{X}}_S$. One can check that the map implicit in this is the inclusion $j_i : \overline{\mathcal{X}}_S(t_i) \to \overline{\mathcal{X}}_{S+}$ composed with the projection $\pi : \overline{\mathcal{X}}_{S+} \to \overline{\mathcal{X}}_S$. So we have a section $\sigma_i = j_i \circ \alpha_i^{-1}$. This is the same as the previous definition since they agree on $\mathcal{X}_S$. One can also do a similar construction with the tree $t_0$.

THEOREM 8.7.8. *The map $\pi : \overline{\mathcal{X}}_{S+} \to \overline{\mathcal{X}}_S$ together with its sections $x : S_+ \times \overline{\mathcal{X}}_S \to \overline{\mathcal{X}}_{S+}$ makes this pair a stable family of S-curves moreover this is canonically the universal curve for $\overline{\mathcal{X}}_S$.*

## 8.8. Functors of coherent sheaves

We will need to show that the function $\theta_S : \overline{\mathcal{X}}_S \to PV_S$ is actually a morphism of varieties. It is well known that morphism to projective spaces are characterized by invertible sheaves over the domain scheme, or in other words line bundles, that are generated by global sections. Our definition and analysis of $\theta_S$ involved cohomology of sheaves. Thus, to define $\theta_S$ for families of stable curves, we need cohomology of families of sheaves, or in other words, push-forward functors and their derivatives. In this section, we recall the basic facts about these functors.

Let $(C, x)$ be a stable $S$-curve over a locally Noetherian scheme $X$. Recall that a sheaf $\mathcal{F}$ of $\mathcal{O}_C$-modules is *coherent* if locally on $C$ one can find a right-exact sequence $\mathcal{O}_C{}^n \to \mathcal{O}_C{}^m \twoheadrightarrow \mathcal{F}$ for some $m$ and $n$. Given such a sheaf $\mathcal{F}$, we can restrict to get a coherent sheaf $\mathcal{F}_a$ over $C_a$ for each $a \in X$, and then compute the cohomology groups $H^i(C_a; \mathcal{F}_a)$. We write $h^i(\mathcal{F}_a)$ for the dimensions of these groups over $\mathbb{C}$ (which are always finite). These numbers can jump discontinuously as $a$ moves, but the Euler characteristic

$$\chi(\mathcal{F}_a) := h^0(\mathcal{F}_a) - h^1(\mathcal{F}_a)$$

is constant on each connected component of $X$. If it happens that $h^i(\mathcal{F}_a)$ is constant, then we would like to assemble the vector spaces $H^i(C_a; \mathcal{F}_a)$ into a vector bundle over $X$. Equivalently, we would like to construct a locally free sheaf $\mathcal{G}$ of $\mathcal{O}_X$-modules with $\mathcal{G}_a = H^i(C_a; \mathcal{F}_a)$ for all $a$. We next explain the obvious candidate for $\mathcal{G}$.

We define a sheaf $\pi_*\mathcal{F}$ of $\mathcal{O}_X$-modules by

$$(\pi_*\mathcal{F})(U) = \mathcal{F}(\pi^{-1}U)$$

for any open subset $U \subseteq X$. This functor is left exact but not right exact in general. However, as our map $\pi$ is proper and flat of relative dimension one, there is only one higher derived functor, denoted by $R^1\pi_*$ (see [**6**, Corollary III.11.2]). Thus, a short exact sequence

$$\mathcal{F}_0 \to \mathcal{F}_1 \to \mathcal{F}_2$$

of sheaves on $C$ gives a six-term exact sequence

$$\pi_*\mathcal{F}_0 \to \pi_*\mathcal{F}_1 \to \pi_*\mathcal{F}_2 \to R^1\pi_*\mathcal{F}_0 \to R^1\pi_*\mathcal{F}_1 \to R^1\pi_*\mathcal{F}_2$$

(and $R^0\pi_*\mathcal{F} = \pi_*\mathcal{F}$). As $\mathcal{F}$ was assumed to be coherent, the sheaves $R^i\pi_*\mathcal{F}$ are also coherent, by [**5**, Théorème III.3.2.1]. Thus, for each point $a \in X$ we have a finite-dimensional complex vector space $(R^i\pi_*\mathcal{F})_a$, whose dimension may vary as $a$ moves in $X$. There is a natural map $\mu^i(a) : (R^i\pi_*\mathcal{F})_a \to H^i(C_a; \mathcal{F}_a)$, which may or may not be an isomorphism.

PROPOSITION 8.8.1. *Let $\mathcal{F}$ be a coherent sheaf of $\mathcal{O}_C$-modules that is flat over $X$.*

(a) *Suppose that $H^1(C_a; \mathcal{F}_a) = 0$ for all $a$. Then $R^1\pi_*\mathcal{F} = 0$, and $\pi_*\mathcal{F}$ is a locally free sheaf, and the map $\mu^0(a) : (\pi_*\mathcal{F})_a \to H^0(C_a; \mathcal{F}_a)$ is an isomorphism for all $a$.*

(b) *Suppose instead that $H^0(C_a; \mathcal{F}_a) = 0$ for all $a$. Then $\pi_*\mathcal{F} = 0$, and $R^1\pi_*\mathcal{F}$ is a locally free sheaf, and the map $\mu^1(a) : (R^1\pi_*\mathcal{F})_a \to H^1(C_a; \mathcal{F}_a)$ is an isomorphism for all $a$.*

We quote the following fact from [**12**, Corollary II.5.3]:

THEOREM 8.8.2. *If $H^{i+1}(C_a; \mathcal{F}_a) = 0$ for all $a$, then the maps $\mu^i(a)$ are isomorphisms.* □

PROOF OF PROPOSITION 8.8.1. For part (a), we certainly have $H^2(C_a; \mathcal{F}_a) = 0$ for all $a$, so $\mu^1(a) : (R^1\pi_*\mathcal{F})_a \to H^1(C_a; \mathcal{F}_a) = 0$ is an isomorphism, so $R^1\pi_*\mathcal{F} = 0$. Next, as $H^1(C_a; \mathcal{F}_a) = 0$ we see that $\mu^0(a)$ is an isomorphism. Recall that the Euler characteristic

$$\chi(\mathcal{F}_a) = \dim_\mathbb{C} H^0(C_a; \mathcal{F}_a) - \dim_\mathbb{C} H^1(C_a; \mathcal{F}_a)$$

is a locally constant function of $a$ (this is also proved in [**12**, Section II.5]). In our case this means that $\dim_\mathbb{C}((\pi_*\mathcal{F})_a)$ is locally constant. As $\pi_*\mathcal{F}$ is coherent and $X$ is locally Noetherian, this is enough to ensure that $\pi_*\mathcal{F}$ is a locally free sheaf.

Part (b) is similar. We argue as before that $\mu^1(a)$ is an isomorphism, and thus (using the Euler characteristic) that $R^1\pi_*\mathcal{F}$ is locally free. Next, suppose we have an open set $U \subseteq X$ and a section $s$ of $\mathcal{F}$ over $\pi^{-1}U$. If $s$ maps to zero in $H^0(C_a; \mathcal{F}_a)$ for all $a \in U$, then it is easy to see that $s = 0$. This means that the maps $(\pi_*\mathcal{F})_a \to H^0(C_a; \mathcal{F}_a)$ are injective, but in our case $H^0(C_a; \mathcal{F}_a) = 0$, so $\pi_*\mathcal{F} = 0$. □

PROPOSITION 8.8.3. *There is a morphism $\psi_S : X \to PV_S$ such that $\psi_S(a) = \theta_S([C_a, x_a])$ for all $a$.*

PROOF. First, to give a morphism $X \to PV_S$ of schemes is the same as to give a subsheaf $\mathcal{L} \leq V_S \otimes_\mathbb{C} \mathcal{O}_X$ such that the quotient $(V_S \otimes_\mathbb{C} \mathcal{O}_X)/\mathcal{L}$ is locally free of rank $n - 2$ (where $n = |S|$). This statement can be recovered from [**6**, Theorem II.7.1], for example. As $V_S = \widetilde{V}_S/\Delta_S$, it is equivalent to give a subsheaf $\mathcal{M} \leq \widetilde{V}_S \otimes_\mathbb{C} \mathcal{O}_X$ such that $\Delta_S \otimes_\mathbb{C} \mathcal{O}_X \leq \mathcal{M}$ and $(\widetilde{V}_S \otimes_\mathbb{C} \mathcal{O}_X)/\mathcal{M}$ is locally free of rank $n - 2$. We obviously want the fibre $\mathcal{M}_a$ to be the two-dimensional space $M = M_a$ in the definition of $\theta_S([C_a, x_a])$; our problem is to fit these together into a sheaf. For this, we simply need to extend our original definition and analysis of the map $\theta_S$ so that it works for families of curves.

We let $D_0$ and $D$ be the images of $\{0\} \times X$ and $S \times X$ under the closed inclusion $x : S_+ \times X \to C$. Let $\mathcal{J}_0$ and $\mathcal{J}$ be the ideal sheaves corresponding to $D_0$ and $D$. As $D \amalg D_0$ is contained in the smooth part of $C$, we see that $\mathcal{J}_0$ and $\mathcal{J}$ are invertible, so we can define $\mathcal{K} = \mathcal{J} \otimes \mathcal{J}_0^{-1}$, and this is again invertible. As $D_0$ and $D$ are disjoint we have $j^*\mathcal{J}_0^{-1} = j^*\mathcal{O}_C = \mathcal{O}_D$ (where $j : D \to C$ is the inclusion). This gives a short exact sequence

$$\mathcal{K} \to \mathcal{J}_0^{-1} \to j_*\mathcal{O}_D.$$

It is easy to see that the three sheaves here are flat over $X$, and that their restrictions to the curves $C_a$ are the sheaves used (under the same names) in our analysis of the map $\phi_S$.

Using Propositions 8.5.8 and 8.8.1, we see that $R^1\pi_*\mathcal{O}_C = (R^1\pi_*)(j_*\mathcal{O}_D) = R^1\pi_*\mathcal{J}_0^{-1} = \pi_*\mathcal{K} = 0$ and that $\pi_*\mathcal{O}_C$, $\pi_*j_*\mathcal{O}_D$, $\pi_*\mathcal{J}_0^{-1}$ and $R^1\pi_*\mathcal{K}$ are locally free sheaves of ranks 1, $n$, 2 and $n-2$. In fact, there is an evident map $\mathcal{O}_X \to \pi_*\mathcal{O}_C$ ("inclusion of constants") which is easily seen to be an isomorphism. Moreover, using the isomorphism $x: S \times X \to D$ we see that $\pi_*j_*\mathcal{O}_D = \prod_{s \in S} \mathcal{O}_X = \widetilde{V}_S \otimes_\mathbb{C} \mathcal{O}_X$. We write $\mathcal{M} = \pi_*\mathcal{J}_0^{-1}$. This is a locally free sheaf of rank two over $X$, and the fibre $\mathcal{M}_a$ is just $H^0(C_a; (\mathcal{J}_0^{-1})_a) = M_a$. The short exact sequence displayed above gives a six-term sequence of sheaves over $X$, but half of the terms are zero and we are left with a short exact sequence

$$\mathcal{M} \to \widetilde{V}_S \otimes_\mathbb{C} \mathcal{O}_X \to R^1\pi_*\mathcal{K},$$

showing that $(\widetilde{V}_S \otimes_\mathbb{C} \mathcal{O}_X)/\mathcal{M}$ is locally free of rank $n-2$. The inclusion $\mathcal{O}_C \to \mathcal{J}_0^{-1}$ gives an inclusion $\mathcal{O}_X = \pi_*\mathcal{O}_C \to \pi_*\mathcal{J}_0^{-1} = \mathcal{M}$, showing that $\mathcal{M}$ contains $\Delta_S \otimes_\mathbb{C} \mathcal{O}_X$. This gives the claimed map $\psi : X \to PV_S$. □

THEOREM 8.8.4. *The map $\theta_S : \overline{\mathcal{X}}_S \to PV_S$ is a morphism of schemes.*

PROOF. It is the map $\psi_S$ for the universal curve $\mathcal{C}_S$ over $\mathcal{X}_S$. □

## 8.9. Contractions

CONSTRUCTION 8.9.1. Now suppose we have a subset $T \subseteq S$ (with $|T| > 1$). If $(C, x)$ is a generic $S$-curve, it is clear that $(C, x|_{T_+})$ is a generic $T$-curve. We can thus define a map $\pi_T^S : \mathcal{X}_S \to \mathcal{X}_T$ by $\pi_T^S[C, x] = [C, x|_{T_+}]$.

EXAMPLE 8.9.2. This illustrates the last construction with 2 points forgotten.

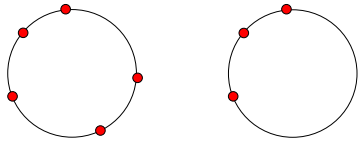

Next, the restriction map $\pi_T^S : V_S \to V_T$, evidently induces an epimorphism $\pi_T^S : U_S \to U_T$. We have the partial map $\rho_T^S : PV_S \to PV_T$ which is defined unless $M_S \leq \ker(\pi_T^S)$ thus we have an induced map $\rho_T^S : PV_S \setminus P\ker(\pi_T^S) \to PV_T$. It is easy to see that $U_S$ is contained in the domain of this map, and that the following diagram commutes

$$\begin{array}{ccc} U_S & \xrightarrow{\phi_S^{-1}} & \mathcal{X}_S \\ \rho_T^S \downarrow & & \downarrow \pi_T^S \\ U_T & \xrightarrow{\phi_T^{-1}} & \mathcal{X}_T \end{array}$$

THEOREM 8.9.3. *For any $S$ and $T$ as above, there is a unique morphism $\pi_T^S : \overline{\mathcal{X}}_S \to \overline{\mathcal{X}}_T$ of schemes extending the map $\pi_T^S : \mathcal{X}_S \to \mathcal{X}_T$ described above. Moreover, when $U \subseteq T \subseteq S$ we have $\pi_U^S = \pi_U^T \circ \pi_T^S$.*

PROOF. We know from [1] that $\overline{\mathcal{X}}_S$ is irreducible and has the same dimension as $\mathcal{X}_S$, and that $\overline{\mathcal{X}}_T$ is separated. It follows that any two maps $\overline{\mathcal{X}}_S \to \overline{\mathcal{X}}_T$ that agree on $\mathcal{X}_S$ are in fact equal. This shows that $\pi_T^S$ is unique if it exists. It follows from [10, Proposition 2.1] that $\pi_T^S$ exists when $|S \setminus T| = 1$. For the general case, just choose a chain $T = T_0 \subset T_1 \subset \ldots \subset T_r = S$ with $|T_i \setminus T_{i-1}| = 1$ and compose the corresponding $\pi$'s to get a map $\overline{\mathcal{X}}_S \to \overline{\mathcal{X}}_T$, showing that $\pi_T^S$ exists. The composition rule $\pi_U^S = \pi_U^T \circ \pi_T^S$ follows immediately from uniqueness. □

## 8.10. The isomorphism from $\overline{\mathcal{X}}_S$ to $\overline{\mathcal{M}}_S$

In this section we define a map $\overline{\theta}_S : \overline{\mathcal{X}}_S \to \overline{\mathcal{M}}_S$ and will prove that this is actually an isomorphism of varieties. This result relies on the construction of the regular map $\theta_S : \overline{\mathcal{X}}_S \to PV_S$ from section 8.5. This map should be compared with Kapranov's map $p_S : \overline{\mathcal{X}}_S \to PV_S$ in [7]. Later we will show that the inverse to this is induced from the map $\pi : \overline{\mathcal{M}}_{S+} \to \overline{\mathcal{M}}_S$ and its structure section $\sigma : S_+ \times \overline{\mathcal{M}}_S \to \overline{\mathcal{M}}_{S+}$ in the usual way. From this it will then follow that $\pi : \overline{\mathcal{M}}_{S+} \to \overline{\mathcal{M}}_S$ is the universal curve of $\overline{\mathcal{M}}_S$. The construction of our map $\overline{\theta}_S$ will be geometrically simple to understand. In order to construct the morphism $\overline{\theta}_S$ we will need to use the results regarding the maps $\pi_T^S : \overline{\mathcal{X}}_S \to \overline{\mathcal{X}}_T$ for any $T$ contained in $S$ that we considered earlier.

CONSTRUCTION 8.10.1. Here we construct a map $\overline{\theta}_S : \overline{\mathcal{X}}_S \to \overline{\mathcal{M}}_S$. For any finite set $S$ let $\theta_S : \overline{\mathcal{X}}_S \to PV_S$ be the map of construction 8.5.5 and for any subset $T \subseteq S$ define a map $\theta_T^S : \overline{\mathcal{X}}_S \to PV_T$ by $\theta_T^S = \theta_T \pi_T^S$ where $\pi_T^S : \overline{\mathcal{X}}_S \to \overline{\mathcal{X}}_T$ is the usual morphism. Then define $\overline{\theta}_S : \overline{\mathcal{X}}_S \to \prod PV_T$ by $\overline{\theta}_S = \prod_{T \subseteq S} \theta_T^S$. We first need to check that image($\overline{\theta}_S$) $\subseteq \overline{\mathcal{M}}_S$. From the commutative diagram of the last section it is easy to see that $\overline{\theta}_S(\mathcal{X}_S) = \mathcal{M}_S$, thus taking closures in the classical topology we see that $\overline{\theta}_S(\overline{\mathcal{X}}_S) = \overline{\theta}_S(\text{cl}\mathcal{X}_S) = \text{cl}\overline{\theta}_S(\mathcal{X}_S) = \text{cl}(\mathcal{M}_S) = \overline{\mathcal{M}}_S$, so that $\overline{\theta}_S : \overline{\mathcal{X}}_S \to \overline{\mathcal{M}}_S$ and the map is surjective.

LEMMA 8.10.2. *Let $[D, x] \in \overline{\mathcal{X}}_S$ and put $\mathcal{T}$ to be the corresponding S-tree. Let $T \in \mathcal{T}$ and put $M_T = \theta_T^S([D, x])$ then $\pi_U^T M_T = 0$ if and only if $U \subseteq W \in M(\mathcal{T}, T)$.*

PROOF. We prove this for $T = S$ the general case then easily follows because restricting the curve $[D, x]$ in $\overline{\mathcal{X}}_T$ is particulary simple for each $T \in \mathcal{T}$. Let $[D, x] \in \overline{\mathcal{X}}_S$ be an isomorphism class of a stable $S$-curve with tree type $t$ say and $f \in M(D, x)$. Let $D_0$ be the irreducible component of $D$ containing $x(0)$ and write $[D_0, y]$ for the corresponding generic curve whose marked point are the special points of $D$ that intersect $D_0$. Write $v_0$ for the internal vertex of $t$ corresponding to $D_0$, $I = C_{v_0}$ and $y : I_+ \to D_0$ for the corresponding markings. For each internal vertex $v$ of $t$ that lies in $I$ let $D_v$ be the curve corresponding to the disconnection of $D$ at the singular points $y(v)$, that is the component not containing $x(0)$. For each such $v$ consider $f$ restricted to $D_v$. Then $f$ is regular on $D_v$ because $x(0) \notin D_v$ and each irreducible component of $D_v$ is a copy of $\mathbb{C}P^1$. Thus $f$ is constant on each of these components. But $D_v$ is connected because $D$ is so that $f$ is constant on $D_v$. Clearly this construction is reversible. By lemmas 8.5.6 and 8.2.3 we see that the restricted map $\theta_I : \mathcal{X}_I \to U_I$ is an isomorphism where $U_I$ consists of the 'injective functions'. Thus consider the map $\theta_S : \overline{\mathcal{X}}_S \to PV_S$ put $\mathcal{T}$ to be the $S$-tree corresponding to $t$ and $M_S = \theta_S([D, x])$ then we have shown that $\pi_T^S(M_S) = 0$ if and only if $T \subseteq W \in M(\mathcal{T}, S)$ as the elements of $M(\mathcal{T}, S)$ are of the form $T_v$ for each intersection point $v$ where $T_v = \{$ elements of $S$ below $v$ away from 0 $\}$

□

COROLLARY 8.10.3. *Let $\overline{\theta}_S : \overline{\mathcal{X}}_S \to \overline{\mathcal{M}}_S$ be the map of the last construction then the following diagram commutes.*

$$\begin{array}{ccc} \overline{\mathcal{X}}_S & \xrightarrow{\overline{\theta}_S} & \overline{\mathcal{M}}_S \\ \text{type} \downarrow & & \downarrow \text{type} \\ \mathbb{U}_S & \xrightarrow{U_3} & \mathbb{T}_S \end{array}$$

PROOF. Let $[C, x] \in \overline{\mathcal{X}}_S$, $t$ the corresponding combinatorial tree, $\underline{M} = \overline{\theta}_S([C, x])$ and $\mathcal{T}$ the corresponding tree in $\mathbb{T}_S$. Let $T \in \mathcal{T}$ and $U \subseteq T$ then by the last lemma we have shown that $\pi_U^T M_T = 0$ if and only if $U \subseteq W \in M(\mathcal{T}, T)$. This is enough to show that $\mathcal{T} = \text{type}(\underline{M})$

□

LEMMA 8.10.4. *We have the following commutative diagram of isomorphisms,*

$$\begin{array}{ccc} \mathcal{X}_S & \xrightarrow{\overline{\theta}_S} & \mathcal{M}_S \\ {\theta_S} \searrow & & \swarrow {\pi_S} \\ & U_S & \end{array}$$

□

LEMMA 8.10.5. *Let $\overline{\theta}_S : \mathcal{X}_S(t) \to \mathcal{M}_S(\mathcal{T})$ be the restricted map. Then $\overline{\theta}_S$ is a bijection such that the following commutes.*

$$\begin{array}{ccc} \prod \mathcal{X}_{C_v} & \longrightarrow & \prod \mathcal{M}_{C_T} \\ \downarrow & & \downarrow \\ \mathcal{X}_S(t) & \longrightarrow & \mathcal{M}_S(\mathcal{T}) \end{array}$$

PROOF. To prove the claim it will be sufficient to show that the elements $M_T$ agree for each $T \in \mathcal{T}$ because then there is only a unique element $\underline{M}$ in $\overline{\mathcal{M}}_S$ with these values. Let $[D_S, x_S]$ be an element of $\mathcal{X}_S(t)$ represented by the curve $(D_S, x_S)$. For each $T \in \mathcal{T}$ write $[D_T, x_T]$ for the image of $[D_S, x_S]$ under the map $\pi_T^S : \overline{\mathcal{X}}_S \to \overline{\mathcal{X}}_T$ and let $v$ be the unique internal vertex of $t$ with $T_v = T$. Let $O_v$ the irreducible component of $D_T$ containing the marked point corresponding to $0$. We consider $O_v$ as a marked curve in the evident way and write $[O_v, y_v]$ for the isomorphism class of this component in $\mathcal{X}_{C_v}$. It is clear that the induced marked curve of $O_v$ in $D_T$ is the same as that of $D_S$. The marked points of $O_v$ correspond to the set $(C_v)_+$ and a non constant function $f \in M(D_T, x_T)$ restricts to a non constant function $g \in M(O_v, y_v)$. Then by construction $y_v(u) = x_T(u)$ for each $u \in (T \cap C_v)_+$ and each other $y(u) \in O_v$ corresponds to the evident intersection points so that $f(y_v(u)) = g(x_T(j))$ for all $j \in T_u$ with $u \in C_v$. Write $M_{C_v} = \theta_{C_v}([O_v, y_v])$ and using the isomorphism of $PV_{C_v}$ with $PV_{C_T}$ induced from evident bijection $b : C_v \to C_T$ put $M_{C_T}$ to be the image of $M_{C_v}$. Finally put $M_T = \theta_T([D_T, x_T])$ then under the map $i_T : PV_{C_T} \to PV_T$ we have $M_T = i_T(M_{C_T})$ thus the diagram commutes. To finish the proof we must also show that the map is bijective but this is clear since the other 3 are bijective. □

COROLLARY 8.10.6. *The map $\overline{\theta}_S : \overline{\mathcal{X}}_S \to \overline{\mathcal{M}}_S$ is a bijective morphism of varieties.*

PROOF. We have $\overline{\mathcal{X}}_S = \coprod \mathcal{X}_S(t)$ and $\overline{\mathcal{M}}_S = \coprod \mathcal{M}_S(\mathcal{T})$ and under the bijection $U_3 : \mathbb{U}_S \to \mathbb{T}_S$ we have the restricted map $\overline{\theta}_S : \mathcal{X}(t) \to \mathcal{M}_S(\mathcal{T})$ which is a bijection by the last result. This proves the claim.

□

I next claim that the map $\overline{\theta}_S : \overline{\mathcal{X}}_S \to \overline{\mathcal{M}}_S$ is actually an isomorphism of varieties this will be immediate from the following proposition. This is a standard result amongst algebraic geometers however I could find no reference so we prove it here for completeness.

LEMMA 8.10.7. *Let $f : X \to Y$ be a bijective morphism of complete smooth varieties then $f$ is an isomorphism.*

PROOF. Consider a point $x \in X$ and put $y = f(x)$. Choose an open affine subscheme $V = \text{spec}(B)$ such that $y \in V$, and an affine open subscheme $U = \text{spec}(A)$ such that $x \in U$ and $f(U)$ is contained in $V$. We then have a ring map $f^* : B \to A$. Let $m$ be the maximal ideal of $A$ that is the ideal of functions that vanish at $x$, and let $n$ be the maximal ideal of $B$ that is the ideal of functions that vanish at $f(x)$. To prove our claim it will be sufficient to check that the induced cotangent map $f^* : n/n^2 \to m/m^2$ is an isomorphism. For this it will be enough to show that the completion $\hat{A}$ of $A$ at $m$ maps isomorphically to the completion $\hat{B}$ of $B$ at $n$. Next consider $U \cap f^{-1}\{y\} = \mathrm{spec}(A/nA)$ as schemes. As $x$ is the only point of this scheme we have $\mathrm{spec}(A/nA) = \mathrm{spec}(A/m)$ as sets. By the nullstellensatz, we have $m^t \leq nA \leq m$ for some $t$. This means that the completion of $A$ at $m$ is the same as the completion of $A$ at $nA$. Thus, it will be enough to show that $\hat{B}$ maps isomorphically to the completion of $A$ at $nA$. Next, as explained in the proof of Hartshorne's Corollary $III$.11.4, we have $f_*(\mathcal{O}_X) = \mathcal{O}_Y$. We now apply Theorem $III$.11.1 ( The theorem on formal functions ) with $\mathcal{F} = \mathcal{O}_X$ and $i = 0$. The left hand side is just $\mathcal{O}_Y$ completed at $y$, or equivalently $B$ completed at $n$. The right hand side is $\mathrm{invlim}_k H^0(X_k, \mathcal{O}_{X_k})$. Here $X_k = \mathrm{spec}(A/(m_y A)^k)$. Thus, the right hand side is the completion of $A$ at $nA$. The theorem therefore tells us that the two completions are the same, as required.

$\square$

PROPOSITION 8.10.8. *Let $\pi : \overline{\mathcal{M}}_{S+} \to \overline{\mathcal{M}}_S$ be the projection and $\sigma : S_+ \times \overline{\mathcal{M}}_S \to \overline{\mathcal{M}}_{S+}$ be the structure sections of construction 7.5.2 then we have an induced map $\overline{\pi} : \overline{\mathcal{M}}_S \to \overline{\mathcal{X}}_S$ and this is the inverse for $\overline{\theta}_S : \overline{\mathcal{X}}_S \to \overline{\mathcal{M}}_S$, moreover the pair $(\pi, \sigma)$ is canonically the universal curve.*

PROOF. Since all varieties are separated and by theorem 8.9.3 the restricted map onto the usual open dense sets commutes we have the following commutative diagram

$$\begin{array}{ccc} \overline{\mathcal{X}}_{S+} \xrightarrow{\overline{\theta}_{S+}} \overline{\mathcal{M}}_{S+} & & \overline{\mathcal{X}}_{S+} \xrightarrow{\overline{\theta}_{S+}} \overline{\mathcal{M}}_{S+} \\ \pi \downarrow \quad \quad \downarrow \pi & & x \uparrow \quad \quad \uparrow \sigma \\ \overline{\mathcal{X}}_S \xrightarrow{\overline{\theta}_S} \overline{\mathcal{M}}_S & & S_+ \times \overline{\mathcal{X}}_S \xrightarrow{\overline{\theta}_S} S_+ \times \overline{\mathcal{M}}_S \end{array}$$

such that the induced diagram on section commutes. This proves our claim.

$\square$

CHAPTER 9

# The cohomology ring of $\overline{\mathcal{M}}_S$

In this chapter we are going to construct a ring $R_S$ that will turn out to be the cohomology ring of $\overline{\mathcal{M}}_S$. The ring $R_S$ will be defined as the quotient of a polynomial ring by certain homogenous relations. These relations will reflect some of the combinatorics of forests that we have already studied and some new combinatorics, that of connected sets, this we explain more precisely further in this chapter. The polynomial ring will contain one generator in degree 2 for each $T \subseteq S$ whose size is at least 3 thus giving $R_S^2$ rank $2^n - n - 1 - \binom{n}{2}$, where $n = |S|$. We compute an additive basis $\overline{B}[S]$ for our ring that uses the combinatorics of forests and restrictions imposed on the exponents that we developed earlier, this basis should be compared with its counterpart in chapter 5. The techniques we use in this section will be mainly developed from those of chapter 5. We also state necessary and sufficient conditions for an element $x$ of $R_S$ to be zero and describe $R_S$ in an equivalent form. In the last part of this chapter we will exhibit an isomorphism from $R_S$ to the Chow ring $A^*(\overline{\mathcal{M}}_S)$. We recall here that we will be using the conventions in chapter 2 regarding rings.

## 9.1. Natural classes for $H^*\overline{\mathcal{M}}_S$

We first introduce certain characteristic classes that will turn out to generate $H^*(\overline{\mathcal{M}}_S)$. This will be proven later in this chapter. We also fix some notation for the various vector bundles that we will need later.

DEFINITION 9.1.1. We have natural projections $\pi_T : \overline{\mathcal{M}}_S \to PV_T$ given by $\pi_T(\underline{M}) = M_T$. Let $\pi_T^* : H^*(PV_T) \to H^*(\overline{\mathcal{M}}_S)$ be the usual induced ring homomorphism. We put $L_T =$ the tautological line bundle over $PV_T$ and $y_T = e(L_T) \in H^2(PV_T)$, where $e$ is the Euler class. We then put $N_T = \pi_T^*(L_T)$, the pullback of $L_T$ over $\overline{\mathcal{M}}_S$ and define $x_T = e(N_T) = \pi_T^*(y_T) \in H^2(\overline{\mathcal{M}}_S)$. It is well known that $H^*(PV_T) = \mathbb{Z}[y_T]/y_T^{|T|-1}$.

## 9.2. Relations for $R_S$

In this section we define a ring $R_S$ as the quotient of a polynomial ring by an ideal generated by homogenous elements. We then show these relations hold in the cohomology ring of $\overline{\mathcal{M}}_S$. This will then gives us a ring map $\overline{\phi}_S : R_S \to H^*(\overline{\mathcal{M}}_S)$. It will turn out that this will be an isomorphism. To explain one of the relations we will need to introduce one definition, that of a connected collection of sets. This name will become apparent after giving the definition.

DEFINITION 9.2.1. Let $\mathcal{L}$ be a collection of subsets of $S$ and put $T = \mathrm{supp}(\mathcal{L})$ we say $\mathcal{L}$ is *connected* if there does not exists a splitting $T = V \coprod W$ of $T$ such that for all $U \in \mathcal{L}$ either $U \subseteq V$ or $U \subseteq W$.

DEFINITION 9.2.2. Let $\mathbb{Z}_S$ be the graded polynomial algebra over $\mathbb{Z}$ with one generator $y_T \in \mathbb{Z}_S^2$ for each $T \subseteq S$ with $|T| > 2$, that is we define $\mathbb{Z}_S = \mathbb{Z}[\, y_T \mid T \subseteq S \text{ and } |T| > 2\,]$.

Let $I_S$ be the ideal generated by the following relations

(1) Let $T \subseteq S$ then $y_T^{|T|-1} = 0$
(2) Let $T \subseteq S$ and $\mathcal{T}$ a $T$-tree of depth 2 then $y_T^{m(\mathcal{T},T)} \prod_{U \in M(\mathcal{T},T)} (y_T - y_U) = 0$
(3) Let $\mathcal{L}$ be a connected set then $\prod_{U \in \mathcal{L}} (y_T - y_U) = 0$ where $T = \mathrm{supp}(\mathcal{L})$

then we define the ring $R_S$ by $R_S = \mathbb{Z}_S / I_S$

PROPOSITION 9.2.3. *Let $\phi_S : \mathbb{Z}_S \to H^*(\overline{\mathcal{M}}_S)$ be the map given by $\phi(y_T) = x_T$, then $\phi(I_S) = 0$ and so there is an induced ring map $\overline{\phi}_S : R_S \to H^*(\overline{\mathcal{M}}_S)$.*

REMARK 9.2.4. If $T \subseteq S$ then we have an inclusion of ideals $I_T \subseteq I_S$ and the inclusion map $i : \mathbb{Z}_T \to \mathbb{Z}_S$ induces a map $r_T^S : R_T \to R_S$ such that the following diagram commutes. We will see later that this map is injective.

$$\begin{array}{ccc} H^*(\overline{\mathcal{M}}_\mathcal{T}) & \longrightarrow & H^*(\overline{\mathcal{M}}_S) \\ \phi_\mathcal{T} \uparrow & & \uparrow \phi_S \\ R_\mathcal{T} & \xrightarrow{r_\mathcal{T}^S} & R_S \end{array}$$

Here we prove relations 1, 2 and 3 are satisfied by our ring. We note here that relation 1 can be considered as a special case of relation 2 by relaxing the condition on the depth of the tree $\mathcal{T}$ to $\mathrm{d}(\mathcal{T}) \leq 2$. Then relation 1 corresponds to the tree $\mathcal{T} = \{T\}$ as $M(\mathcal{T}, T)$ is empty and $m(\mathcal{T}, T) = |T| - 1$. For the proof of relations 1 and 2 we use results obtained in chapter 5.

LEMMA 9.2.5. *Relations 1 and 2 hold in the ring $H^*(\overline{\mathcal{M}}_S)$.*

PROOF. Let $\mathcal{T}$ be a $T$-tree of depth 2 or less and $\pi : \overline{\mathcal{M}}_S \to \overline{\mathcal{M}}_\mathcal{T}$ be the projection map. Then we claim $x_T^{m(\mathcal{T},T)} \prod_{U \in M(\mathcal{T},T)} (x_T - x_U) = 0$ in $H^*(\overline{\mathcal{M}}_S)$. Then by lemma 5.1.10 we have the equivalent relation in the ring $H^*(\overline{\mathcal{M}}_\mathcal{T})$ thus applying the map $\pi^* : H^*(\overline{\mathcal{M}}_\mathcal{T}) \to H^*(\overline{\mathcal{M}}_S)$ we obtain the required relation. □

Before we can prove relation 3 we need a lemma regarding connected sets. To prove this relation we then show that there is an injective map of vector bundles $N_T \rightarrowtail \bigoplus_{U \in \mathcal{L}} N_U$.

LEMMA 9.2.6. *Let $\mathcal{L}$ be a collection of subsets of $S$ and put $T = \mathrm{supp}(\mathcal{L})$, then the map $\pi_\mathcal{L} : V_T \to \bigoplus_{U \in \mathcal{L}} V_U$ is injective if and only if $\mathcal{L}$ is connected.*

PROOF. Suppose $\mathcal{L}$ is connected and let $\overline{x} \in \ker(\pi_\mathcal{L})$ then $x|_U$ is constant for each $U \in \mathcal{L}$. Choose $U \in \mathcal{L}$ and put $c = x(u)$ for any $u \in U$. Put $T_0 = \{\, t \in T \mid x(t) = c \,\}$ and $T_1 = \{\, t \in T \mid x(t) \neq c \,\}$ then $T = T_0 \amalg T_1$ and for every $W \in \mathcal{L}$ either $W \subseteq T_0$ or $W \subseteq T_1$ therefore by connectivity either $T_0 = \emptyset$ or $T_1 = \emptyset$. But $U \subseteq T_0$ therefore $T_1 = \emptyset$, $T_0 = T$ and $x = c$ thus $\overline{x} = 0$ and $\pi_\mathcal{L}$ is injective. Next suppose $\mathcal{L}$ is not connected and let $T = T_0 \amalg T_1$ be such a splitting. Consider the map $\pi : V_T \to V_{T_0} \oplus V_{T_1}$ then one readily see that $\ker(\pi) \subseteq \ker(\pi_\mathcal{L})$ and by lemma 4.1.5 we see $\dim(\ker(\pi_\mathcal{L})) > 0$. □

LEMMA 9.2.7. *Let $\mathcal{L}$ be a connected set and $T = \mathrm{supp}(\mathcal{L})$ then $\prod_{U \in \mathcal{L}} (x_T - x_U) = 0$ in the cohomology ring of $\overline{\mathcal{M}}_S$.*

PROOF. Let $\underline{M}$ be an element of $\overline{\mathcal{M}}_S$ then by definition $M_T \leq \bigcap_{U \in \mathcal{L}} (\pi_U^T)^{-1} M_U = \pi_\mathcal{L}^{-1}(\bigoplus_{U \in \mathcal{L}} M_U)$ and we have the following diagram.

$$\begin{array}{ccccc}
M_T & \xrightarrow{i} & \pi_\mathcal{L}^{-1}(\oplus M_U) & \xrightarrow{\pi_\mathcal{L}} & \oplus M_U \\
\Big\| & & \Big\downarrow i & & \Big\downarrow i \\
M_T & \xrightarrow{i} & V_T & \xrightarrow{\pi_\mathcal{L}} & \oplus V_U
\end{array}$$

where $\pi_\mathcal{L}$ is injective. Thus

$$\begin{array}{rcl}
\pi_\mathcal{L}^{-1}(\oplus M_U) & \rightarrowtail & \bigoplus_{U \in \mathcal{L}} M_U \\
M_T & \rightarrowtail & \bigoplus_{U \in \mathcal{L}} M_U \\
N_T & \rightarrowtail & \bigoplus_{U \in \mathcal{L}} N_U \text{ as vector bundles} \\
\mathbb{C}[1] & \rightarrowtail & \bigoplus_{U \in \mathcal{L}} (N_T^* \otimes N_U) \text{ where } \mathbb{C}[1] \text{ is trivial} \\
0 & = & e(\bigoplus_{U \in \mathcal{L}} (N_T^* \otimes N_U)) \text{ on taking the Euler class} \\
& = & \prod_{U \in \mathcal{L}} (e(N_T) - e(N_U)) \\
& = & \prod_{U \in \mathcal{L}} (x_T - x_U)
\end{array}$$

$\square$

## 9.3. Construction of a basis for $R_S$

Our next problem is to understand in more detail the structure of the ring $R_S$. Using certain combinatorial conditions, we will describe a set $\overline{B}[S]$ of monomials in the generators $x_T$, which will turn out to be a basis for $R_S$ over $\mathbb{Z}$.

DEFINITION 9.3.1. We put $B_S = \{$ monomial basis for $\mathbb{Z}_S\}$ and say that a monomial $x = \prod_{T \in \mathcal{F}} y_T^{n_T}$ is $forest - like$ if $\mathcal{F}$ is a forest.

DEFINITION 9.3.2. For any monomial $y = \prod_{T \in \mathcal{L}} y_T^{n_T} \in B_S$ we define functions called the $shape$, shape $: B_S \to P^2(S)$ and $support$, supp $: B_S \to P(S)$ by shape$(y) = \mathcal{L}$ and supp$(y) = \bigcup_{T \in \mathcal{L}} T$.

REMARK 9.3.3. Clearly a monomial is forest-like if and only if its shape is a forest.

DEFINITION 9.3.4. For any monomial $y$ with $y = \prod_{T \in \mathcal{L}} y_T^{n_T}$ so that shape$(y) = \mathcal{L}$ and for any $\mathcal{U} \subseteq \mathcal{L}$ we write $y|_\mathcal{U} = \prod_{U \in \mathcal{U}} y_U^{n_U}$ called the $restriction$ of $y$ to $\mathcal{U}$.

REMARK 9.3.5. The numbers $m(\mathcal{F}, T)$ defined in chapter 3 will turn out to give us conditions on the exponents of generators of forest-like monomials $x$ which will give us a basis for $R_S$. These numbers come naturally from condition 2 in definition 9.2.2. We next give a precise description of a basis $\overline{B}[S]$ of $R_S$.

DEFINITION 9.3.6. For each forest $\mathcal{F}$ of $S$ we define

$$B[S][\mathcal{F}] = \left\{ \prod_{T \in \mathcal{F}} y_T^{n_T} \mid 1 \leq n_T < m(\mathcal{F}, T) \text{ for all } T \in \mathcal{F} \right\}$$
$$B[S] = \coprod_{forests\ \mathcal{F}} B[S][\mathcal{F}]$$
$$\overline{B}[S] = q_S(B[S]) \text{ where } q_S : \mathbb{Z}_S \to R_S \text{ is the quotient map.}$$

REMARK 9.3.7. *For any $T \subseteq S$ it is clear that $B[T] \subseteq B[S]$*

LEMMA 9.3.8. *If $y = \prod_{T \in \mathcal{F}} y_T^{n_T} \in B[S]$ then*

(1) *For each $U \in \mathcal{F}$ we have $\sum_{T \in \mathcal{F}|_U} n_T \leq |U| - 2$ with equality if and only if $\mathcal{F}|_U = \{U\}$ and $n_U = |U| - 2$.*

(2) *In particular* $\sum_{T \in \mathcal{F}} n_T \leq |S| - 2$ *with equality if and only if* $\mathcal{F} = \{S\}$ *and* $n_S = |S| - 2$.

PROOF. We prove the second case. The first then follows by considering $\mathcal{F}|_U$ and applying the second case and observing $M(\mathcal{F}|_U) = \{U\}$ for each $U \in \mathcal{F}$.

$$
\begin{aligned}
\sum_{T \in \mathcal{F}} n_T &\leq \sum_{T \in \mathcal{F}} (m(\mathcal{F}, T) - 1) \\
&= \sum_{T \in \mathcal{F}} m(\mathcal{F}, T) - |\mathcal{F}| \\
&= \sum_{T \in M(\mathcal{F})} |T| - |M(\mathcal{F})| - |\mathcal{F}| \text{ by lemma 3.3.2} \\
&= |\coprod_{T \in M(\mathcal{F})} T| - |M(\mathcal{F})| - |\mathcal{F}| \\
&\leq |S| - 2 \text{ as } |M(\mathcal{F})|, |\mathcal{F}| \geq 1
\end{aligned}
$$

It is now clear that we have equality if and only if $\mathcal{F} = \{S\}$ and $n_S = |S| - 2$　　□

## 9.4. A filtration for $R_S$

We will next introduce the concept of the weight of a monomial which will induce a filtration on $R_S$ that will enable us to show that the set $\overline{B}[S]$ does indeed span $R_S$. Most of these constructions will be analogous to those in section 3 of chapter 5 so we do not repeat proofs of those statements here.

DEFINITION 9.4.1. Let $y \in B_S$, so $y = \prod_{T \in \mathcal{L}} y_T^{n_T}$ say. We define a function $\text{wt} : B_S \to \mathbb{N}$ called the *weight* by $\text{wt}(y) = \sum_{T \in \mathcal{L}} n_T |T|$. In particular this gives a monomial $y_T$ weight $|T|$. We also put $\deg : R_S \to \mathbb{N}$ to be the cohomological degree function, so $\deg(y_T) = 2$.

LEMMA 9.4.2. $\text{wt}(xy) = \text{wt}(x) + \text{wt}(y)$ *and* $\deg(x) \leq \text{wt}(x) \leq \frac{n}{2} \deg(x)$　　□

DEFINITION 9.4.3. We define filtration's on $\mathbb{Z}_S$ by

$$F_k \mathbb{Z}_S = \text{span}\{\, y \in B_S \mid \text{wt}(y) \geq k \,\}$$
$$G_k \mathbb{Z}_S = \text{span}\{\, y \in B_S \mid \deg(y) \geq k \,\}$$

LEMMA 9.4.4. *Put $H = F$ or $G$ then $H_k \mathbb{Z}_S$ is a convergent decreasing filtration of $\mathbb{Z}_S$, that is,*

$$\mathbb{Z}_S = H_0 \mathbb{Z}_S \supseteq H_1 \mathbb{Z}_S \supseteq H_2 \mathbb{Z}_S \supseteq, \dots \text{ and } \bigcap_{k \geq 0} H_k \mathbb{Z}_S = \{0\}$$

□

DEFINITION 9.4.5. We define a function wt : $\mathbb{Z}_S \to \mathbb{N} \cup \{\infty\}$ called the *weight* that extends the last function as follows. For every non-zero $y \in \mathbb{Z}_S$ we know by the previous lemma that there exists a largest k such that $y \in F_k \mathbb{Z}_S$ but $y \notin F_{k+1} \mathbb{Z}_S$. We put $\text{wt}(y) = k$ and define $\text{wt}(0) = \infty$.

LEMMA 9.4.6. *The weight function has the following properties*

(1) *If $y = \sum_{i \in I} a_i y_i$ with $y_i \in B_S$ then, $\text{wt}(y) = \min\{\, \text{wt}(y_i) \mid i \in I \,\}$*
(2) $\text{wt}(xy) = \text{wt}(x) + \text{wt}(y)$

□

DEFINITION 9.4.7. We define filtration's on $R_S$ by $F_k R_S = q_S(F_k \mathbb{Z}_S)$ and $G_k R_S = q_S(G_k \mathbb{Z}_S)$ where $q_S : \mathbb{Z}_S \to R_S$ is the quotient map.

LEMMA 9.4.8. *Put $H = F$ or $G$ then $H_k R_S$ is a convergent decreasing filtration of $R_S$, that is,*

$$R_S = H_0 R_S \supseteq H_1 R_S \supseteq H_2 R_S \supseteq, \dots \quad \text{and } \bigcap_{k \geq 0} H_k R_S = \{0\}$$

□

PROOF. The first part of the claim is clear. By lemma 9.4.2 it is clear that $F_k\mathbb{Z}_S \subseteq G_l\mathbb{Z}_S$ where $l = [2k/n]$ and $[a]$ is the integer part of $a$ thus $F_k R_S = q_S(F_k\mathbb{Z}_S) \subseteq q_S(G_l\mathbb{Z}_S) = G_l R_S$. Therefore $\bigcap_k F_k R_S \subseteq \bigcap_l G_l R_S = \{0\}$ as $R_S$ is graded by degree. □

DEFINITION 9.4.9. We next define a function wt : $R_S \to \mathbb{N} \cup \{\infty\}$ called the *weight* as follows. For every non-zero $x \in R_S$ we know by the previous lemma that there exists a largest $k \in \mathbb{N}$ such that $x \in F_k R_S$, but $x \notin F_{k+1} R_S$. We put wt$(x) = k$ and define wt$(0) = \infty$.

LEMMA 9.4.10. *The function* wt *has the following properties.*

(1) *If* $x = \sum_{i \in I} a_i x_i$ *with* $x_i \in R_S$ *then* wt$(x) \geq \min\{$wt$(x_i) \mid i \in I\}$
(2) *There is an* $m \in \mathbb{N}$ *such that* wt$(x) \leq m$ *for all non zero* $x \in R_S$.
(3) wt$(xy) \geq$ wt$(x) +$ wt$(y)$
(4) wt$(q_S(x)) \geq$ wt$(x)$

□

DEFINITION 9.4.11. We define $F_k B[S] = \{\, x \in B[S] \mid$ wt$(x) \geq k\,\}$ and the induced set $F_k \overline{B}[S] = q_S(F_k B[S])$.

DEFINITION 9.4.12. Let $y \in \mathbb{Z}_S$ be a monomial. We say $y$ is admissible if $y \in B[S]$ and inadmissible if $y \notin B[S]$. We also say $y \in \mathbb{Z}_S$ is minimally inadmissible if it has the form

(1) $y_T y_U$ with $T \cap U$ non-empty, $T \not\subseteq U$ and $U \not\subseteq T$
(2) $y_T^{m(\mathcal{T},T)} \prod_{U \in M(\mathcal{T},T)} y_U$ with $\mathcal{T}$ a $T$-tree of depth 2 .
(3) $y_T^{|T|-1}$ for some $T \subseteq S$.

LEMMA 9.4.13. *Let* $z \in \mathbb{Z}_S$ *be a non-zero monomial, then* $z$ *is inadmissible if and only if* $z = xy$ *with* $x$ *minimally inadmissible.*

PROOF. The result is immediate from the definitions. □

LEMMA 9.4.14. *If $y \in \mathbb{Z}_S$ is minimally inadmissible then $\operatorname{wt}(q_S(y)) > \operatorname{wt}(y)$.*

PROOF. In case one we have sets $T$ and $U$ with $T \cap U$ non-empty $T \not\subseteq U$ and $U \not\subseteq T$. Put $W = U \cup V$, then by relation 3 we have $(x_W - x_T)(x_W - x_U) = 0$ and so $x_T x_U = x_W x_T + x_W x_U - x_W^2$ therefore using lemma 9.4.10 we see that

$$\begin{aligned}
\operatorname{wt}(q_S(y_T y_U)) &= \operatorname{wt}(x_T x_U) \\
&= \operatorname{wt}(q_S(y_W y_T + y_W y_U - y_W^2)) \\
&= \operatorname{wt}(q_S(y_W y_T) + q_S(y_W y_U) - q_S(y_W^2)) \\
&\geq \min\{\operatorname{wt}(q_S(y_W y_T)), \operatorname{wt}(q_S(y_W y_U)), \operatorname{wt}(q_S(y_W^2))\} \\
&\geq \min\{\operatorname{wt}(y_W y_T), \operatorname{wt}(y_W y_U), \operatorname{wt}(y_W^2)\} \\
&= \min\{|T \cup U| + |T|, |T \cup U| + |U|, 2|T \cup U|\} \\
&> |T| + |U| \\
&= \operatorname{wt}(y_T y_U)
\end{aligned}$$

The second and third case are essentially lemma 5.3.14 of chapter 5. □

COROLLARY 9.4.15. *Let $y \in \mathbb{Z}_S$ be a monomial. Then if $y$ is inadmissible that is $y \notin B[S]$ we have $\operatorname{wt}(q_S(y)) > \operatorname{wt}(y)$.*

PROOF. The proof is essentially that of lemma 5.3.15 so we do not prove it again. □

LEMMA 9.4.16. $R_S = \operatorname{span}\overline{B}[S]$

PROOF. The proof will follow from a downward induction on the weight $w$ of the statement that for every $k \in \mathbb{N}$  $F_k R_S = \operatorname{span} F_k \overline{B}[S]$. For $k \gg 0$ we know by lemmas 9.4.10 part 2 and 9.3.8 that $F_k R_S = 0 = \operatorname{span} F_k \overline{B}[S]$. Suppose it is true for $k > w$. Let $x \in F_w R_S = q_S(F_w \mathbb{Z}_S)$, then $x = q_S(y)$ where $y = \sum_{i \in I} a_i y_i \in F_w \mathbb{Z}_S$ thus $\operatorname{wt}(y_i) \geq w$. So it is enough to show that for every $i \in I$ we have $q_S(y_i) \in \operatorname{span} F_w \overline{B}[S]$. If $y_i \in B[S]$ then we are ok so we may assume that $y_i \notin B[S]$. In this case we know that $\operatorname{wt}(q_S(y_i)) > \operatorname{wt}(y_i) \geq w$ and so by induction $q_S(y_i) \in F_{k+1}R_S = \text{span}F_{k+1}\overline{B}[S] \subseteq \text{span}F_k\overline{B}[S]$. Therefore $R_S = F_0R_S = \text{span}F_0\overline{B}[S] = \text{span}\overline{B}[S]$ and we are done.

□

COROLLARY 9.4.17. *For every $i > |S| - 2$ we have $R_S^{2i} = 0$ and $R_S^{2(|S|-2)} = \mathbb{Z}[x_S^{|S|-2}]$*

PROOF. For every $i > |S| - 2$ we know by lemma 9.3.8 that $B^{2i}[S] = \emptyset$. Therefore $R_S^{2i} = \text{span}B^{2i}[S] = 0$. We have shown by lemma 6.4.18 part 7 that $H^{2(|S|-2)}(\overline{\mathcal{M}}_S) = \mathbb{Z}[x_S^{|S|-2}]$ thus the claim follows by lemma 9.3.8 using the map $\overline{\phi}_S : R_S \to H^*(\overline{\mathcal{M}}_S)$.

□

COROLLARY 9.4.18. *Let $y \in \mathbb{Z}_S$ with $\text{shape}(y) = \mathcal{V}$ say and put $x = q_S(y)$ where $q_S : \mathbb{Z}_S \to R_S$ is the usual quotient map. Suppose there exists a subset $\mathcal{U}$ of $\mathcal{V}$ such that $\deg(x|_\mathcal{U}) \geq 2(|\text{supp}(\mathcal{U})| - 1)$ then $x = 0$.*

PROOF. Put $T = \text{supp}(\mathcal{U})$ by considering $x|_\mathcal{U}$ as an element in $R_T$ we see by the previous result that $x|_\mathcal{U} = 0$ in $R_T$. Let $r_T^S : R_T \to R_S$ be the ring map induced by the inclusion $i : \mathbb{Z}_T \to \mathbb{Z}_S$ then we see that $x_\mathcal{U}$ is zero in $R_S$. Therefore $x = x|_\mathcal{U} \cdot x'$ is zero in $R_S$.

□

LEMMA 9.4.19. *$\overline{B}[S]$ is a basis for $R_S$ and $\overline{\phi}_S : R_S \to H^*(\overline{\mathcal{M}}_S)$ is an isomorphism.*

PROOF. We have shown in lemma 6.4.18 part 5 that $H^*(\overline{\mathcal{M}}_S)$ is generated by $H^2(\overline{\mathcal{M}}_S)$. Part 3 tells us that $\overline{\phi}_S : R_S^2 \to H^2(\overline{\mathcal{M}}_S)$ is a surjection. Thus the map $\overline{\phi}_S$ is surjective. We have also shown in lemma 6.4.18 that $H^*(\overline{\mathcal{M}}_S)$ is a finite free module of rank $d$ where $d = P_S(1)$ and $P_S$ is the Poincaré series of $\overline{\mathcal{M}}_S$ constructed in definition 6.4.15. It is then clear by the construction of $B[S]$ that $|B[S]| = d$. By lemma 9.4.16 we have that $\overline{B}[S]$ spans $R_S$ and we have the map $\overline{\phi}_S : R_S \to H^*(\overline{\mathcal{M}}_S)$ is surjective. Thus $\overline{\phi}_S(\overline{B}[S])$ spans $H^*(\overline{\mathcal{M}}_S)$. Then $d \leq |\overline{\phi}(\overline{B}[S])| \leq |\overline{B}[S]| \leq |B[S]| = d$. Thus $|\overline{\phi}(\overline{B}[S])| = \text{rank}(H^*(\overline{\mathcal{M}}_S))$ and $\overline{\phi}(\overline{B}[S])$ is a basis. Therefore $\overline{B}[S]$ is a basis for $R_S$ and the map $\overline{\phi}_S : R_S \to H^*(\overline{\mathcal{M}}_S)$ is an isomorphism.

□

LEMMA 9.4.20. *The map $r_T^S : R_T \to R_S$ is injective.*

PROOF. Topologically this is clear since we have maps $\overline{\mathcal{M}}_T \to \overline{\mathcal{M}}_S \to \overline{\mathcal{M}}_T$ such that the composition is the identity. To see this algebraically we have an inclusion of sets $B[T] \subseteq B[S]$ in the evident way.

□

## 9.5. The zero conditions for monomials of $R_S$

The aim of this section is to compute precisely the set of non-zero monomials of $R_S$. This set is given in the next definition. Corollary 9.4.18 gives us a condition for a monomial to be zero. We will show that this is condition is also necessary. In this section we will be using the ordinary degree of a monomial.

DEFINITION 9.5.1. We define $N[S]$ to be the set of monomials of the form $y = \prod_{T \in \mathcal{L}} y_T^{n_T}$ where $\mathcal{L}$ is any collection of subsets of $S$ such that for every $U \in \mathcal{L}$ we have $n_U \geq 1$ and for every $T \subseteq S$ we have $\sum_{\substack{U \in \mathcal{L} \\ U \subseteq T}} n_U \leq |T| - 2$. We also write $N_F[S]$ for the subset of $N[S]$ consisting of the monomials whose shapes are forests. We then define $\overline{N}[S] = q_S(N[S])$ and $\overline{N}_F[S] = q_S(N_F[S])$ where $q_S : \mathbb{Z}_S \to R_S$ is the quotient map.

LEMMA 9.5.2. $\overline{N}_F[S]$ *is precisely the set of non zero monomials whose shapes are forests moreover for any $x \in \overline{N}_F[S]$ we have $x_S^{|S|-2-\deg(x)} x = x_S^{|S|-2}$*

PROOF. let $y$ be a monomial in $N_F[S]$ with shape $\mathcal{F}$ say. Put $x = q_S(y)$, then we can consider $x$ as a monomial in $R_\mathcal{F}$ in the obvious way. Let $i : \mathbb{Z}_\mathcal{F} \to \mathbb{Z}_S$ be the inclusion map. Then $I_\mathcal{F} \subseteq I_S$ so that we have an induced ring map $r : R_\mathcal{F} \to R_S$. One then readily verifies that $\overline{B}[\mathcal{F}] \subseteq \overline{B}[S]$ so that the map is injective. Then as $x \in N[\mathcal{F}]$ it is non-zero in $R_\mathcal{F}$ by lemma 5.4.3. Then applying $r$ we see that $x$ is non-zero in $R_S$. In particular $x_S^{|S|-2-\deg(x)} x = x_S^{|S|-2}$ in $R_\mathcal{F}$ thus applying $r$ again we obtain our result.

□

DEFINITION 9.5.3. Let $\mathcal{L}$ be a collection of subsets of $S$. A *chain* in $\mathcal{L}$ is a subset $\mathcal{C} = \{T_1, ..., T_n\}$ of $\mathcal{L}$ such that $T_i \not\subseteq T_{i+1}$, $T_{i+1} \not\subseteq T_i$ and $T_i \cap T_{i+1}$ is non-empty. A $sub-chain$ $\mathcal{L}'$ of a *chain* $\mathcal{L}$ is a subset of $\mathcal{L}$ that is itself a *chain*. Note that any chain is a connected set.

DEFINITION 9.5.4. Let $\mathcal{L}$ be a collection of subsets of $S$, we define an equivalence relation $\sim$ on $\mathcal{L}$ as follows. Given $U$ and $V$ in $\mathcal{L}$ we say $U \sim V$ if and only if $U = V$ or there is a chain $(T_1, ..., T_n)$ of $\mathcal{L}$ with $U = T_1$ and $V = T_n$. Note we allow repetition in our chain. We write $\overline{U}$ for the equivalence class of $U \in \mathcal{L}$ and it is clear that $\overline{U}$ is itself a chain. It is the maximal chain containing $U$.

DEFINITION 9.5.5. Let $\mathcal{L}$ be a collection of subsets of $S$. We define the set

$$\mathcal{L}^* = \{ \operatorname{supp}(\mathcal{U}) \mid \mathcal{U} \subseteq \mathcal{L} \text{ and } \mathcal{U} \text{ is a chain } \}$$

and call this the *completion* of $\mathcal{L}$. Note that any singleton is a chain, so $\mathcal{L} \subseteq \mathcal{L}^*$

LEMMA 9.5.6. *Let $\mathcal{L}$ be a collection of subsets of $S$ and $\mathcal{C}, \mathcal{D}$ be chains of $\mathcal{L}$ that are subsets of different equivalence classes. Put $T = supp(\mathcal{C})$ and $W = supp(\mathcal{D})$ then $\{T, W\}$ is a forest.*

PROOF. Suppose for a contradiction that $\{T, W\}$ is not a forest. Then we can find $C \in \mathcal{C}$ and $D \in \mathcal{D}$ with $C \cap D$ non-empty. Now the pair $\{C, D\}$ must be a tree for otherwise $\mathcal{C}$ and $\mathcal{D}$ would belong to the same equivalence class. We may suppose then without loss of generality that $C \subseteq D$. Now not every $E$ in $\mathcal{C}$ can be contained in $D$ else $\{T, W\}$ would be a tree, thus there is some $E \in \mathcal{C}$ not contained in $D$. Because $C, E \in \mathcal{C}$ there must be a chain $\{C_1, ..., C_r\}$ with $C_1 = C$ and $C_r = E$. We have $C_1 \subseteq D$ and $C_r \not\subseteq D$ so for some $k$ we must have $C_k \subseteq D$ and $C_{k+1} \not\subseteq D$. One then checks that $\{C_{k+1}, D\}$ is a chain and therefore $C_{k+1} \sim D$, contrary to the assumption on $\mathcal{C}$ and $\mathcal{D}$. This completes our proof.

□

COROLLARY 9.5.7. *Let $\mathcal{L}$ be a collection of subsets of $S$ and $U, V$ be elements of $\mathcal{L}$ that are not related under $\sim$. Let $\mathcal{U}$ be a forest on $(\overline{U})^*$ and $\mathcal{V}$ be a forest on $(\overline{V})^*$ then $\mathcal{F} = \mathcal{U} \cup \mathcal{V}$ is a forest.*

PROOF. The proof is immediate given the last lemma.
□

COROLLARY 9.5.8. *For any collections $\mathcal{U}, \mathcal{V}$ of subsets of $S$ we have $\mathcal{U} \subseteq \mathcal{U}^*$, $\mathcal{U}^{**} = \mathcal{U}^*$, $\mathcal{U}^* \cup \mathcal{V}^* \subseteq (\mathcal{U} \cup \mathcal{V})^*$ and if $\mathcal{U} \subseteq \mathcal{V}$ then $\mathcal{U}^* \subseteq \mathcal{V}^*$.*

PROOF. Only the second part requires some comment. We already know that $\mathcal{L}^* \subseteq \mathcal{L}^{**}$ so it suffices to prove the reverse inclusion. Let $W \in \mathcal{L}^{**}$ then $W = \text{supp}\{V_1, ..., V_m\}$ for some chain $\{V_1, ..., V_m\}$ of $\mathcal{L}^*$. We then have each $V_i = \text{supp}(\mathcal{C}_i)$ for some chain $\mathcal{C}_i$ of $\mathcal{L}$. Now consider the pair $\mathcal{C}_j, \mathcal{C}_{j+1}$ for some $j$ and put $\mathcal{D}_j = \mathcal{C}_j \cup \mathcal{C}_{j+1}$. Then apply the relation $\sim$ to $\mathcal{D}_j$. Either there union form a chain or each is a separate equivalence class. But the latter case cannot happen since then by lemma 9.5.6 we would have $\{V_j, V_{j+1}\}$ forming a tree. Thus $\mathcal{D} = \bigcup_i \mathcal{C}_i$ is a chain on $\mathcal{L}$ and $W = \text{supp}(\mathcal{D})$.
□

DEFINITION 9.5.9. *Given two collections of elements $\mathcal{U}, \mathcal{V}$ of S we say $\mathcal{U} \leq \mathcal{V}$ if $\mathcal{U} \subseteq \mathcal{V}^*$*

LEMMA 9.5.10. *The above relation is reflexive and transitive.*

PROOF. Clearly $\mathcal{U} \leq \mathcal{U}$ as $\mathcal{U} \subseteq \mathcal{U}^*$. If $\mathcal{U} \leq \mathcal{V}$ and $\mathcal{V} \leq \mathcal{W}$ then $\mathcal{U} \subseteq \mathcal{V}^*$ and $\mathcal{V} \subseteq \mathcal{W}^*$ thus $\mathcal{V}^* \subseteq \mathcal{W}^{**} = \mathcal{W}^*$ and $\mathcal{U} \subseteq \mathcal{W}^*$ therefore $\mathcal{U} \leq \mathcal{W}$.
□

DEFINITION 9.5.11. *Let $x, y \in N[S]$ with $x = \prod_{U \in \mathcal{U}} y_U^{n_U}$ and $y = \prod_{V \in \mathcal{V}} y_V^{m_V}$ then we say $x \leq y$ if $\mathcal{U} \leq \mathcal{V}$ and for each $T \subseteq S$ we have $\sum_{\substack{U \in \mathcal{U} \\ U \subseteq T}} n_U \leq \sum_{\substack{V \in \mathcal{V} \\ V \subseteq T}} m_U$*

LEMMA 9.5.12. *The above relation is a reflexive and transitive ordering on monomials $x$ of $N[S]$ and if $x \leq y$ and $x' \leq y'$ and $yy' \in N[S]$ then $xx' \leq yy'$*

PROOF. The proofs of these statements are clear though long and tedious so we do not supply any proofs.

□

DEFINITION 9.5.13. Let $y = \prod_{U \in \mathcal{L}} y_T^{n_T}$ be a monomial in $N[S]$ we say $y$ is a proper monomial if we can write $x = q_S(y)$ as $x = \sum_{i \in I} \epsilon(i) x_i$ with $x_i = q_S(y_i)$ for some monomials $y_i$ in $N_F[S]$ such that $\epsilon(i) = \pm 1$, $\sum_{i \in I} \epsilon(i) = 1$ and each $y_i \leq y$

LEMMA 9.5.14. *If $y$ is a proper monomial then $x = q_S(y)$ is non-zero.*

PROOF. Since $y$ is proper we may write $x = \sum_{i \in I} \epsilon(i) x_i$ with the $x_i = q_S(y_i)$ and $y_i$ monomials in $N_F[S]$ such that $\epsilon(i) = \pm 1$, $\sum_{i \in I} \epsilon(i) = 1$ and each $x_i \leq x$ then

$$\begin{aligned} x_S^{|S|-2-\deg(x)} x &= x_S^{|S|-2-\deg(x)} \sum_{i \in I} \epsilon(i) x_i \\ &= \sum_{i \in I} \epsilon(i) x_S^{|S|-2-\deg(x)} x_i \\ &= \sum_{i \in I} \epsilon(i) x_S^{|S|-2} \text{ by lemma 9.5.2} \\ &= x_S^{|S|-2} \sum_{i \in I} \epsilon(i) \\ &= x_S^{|S|-2} \end{aligned}$$

□

LEMMA 9.5.15. *Let $\mathcal{L}$ be a collection of subsets of $S$, $y = \prod_{U \in \mathcal{L}} y_U^{n_U}$ in $N[S]$ be a proper monomial and $T \supseteq \mathrm{supp}(\mathcal{L})$. Suppose that $n_T + \deg(y) \leq |T| - 2$ then $y_T^{n_T} y$ is a proper monomial.*

PROOF. As $y$ is proper we may write $q_S(y) = \sum_{i \in I} \epsilon(i) q_S(y_i)$ with $\sum_{i \in I} \epsilon(i) = 1$ and $y_i \leq y$ then

$$q_S(y_T^{n_T} y) = \sum_{i \in I} \epsilon(i) q_S(y_T^{n_T} y_i)$$

Because $T \supseteq \text{supp}(\mathcal{L})$ and $n_T + \deg(y) \leq |T| - 2$ it is clear that $y_T^{n_T} y \in N[S]$ and each $y_T^{n_T} y_i \in N_F[S]$. We also see that $y_T^{n_T} y_i \leq y_T^{n_T} y$, thus $y_T^{n_T} y$ is proper.

$\square$

LEMMA 9.5.16. *Let $\mathcal{L}$ be a collection of subsets of $S$ and $y \in N[S]$ a monomial on $\mathcal{L}$. Suppose $q_S(y) = \sum_{i \in I} \epsilon(i) q_S(y_i)$ with $\epsilon(i) = \pm 1$ and $\sum_{i \in I} \epsilon(i) = 1$. Then if each $y_i$ is proper and $y_i \leq y$ then $y$ is proper.*

PROOF. As each $y_i$ is proper then we may write $q_S(y_i) = \sum_{j \in I_i} \epsilon(i,j) q_S(y_{i,j})$ with $\sum_{j \in I_i} \epsilon(i,j) = 1$ and $y_{i,j} \leq y_i$ then

$$q_S(y) = \sum_{i \in I} \sum_{j \in I_i} \epsilon(i) \epsilon(i,j) q_S(y_{i,j})$$

then since each $y_{i,j} \leq y_i$ and $y_i \leq y$ we have $y_{i,j} \leq y$ and

$$\begin{aligned} \sum_{i \in I} \sum_{j \in I_i} \epsilon(i)\epsilon(i,j) &= \sum_{i \in I} \epsilon(i) \sum_{j \in I_i} \epsilon(i,j) \\ &= \sum_{i \in I} 1 \epsilon(i) \\ &= 1 \end{aligned}$$

thus $y$ is proper

$\square$

LEMMA 9.5.17. *Suppose $x$ and $y$ in $N[S]$ are proper so that $q_S(x) = \sum_{i \in I} \epsilon(i) q_S(x_i)$ and $q_S(y) = \sum_{j \in J} \epsilon(j) q_S(y_j)$ and for each $i \in I$, $j \in J$ $\text{shape}(x_i y_j)$ is a tree then if $xy \in N[S]$ we have that $xy$ is proper.*

PROOF.

$$q_S(xy) = \sum_{i \in I} \sum_{j \in J} \epsilon(i) \epsilon(j) q_S(x_i y_j)$$

as $x_i \leq x$ and $y_i \leq y$ then by lemma 9.5.12 we have $x_i y_j \leq xy$ and shape($x_i y_j$) being a tree is given thus it remains to prove the epsilon condition

$$\sum_{i \in I}\sum_{j \in J} \epsilon(i)\epsilon(j) = \sum_{i \in I}\epsilon(i) \sum_{j \in J}\epsilon(j)$$
$$= 1 \times 1$$
$$= 1$$

□

LEMMA 9.5.18. *Let $\mathcal{L}$ be a collection of subsets of $S$ and $y = \prod_{U \in \mathcal{L}} y_U^{n_U}$ be a monomial in $N[S]$. Put $\mathcal{C}(\mathcal{L})$ to be the set of equivalence classes of $\mathcal{L}$ and for each $\mathcal{C} \in \mathcal{C}(\mathcal{L})$ let $y(\mathcal{C})$ be the restriction of the monomial $y$ to $\mathcal{C}$ then if for every $\mathcal{C} \in \mathcal{C}(\mathcal{L})$ the elements $y(\mathcal{C})$ are proper monomials associated to $\mathcal{C}$ then $y$ is a proper monomial associated to $\mathcal{L}$.*

PROOF. We first note that as $\mathcal{L} = \coprod_{\mathcal{C} \in \mathcal{C}(\mathcal{L})} \mathcal{C}$ we have $y = \prod_{\mathcal{C} \in \mathcal{C}(\mathcal{L})} y(\mathcal{C})$. We are given that each $y(\mathcal{C})$ is a proper monomial of $\mathcal{C}$, we may then use corollary 9.5.7 and apply the obvious extension of lemma 9.5.17 to deduce our result.

□

LEMMA 9.5.19. *Let $\mathcal{L}$ be a collection of subsets of $S$ and $\mathcal{C} = \{T_1, ..., T_n\}$ be a chain on $\mathcal{L}$. Consider the monomial $y = \prod_{T \in \mathcal{C}} y_T^{n_T}$ then if $y \in N[S]$ it is proper so that $x = q_S(y) = \sum_{i=1}^{n} \epsilon(i) x_i$ and the top of each tree is $T = \text{supp}(\mathcal{C})$.*

PROOF. First put $T = \text{supp}(\mathcal{C}) = \bigcup_{U \in \mathcal{C}} U$ and $n = \deg(y)$. We know that

$$\prod_{U \in \mathcal{C}}(x_U - x_T) = 0$$
$$\prod_{U \in \mathcal{C}} x_U^{n_U - 1} \prod_{U \in \mathcal{C}}(x_U - x_T) = 0$$
$$\prod_{U \in \mathcal{C}}(x_U^{n_U} - x_T x_V^{n_U - 1}) = 0$$

So on expanding this we get

$$x = q_S(y) = \sum_{i \in I} \epsilon(i) x_T^{n-\deg(y_i)} q_S(y_i) \text{ clearly } \epsilon(i) = \pm 1 \text{ and } \sum_{i \in I} \epsilon(i) = 1$$

Put $\mathcal{C}_i = \text{shape}(y_i)$ then we observe that $\mathcal{C}_i \subseteq \mathcal{C}$ is a possibly disjoint union of chains, $y_i$ is a proper divisor of $y$ thus $\deg(y_i) < \deg(y)$ and $y_i \leq y$. Thus by induction we may suppose the claim is true for each equivalence class of $\mathcal{C}_i$ and by lemma 9.5.18 we can deduce that $y_i$ is a proper monomial associated to $\mathcal{C}_i$. It now remains to check that $y_T^{n-\deg(y_i)} y_i$ is a proper monomial and $y_T^{n-\deg(y_i)} y_i \leq y$, because $T = \text{supp}(\mathcal{C})$ applying lemma 9.5.15 we see this is clear and then apply lemma 9.5.16 to obtain the result. □

COROLLARY 9.5.20. *Every monomial $y \in N[S]$ is a proper monomial.*

PROOF. The proof is immediate from lemma 9.5.18 then 9.5.19.

□

LEMMA 9.5.21. *Let $\mathcal{L}$ be a collection of subsets of $S$ and consider the monomial $x = \prod_{T \in \mathcal{L}} x_T^{n_T}$ then $x = 0$ if and only if $x \notin \overline{N}[S]$ further if $x$ is non zero and $\deg(x) = |T| - 2$ that is has maximal degree then $x = x_T^{|T|-2}$ where $T = \text{supp}(\mathcal{U})$.*

PROOF. If $y \notin N[S]$ we have already seen that $x = q_S(y) = 0$. Given any $y \in N[S]$ we have show that $y$ is proper and thus $x = q_S(y)$ is non-zero. As $\overline{N}[S] = q_S(N[S])$ we see our claim is true.

□

## 9.6. Minimal relations

In this short section we consider the relations for connected sets and explain how these can be deduced from a minimal set. We will prove this explicitly, however this can already be seen from the work we have done.

DEFINITION 9.6.1. Let $\mathbb{Z}_S$ be the graded polynomial algebra over $\mathbb{Z}$ with one generator $y_T \in \mathbb{Z}_S^2$ for each $T \subseteq S$ with $|T| > 2$, that is we define $\mathbb{Z}_S = \mathbb{Z}[\, y_T \mid T \subseteq S \text{ and } |T| > 2\,]$.

Let $J_S$ be the ideal generated by the following relations

(1) Let $T \subseteq S$ then $y_T^{|T|-1} = 0$
(2) Let $T \subseteq S$ and $\mathcal{T}$ a $T$-tree of depth 2 then $y_T^{m(\mathcal{T},T)} \prod_{U \in M(\mathcal{T},T)} (y_T - y_U) = 0$
(3) For all $U, V \subseteq S$ with $U \cap V$ non-empty we have $(y_{U \cup V} - y_U)(y_{U \cup V} - y_V) = 0$

then we define the ring $Q_S$ by $Q_S = \mathbb{Z}_S/J_S$

REMARK 9.6.2. Clearly we have the inclusion of ideals $J_S \subseteq I_S$. By considering our analysis of the filtration in section 9.4 in particular the admissability of monomials we easily see that $Q_S$ is also the cohomology ring of $\overline{\mathcal{M}}_S$. Thus $R_S = Q_S$ as claimed. It would be better to see this algebraically. This we do next as it is straightforward.

LEMMA 9.6.3. *Let $\mathcal{L}$ be a connected set then $\prod_{U \in \mathcal{L}} (x_T - x_U) = 0$ in $Q_S$ where $T = \mathrm{supp}(\mathcal{L})$.*

PROOF. The proof is by induction on $\mathcal{L}$. Clearly we may suppose that $\mathcal{L}$ is not a forest. When $|\mathcal{L}| = 2$ the claim is immediate. Suppose the case is true for $|\mathcal{L}| = n$. Let $\mathcal{L}$ be connected with $|\mathcal{L}| = n + 1$. Then we can find $U, V \in \mathcal{L}$ with $\{U, V\}$ not a forest. Put $W = U \cup V$ so that $(x_W - x_U)(x_W - x_V) = 0$. This implies $(x_T - x_U)(x_T - x_V) = f(x_T - x_W)$. Put $\mathcal{L}'$ to be $\mathcal{L}$ with $U, V$ removed and $W$ added. Then $\prod_{U \in \mathcal{L}} (x_T - x_U) = f \prod_{V \in \mathcal{L}'} (x_T - x_V)$, $|\mathcal{L}'| = n$ and clearly $\mathcal{L}'$ is connected, $\mathrm{supp}(\mathcal{L}') = T$ thus inductively we obtain our relation from the last equation. □

COROLLARY 9.6.4. *Let $S$ be a finite set then we have an equality $R_S = Q_S$* □

## 9.7. An isomorphism from $A^*(\overline{\mathcal{M}}_S)$ to $H^*(\overline{\mathcal{M}}_S)$.

In this section we consider the natural map $\mathrm{cl} : A^*(\overline{\mathcal{M}}_S) \to H^*(\overline{\mathcal{M}}_S)$ that we proved in chapter 6 was an isomorphism. We will state a presentation for $A^*(\overline{\mathcal{M}}_S)$ given by Keel in [9] and compute explicitly this map and also its inverse.

DEFINITION 9.7.1. Let $S$ be a finite set with $|S| = n$ and $T \subseteq S$ then for any $i, j \in T$ with $i \neq j$ we define the following spaces.

$$\begin{aligned} F^{i,j}(T, \mathbb{C}) &= \{\, f : T \to \mathbb{C} \mid f(i) = f(j) \,\} \\ V_T^{i,j} &= F^{i,j}(T, \mathbb{C})/C_T \\ PV_T^{i,j} &\subseteq PV_T \end{aligned}$$

We have the maps $\pi_T^S : \overline{\mathcal{M}}_S \to PV_T$ and define $X_T^{i,j} = (\pi_T^S)^{-1}(PV_T^{i,j}) \subseteq \overline{\mathcal{M}}_S$. For any $T \subset S$ we also define the divisors $D^T \subseteq \overline{\mathcal{M}}_S$ by $D^T = \overline{\mathcal{M}}_{\{S,T\}}$.

REMARK 9.7.2. Note that these are all codimension 1 subspaces of their respective sets and Keel proves that the $D^U$ generate the Chow ring.

LEMMA 9.7.3. $X_T^{i,j} = \sum \{\, D^U \mid U \supseteq \{i, j\} \text{ and } U \not\supseteq T \,\}$.

PROOF. Choose a set $U \subseteq S$ with $i, j \in U$ and $U \not\supseteq T$. Let $\underline{M} \in D^U$ and put $\mathcal{T} = \text{type}(\underline{M})$, so $\text{type}(\underline{M}) \supseteq \{S, U\}$. Put $V = \text{root}(T)$ then $U \subset V$ since $\text{type}(\underline{M})$ is a tree, $U \cap V$ is non-empty and $U \not\supseteq T$. Then $\{i, j\} \subseteq T \cap U \subseteq U$ and $\pi_{T \cap U}^T M_T = \pi_{T \cap U}^T \pi_T^V M_V = \pi_{T \cap U}^V M_V = 0$ thus $M_T \in X_T^{i,j}$. Now suppose $\underline{M} \in X_T^{i,j}$, put $\mathcal{T} = \text{type}(\underline{M})$ and $V = \text{root}(T)$ then there exists a $U \in M(\mathcal{T}, V)$ with $i, j \in U$ and we must have $U \not\supseteq T$ for otherwise $U \supseteq T$ and $\text{root}(T) \subseteq U$ a contradiction. Next put $\mathcal{U} = \{S, U\}$ then $\mathcal{T} \supseteq \mathcal{U}$ and $\underline{M} \in D^U$. This proves the claim. □

DEFINITION 9.7.4. Let $L_T$ be the tautological line bundle over $PV_T$. Put $N_T = (\pi_T^S)^*(L_T)$ and $N_T^*$ the dual bundle over $\overline{\mathcal{M}}_S$ where $\pi_T^S : \overline{\mathcal{M}}_S \to PV_T$ is the usual map. Then given a linear map $\alpha \in \hom(V_T, \mathbb{C})$ we write $s_\alpha$ for the induced section on $N_T^*$ given by $s_\alpha(\underline{M}) = \alpha|_{M_T} \in \hom(M_T, \mathbb{C})$. The zero set of this section is clearly given by the space $\{\, \underline{M} \in \overline{\mathcal{M}}_S \mid M_T \leq \ker(\alpha) \,\}$

LEMMA 9.7.5. *For every $T \subseteq S$ and any $i, j \in T$ with $i \neq j$, there is an algebraic section $s_T^{i,j}$ of the usual line bundle $N_T^*$ over $\overline{\mathcal{M}}_S$ whose zero set is $X_T^{i,j}$.*

PROOF. Let $N_T$ be the usual bundle over $\overline{\mathcal{M}}_S$ and $N_T^*$ be the dual bundle. Define $s_\alpha : \overline{\mathcal{M}}_S \rightarrowtail EN_T^*$ as follows. Define $\alpha : F(T, \mathbb{C}) \to \mathbb{C}$ by $\alpha(f) = f(i) - f(j)$ and the induced map $\overline{\alpha} : V_T \to \mathbb{C}$ by $\overline{\alpha}(f + C_T) = \alpha(f)$ then clearly,

$$\begin{aligned} s_\alpha(\underline{M}) = 0 &\iff M_T \leq V_T^{i,j} = \ker(\overline{\alpha}) \\ &\iff M_T \in PV_T^{i,j} \\ &\iff \underline{M} \in (\pi_T^S)^{-1}(PV_T^{i,j}) = X_T^{i,j} \end{aligned}$$

□

COROLLARY 9.7.6. *The cohomology class $[X_T^{i,j}] \in H^2(\overline{\mathcal{M}}_S)$ is independent of $i, j \in T$ with $i \neq j$ and is $x_T = -e(N_T^*)$ the euler class of $N_T^*$*

COROLLARY 9.7.7. *We have the following equations*

$$\begin{aligned} x_S - x_T &= \sum_{T \subseteq U \subset S} D^U \\ x_S &= \sum_{\{i,j\} \subseteq U \subset S} D^U \\ D^U &= \sum_{U \subseteq V \subseteq S} (-1)^{|V/U|} x_V \end{aligned}$$

PROOF. Put $E^U = \sum_{U \subseteq V \subseteq S} (-1)^{|V/U|} x_V$. First note that $D^S$ is not a divisor but we can make the arbitrary definition $D^S := E^S = -x_S$ without affecting the following sum. Then by lemma 9.7.3 have shown that $x_T = \sum_{\{i,j\} \subseteq U \subseteq S} D^U - \sum_{T \subseteq U \subseteq S} D^U$, we first show that $x_T = \sum_{\{i,j\} \subseteq U \subseteq S} E^U - \sum_{T \subseteq U \subseteq S} E^U$. Let $W \subseteq S$, then we consider the number of occurrences of $x_W$ in the sum $\sum_{T \subseteq U \subseteq S} E^U$. That is, we require the number of sets $U$ such that $T \subseteq U \subseteq W$. If $T \not\subseteq W$ then there are no such occurrences otherwise there are $2^{|W \setminus T|}$ sets of this type. This is even if $W \neq T$ and takes the value 1 otherwise. Because the $x_W$ alternate in sign we see that there is no total contribution of $x_W$ in the sum if $W \neq T$ and a total contribution of $-x_T$ in the case $W = T$. We next prove that the sum $\sum_{\{i,j\} \subseteq U \subseteq S} E^U$ is zero. Let $W \subseteq S$, then as before we consider the number of occurrences of $x_W$ in the second sum. That is, we require the number of sets $U$ such that $\{i,j\} \subseteq U \subseteq W$. There are no such occurrences if $\{i,j\} \not\subseteq W$ and $2^{|W|-2}$ such sets otherwise. Because $x_W = 0$ if $W = \{i,j\}$ we may suppose $|W| > 2$. But then the number of occurrences of $x_W$ is even and thus the total contribution of $x_W$ is zero by there alternating sum. Thus we have proven that the only contribution is by $x_T$ and the sum is as claimed.

Next put $Q^T = D^T - E^T$ and subtract the equations

$$x_T = \sum_{\{i,j\} \subseteq U \subseteq S} E^U - \sum_{T \subseteq U \subseteq S} E^U$$

$$x_T = \sum_{\{i,j\} \subseteq U \subseteq S} D^U - \sum_{T \subseteq U \subseteq S} D^U$$

then we obtain $\sum_{T \subseteq U \subseteq S} Q^U = \sum_{\{i,j\} \subseteq U \subseteq S} Q^U$ for all $T$ and all $\{i,j\} \subseteq T$. Putting $T = S$ we obtain $\sum_{\{i,j\} \subseteq U \subseteq S} Q^U = 0$ because $Q^S = 0$ thus substituting this back we obtain $\sum_{T \subseteq U \subseteq S} Q^U = 0$ for all $T$. We now use a downward induction on $|T|$ to deduce that $Q^T = 0$ for all $T$. For $|T| = n = |S|$ the claim is clear. Assume the claim is true for $m > k$ for some $k < n$ then for a set $T$ with $|T| = k$ we see immediately from $\sum_{T \subseteq U \subseteq S} Q^U = 0$ that $Q^T = 0$ thus $D^T = E^T$ and $\sum_{\{i,j\} \subseteq U \subseteq S} D^U = 0$ as claimed.

□

DEFINITION 9.7.8. Let $\mathbb{Z}_S$ be the polynomial ring $\mathbb{Z}[\, D^U \mid U \subset S\,]$ and $I_S$ the ideal generated by

$$\sum_{\{i,j\} \subseteq U \subseteq S} D^U = \sum_{\{k,l\} \subseteq U \subseteq S} D^U \text{ for every } i,j,k,l$$

$$D^T D^U = 0 \text{ unless } T \subset U,\ U \subset T,\ T \subset U^c,\ U \subset T^c$$

then we define $T_S = \mathbb{Z}_S / I_S$ and we have the following theorem due to Keel [9]

THEOREM 9.7.9. *For each finite set $S$ the Chow ring of $\overline{\mathcal{M}}_S$ is $T_S$*

## 9.8. The ring $R_\mathcal{L}$

Having computed the cohomology ring of $\overline{\mathcal{M}}_S$ we observe in this section that in practice everything we have done for $\overline{\mathcal{M}}_S$ could equally well have been applied to the space $\overline{\mathcal{M}}_\mathcal{L}$ for a thicket $\mathcal{L}$. We could define in the evident way the ring $R_\mathcal{L}$ and prove that $R_\mathcal{L}$ is the cohomology ring of $\overline{\mathcal{M}}_\mathcal{L}$, should one be interested in such an object. We also observe that given our approach in section 4 of chapter 6 we could without much further effort deduce the following commutative diagram in even degrees

$$\begin{array}{ccc} R_{\mathcal{L}_+} & \longrightarrow & R_{\overline{\mathcal{L}}} \otimes R_{\{T\}} \\ \uparrow & & \uparrow \\ R_\mathcal{L} & \longrightarrow & R_{\overline{\mathcal{L}}} \end{array}$$

and presumably one could deduce the following short exact sequence in even degrees

$$R_\mathcal{L} \to R_{\mathcal{L}_+} \oplus R_{\overline{\mathcal{L}}} \to R_{\overline{\mathcal{L}}} \otimes R_{\{T\}}$$

This would offer an alternative approach to the analysis of $R_S$ and in particular would enable us to deduce most of our results without comparing them to the cohomology ring.

# Bibliography

[1] P. Deligne and D. Mumford. Irreducibility of the space of curves of given genus. *Inst. Hautes Études Sci. Publ. Math.*, 36:75–109, 1969.

[2] D. Eisenbud and J. Harris. *The Geometry of Schemes*, volume 197 of *Graduate Texts in Mathematics*. Springer–Verlag, 1991.

[3] W. Fulton. *Intersection Theory*. Springer, 1998

[4] P. Griffiths and J. Harris. *Principles of Algebraic Geometry*. Wiley, 1994

[5] A. Grothendieck and J. Dieudonné. Eléments de Géométrie Algébrique III (i). *Publ. Math. IHES*, 11, 1961.

[6] R. H. Hartshorne. *Algebraic Geometry*, volume 52 of *Graduate Texts in Mathematics*. Springer–Verlag, 1977.

[7] M. M. Kapranov. *Veronese curves and the Grothendieck-Knudsen moduli space $M(0,n)$*, Journal of Algebraic Geometry, 2, (1993), p. 239–262.

[8] M. M. Kapranov. *Chow quotient of Grassmannians I*, Adv. Soviet. Math. 16, part 2,(October 1992), p. 29-110.

[9] S. Keel. *Intersection theory of moduli space of stable N-pointed curves of genus zero*, Transactions of the American mathematical society. Volume 330, issue 2 (April 1992), p. 545-574.

[10] F. F. Knudsen. The projectivity of the moduli space of stable curves. II. The stacks $M_{g,n}$. *Math. Scand.*, 52(2):161–199, 1983.

[11] J. W. Milnor and J. D. Stashheff. *Characteristic Classes*. Princeton University Press, 1974.

[12] D. Mumford. *Abelian Varieties*, volume 5 of *Tata institute of fundamental research series in mathematics*. Oxford University Press, 1970.

[13] I. R. Shafarevich. *Basic Algebraic Geometry* 2 Springer-Verlag, 1994.